\documentclass[11pt,reqno]{amsart}
\usepackage[T1]{fontenc}
\usepackage{lmodern}
\usepackage[letterpaper,margin=1.08in]{geometry}
\usepackage{microtype}
\usepackage{amsmath,amssymb,amsthm,mathtools}
\usepackage{booktabs,array,longtable}
\usepackage{enumitem}
\usepackage{needspace}
\usepackage{xcolor}
\usepackage[hidelinks]{hyperref}
\usepackage[nameinlink,noabbrev]{cleveref}
\usepackage{mathrsfs}

\numberwithin{equation}{section}

\newtheorem{theorem}{Theorem}[section]
\newtheorem{mainthm}{Theorem}

\newtheorem{proposition}[theorem]{Proposition}
\newtheorem{lemma}[theorem]{Lemma}
\newtheorem{corollary}[theorem]{Corollary}
\theoremstyle{definition}

\theoremstyle{remark}
\newtheorem{remark}[theorem]{Remark}

\DeclareMathOperator{\GU}{GU}
\DeclareMathOperator{\U}{U}
\DeclareMathOperator{\GL}{GL}

\DeclareMathOperator{\Res}{Res}
\DeclareMathOperator{\Hom}{Hom}
\DeclareMathOperator{\Ad}{Ad}
\DeclareMathOperator{\Lie}{Lie}

\DeclareMathOperator{\tr}{tr}
\DeclareMathOperator{\im}{im}
\DeclareMathOperator{\rk}{rk}
\DeclareMathOperator{\diag}{diag}
\DeclareMathOperator{\N}{N}
\DeclareMathOperator{\Gal}{Gal}
\DeclareMathOperator{\supp}{supp}
\DeclareMathOperator{\End}{End}
\DeclareMathOperator{\Tr}{Tr}

\DeclareMathOperator{\Aut}{Aut}

\newcommand{\Kt}{\widetilde K}

\newcommand{\Ddist}{\mathcal D}

\newcommand{\triv}{\mathrm{triv}}
\newcommand{\eps}{\epsilon}
\newcommand{\lam}{\lambda}
\newcommand{\barF}{\overline F}

\newcommand{\Mat}{\operatorname{M}}
\newcommand{\Gm}{\mathbb G_m}

\setlist{nosep}
\allowdisplaybreaks
\newcolumntype{L}[1]{>{\raggedright\arraybackslash}p{#1}}
\DeclareMathOperator{\Ext}{Ext}
\newcommand{\cX}{\mathcal X}

\newcommand{\cAZero}{\mathcal A_0}
\newcommand{\Fb}{\overline F}
\newcommand{\Co}{\mathbf C}
\newcommand{\Vlegs}{\mathcal V}
\newcommand{\Efac}{\mathsf E}
\newcommand{\Rfac}{\mathsf R}
\newcommand{\Qfac}{\mathsf Q}
\newcommand{\Hfac}{\mathsf H}
\newcommand{\qvr}{\mathscr Q}
\DeclareMathOperator{\Irr}{Irr}
\DeclareMathOperator{\Ind}{Ind}
\newcommand{\cInd}{\operatorname{c\!-Ind}}

\crefname{mainthm}{Theorem}{Theorems}
\Crefname{mainthm}{Theorem}{Theorems}
\makeatletter
\renewcommand{\paragraph}{%
  \@startsection{paragraph}{4}{0pt}%
  {0pt}{-0.8em}{\normalfont\bfseries}%
}
\makeatother
\hypersetup{
 pdftitle={Multiplicity one for a linear spherical pair and its unitary analogue},
 pdfauthor={Yuxuan Ge},
 pdfsubject={Multiplicity one for two spherical pairs and unitary similitude Shalika models},
 pdfkeywords={spherical pairs, multiplicity one, Shalika models, invariant distributions, theta correspondence, unitary similitude groups}
}
\AtBeginEnvironment{thebibliography}{\raggedright\microtypesetup{expansion=false}}

\newcommand{\C}{\mathbb C}
\newcommand{\cS}{\mathcal S}
\newcommand{\abs}[1]{\lvert#1\rvert}

\title[Multiplicity one for a linear spherical pair and its unitary analogue]{Multiplicity one for a linear spherical pair and its unitary analogue}
\author{Yuxuan Ge}
\address{School of Mathematics, University of Minnesota, Minneapolis, MN 55455, USA}
\email{ge005678@umn.edu}
\date{}
\begin{document}
\begin{abstract}
We prove multiplicity one for a linear spherical pair associated with \(\mathrm{GL}_4(F)\times\mathrm{GL}_2(F)\) and its unitary similitude analogue associated with \(\mathrm{GU}(2,2)(F)\times\mathrm{GU}(1,1)(F)\). Both results hold for all irreducible admissible representations over non-Archimedean local fields of characteristic zero. Generalized Harish–Chandra descent reduces the proofs to the analysis of closed orbits and their normal representations. The unitary proof combines reductions to the linear case with a rational classification of closed orbits and a separate analysis of the remaining normal representations.
We also establish multiplicity one for unitary Shalika periods of \(\mathrm{GU}(2,2)(F)\) with arbitrary irreducible admissible coefficients. A theta-period identity relates these periods to the unitary spherical pair. By varying the splitting character, we prove that every irreducible admissible representation of \(\mathrm{GU}(2,2)(F)\) occurs as a quotient of a suitable big theta lift, allowing the multiplicity-one bound to pass to the Shalika periods. The linear argument also recovers Rallis's multiplicity-one theorem for twisted Jacquet modules of \(\mathrm{GL}_4(F)\).
\end{abstract}
\maketitle
\setcounter{tocdepth}{2}
\tableofcontents
\section{Introduction and main theorems}\label{sec:intro}

\subsection{The spherical pairs}
Let $F$ be a non-Archimedean local field of characteristic zero.
We first consider the linear pair
\begin{equation}\label{eq:main-linear-subgroup}
 G_{\mathrm L}=\GL_4(F)\times\GL_2(F),\qquad
 H_{\mathrm L}=\{(\diag(a,b),b):a,b\in\GL_2(F)\}.
\end{equation}
\begin{mainthm}[Linear spherical pair]\label{thm:main-linear}\label{m1:thm:main}
For every irreducible admissible representation $\pi$ of $G_{\mathrm L}$,
\[
 \dim\Hom_{H_{\mathrm L}}(\pi,\mathbb C)\leq1.
\]
\end{mainthm}

For the unitary statement let $E/F$ be a quadratic field extension,
write $z\mapsto\bar z$ for its nontrivial automorphism, and put
$E^1=\{z\in E^\times:z\bar z=1\}$. Forms are linear in the first
variable. For a Hermitian space $Z$, write $\U(Z)$ and $\GU(Z)$
for its groups of $F$-rational isometries and similitudes.
Let $V_1,V_2$ be hyperbolic Hermitian planes. Write $\lambda$ for
every similitude multiplier. With $V=V_1\perp V_2$, define
\begin{equation}\label{eq:main-unitary-subgroup}
 \begin{split}
 G_{\mathrm U}&=\GU(V)\times\GU(V_2),\\
 H_{\mathrm U}&=\{(u\oplus h,h):u\in\GU(V_1),\ h\in\GU(V_2),
                              \ \lambda(u)=\lambda(h)\}.
 \end{split}
\end{equation}
Since $V$ is hyperbolic,
$G_{\mathrm U}\simeq\GU(2,2)(F)\times\GU(1,1)(F)$.
In the antidiagonal realization the two planes are the outer and middle
coordinate planes; the fixed permutation in
Section~\ref{m2:red:pair} identifies that realization with
\eqref{eq:main-unitary-subgroup}.
\begin{mainthm}[Unitary similitude spherical pair]\label{thm:main-unitary}\label{m2:thm:main}
For every irreducible admissible representation $\pi$ of $G_{\mathrm U}$,
\[
 \dim\Hom_{H_{\mathrm U}}(\pi,\mathbb C)\leq1.
\]
\end{mainthm}
These are Models~1 and~2 among the strongly tempered spherical pairs
studied by Wan--Zhang~\cite{WanZhang25}. They computed the unramified local relative characters and established multiplicity formulas in the tempered setting. For the linear pair, global questions have also been
studied by Pollack--Wan--Zydor~\cite{PWZ21} and by
Cauchi--Gutierrez Terradillos~\cite{CGTLinear}, including connections
with central $L$-values.

\subsection{Unitary similitude Shalika uniqueness}
Fix a nontrivial additive character $\psi:F\to\mathbb C^\times$.
Let $W=X\oplus Y$ be the split skew-Hermitian four-space with
$X=Y=E^2$ and matrix
$J_W=\left(\begin{smallmatrix}0&I_2\\-I_2&0\end{smallmatrix}\right)$.
Let $B$ be a hyperbolic Hermitian form on $X$, with matrix
$\mathsf B=\mathsf B^\dagger$, where $a^\dagger=\bar a^{\,t}$.
Write $\lambda_B$ for the multiplier on $\GU(B)$ and $\tr_E$ for
the matrix trace over $E$. The embedding $m_B:\GU(B)\to\GU(W)$,
the Shalika subgroup, and its character are
\begin{equation}\label{eq:intro-shalika}
 \begin{gathered}
 N(X)=\left\{n(u)=\begin{pmatrix}I_2&u\\0&I_2\end{pmatrix}:
                                    u^\dagger=u\right\},\\
 m_B(a)=\diag\bigl(a,\lambda_B(a)\bar a^{-t}\bigr),\qquad
 S_B^{\GU}=m_B(\GU(B))N(X),\\
 \psi_B(n(u))=\psi\bigl(-\tfrac12\tr_E(u\mathsf B)\bigr).
 \end{gathered}
\end{equation}
The trace $\tr_E(u\mathsf B)$ belongs to $F$. Representations of
$\GU(B)$ are viewed as representations of $m_B(\GU(B))$ through
this embedding.
\begin{mainthm}[Unitary similitude Shalika models]\label{thm:main-gu}\label{gu:thm:main}
For every irreducible admissible representation $\Pi$ of $\GU(W)$
and every irreducible admissible representation $\sigma$ of $\GU(B)$,
\begin{equation}\label{gu:eq:main}
 \dim\Hom_{S_B^{\GU}}(\Pi,\sigma\boxtimes\psi_B)\leq1.
\end{equation}
\end{mainthm}

For ordinary unitary groups, Gan~\cite{GanPeriods} and Chen--Gan
\cite[Proposition~3.4 and Theorem~3.5]{ChenGan} establish the theta-period
identity for Shalika periods with arbitrary coefficient representations. The distinction between
isometry and similitude groups matters: Shalika periods for ordinary unitary groups can have multiplicity greater than one
\cite[Theorem~4.3]{ChenGan}. Here we establish the corresponding period identity
for unitary similitudes and use it to prove Theorem~\ref{thm:main-gu}.

The twisted Jacquet module of the similitude oscillator gives a
Hom-space identity (Theorem~\ref{gu:thm:hom}) relating the
Shalika model to the unitary spherical pair. Combined with
Theorem~\ref{thm:main-unitary}, it bounds the Shalika Hom space of a
big theta lift by one. To pass to an arbitrary irreducible admissible
$\GU(W)$-representation $\Pi$, Theorem~\ref{gu:asm:coverage} provides a
unitary character $\chi$ of $E^\times/F^\times$ and
$\tau\in\Irr(\GU(V))$ with
\begin{equation}\label{gu:eq:intro-coverage}
 \Theta^{\GU}_{\chi,1,1}(\tau)\twoheadrightarrow\Pi.
\end{equation}
The surjection in \eqref{gu:eq:intro-coverage} induces,
by precomposition, an injection from the Shalika Hom space
of $\Pi$ into that of the big theta lift.
Theorem~\ref{thm:main-gu} follows from the
multiplicity-one bound for the latter.
The coverage theorem uses analytic equal-rank theta
dichotomy and highly ramified stability to obtain
nonvanishing in the prescribed quasi-split tower.

\subsection{Outline of the proofs}
Section~\ref{sec:quiver} develops the determinant-open quiver geometry
shared by the two pairs: closed orbits, stabilizers, and normal
representations. Section~\ref{sec:linear} establishes the linear
two-leg theorem and treats the identity descendant.
These results, together with generalized
Harish--Chandra descent, yield Theorem~\ref{thm:main-linear}.

For the unitary pair, this common geometry is only the starting point.
Section~\ref{sec:unitary-geometry} determines the Hermitian domains
admitting rational points, proves uniqueness of the rational orbit
inside each admitted closed geometric orbit, and deduces rational
transpose stability. It then determines the corresponding \(F\)-rational normal representations and their corrected odd operators. Section~\ref{sec:unitary-descendants}
shows that stabilizer-invariant distributions on these
normal spaces are fixed by the corrected
involutions. The split-adjoint regular case
reduces to the linear two-leg theorem, while
the anisotropic case is treated using Hermitian
Gram fibers. The scalar case requires a separate
identity-descendant argument. This proves
Theorem~\ref{thm:main-unitary}.

Section~\ref{sec:gl-shalika} recovers Rallis's multiplicity-one
theorem for twisted Jacquet modules of $\GL_4(F)$, as recorded by
Prasad~\cite[Theorem~1]{PrasadShalika}, using the linear spherical-pair
theorem and Gan's theta-period identity~\cite[Theorem~3.1]{GanPeriods}.
Section~\ref{sec:gu-transfer} establishes the corresponding Shalika
period identity for unitary similitude theta lifts. Finally,
Section~\ref{sec:gu-coverage} shows that, after varying the splitting
character, every irreducible admissible representation of $\GU(W)$ is
a quotient of a suitable big theta lift. This completes the proof of
Theorem~\ref{thm:main-gu}.

Appendices~A and~B supply the functional-analytic and normalization
details for the theta-period arguments. Appendix~C treats general
splitting data, while Appendix~D records the rational transpose
transporter constructions used in the unitary orbit analysis.

\subsection{Notation}
Unless otherwise stated, representations of $p$-adic groups
are smooth and defined over $\mathbb C$.
We write $\pi^\vee$ for the smooth contragredient of $\pi$
and $\omega_\pi$ for its central character when $\pi$
is irreducible admissible.
We denote by $\Irr(G)$ the set of isomorphism classes
of irreducible admissible representations of $G$. The symbols $G_{\mathrm L},H_{\mathrm L}$ and
$G_{\mathrm U},H_{\mathrm U}$ denote the two spherical pairs.
In the theta part $G_\theta=\GU(V)$ and
$G'_\theta=\GU(W)$ denote the two members of the dual pair;
$K=\GU(B)$ denotes the Shalika coefficient group.
All similitudes scale their forms by $\lambda$ in the direction
$h(gv,gw)=\lambda(g)h(v,w)$. We distinguish ordinary transpose $t$, conjugate
transpose $\dagger$, and the form adjoint $*$.
Absolute values satisfy $|z|_E=|N_{E/F}z|_F$ and
$|s|_E=|s|_F^2$ for $s\in F^\times$.
Schwartz spaces consist of compactly supported locally constant
functions, and tensor products are algebraic. Ordinary and compact induction are written
$\Ind$ and $\cInd$; their normalization and duality are specified
where used.

\subsection*{Acknowledgments}
The author thanks Professor Dihua Jiang for helpful suggestions on an earlier
draft of this paper.

\subsection*{Disclosure of AI Use}
ChatGPT (including GPT-5.6 Sol and GPT-6 Astra) and Claude Opus 5.5 were used for proof exploration and drafting this
paper. The author takes responsibility for the
mathematical content.
\section{Distributional tools and multiplicity criteria}\label{sec:tools}
\label{m1:sec:tools}\label{m2:sec:tools}
This section collects the distributional tools used in both
multiplicity-one proofs. We recall the Gelfand--Kazhdan criterion
and generalized Harish--Chandra descent, together with the
localization and homogeneity arguments needed for the normal spaces.

Throughout the paper $F$ is a non-Archimedean local field of characteristic
zero, and $\Fb$ is an algebraic closure of $F$. An $l$-space is a locally
compact totally disconnected Hausdorff space; an $l$-group is a topological
group with this property. Write $\mathcal S(X)=C_c^\infty(X,\mathbf C)$ and
$\mathcal D(X)=\mathcal S(X)^*$ for its Schwartz functions and distributions. For an
$l$-group $L$ acting on $X$, our conventions are
\begin{equation}\label{m1:eq:distribution-conventions}
 (l\cdot f)(x)=f(l^{-1}x),\qquad (lD)(f)=D(f\circ l),\qquad
 \mathcal D(X)^{L,\chi}=\{D:lD=\chi(l)D\}.
\end{equation}
In particular, an invertible linear map $l$ acts on additive Haar measure
by $|\det l|^{-1}$. We omit $\chi$ when it is trivial.

\begin{proposition}[Open and closed support]\label{m2:tools:support}
For a closed subset $Z\subset X$ and its open complement $U$, there
is an exact sequence
\[
 0\longrightarrow\mathcal D(Z)\longrightarrow\mathcal D(X)
   \longrightarrow\mathcal D(U)\longrightarrow0.
\]
The first arrow identifies $\mathcal D(Z)$ with distributions supported
in $Z$, and the second arrow is restriction to the open set $U$.
If a group preserves $Z$, these maps are equivariant. In particular,
vanishing of a fixed equivariance space on both $U$ and $Z$ implies
its vanishing on $X$.
\end{proposition}

\begin{proof}
Compact-open extension of locally constant functions gives the exact
sequence
\[
 0\longrightarrow\mathcal S(U)\longrightarrow\mathcal S(X)
  \longrightarrow\mathcal S(Z)\longrightarrow0.
\]
Dualizing proves the assertion; taking equivariants preserves its
injection and kernel. In particular, a distribution supported at one
point is a multiple of its Dirac distribution.
\end{proof}

\subsection{The Gelfand--Kazhdan criterion and descent}
\label{m1:subsec:GK}
Let $G$ be reductive and $H\subset G$ a closed reductive subgroup. Suppose
that $\tau$ is an involutive anti-automorphism with $\tau(H)=H$. The action
of $H\times H$ is $(h_1,h_2)g=h_1gh_2^{-1}$.

\begin{proposition}[Gelfand--Kazhdan criterion]\label{m1:prop:GK}\label{m2:tools:gk}
If every $H\times H$-invariant distribution on $G$ is $\tau$-invariant,
then every irreducible admissible smooth complex representation $\pi$ of
$G$ satisfies
\[
 \dim\Hom_H(\pi,\Co)\,\dim\Hom_H(\pi^\vee,\Co)\leq1.
\]
\end{proposition}
This is \cite[Lemma~4.2]{Prasad90}; see also
\cite[Remark~2.4.6]{AGSGelfand}. The invariant measure on $G/H$ required in
Prasad's formulation exists because $G$ and $H$ are unimodular.
Contragredient duality will identify the two multiplicities for each block pair.

Whenever an involution $\tau$ normalizes an action of $L$, write
$\widetilde L=L\rtimes\langle\tau\rangle$ and let $\varepsilon$ be its
sign character: $\varepsilon|_L=1$ and $\varepsilon(\tau)=-1$.
Write $\mathcal D(X)^{L,\tau}$ for distributions fixed by both $L$ and $\tau$.
Transpose invariance of $\mathcal D(X)^L$ is equivalent to
$\mathcal D(X)^{\widetilde L,\varepsilon}=0$, by taking the odd part of a
distribution. We retain this convention for all extended groups.

For an algebraic action of a reductive group $L$ on a smooth affine
variety $X$, write $L_x$ for the stabilizer of $x$ and put
$N_x=T_xX/T_x(Lx)$, the normal space to the orbit in the ambient variety $X$.
We write $\mathcal D(Y)$ for $\mathcal D(Y(F))$ when $Y$ is an $F$-variety.
A point $x\in X(F)$ is semisimple if its $L(F)$-orbit
is analytically closed, equivalently its geometric orbit is Zariski
closed \cite[Theorem~2.3.8]{AGDescent}.

\begin{theorem}[Generalized Harish--Chandra descent]
\label{m1:thm:HC-descent}\label{m2:tools:descent}\label{m1:subsec:HC-descent}
Let $L$ be reductive, possibly disconnected, acting on a smooth affine
$F$-variety $X$, and let $\chi$ be a character of $L(F)$. If
$\mathcal D(N_x(F))^{L_x(F),\chi}=0$ for every semisimple $x\in X(F)$,
then $\mathcal D(X(F))^{L(F),\chi}=0$.
\end{theorem}
We use \cite[Theorem~3.1.1]{AGDescent}, with the convention on disconnected
reductive groups of that paper.

Suppose now that $\tau$ and the automorphism $\alpha$ of $L$ that it
induces are algebraic over $F$, and that every closed rational orbit is
$\tau$-stable. Choose $g_x\in L(F)$ with $g_x\tau(x)=x$ and set
$\tau_x=g_x\tau$. Then
\begin{equation}\label{m1:eq:corrected-extension}
 \widetilde L_x=L_x(F)\sqcup L_x(F)\tau_x,
 \qquad \tau_x^2=g_x\alpha(g_x)\in L_x(F).
\end{equation}
Changing $g_x$ multiplies $\tau_x$ by an element of $L_x(F)$,
so its action on invariant distributions is independent of this choice.
Since $\tau_x^2\in L_x(F)$, the induced action on
$L_x(F)$-invariant distributions on $N_x$ is an involution. Denote the
restricted sign character by $\varepsilon_x$.

\begin{corollary}\label{m1:cor:HC-invariance}
Under these hypotheses, if
$\mathcal D(N_x)^{\widetilde L_x,\varepsilon_x}=0$ at every closed rational
orbit, then every $L(F)$-invariant distribution on $X(F)$ is
$\tau$-invariant.
\end{corollary}
\begin{proof}
We apply Theorem~\ref{m1:thm:HC-descent} to the reductive algebraic finite
extension $\mathbf{\widetilde L}=L\rtimes_\alpha C_2$. Over $\Fb$ an
extended orbit is a union of at most two $L$-orbits of equal dimension.
If that union is closed, the boundary of either constituent lies in the
union. Boundary orbits have strictly smaller dimension, so this boundary
is empty. Conversely, a finite union of closed constituent orbits is
closed. The geometric--analytic comparison for both groups therefore
shows that the finite extension has exactly the same semisimple points.
Its orbit tangent equals the ordinary orbit tangent, since the extension
is finite. Hence its normal space is $N_x$, and its stabilizer and character
are precisely \eqref{m1:eq:corrected-extension}. The descent theorem applies.
\end{proof}

\subsection{Frobenius descent and orbital measures}\label{m1:subsec:Frobenius}
For a locally compact group $M$ our modular convention is
\begin{equation}\label{m1:eq:modulus}
 \int_M f(mh)\,dm=\Delta_M(h)\int_M f(m)\,dm,
 \qquad \mu_{L,H}=\Delta_L|_H\Delta_H^{-1}.
\end{equation}
For an algebraic group over $F$, $\Delta_M(h)$ is
$|\det\operatorname{Ad}(h)|$ on its Lie algebra.

\begin{proposition}[Frobenius descent {\cite[\S1.5]{Bernstein84}}]
\label{m1:prop:Frobenius-descent}\label{m2:tools:frobenius}
Let $p:X\to Y$ be a continuous equivariant map of $l$-spaces, where $L$
acts transitively on $Y\simeq L/L_{y_0}$ with its quotient topology. For
any character $\chi$ of $L$ there is a natural isomorphism
\begin{equation}\label{m2:tools:frobformula}
 \mathcal D(X)^{L,\chi}\simeq
 \mathcal D(p^{-1}(y_0))^{L_{y_0},\,\chi|_{L_{y_0}}\mu_{L,L_{y_0}}}.
\end{equation}
\end{proposition}
All homogeneous bases used below are rational algebraic orbits with this
topology. At $y_0$ the base density transforms by
$\Delta_{L_{y_0}}/\Delta_L=\mu_{L,L_{y_0}}^{-1}$, explaining the direction
of the fiber character. For example, projecting $F^2\setminus\{0\}$ to
$\mathbf P^1(F)$ gives the upper triangular stabilizer $P$ at $Fe_1$,
with $\Delta_P(\diag(\alpha,\beta))=|\alpha/\beta|$. Additive measure has
character $|\alpha\beta|^{-1}$, and its radial measure $|c|\,dc$ has
character $|\alpha|^{-2}=|\alpha\beta|^{-1}|\beta/\alpha|$.

\begin{lemma}[Homogeneous spaces and positive measures]
\label{m1:lem:orbital-measures}\label{m2:tools:orbit}\label{m1:cor:homogeneous-transpose}
For $X=L/H$, the space $\mathcal D(X)^{L,\chi}$ is zero unless
$\chi|_H\mu_{L,H}=1$, and is one-dimensional when this condition holds.
If $\chi$ is positive, the latter space has a positive measure generator.
Suppose an operator $\sigma$ preserves this space and the orbit, and its
square acts as the identity on the space. Then $\sigma$ fixes that
positive generator. In particular there is no distribution in the space
on which $\sigma$ acts by $-1$.
\end{lemma}
\begin{proof}
Frobenius descent reduces the assertion to a one-point fiber, whose
distribution space is one-dimensional with the stated character condition. When $\chi$ is positive, the
usual Haar-density construction gives a positive equivariant measure.
Pushforward by $\sigma$ acts on this ray by a positive scalar. Its square
is one, so that scalar is one.
\end{proof}
For $\sigma^2=h\in L$, the square condition is equivalent to $\chi(h)=1$.

\begin{proposition}[Descent with an involution]\label{m1:prop:odd-Frobenius}
Suppose $p:X\to Y$ is $\widetilde L$-equivariant and $L$ is transitive on
$Y$. Choose $g_0\tau(y_0)=y_0$, put $\tau_{y_0}=g_0\tau$, and let
$\widetilde L_{y_0}=L_{y_0}\sqcup L_{y_0}\tau_{y_0}$. Then
\begin{equation}\label{m1:eq:odd-Frobenius}
 \mathcal D(X)^{\widetilde L,\varepsilon}
 \simeq\mathcal D(p^{-1}(y_0))^{\widetilde L_{y_0},\,
 \varepsilon|_{\widetilde L_{y_0}}\widetilde\mu_{y_0}},\qquad
 \widetilde\mu_{y_0}=
 \Delta_{\widetilde L}|_{\widetilde L_{y_0}}
 \Delta_{\widetilde L_{y_0}}^{-1}.
\end{equation}
\end{proposition}
Proposition~\ref{m1:prop:Frobenius-descent}, applied to the extended
group, gives the claimed identification, with the relative modulus
evaluated on both components. An order-two
correcting operator has modulus one. In the non-unimodular applications below we check that the induced operator
preserves the positive ordinary fiber character, as required by
Lemma~\ref{m1:lem:orbital-measures}.

\subsection{Localization and opposite weights}\label{m1:subsec:localization}
\begin{proposition}[Bernstein localization {\cite[\S1.4]{Bernstein84}}]
\label{m1:prp:localization}\label{m2:tools:localization}
Let an $l$-group $L$ act on an $l$-space $X$, and let $p:X\to Y$ be a
continuous map to a Hausdorff $l$-space, constant on $L$-orbits. If
$\mathcal D(p^{-1}(y))^{L,\chi}=0$ for every $y\in Y$, then
$\mathcal D(X)^{L,\chi}=0$.
\end{proposition}
For the sign character, $p$ is required to be invariant under the extended group.

\begin{lemma}[Hyperbolas and the two-axis node]
\label{m1:lem:nonzero-xy-fiber}\label{m1:lem:Prasad-two-axis}\label{m1:subsec:two-axis}
For $c\in F$, let $V_c=\{(x,y):xy=c\}$, with
$t(x,y)=(tx,t^{-1}y)$ and $\sigma(x,y)=(y,x)$. Every
$F^\times$-invariant distribution on $V_c$ is $\sigma$-invariant.
\end{lemma}
\begin{proof}
For $c\neq0$, the equivariant parametrization
$a\mapsto(a,c/a)$ identifies $V_c$ with $F^\times$ and $\sigma$ with
$a\mapsto c/a$. Multiplicative Haar measure is invariant under translation
and inversion. For $c=0$, consider the closed fiber
\[
 Z=V_0=(F\times\{0\})\cup(\{0\}\times F).
\]
Let $T\in\mathcal D(Z)^{F^\times}$. Its restriction to the open
subset $Z\setminus\{0\}$ is
$c_1\nu_1+c_2\nu_2$, where $\nu_i$ is multiplicative Haar
measure on the corresponding punctured axis, normalized so that
the units have measure one. This follows from Frobenius descent
on each of its two open orbits, each with trivial stabilizer.
Let $\mathfrak o$ be the ring of integers and $\varpi$ a
uniformizer. The set $C=Z\cap\mathfrak o^2$ is compact open
in $Z$. For $t=\varpi$, the first-axis part of $tC$ is
$\varpi\mathfrak o$, and its second-axis part is
$\varpi^{-1}\mathfrak o$. Therefore
\[
 \mathbf1_C-\mathbf1_{tC}
 =\mathbf1_{\mathfrak o^\times\times\{0\}}
  -\mathbf1_{\{0\}\times\varpi^{-1}\mathfrak o^\times}.
\]
The difference vanishes near the origin. Invariance and the
normalization of the two Haar measures give
$0=T(\mathbf1_C-\mathbf1_{tC})=c_1-c_2$.
Thus $T-\sigma_*T$ is supported at the origin. By
\cref{m2:tools:support} it is a multiple of $\delta_0$; it is also
anti-invariant under $\sigma$, whereas $\delta_0$ is fixed.
It vanishes; see also \cite[Lemma~4.6]{Prasad90}.

\end{proof}

\begin{proposition}[Opposite-weight representations]
\label{m1:prop:opposite-weight-plane}\label{m2:elem:split}\label{m1:cor:opposite-weight-character}
For the same action and swap, every invariant distribution on $F^2$ is
fixed by $\sigma$. The assertion remains true when the action is pulled
back along a surjective character $\alpha:S\to F^\times$.
\end{proposition}
\begin{proof}
Proposition~\ref{m1:prp:localization}, applied to $(x,y)\mapsto xy$,
and Lemma~\ref{m1:lem:nonzero-xy-fiber} give the conclusion. Surjectivity makes $S$-invariance exactly
the same as $F^\times$-invariance.
\end{proof}

\begin{corollary}[Branch coefficients]\label{m1:cor:two-axis-coefficients}
Let $\nu_x,\nu_y$ be positive invariant measures on the punctured
branches of $V_0$, normalized by $\sigma_*\nu_x=\nu_y$. If an invariant distribution
on the node has restrictions $c_x\nu_x,c_y\nu_y$, then $c_x=c_y$.
\end{corollary}
\begin{proof}
The equality $\sigma D=D$ from Lemma~\ref{m1:lem:Prasad-two-axis},
restricted to the open punctured $y$-branch, gives $c_x=c_y$.
\end{proof}

\begin{lemma}[Fixed parameters and independent products]
\label{m1:lem:passive-products}\label{m2:tools:passive}\label{m1:cor:opposite-weight-passive}
Suppose every $S$-invariant distribution on $V$ is fixed by $\sigma$.
The same holds on $V\times W$ for the product actions
$s(v,w)=(sv,w)$ and $\sigma(v,w)=(\sigma v,w)$.
If $(\alpha_1,\ldots,\alpha_m):S\to(F^\times)^m$ is jointly
surjective, the independent opposite-weight action on a product of
planes, hyperbolas, or two-axis nodes fixes every invariant distribution
under each individual coordinate swap, and hence under their product.
\end{lemma}
\begin{proof}
Fixing a test function in $W$ gives an invariant distribution on $V$.
Compactly supported locally constant functions on a product are finite
sums of product test functions: partition a compact support into finitely
many compact open rectangles on which the function is constant.
This proves the first assertion. For the second, joint surjectivity
provides each factor's scaling group while fixing all other factors.
The result follows by applying the first assertion successively to
each factor, together with Lemma~\ref{m1:lem:nonzero-xy-fiber} or
Proposition~\ref{m1:prop:opposite-weight-plane} as appropriate.
\end{proof}

\subsection{Quadratic norm fibers and matrix pairs}
\begin{lemma}[Quadratic norm fibers]\label{m2:elem:norm}
Let $L/F$ be a quadratic field extension, with nontrivial
automorphism $w\mapsto\bar w$, and let
$L^1=\{q\in L^\times:N_{L/F}(q)=1\}$. For the multiplication
action of $L^1$ on $L$, every invariant distribution is fixed by
conjugation.
\end{lemma}

\begin{proof}
The $L^1$-orbits lie in fibers of the norm map $N_{L/F}:L\to F$.
A nonempty nonzero fiber
is $w_0L^1$ and is simply transitive under $L^1$. In the orbit
coordinate $w_0q$, conjugation sends
\[
 q\longmapsto(\bar w_0/w_0)q^{-1}.
\]
The multiplier lies in $L^1$, and both inversion and translation
preserve Haar measure on this abelian group. The stabilizer is
trivial, so its modulus and that of $L^1$ are both one. The zero
fiber is $\{0\}$ because $L$ is a field; its Dirac
distribution is fixed. All other fibers are empty, so
\cref{m2:tools:localization} completes the proof.
\end{proof}

\begin{lemma}[Independent norm coordinates]\label{m2:elem:products}
Let $E/F$ be a quadratic field extension, and let $(E^1)^r$ act
coordinatewise on $E^r$. Let the odd
operator be coordinatewise conjugation. Then its sign
distribution space is zero. The same conclusion holds after
pulling back the action through a surjective homomorphism
$S(F)\to(E^1)^r$, and after adjoining any passive parameter
space fixed by the extended group.
\end{lemma}

\begin{proof}
By \cref{m2:tools:passive}, the other coordinates may be treated as
fixed parameters when applying \cref{m2:elem:norm} to each coordinate. Each
coordinate conjugation fixes every invariant distribution, so their
product does also. Pullback through a surjective homomorphism gives
the same transformations and invariant distributions. The last
assertion follows again from \cref{m2:tools:passive}.
\end{proof}

\begin{lemma}[Simultaneous conjugation of matrix pairs]
\label{m1:lem:matrix-pair-transpose}
Every distribution on $M_2(F)^2$ invariant under simultaneous
$\GL_2(F)$-conjugation is fixed by componentwise transpose.
\end{lemma}
\begin{proof}
For $2\times2$ matrices put
$\bar P=(\tr P)I-P$ and
\[
 J_0=\begin{pmatrix}0&1\\-1&0\end{pmatrix},\qquad
 \bar P=J_0P^tJ_0^{-1},\qquad
 \jmath(P,Q)=(\bar P,\bar Q).
\]
For $s,t\in F$ the invariant open sets
\begin{equation}\label{m1:eq:scalar-translation-cover}
 U_{s,t}=\{(P,Q):P+sI,\ Q+tI\in\GL_2(F)\}
\end{equation}
cover $M_2(F)^2$, since each determinant excludes at most two scalar
shifts and $F$ is infinite. Translation
$T_{s,t}(P,Q)=(P+sI,Q+tI)$ identifies $U_{s,t}$ equivariantly with
$\GL_2(F)^2$ and commutes with $\jmath$, because
$\overline{P+sI}=\bar P+sI$. These opens are also $\jmath$-stable.
Prasad's simultaneous-conjugation theorem
\cite[Proposition~4.5]{Prasad90} applies on the invertible pair space and
gives $(\jmath D-D)|_{U_{s,t}}=0$ for every invariant distribution $D$.
The open cover gives $\jmath D=D$ on all matrix pairs. Since $\jmath$
is componentwise transpose followed by simultaneous conjugation by
$J_0$, this is the required transpose invariance.
\end{proof}

\subsection{Fourier support and homogeneity}
Fix a nontrivial continuous unitary additive character $\psi:F\to\mathbb C^\times$.
For a nondegenerate symmetric bilinear form $B_V$ on a finite
dimensional $F$-space $V$, use its self-dual Haar measure to define
\[
 (\mathcal F_V f)(v)=\int_V f(w)\psi(B_V(v,w))\,dw.
\]
The same formula in the $V$ variable defines partial Fourier
transform on $V\times M$. On distributions, set
$(\mathcal F_VT)(f)=T(\mathcal F_Vf)$.

\begin{theorem}[Odd-dimensional homogeneity]\label{m2:tools:homogeneity}
Let $V$ have odd dimension and nondegenerate symmetric form $B_V$,
and let $M$ be an analytic manifold over $F$. Put
$\mathcal N_V=\{v:B_V(v,v)=0\}$. A distribution on $V\times M$
whose support and partial Fourier support are both contained in
$\mathcal N_V\times M$ is zero.
\end{theorem}

This is the odd-dimensional conclusion of
\cite[Theorem~5.1.5]{AGDescent}, including its analytic parameter $M$.
In the non-Archimedean category its Schwartz distributions are the
distributions used here. In both identity descendants, $\dim_FV=3$.

\section{Determinant-open quiver geometry and normal representations}\label{sec:quiver}
The same determinant-open quiver governs the linear pair
and the scalar extension of its unitary counterpart. We first establish
its closed-orbit geometry and normal spaces over $F$. The
unitary application will retain its rational descent data and the
additional multiplier direction.
The algebraic statements and constructions of this section hold over
every field of characteristic zero; in Section~\ref{sec:unitary-geometry}
we also apply them over $E$ and $\Fb$.

\subsection{The determinant open and its rational closed orbits}\label{m1:sec:geometry}
\subsubsection{Determinant-relative stability}\label{m1:subsec:determinant-stability}
Let $\qvr$ be the quiver with vertices $U,Q,P$ and arrows
$A:U\to P$, $B:Q\to P$, $C:U\to Q$, and $D:Q\to Q$. Dimension vectors
are ordered as $(\dim U,\dim Q,\dim P)$. Its representation space in
dimension $(2,2,2)$ is
\[
 V=\Hom(U,P)\oplus\Hom(Q,P)\oplus\Hom(U,Q)\oplus\End(Q),
\]
with base-change group $L=\GL(P)\times\GL(Q)\times\GL(U)$.
For $x=(A,B,C,D)$ its action, the transpose involution, and the
induced automorphism of $L$ are
\begin{align}
 (p,q,u)x&=(pAu^{-1},pBq^{-1},qCu^{-1},qDq^{-1}),
 \label{m1:eq:residual-action}\\*
 \tau x&=(A^t,C^t,B^t,D^t),\qquad
 \alpha(p,q,u)=(u^{-t},q^{-t},p^{-t}).\label{m1:eq:tau}
\end{align}
 For $E=(U,Q,P;A,B,C,D)$ put
\begin{equation}\label{m1:eq:block-map}
 \Phi_E:U\oplus Q\longrightarrow P\oplus Q,
 \qquad\Phi_E(u,q)=(Au+Bq,Cu+Dq).
\end{equation}
The residual variety is the principal open
\begin{equation}\label{m1:eq:X}
 \cX=\{E\in V:f(E)\neq0\},\qquad f(E)=\det\Phi_E.
\end{equation}
The determinant is a semi-invariant:
\begin{equation}\label{m1:eq:det-character}
 f((p,q,u)E)=\chi_{\det}(p,q,u)f(E),
 \qquad\chi_{\det}(p,q,u)=\det p/\det u.
\end{equation}
For a field $k$ and a character $\gamma$ with values in $k^\times$,
write $k_\gamma$ for the one-dimensional representation $g\cdot z=\gamma(g)z$,
and write $k_{\triv}$ for the trivial representation.

A subrepresentation $N=(U',Q',P')$ satisfies
$AU'\subset P'$, $BQ'\subset P'$, $CU'\subset Q'$, and $DQ'\subset Q'$.
Define its balance by
\begin{equation}\label{m1:eq:theta}
 \theta(N)=\dim P'-\dim U'.
\end{equation}
We call $N$ balanced when $\theta(N)=0$.

\begin{lemma}\label{m1:lem:theta-nonnegative}
If $E\in\cX$, then $\theta(E)=0$ and $\theta(N)\geq0$ for every
subrepresentation $N$. Every balanced subrepresentation and its quotient
have invertible induced block maps.
\end{lemma}
\begin{proof}
The restriction of $\Phi_E$ is an injection
$U'\oplus Q'\to P'\oplus Q'$, so $\dim U'\leq\dim P'$. If equality
holds, this restriction is an isomorphism, and so is the induced quotient
map. The argument applies successively to balanced filtration factors.
\end{proof}

Let $\cAZero$ be the full subcategory of finite-dimensional quiver
representations with $\theta(E)=0$ and $\theta(N)\geq0$ for all
subrepresentations. A nonzero object is stable if the latter inequality
is strict for every nonzero proper subrepresentation, and polystable if
it is a direct sum of stable objects. This is the zero-slope semistable
category, with slope $\theta(E)/(\dim U+\dim Q+\dim P)$ for nonzero $E$.

\begin{lemma}\label{m1:lem:A0-abelian}
The category $\cAZero$ is abelian of finite length. Its simple objects
are precisely its stable objects.
\end{lemma}
\begin{proof}
For a morphism $h:E\to E'$ in $\cAZero$, its kernel and image have
nonnegative balance. Additivity and $\theta(E)=0$ force both balances
to be zero; the cokernel is then balanced as well. Nonnegativity passes
to subobjects of the kernel and image, and to the cokernel by lifting
subobjects to $E'$. These are therefore the kernel, image, and cokernel
in the full subcategory. Direct sums satisfy the same condition by
intersecting a subobject with one summand and projecting to the other.
Every balanced subrepresentation of an object of $\cAZero$ belongs to
$\cAZero$, so simplicity is exactly stability. Total vector-space
dimension strictly decreases in a proper subobject or quotient, proving
finite length.
\end{proof}

Over an algebraically closed field, these definitions agree with
character-GIT semistability and stability for $\chi_{\det}$ by
\cite[Proposition~3.1]{King94}. In the semistable locus, closed orbits
correspond to polystable objects by \cite[Proposition~3.2]{King94}.
The semistable locus also contains points with singular block map,
such as an $A$-only stable factor together with a $Q$-supported
zero-loop factor. We prove the closed-orbit criterion on $\cX$ directly.

\begin{proposition}[Polystability on the determinant open]
\label{m1:prop:polystability}
For $E\in\cX$, its geometric $L$-orbit is closed in $\cX$ if and only if
$E_{\Fb}$ is semisimple in $\cAZero$, equivalently polystable.
\end{proposition}
\begin{proof}
Work over $\Fb$. The map $v\mapsto(v,f(v)^{-1})$ identifies $\cX$
equivariantly with the closed invariant graph
\begin{equation}\label{m1:eq:closed-graph}
 \widehat\cX=\{(v,z)\in V\oplus\Fb_{\chi_{\det}^{-1}}:zf(v)=1\}.
\end{equation}
We apply \cite[Proposition~2.6(i)]{King94}, with the trivial character,
to this ambient linear representation of the reductive group $L$.
For that character the semistable locus is the entire representation
space, and its pairing with every one-parameter subgroup is zero.
Thus an orbit is closed if and only if every existing one-parameter
limit stays in that orbit. Since the graph is closed, any such limit
from the graph remains in the graph.

If $E$ is not semisimple, choose a Jordan--H\"older filtration
$0=E_0\subset E_1\subset\cdots\subset E_r=E$ in $\cAZero$, with
$S_i=E_i/E_{i-1}$. Split the vertex filtrations as vector spaces, and
give the $i$th graded part at every vertex a weight $w_i$, where
$w_1>\cdots>w_r$. For any arrow, its block from the $j$th source part
to the $i$th target part is zero unless $i\leq j$, and is multiplied by
$t^{w_i-w_j}$. The one-parameter limit is therefore $\bigoplus_i S_i$.
Its determinant weight is
$\sum_i w_i\theta(S_i)=0$, so the graph coordinate stays fixed and the
limit lies in $\cX$. It is not isomorphic to $E$, proving nonclosedness.

Conversely suppose $E$ is semisimple in $\cAZero$, and a one-parameter
subgroup has a limit in \eqref{m1:eq:closed-graph}. Its weights decompose
every vertex space; let $F^cE$ consist of the weights at least $c$.
An arrow block from weight $b$ to weight $a$ is multiplied by $t^{a-b}$;
existence of the limit forces it to vanish for $a<b$. Thus every $F^cE$
is a subrepresentation. If the distinct weights are
$m_1<\cdots<m_s$, put $\delta_i=\dim P_{m_i}-\dim U_{m_i}$. Then
$\sum_i\delta_i=0$ and $\theta(F^{m_i}E)=\sum_{j\geq i}\delta_j\geq0$.
The graph equation forces the determinant weight to be zero: both
$f(\lambda(t)E)=t^k f(E)$ and its inverse have finite nonzero limits,
so $k=0$. Summation by parts gives
\begin{equation}\label{m1:eq:balanced-weight-identity}
 0=\sum_i m_i\delta_i
   =\sum_{i=2}^s(m_i-m_{i-1})\theta(F^{m_i}E).
\end{equation}
Each term is nonnegative, hence every threshold subobject is balanced
and belongs to $\cAZero$. This filtration splits because $E$ is
semisimple there. Its associated graded, which is the one-parameter
limit, is therefore isomorphic to $E$. The affine criterion proves
closedness.
\end{proof}
For $E$ defined over $F$, geometric closedness is
equivalent to analytic closedness of $L(F)E$ by
\cite[Theorem~2.3.8]{AGDescent}; all rational orbit identifications below
are proved over $F$.

\subsubsection{Stable factors and arrow ranks}\label{m1:subsec:stable-factors}
Duality interchanging $U$ and $P$ preserves stability. Indeed a proper
nonzero subobject of the dual is dual to a quotient $E/N$, and its
balance is $-\theta(E/N)=\theta(N)>0$ for stable $E$.

\begin{lemma}\label{m1:lem:stable-factors}
Over $\Fb$, a stable balanced factor of a point of $\cX$, of dimension
$(r,s,r)$, has $s=1$ if $r=0$, and $r=1$ if $s=0$. If $r,s>0$, its
$C$ is injective and its $B$ is surjective. The only proper stable
dimension vectors in ambient dimension $(2,2,2)$ are
\[
 (1,0,1),\quad(0,1,0),\quad(1,1,1),\quad(1,2,1).
\]
\end{lemma}
\begin{proof}
For $r=0$ only the loop remains, and a loop on a space of dimension
greater than one over $\Fb$ has a proper invariant subspace. For $s=0$,
the block map is the isomorphism $A:U\to P$; a proper subspace $U'$
would give a balanced subobject $(U',0,AU')$. These prove the first two
assertions. If $r,s>0$ and $U_0=\ker C\neq0$, injectivity of $\Phi_E$
makes $A|_{U_0}$ injective. Then $(U_0,0,AU_0)$ is a proper nonzero
balanced subobject, contradicting stability. Duality gives surjectivity
of $B$. Hence $r\leq s$, which gives the list under the ambient
dimension bounds. The same argument applies to a full-dimensional
stable object.
\end{proof}

Denote the first two factors by $\Efac$ and $\Qfac_e$. On $\Efac$ we
normalize $A=1$, and on $\Qfac_e$ the loop is $e\neq0$. A stable
$(1,1,1)$-factor is normalized as
\begin{equation}\label{m1:eq:R-factor}
 \Rfac_{x,d}:\qquad \Phi=\begin{pmatrix}x&1\\1&d\end{pmatrix},
 \qquad xd\neq1.
\end{equation}
It is stable because $B$ and $C$ are nonzero, leaving no proper balanced
subobject. We denote a stable $(1,2,1)$-factor by $\Hfac$.

\begin{corollary}\label{m1:prop:rank-exclusion}
Every closed orbit in $\cX$ has $\rk B=\rk C$. Its rank type is
$(2,2)$, $(1,1)$, or $(0,0)$.
\end{corollary}
\begin{proof}
On each stable factor the two ranks agree: they are zero on $\Efac$
and $\Qfac_e$, one on $\Rfac$ and $\Hfac$, and two on a full stable
factor. Ranks add under the polystable decomposition of
Proposition~\ref{m1:prop:polystability}. This excludes all six unequal-rank
types, including those with only one nonzero arrow.
\end{proof}
Geometrically the rank-zero factors are $\Efac^{\oplus2}$ and a
semisimple two-dimensional loop object. Rank one has either
$\Efac\oplus\Rfac\oplus\Qfac$ or $\Efac\oplus\Hfac$. Rank two is
stable of full dimension or a sum of two $(1,1,1)$-factors, possibly
isomorphic. This is a geometric list; a nonsplit rational simple may
split into geometric factors and need not have one of the proper
dimension vectors above.

\subsubsection{Dense pairs and rational semisimplicity}\label{m1:subsec:dense-geometry}
The regular $L$-invariant function
$\rho(E)=\det B\det C/\det\Phi_E$ shows that an orbit in the dense
rank-$(2,2)$ locus cannot degenerate to a lower rank in $\cX$: its
closure lies in the nonzero fiber $\rho^{-1}(\rho(E))$.
Taking $p=qB^{-1}$ and $u=qC$ gives the gauge $B=C=I$ with
\begin{equation}\label{m1:eq:dense-gauge}
 X=B^{-1}AC^{-1},\qquad Y=D,
 \qquad (X,Y)\mapsto(qXq^{-1},qYq^{-1}).
\end{equation}
The condition is $\Delta(X,Y)=\det\left(\begin{smallmatrix}X&I\\I&Y\end{smallmatrix}\right)\neq0$.
As $\Delta$ is conjugation-invariant, a pair orbit's closure in $M_2^2$
also remains in this open.

\begin{proposition}\label{m1:prop:dense-closed}
A dense rational orbit is closed exactly when
$M=(F^2;X,Y)$ is semisimple as an $F\langle x,y\rangle$-module.
Its four rational types and stabilizers are: an absolutely irreducible
module, with stabilizer $F^\times$; two distinct rational characters,
with stabilizer $(F^\times)^2$; a rational irreducible module splitting
over a quadratic field $K/F$, with stabilizer $K^\times$; and two equal
characters, with stabilizer $\GL_2(F)$.
\end{proposition}
\begin{proof}
In the gauge, a balanced subobject has $\dim U'=\dim P'$ and
$U'\subset Q'\subset P'$, because $B=C=I$. Thus all three spaces agree.
Such subobjects are exactly the common $X,Y$-invariant subspaces.
This identifies geometric polystability with semisimplicity of the
free-algebra module by Proposition~\ref{m1:prop:polystability}.

The rational comparison requires both directions. If $M$ is semisimple
over $F$, let $A_M\subset\End_F(M)$ be the finite-dimensional image of
$F\langle x,y\rangle$. Its radical annihilates the faithful semisimple
module $M$, hence is zero. In characteristic zero, a finite-dimensional
semisimple algebra is separable: its centers are finite separable
extensions and its central simple factors split over $\Fb$. Therefore
$A_M\otimes_F\Fb$ and $M\otimes_F\Fb$ are semisimple. Conversely, for
an $F$-defined submodule $W\subset M$, a retraction over $\Fb$ satisfies
linear equations over $F$: it commutes with $X,Y$ and is the identity
on $W$. Row reduction gives a solution over $F$. Thus geometric
semisimplicity implies rational semisimplicity.

The stabilizer is $\End_{A_M}(M)^\times$. Absolute irreducibility gives
$\End_{A_M}(M)=F$. Two distinct rational characters give $F\oplus F$,
and two equal characters give $M_2(F)$. In the remaining irreducible
case, $D_0=\End_{A_M}(M)$ is a division algebra, and $\dim_FD_0$ divides
$\dim_FM=2$. If $D_0=F$, the finite-dimensional density theorem gives
$A_M=M_2(F)$, contrary to failure of absolute irreducibility. Hence
$\dim_FD_0=2$; any nonscalar element generates it, so $D_0$ is a
quadratic field $K$. Identifying $M$ with $K$, the pair consists of
multiplication operators $m_x,m_y$, at least one nonscalar. These exhaust
all semisimple modules of dimension two.
\end{proof}

\subsubsection{Rank zero and rank one over \texorpdfstring{$F$}{F}}
\label{m1:subsec:rational-boundary}
If $B=C=0$, invertibility gives $A,D\in\GL_2(F)$ and base change
normalizes $A=I$. The object is $\Efac^{\oplus2}\oplus\Qfac_D$, where
$\Qfac_D$ is supported on the loop vertex. Polystability says exactly
that $D$ is semisimple. Rational canonical form yields the following.

\begin{proposition}\label{m1:prop:rankzero}\label{m1:prop:rank00-classification}
Every closed rank-zero orbit has a representative
$x_D=(I,0,0,D)$ with $D\in\GL_2(F)$ semisimple. Its rational type is
scalar $D=dI$ with $d\neq0$, split regular with two distinct nonzero
eigenvalues, or nonsplit regular with irreducible quadratic characteristic
polynomial.
\end{proposition}

For rank one, $U_0=\ker C$ and $P_0=AU_0$ are rational lines, and
$\Efac=(U_0,0,P_0)$ is a balanced subrepresentation.

\begin{lemma}\label{m1:lem:rank11-split}
If $E$ has closed orbit and rank type $(1,1)$, there is an $F$-rational
splitting $E\simeq\Efac\oplus\mathcal H$, with
$\dim\mathcal H=(1,2,1)$.
\end{lemma}
\begin{proof}
The inclusion of $\Efac$ splits over $\Fb$ by polystability. A retraction
is a triple of linear maps satisfying the quiver morphism equations and
restricting to the identity on $\Efac$. All these equations are linear
over $F$, so solvability over $\Fb$ implies solvability over $F$. Its
kernel is the required complement, geometrically semisimple in
$\cAZero$.
\end{proof}
Write the complement's arrows as $c:F\to Q$, $b:Q\to F$, $D:Q\to Q$,
and its remaining scalar arrow as $a$. Both $b$ and $c$ are nonzero.
Put $L_c=\operatorname{im}c$ and $K_b=\ker b$.

\begin{lemma}\label{m1:lem:rank11-balanced}
The only proper nonzero balanced subobjects of $\mathcal H$ are
$(0,K_b,0)$, present exactly when $K_b$ is $D$-stable, and
$(F,L_c,F)$, present exactly when $L_c$ is $D$-stable. Thus $\mathcal H$
is geometrically stable exactly when neither line is $D$-stable.
\end{lemma}
\begin{proof}
The two end dimensions of a balanced subobject are either both zero or
both one. In the first case its middle space is a nonzero subspace of
$K_b$, hence is $K_b$. In the second it must contain $L_c$, and
properness forces it to be $L_c$. Each is a subobject precisely when
its middle line is $D$-stable. This argument holds over every scalar
extension, and both candidate lines are defined over $F$.
\end{proof}

\begin{proposition}[Rank-one dichotomy]\label{m1:prop:rank11}
A rank-one orbit is closed if and only if it admits the rational
splitting above and either neither of $L_c,K_b$ is $D$-stable, or both
are $D$-stable and distinct. The first case is
$E=\Efac\oplus\Hfac$, with $\Hfac$ stable. The second is
$E=\Efac\oplus\Rfac\oplus\Qfac$.
\end{proposition}
\begin{proof}
The splitting is necessary by Lemma~\ref{m1:lem:rank11-split}. If neither
line is stable, Lemma~\ref{m1:lem:rank11-balanced} proves stability of the
complement over $\Fb$. If both are stable and distinct, then
$Q=L_c\oplus K_b$ is a rational $D$-stable decomposition, with $b$
nonzero on $L_c$. Its two summands have stable dimension vectors
$(1,1,1)$ and $(0,1,0)$, proving polystability.
If exactly one line is stable, or both are stable and coincide, the
balanced subobject supplied by the lemma has no complement: a complement
to $(0,K_b,0)$ must use $L_c$, and a complement to $(F,L_c,F)$ must use
$K_b$. Either splitting forces the lines to be stable and distinct.
Thus these remaining cases are not semisimple, and their orbits are
not closed by Proposition~\ref{m1:prop:polystability}.
\end{proof}

The decomposable family has the representatives
\begin{equation}\label{m1:eq:decomposable-rep}
 A=\diag(x,1),\quad B=C=E_{11},\quad D=\diag(d,e),
 \qquad e\neq0,\quad xd\neq1.
\end{equation}
These are exactly the invertibility conditions on
$\Efac\oplus\Rfac_{x,d}\oplus\Qfac_e$; they allow $x=0$ and $d=e$.

For the stable-system family choose the rational cyclic basis $(c,Dc)$
of $Q$, possible because $L_c$ is not $D$-stable. Then $c=e_1$ and
$D=\left(\begin{smallmatrix}0&-d\\1&s\end{smallmatrix}\right)$.
Adjoining $\Efac$ gives
\begin{equation}\label{m1:eq:stable-rep-general}
 A=\diag(a,1),\quad
 B=\begin{pmatrix}b_0&b_1\\0&0\end{pmatrix},\quad C=E_{11},\quad
 D=\begin{pmatrix}0&-d\\1&s\end{pmatrix}.
\end{equation}
The determinant on the $(1,2,1)$-factor is
$\Delta=ad-sb_0+b_1\neq0$. Scaling the first $P$-coordinate by
$\Delta^{-1}$ scales $a,b_0,b_1,\Delta$ together and fixes $C,D$ and
the complementary identity block. Relabeling the scaled parameters gives
\begin{equation}\label{m1:eq:b1-normalized}
 b_1=1-ad+sb_0.
\end{equation}
Writing $b=(b_0,b_1)$, set
\begin{equation}\label{m1:eq:H}
 H=\begin{pmatrix}b_0&b_1\\b_1&sb_1-db_0\end{pmatrix}
   =\begin{pmatrix}b\\bD\end{pmatrix},\qquad\det H\neq0.
\end{equation}
The rows $b,bD$ are dependent exactly when $\ker b$ is $D$-stable:
dependence means $bD$ vanishes on $\ker b$, and conversely any functional
vanishing there is a multiple of the nonzero $b$. Thus \eqref{m1:eq:H} is
precisely the remaining stability condition. These normalizations also apply to nilpotent and repeated-root values of $D$.

\begin{theorem}[Rational closed-orbit classification]\label{m1:thm:closed-orbits}
Every closed $L(F)$-orbit in $\cX(F)$ is represented by exactly one of
the following families of types: a dense semisimple pair from
Proposition~\ref{m1:prop:dense-closed}; a rank-zero point from
Proposition~\ref{m1:prop:rankzero}; a decomposable representative
\eqref{m1:eq:decomposable-rep}; or a stable-system representative
\eqref{m1:eq:stable-rep-general} satisfying
\eqref{m1:eq:b1-normalized}--\eqref{m1:eq:H}. Each representative above has a closed orbit. The nine families are disjoint; within each family,
the parameters may represent the same orbit.
\end{theorem}
\begin{proof}
Corollary~\ref{m1:prop:rank-exclusion} gives the three ranks, and the dense,
rank-zero, and rational rank-one classifications exhaust them. Their
proofs establish both closedness and the exhaustion of the rational orbits.
\end{proof}

\subsection{Normal representations and rational transpose}\label{m1:sec:normals}
The classification supplies the objects at which descent is applied.
We now identify their stabilizers, normal representations, and
corrected transpose operators over $F$.

\subsubsection{The Hom--Ext complex}\label{m1:subsec:hom-ext}
For representations $S,T$ of $\qvr$ over $F$, let
\[
 \begin{gathered}
 C^0(S,T)=\bigoplus_v\Hom_F(S_v,T_v),\qquad
 C^1(S,T)=\bigoplus_{a:v\to w}\Hom_F(S_v,T_w),\\
 d(h)_a=h_wS_a-T_ah_v.
 \end{gathered}
\]

\begin{lemma}\label{m1:lem:hom-ext}
There is a natural exact sequence
\begin{equation}\label{m1:eq:hom-ext-complex}
 0\longrightarrow\Hom_{\qvr}(S,T)\longrightarrow C^0(S,T)
 \xrightarrow{d}C^1(S,T)\longrightarrow\Ext^1_{\qvr}(S,T)
 \longrightarrow0.
\end{equation}
For $E\in\cX$ this identifies, equivariantly under $L_E=\Aut_{\qvr}(E)$,
\begin{equation}\label{m1:eq:full-normal-ext}
 N_E=T_E\cX/T_E(LE)\simeq\Ext^1_{\qvr}(E,E),
 \qquad\dim_F N_E=4+\dim_F\End_{\qvr}(E).
\end{equation}
\end{lemma}
\begin{proof}
The kernel consists exactly of quiver morphisms. Split an extension
$0\to T\to Z\to S\to0$ at each vertex as vector spaces. Its arrow
matrices are $\left(\begin{smallmatrix}T_a&c_a\\0&S_a\end{smallmatrix}\right)$.
There are no relations on the arrows, including the loop, so $c_a$ is
arbitrary. Changing the vertex splittings changes $(c_a)$ by an element
of $\operatorname{im}d$, up to sign. Conversely this change gives an
equivalent extension. Addition and scalar multiplication of the upper
right blocks give the corresponding operations on extension classes,
proving the cokernel assertion.

Since $\cX$ is open, $T_E\cX=C^1(E,E)$. Its orbit differential is
\begin{equation}\label{m1:eq:orbit-differential}
 d_E(P_0,Q_0,U_0)
 =(P_0A-AU_0,\ P_0B-BQ_0,\ Q_0C-CU_0,\ Q_0D-DQ_0).
\end{equation}
Its kernel is $\End_{\qvr}(E)$ and its image is $T_E(LE)$.
Taking the quotient proves \eqref{m1:eq:full-normal-ext}, including the
dimension formula $16-12+\dim\End_{\qvr}(E)$. All constructions commute
with automorphisms of $E$.
\end{proof}
For dimension vectors $\mathbf d=(u,q,p)$ and
$\mathbf e=(u',q',p')$, the Euler form from
\eqref{m1:eq:hom-ext-complex} is
\[
 \langle\mathbf d,\mathbf e\rangle
 =uu'+qq'+pp'-up'-qp'-uq'-qq'.
\]
On balanced vectors it is symmetric:
\begin{equation}\label{m1:eq:balanced-euler}
 \langle(r,s,r),(r',s',r')\rangle=rr'-sr'-rs',\qquad
 \dim\Ext^1(S,T)=\dim\Hom(S,T)-\langle\dim S,\dim T\rangle.
\end{equation}
Distinct stable objects have no Hom between them, since they are
nonisomorphic simples in the full subcategory $\cAZero$. For an
absolutely stable object defined over $F$, its endomorphism algebra is
$F$, by Schur's lemma after scalar extension and descent of the linear
morphism equations. Thus for pairwise nonisomorphic absolutely stable
$S_i$ defined over $F$ and finite-dimensional $F$-spaces $M_i$,
\begin{equation}\label{m1:eq:multiplicity-spaces}
 \begin{split}
 E=\bigoplus_i S_i\otimes_F M_i\quad&\Longrightarrow\quad
 \Aut(E)=\prod_i\GL(M_i),\\
 \Ext^1(E,E)&=\bigoplus_{i,j}
 \Hom_F(M_i,M_j)\otimes_F\Ext^1(S_i,S_j).
 \end{split}
\end{equation}
These identities follow by decomposing the complex into its $(i,j)$
blocks. On the factor $\Hom_F(M_i,M_j)$ in \eqref{m1:eq:multiplicity-spaces}, the action is $f\mapsto g_jfg_i^{-1}$.
The nonsplit rational cases below are computed in their quadratic-field models.

\subsubsection{Normal representations for dense orbits}\label{m1:subsec:dense-normals}
At the gauge $B=C=I$ the differential of \eqref{m1:eq:dense-gauge} is
\[
 (\dot A,\dot B,\dot C,\dot D)\longmapsto
 (\dot A-\dot B X-X\dot C,\dot D).
\]
It is onto $M_2^2$. Its kernel consists of orbit directions: use
$Q_0=0$, $P_0=\dot B$, $U_0=-\dot C$ in
\eqref{m1:eq:orbit-differential}. The image of a general orbit direction
is $([Q_0,X],[Q_0,Y])$. Hence
\begin{equation}\label{m1:eq:dense-normal}
 N_{X,Y}\simeq\Ext^1_{F\langle x,y\rangle}(M,M)
 =\frac{M_2(F)\oplus M_2(F)}
 {\{([Z,X],[Z,Y]):Z\in M_2(F)\}}.
\end{equation}
The Ext group is taken over the free algebra, so both generators may
vary. The two-loop Euler form is
$\langle m,n\rangle=-mn$.

\begin{lemma}[Dense rational transporters]\label{m1:lem:dense-transpose}
Every closed dense pair has an $F$-rational symmetric invertible matrix
$H_0$ with $H_0X=X^tH_0$ and $H_0Y=Y^tH_0$. The transporter in $L$ is
$(H_0^{-1},H_0^{-1},H_0^{-1})$, and the corrected operator on
\eqref{m1:eq:dense-normal} is induced by
\begin{equation}\label{m1:eq:dense-transpose}
 (P_1,Q_1)\longmapsto(H_0^{-1}P_1^tH_0,H_0^{-1}Q_1^tH_0).
\end{equation}
\end{lemma}
\begin{proof}
For an absolutely irreducible pair, the space of symmetric matrices has
dimension three over $F$. Each of $H_0X-X^tH_0$ and $H_0Y-Y^tH_0$ is
skew-symmetric, imposing at most one linear condition. There is a
nonzero rational symmetric solution. If it were singular, its kernel
would be a common invariant line, contradicting irreducibility. For
split and scalar pairs choose diagonal representatives and $H_0=I$.
For a nonsplit pair on $K$, take the Gram matrix of
$\operatorname{Tr}_{K/F}(uv)$; multiplication operators are self-adjoint
for this nondegenerate symmetric form. To return from a normal-form pair
$(X_0,Y_0)$ to $(X,Y)=g(X_0,Y_0)g^{-1}$, replace its intertwiner
$H_{\mathrm{rep}}$ by $H_0=g^{-t}H_{\mathrm{rep}}g^{-1}$.
This matrix is symmetric and invertible and satisfies
$H_0X=X^tH_0$ and $H_0Y=Y^tH_0$.
The transporter $(H_0^{-1},H_0^{-1},H_0^{-1})$ carries the transposed point back to the original one. Differentiation then gives \eqref{m1:eq:dense-transpose}.
Symmetry makes the corrected operator square to one.
\end{proof}

\begin{proposition}[Normal space for an absolutely irreducible pair]\label{m1:prop:abs-irred-normal}
For an absolutely irreducible pair, the differentials of
\[
 f_1=\tr X,\quad f_2=\det X,\quad f_3=\tr Y,\quad
 f_4=\det Y,\quad f_5=\tr(XY)
\]
identify the normal space with $F^5$. Both the scalar stabilizer and the
corrected transpose act trivially on it.
\end{proposition}
\begin{proof}
For a variation $(P_1,Q_1)$, the trace and product differentials are
$\tr P_1$, $\tr Q_1$, and $\tr(P_1Y+XQ_1)$. Since
$\operatorname{adj}(X)=(\tr X)I-X$, the first two differentials span
$\tr P_1$ and $\tr(XP_1)$, and likewise the third and fourth span
$\tr Q_1$ and $\tr(YQ_1)$. A linear relation among all five would give
\[
 a\tr P_1+b\tr(XP_1)+c\tr Q_1+d\tr(YQ_1)
 +e\tr(P_1Y+XQ_1)=0.
\]
Nondegeneracy of the trace pairing, after separating the variables,
implies $aI+bX+eY=0$ and $cI+dY+eX=0$. If $e\neq0$, then
$Y\in F[I,X]$, so an eigenline of $X$ over $\Fb$ is also $Y$-stable,
contradicting absolute irreducibility. Thus $e=0$. Neither matrix can
be scalar, since the other has an eigenline over $\Fb$; hence
$a=b=c=d=0$. The differentials are independent.

They vanish on the orbit tangent, whose dimension is $4-1=3$. Their
common kernel has dimension $8-5=3$, and is therefore exactly that
tangent. This agrees with the self-Ext dimension $1+2^2=5$ and gives
explicit coordinates on the normal space. The map
\eqref{m1:eq:dense-transpose} fixes the base point and all five functions.
Differentiation shows that it fixes every normal coordinate. The
stabilizer is central and acts trivially.
\end{proof}

For the split type write
$X=\diag(x_1,x_2)$, $Y=\diag(y_1,y_2)$, with
$(x_1,y_1)\neq(x_2,y_2)$. The two character factors have self-Ext
dimension two and cross-Ext dimension one in each direction. Thus, with
$\alpha(t_1,t_2)=t_1/t_2$,
\begin{equation}\label{m1:eq:dense-split-normal}
 S=(F^\times)^2,\qquad N_{X,Y}=F^4_{\triv}\oplus F_\alpha
 \oplus F_{\alpha^{-1}}.
\end{equation}
The four diagonal parameter variations give the self-Ext coordinates and
are fixed by transpose. On the off-diagonal quotient transpose
interchanges the two one-dimensional spaces. Its square is one, so a
basis in one and its image in the other make this a literal swap over
$F$. The character is rationally surjective, since $\alpha(t,1)=t$.

For the nonsplit type, identify $M=K$, write $m_z$ for multiplication
by $z\in K$, and write $z\mapsto\bar z$ for the nontrivial automorphism
of $K/F$. Then $X=m_x,Y=m_y$, and with $\iota(w)=\bar w$ there is a rational decomposition
$\End_F(K)=K\oplus K\iota$. The two $K$-summands in the pair of
endomorphism spaces give four trivial directions. On the remaining
part,
\[
 [m_d\iota,m_x]=m_{d(\bar x-x)}\iota,\qquad
 [m_d\iota,m_y]=m_{d(\bar y-y)}\iota.
\]
Put $a=\bar x-x$, $b=\bar y-y$. At least one is nonzero, and hence
\begin{equation}\label{m1:eq:dense-nonsplit-quotient}
 N_{X,Y}=F^4_{\triv}\oplus V_{\mathrm{ns}},\qquad
 V_{\mathrm{ns}}=K^2/K(a,b),\qquad
 \ell(d_1,d_2)=bd_1-ad_2.
\end{equation}
The kernel of $\ell$ is $K(a,b)$, giving $V_{\mathrm{ns}}\simeq K$.
Conjugation by $t\in K^\times$ multiplies both root coordinates by
$t/\bar t$. For the trace-form adjoint,
\[
 H_0^{-1}(m_c+m_d\iota)^tH_0=m_c+m_{\bar d}\iota,
\]
because $\operatorname{Tr}(d\bar u v)=\operatorname{Tr}(\bar d u\bar v)$.
Thus the adjoint fixes the four multiplication directions, conjugates
the two ambient root coordinates, and preserves $K(a,b)$. Since
$\bar a=-a$ and $\bar b=-b$, its action after $\ell$ is
\begin{equation}\label{m1:eq:dense-nonsplit-action}
 t:d\longmapsto(t/\bar t)d,\qquad \tau_E:d\longmapsto-\bar d.
\end{equation}
The effective group is $K^\times/F^\times\simeq K^1$, by Hilbert 90.

For the scalar pair $X=xI,Y=yI$, let $\chi$ be the one-dimensional
module on which the two generators act by $x$ and $y$. The orbit tangent in
\eqref{m1:eq:dense-normal} is zero. Thus $S=\GL_2(F)$ acts by simultaneous
conjugation on $N_{X,Y}=M_2(F)^2$, and transpose is componentwise. In
multiplicity-space terms this is $M_2(F)\otimes_F\Ext^1(\chi,\chi)$,
where $\dim\Ext^1(\chi,\chi)=2$.

\subsubsection{Normal representations for rank-zero orbits}\label{m1:subsec:rankzero-normal}
At $x_D=(I,0,0,D)$, the stabilizer is
\[
 S_D=\{(g,t,g):g\in\GL_2(F),\ t\in Z_{\GL_2(F)}(D)\}
 \simeq\GL_2(F)\times Z_{\GL_2(F)}(D).
\]
This also follows from $\Efac^{\oplus2}\oplus\Qfac_D$: there are no
Hom terms between the two supports. The loop centralizer is
$(F^\times)^2$, $K^\times$, or $\GL_2(F)$ for split regular, nonsplit
regular, or scalar $D$, respectively. Put
$\mathfrak z(D)=\{Z\in M_2(F):ZD=DZ\}$. The differential specializes to
$(P_0-U_0,0,0,[Q_0,D])$, so
\begin{equation}\label{m1:eq:rankzero-normal}
 N_D=M_2(F)\oplus M_2(F)\oplus M_2(F)/[M_2(F),D]
 \simeq M_2(F)^2\oplus\mathfrak z(D).
\end{equation}
Semisimplicity of $\operatorname{ad}D$ gives the direct sum
$M_2(F)=[M_2(F),D]\oplus\mathfrak z(D)$ over $F$. Equivalently, the
trace pairing is nondegenerate on the centralizer. The action is
\begin{equation}\label{m1:eq:rankzero-action}
 (g,t)(X,Y,Z)=(gXt^{-1},tYg^{-1},tZt^{-1}).
\end{equation}

For regular $D$, choose a symmetric matrix $S$ such that $SD=D^tS$
and $St=t^tS$ for every $t\in F[D]^\times$. In the split eigenbasis use
$S=I$; in the nonsplit model $K=F[D]$, use the trace form. The rational
transporter $(I,S^{-1},I)$ corrects the point and induces
\begin{equation}\label{m1:eq:rankzero-transpose}
 \tau_D(X,Y,Z)=(Y^tS,S^{-1}X^t,S^{-1}Z^tS).
\end{equation}
This squares to one and conjugates $(g,t)$ to $(g^{-t},t^{-1})$.
Every $Z=aI+bD$ is fixed by its last component, since
$S^{-1}Z^tS=aI+bD$. The torus acts trivially there as well. Thus the
centralizer directions are $F^2$ fixed by both actions. As a concrete
alternative, for the companion matrix
$D=\left(\begin{smallmatrix}0&-d\\1&s\end{smallmatrix}\right)$,
$j=\left(\begin{smallmatrix}-s&1\\1&0\end{smallmatrix}\right)$ satisfies
$jD^t=Dj$; the trace-form construction above also determines the action
on both legs.

For scalar $D=dI$, with $d\neq0$, no correction is needed:
\begin{equation}\label{m1:eq:rankzero-scalar}
 S_D=\GL_2(F)^2,\quad N_D=M_2(F)^2\oplus\mathfrak{gl}_2(F),\quad
 \tau_D(X,Y,Z)=(Y^t,X^t,Z^t).
\end{equation}
The last summand carries the adjoint action. The
regular cases will use the Two-Leg Lemma with two fixed parameters;
the scalar case uses the Identity Descendant on all three summands.

\subsubsection{The decomposable case}\label{m1:subsec:decomposable-normal}
At \eqref{m1:eq:decomposable-rep}, the three factors
$\Efac,\Rfac_{x,d},\Qfac_e$ are pairwise nonisomorphic and absolutely
stable. Their scalar automorphisms give the embedding
\begin{equation}\label{m1:eq:decomposable-stabilizer}
 S=\{(\diag(r,b),\diag(r,s),\diag(r,b)):r,b,s\in F^\times\}.
\end{equation}
Their dimension vectors differ, so neither $x=0$ nor $d=e$ merges
factors or enlarges this stabilizer. Put
\begin{equation}\label{m1:eq:decomposable-characters}
 \alpha(r,b,s)=r/s,\qquad\beta(r,b,s)=b/s.
\end{equation}
The common scalar acts trivially, and $(\alpha,\beta)$ is jointly
surjective on rational points: $(r,b,s)=(a_0,b_0,1)$ maps to
$(a_0,b_0)$.

By \eqref{m1:eq:balanced-euler}, the Ext dimensions, with first argument
in the row and second in the column, are
\begin{equation}\label{m1:eq:decomposable-ext-table}
 \begin{array}{c|ccc}
 \dim\Ext^1&\Efac&\Rfac&\Qfac\\\hline
 \Efac&0&0&1\\
 \Rfac&0&2&1\\
 \Qfac&1&1&1
 \end{array}
\end{equation}
For example the $\Efac$--$\Rfac$ Euler value is zero, the
$\Efac$--$\Qfac$ and $\Rfac$--$\Qfac$ values are $-1$, and the
self-values are $1,-1,0$, respectively. The Hom dimensions are one
on the diagonal and zero elsewhere. Target scalar divided by source
scalar is the action on each cross-Ext. Therefore
\begin{equation}\label{m1:eq:decomposable-normal}
 N_E=F^3_{\triv}\oplus F_\alpha\oplus F_{\alpha^{-1}}
                  \oplus F_\beta\oplus F_{\beta^{-1}}.
\end{equation}
The normal space has dimension $7=16-(12-3)$.

The representative is fixed by transpose, whose ambient differential is
$(\dot A,\dot B,\dot C,\dot D)\mapsto
(\dot A^t,\dot C^t,\dot B^t,\dot D^t)$. To identify its action on
self-Ext, write the trivial ambient coordinates as
$a_1=\dot A_{11}$, $a_2=\dot A_{22}$, $b_1'=\dot B_{11}$,
$c_1'=\dot C_{11}$, $d_1=\dot D_{11}$, and $d_2=\dot D_{22}$.
For diagonal infinitesimal changes of basis,
$a_1=x(b_1'+c_1')$, $a_2$ is an independent orbit direction, and
the two $D$-coordinates are zero. The three quotient coordinates are
therefore
\begin{equation}\label{m1:eq:decomposable-fixed-coordinates}
 a_1-x(b_1'+c_1'),\qquad d_1,\qquad d_2.
\end{equation}
They are independent on the quotient by
\eqref{m1:eq:decomposable-ext-table}. Transpose interchanges $b_1',c_1'$
and fixes $a_1,d_1,d_2$, so it fixes all three.

On the stabilizer, transpose inverts $(r,b,s)$, hence both effective
characters. It exchanges each pair of cross-Ext lines. Its square is
one; choosing a rational nonzero vector in each positive line and its
image in the inverse line gives
\begin{equation}\label{m1:eq:decomposable-transpose}
 \tau_E(z_1,z_2,z_3,x,y,u,v)=(z_1,z_2,z_3,y,x,v,u).
\end{equation}
Thus the normal representation is a sum of two independently scaled
opposite-weight planes and three fixed directions, for $e\neq0$ and $xd\neq1$.

\subsubsection{The stable-system case}\label{m1:subsec:stable-normal}
At \eqref{m1:eq:stable-rep-general}--\eqref{m1:eq:H}, the decomposition is
$E=\Efac\oplus\Hfac$. Absolute stability and the distinct dimension
vectors give
\begin{equation}\label{m1:eq:stable-stabilizer}
 S=\{s(r,h)=(\diag(r,h),rI,\diag(r,h)):r,h\in F^\times\},
 \qquad\eta(r,h)=h/r.
\end{equation}
The effective character is surjective, since $\eta(1,t)=t$. This also
follows directly from cyclicity: a matrix commuting with $D$ and
satisfying $qe_1=re_1$ has $qDc=rDc$, hence is $rI$.

The balanced Euler formula now gives
\begin{equation}\label{m1:eq:stable-ext-table}
 \begin{array}{c|cc}
 \dim\Ext^1&\Efac&\Hfac\\\hline
 \Efac&0&1\\\Hfac&1&4
 \end{array},\qquad
 N_E=F^4_{\triv}\oplus F_\eta\oplus F_{\eta^{-1}}.
\end{equation}
Indeed the self-Euler values are $1,-3$ and the cross-Euler value is
$-1$. The four trivial directions are $\Ext^1(\Hfac,\Hfac)$, while the
two moving lines are cross-Ext. The total dimension is $6=16-(12-2)$.

The symmetric matrix $H$ from \eqref{m1:eq:H} satisfies
\begin{equation}\label{m1:eq:H-identities}
 HD=D^tH,\qquad He_1=b^t,\qquad E_{11}H=B,
 \qquad H^{-1}B^t=E_{11}.
\end{equation}
For instance both $HD$ and $D^tH$ equal the symmetric matrix with
entries $b_1$, $sb_1-db_0$, and $-db_1+s(sb_1-db_0)$. Consequently
\begin{equation}\label{m1:eq:stable-transporter}
 g(E)\tau(E)=E,\qquad g(E)=(I,H^{-1},I).
\end{equation}
The induced ambient operator is
\begin{equation}\label{m1:eq:stable-ambient-transpose}
 T_H(\dot A,\dot B,\dot C,\dot D)
 =(\dot A^t,\dot C^tH,H^{-1}\dot B^t,H^{-1}\dot D^tH).
\end{equation}
It preserves the orbit tangent because it fixes the base point and
normalizes $L$. Symmetry of $H$ gives $T_H^2=1$. Moreover
$\tau_Es(r,h)\tau_E^{-1}=s(r^{-1},h^{-1})$: the middle component of
$s$ is scalar and commutes with $H^{-1}$. Thus $\eta$ is inverted.
Choosing $e_{\eta^{-1}}=\tau_Ee_\eta$ makes the two cross-Ext lines an
ordinary rational coordinate swap.

It remains to identify the operator on all four self-Ext directions.
For the normalized family $E=E(a,b_0,d,s)$ its parameter derivatives are
\begin{equation}\label{m1:eq:stable-parameter-derivatives}
 \begin{aligned}
 v_a&=(E_{11},-dE_{12},0,0),&
 v_{b_0}&=(0,E_{11}+sE_{12},0,0),\\
 v_d&=(0,-aE_{12},0,-E_{12}),&
 v_s&=(0,b_0E_{12},0,E_{22}).
 \end{aligned}
\end{equation}
All have trivial weight.

\begin{lemma}\label{m1:lem:stable-parameter-basis}
Their classes form a basis of the trivial summand of the normal space,
and the corrected transpose fixes this summand pointwise.
\end{lemma}
\begin{proof}
The split stabilizer torus decomposes the complex into weight spaces,
so a weight-zero orbit tangent has a weight-zero preimage. On the
$\Hfac$-factor such a preimage consists of scalars $p,u$ at the ends
and $Q_0=(q_{ij})$ in the middle. The $\Efac$-factor has zero self-Ext.
Suppose a parameter variation is an orbit tangent. Since $c=e_1$ is
fixed, $Q_0e_1=ue_1$, so $q_{11}=u$ and $q_{21}=0$. Comparing
\[
 [Q_0,D]=\begin{pmatrix}0&-\dot d\\0&\dot s\end{pmatrix}
\]
in entries $(1,1)$ and $(2,1)$ gives $q_{12}=0$, $q_{22}=u$. Thus
$Q_0=uI$ and $\dot d=\dot s=0$. With $\delta=p-u$, the remaining
equations are
\[
 \dot a=a\delta,\qquad\dot b_0=b_0\delta,\qquad
 -d\dot a+s\dot b_0=b_1\delta.
\]
The normalization $b_1=1-ad+sb_0$ forces
\begin{equation}\label{m1:eq:stable-independence-identity}
 (-ad+sb_0-b_1)\delta=-\delta=0.
\end{equation}
All four variations vanish. The classes are independent and, by
\eqref{m1:eq:stable-ext-table}, span the four-dimensional trivial quotient.

Differentiate \eqref{m1:eq:stable-transporter} along any of these parameter
directions $v$. It gives
\[
 (dg)_E(v)\cdot\tau(E)+g(E)\tau(v)=v.
\]
The first term is an orbit tangent: with $\mathfrak l=\Lie(L)$,
$\xi=(dg)_E(v)g(E)^{-1}\in\mathfrak l$ and the term equals $\xi\cdot E$.
Hence $[g(E)\tau(v)]=[v]$ in the normal space. These four classes span
its trivial summand, proving pointwise fixation.
\end{proof}
Thus \eqref{m1:eq:stable-ext-table} has four fixed coordinates and a swap
on $F_\eta\oplus F_{\eta^{-1}}$. The argument is valid at every
$\det H\neq0$, including nilpotent and repeated-root loop values.

\subsubsection{Normal representations and transpose}\label{m1:subsec:normal-interface}
\begin{corollary}[Rational transpose stability]\label{m1:cor:transpose-stable}
Every closed $L(F)$-orbit in $\cX(F)$ is stable under $\tau$. At every
point its extended stabilizer and normal representation are as in
\eqref{m1:eq:corrected-extension} and the calculations above.
\end{corollary}
\begin{proof}
The transporters for dense and rank-zero representatives were given in
Lemma~\ref{m1:lem:dense-transpose} and \eqref{m1:eq:rankzero-transpose}.
Decomposable representatives are fixed by transpose; stable-system
representatives use \eqref{m1:eq:stable-transporter}. For $y=lE$ a
transporter is $l g(E)\alpha(l)^{-1}$, so the assertion holds on each
rational orbit. The classification is exhaustive.
\end{proof}
For the untwisted stabilizer action, two choices of transporter differ
by $L_E(F)$ and hence agree on $L_E(F)$-invariant normal distributions.

The following table summarizes the nine rational normal
representations. The last column indicates the method used in
Section~\ref{m1:sec:descendants} to prove invariance under the
corrected transpose. We write
$P_\gamma=F_\gamma\oplus F_{\gamma^{-1}}$, with its coordinate swap,
and $\mathcal V_T=\Hom(F^2,F^2)\oplus\Hom(F^2,F^2)$ with action
$(g,t)(X,Y)=(gXt^{-1},tYg^{-1})$ and the adjoint exchange
\eqref{m1:eq:rankzero-transpose}. Each $F^j$-summand in the classification is fixed by both
the stabilizer and corrected transpose. The nonsplit $K$ in the dense
row has action and odd operator \eqref{m1:eq:dense-nonsplit-action}.
\begin{center}
\small
\begin{tabular}{@{}llll@{}}
\toprule
Family & Stabilizer & Normal space & Vanishing method\\
\midrule
Dense absolute & $F^\times$ & $F^5$ & fixed operator\\
Dense split & $(F^\times)^2$ & $F^4\oplus P_\alpha$ & opposite weights\\
Dense nonsplit & $K^\times$ & $F^4\oplus K$ & norm fibers\\
Dense scalar & $\GL_2(F)$ & $M_2(F)^2$ & Prasad, translations\\
Zero split regular & $\GL_2(F)\times(F^\times)^2$ & $\mathcal V_T\oplus F^2$ & Two-Leg Lemma\\
Zero nonsplit regular & $\GL_2(F)\times K^\times$ & $\mathcal V_T\oplus F^2$ & Two-Leg Lemma\\
Zero scalar & $\GL_2(F)^2$ & $M_2(F)^2\oplus\mathfrak{gl}_2(F)$ & Identity Descendant\\
One decomposable & $(F^\times)^3$ & $F^3\oplus P_\alpha\oplus P_\beta$ & independent planes\\
One stable system & $(F^\times)^2$ & $F^4\oplus P_\eta$ & opposite weights\\
\bottomrule
\end{tabular}
\end{center}
The scalar matrix-pair row has componentwise transpose, and the scalar
rank-zero row has \eqref{m1:eq:rankzero-scalar}, including its adjoint
summand. The rational character maps are
\eqref{m1:eq:dense-split-normal}, \eqref{m1:eq:dense-nonsplit-action},
\eqref{m1:eq:decomposable-characters}, and \eqref{m1:eq:stable-stabilizer};
the remaining effective actions are given by the matrix formulas above.
In every row $\dim N_E=4+\dim\Aut(E)$: respectively
$5,6,6,8,10,10,12,7,6$. The finite extension acts on the same normal space through the corresponding
odd operator and sign character.

\section{Linear descendants and multiplicity one}\label{sec:linear}
The linear multiplicity-one problem reduces by Frobenius descent to the quiver
action of Section~\ref{sec:quiver}. The regular rank-zero cases
are handled by the two-leg theorem, while the scalar rank-zero
case requires a separate identity descendant. We prove these statements before treating all nine rational families
and completing the linear multiplicity-one argument.

\subsection{The residual action}\label{m1:subsec:reduction-to-GL4}
For $G_{\mathrm L}=\GL_4(F)\times\GL_2(F)$ and the subgroup $H_{\mathrm L}$
of Theorem~\ref{m1:thm:main}, the group $J=H_{\mathrm L}\times H_{\mathrm L}$ acts by
\[
 ((a_1,b_1),(a_2,b_2))(g,z)
 =\bigl(\diag(a_1,b_1)g\diag(a_2,b_2)^{-1},\ b_1zb_2^{-1}\bigr).
\]
Componentwise transpose preserves $H_{\mathrm L}$ and normalizes this action, sending
the two acting factors to their inverse transposes in reversed order.
The projection $(g,z)\mapsto z$ is equivariant with transitive base
$\GL_2(F)$. At $z=I_2$ its stabilizer is
\[
 J_{I_2}=\{((p,q),(u,q)):p,q,u\in\GL_2(F)\}
 \simeq L(F),\qquad L=\GL(P)\times\GL(Q)\times\GL(U),
\]
where $P,Q,U$ have dimension two. Both $J$ and $J_{I_2}$ are unimodular.
The fiber is $\cX\simeq\GL_4(F)$ and transpose preserves its base point.
The finite extensions are unimodular as well: a modulus is trivial on
the ordinary group and factors through a finite group. Thus
\begin{equation}\label{m1:eq:residual-Frobenius}
 \mathcal D(G_{\mathrm L})^J\simeq\mathcal D(\cX)^{L(F)}
\end{equation}
is transpose-compatible, including its odd component. In block
coordinates $x=(A,B,C,D)$, the resulting action and involution are
exactly \eqref{m1:eq:residual-action}--\eqref{m1:eq:tau}.
The variety $\cX$ is smooth affine and $\alpha$ is algebraic over $F$.
It remains to prove $\mathcal D(\cX)^{\widetilde L,\varepsilon}=0$.
By Corollary~\ref{m1:cor:HC-invariance}, this follows from rational transpose
stability of every closed orbit and the absence of sign-equivariant
distributions on its normal space.
Only one representative per rational orbit is needed: $y=l x$ identifies
both stabilizers and normal spaces by conjugation and the differential of
$l$. The geometric data were established in Section~\ref{sec:quiver};
it remains to rule out sign-equivariant distributions on these normal spaces.

\subsection{The two-leg distribution theorem}\label{m1:sec:two-leg}
\subsubsection{Adjoints and the statement}\label{m1:subsec:two-leg-statement}

Let $V$ and $W$ be two-dimensional $F$-vector spaces, and let
$T\subset \GL(W)$ be a maximal torus.  Put \(\Gamma=\GL(V)\times T\) and \(\Vlegs=\Hom(W,V)\oplus
\Hom(V,W)\). We write a point of $\Vlegs$ as $(X,Y)$, with \(X:W\longrightarrow V\) and
\(Y:V\longrightarrow W\), and let $\Gamma$ act by
\[
    (a,t)\cdot(X,Y)=\bigl(aXt^{-1},\,tYa^{-1}\bigr).
\]
Fix a basis of $V$ with its standard symmetric form. Choose a
nondegenerate symmetric matrix $S$ on $W$ such that
\begin{equation}\label{m1:eq:torus-self-adjoint}
 S t=t^t S\qquad(t\in T).
\end{equation}
Define the adjoint exchange
\begin{equation}\label{m1:eq:two-leg-adjoint}
 \tau=\tau_S,\qquad \tau_S(X,Y)=(Y^tS,S^{-1}X^t).
\end{equation}
It satisfies $\tau_S^2=1$ and, by direct multiplication,
\[
 \tau_S\bigl((a,t)\cdot(X,Y)\bigr)
 =(a^{-t},t^{-1})\cdot\tau_S(X,Y).
\]
In particular it normalizes the given action. Such a form exists over $F$:
in a diagonal realization of a split torus take $S=I_2$; for a nonsplit
torus identify $W$ with a quadratic field $K$ and $T$ with $K^\times$
acting by multiplication, and take the Gram matrix of
$\operatorname{Tr}_{K/F}(ww')$. The trace form is symmetric and
nondegenerate because $K/F$ is separable, and multiplication is self-adjoint.

A change of basis $w=Cw'$ replaces $(X,Y,S)$ by
$(XC,C^{-1}Y,C^tSC)$ and conjugates $T$ accordingly, preserving \eqref{m1:eq:torus-self-adjoint} and
\eqref{m1:eq:two-leg-adjoint}. In a fixed realization, any other such form is
$S=S_0c$ with $c\in T$: indeed $S_0^{-1}S$ commutes with $T$, whose
centralizer in $\operatorname{End}_F(W)$ is its quadratic \'{e}tale algebra.
Then $\tau_S=(1,c^{-1})\tau_{S_0}$, so these forms induce the same
operator on $\Gamma$-invariant distributions. It therefore suffices to prove
the split and trace-form realizations just specified. In this section,
the adjoint exchange on the two legs always means $\tau_S$; on
$\operatorname{End}_F(W)$ it means $A^*=S^{-1}A^tS$.

Let \(\widetilde \Gamma=\Gamma\rtimes\langle\tau\rangle\) and let $\varepsilon$ be the sign character,
\(\varepsilon|_\Gamma=1\) and \(\varepsilon(\tau)=-1\).

\begin{theorem}[Two-Leg Lemma]\label{m1:thm:two-leg}
For every maximal torus $T\subset \GL_2(F)$ and every nondegenerate symmetric form $S$ satisfying
\eqref{m1:eq:torus-self-adjoint},
\[
    \mathcal D(\Vlegs)^{\widetilde \Gamma,\varepsilon}=0.
\]
Equivalently, \(\mathcal D(\Vlegs)^\Gamma=\mathcal D(\Vlegs)^{\Gamma,\tau}\). \end{theorem}

We first establish the matrix-axis boundary statement and
the rank-one factorization result common to both rational torus types.

\subsubsection{Matrix axes and local zeta distributions}\label{m1:subsec:matrix-axes}

Let
\[
    B=
    \bigl(\Hom(W,V)\oplus\{0\}\bigr)
    \cup
    \bigl(\{0\}\oplus\Hom(V,W)\bigr).
\]
Thus $B$ is the union of the two matrix axes.

\begin{lemma}[Matrix two-axis boundary]\label{m1:lem:matrix-boundary}
For every maximal torus $T\subset \GL(W)$,
\[
    \mathcal D(B)^\Gamma=\mathcal D(B)^{\Gamma,\tau}.
\]
Moreover, if \(\lambda(a,t)=\left|\frac{\det a}{\det t}\right|\), then the invariant distributions on the
two rank-two open axes admit
finite-part extensions $P_X$ and $P_Y$ satisfying
\begin{equation}\label{m1:eq:matrix-finite-parts}
    \tau P_X=P_Y,
    \qquad
    gP_X-P_X=-\log\lambda(g)\,R,
    \qquad
    gP_Y-P_Y=+\log\lambda(g)\,R,
\end{equation}
where $R=c\delta_0$ with $c\neq0$.
\end{lemma}

\begin{proof}
\emph{Rank-one axes.} On the rank-one stratum of the $X$-axis, the map \(X\longmapsto \operatorname{im}X\in
\mathbf P(V)\) is $\Gamma$-equivariant, and $\GL(V)$ acts transitively on $\mathbf P(V)$.
At the line $Fe_1$, Frobenius descent leaves the group
\[
    P\times T,
    \qquad
    P=
    \left\{
    \begin{pmatrix}
        \alpha&b\\0&\beta
    \end{pmatrix}:
    \alpha,\beta\in F^\times,\ b\in F
    \right\},
\]
acting on $W^\vee\setminus\{0\}$ by \(c\longmapsto \alpha c t^{-1}\). With our convention for
modular characters,
\[
    \Delta_P
    \begin{pmatrix}
        \alpha&b\\0&\beta
    \end{pmatrix}
    =\left|\frac{\alpha}{\beta}\right|.
\]
Since $\Gamma$ is unimodular, Proposition~\ref{m1:prop:Frobenius-descent} identifies $\Gamma$-invariant
 distributions on the rank-one stratum with distributions on the descended
fiber transforming by the character $\Delta_P^{-1}=|\beta/\alpha|$.  Now take
\[
    p_\beta=\begin{pmatrix}1&0\\0&\beta\end{pmatrix},
    \qquad t=1.
\]
This element acts trivially on the descended geometric fiber, whereas the
required equivariance character is $|\beta|$.  Choosing
$\beta\in F^\times$ with $|\beta|\neq1$ forces the descended distribution to
vanish.  Hence the rank-one $X$-stratum carries no invariant distribution of
the required Frobenius type.  The same argument, after transpose, treats the
rank-one $Y$-stratum.

\emph{Determinant zeta distributions.}\label{m1:para:det-zeta}
We normalize additive Haar measure by
$\operatorname{vol}(M_2(\mathcal O_F))=1$ and, initially for
$\Re(s)\gg0$, define
\[
 Z_X(s)(f)=\int_{M_2(F)}f(X,0)|\det X|^s\,dX,
 \qquad Z_Y(s)=\tau_S Z_X(s).
\]
With the same coordinate Haar measure on the two matrix spaces,
\[
 Z_Y(s)(f)=|\det S|^{s+2}
       \int_{M_2(F)}f(0,Y)|\det Y|^s\,dY,
\]
since $Y=S^{-1}X^t$ has absolute Jacobian $|\det S|^{-2}$.

Theorem~5.13 and Section~5.5 of \cite{HongSun2016}, applied to additive
Haar measure on $\{\det\ne0\}$, give meromorphic continuation and
distribution-valued Laurent coefficients. More precisely, with
$u=s+2$, the family $Z_X(u-2)$ is the case $n=2$, $\chi=1$ of the
matrix zeta integral in \cite[Section~7]{HongSun2016}.
Proposition~7.2(a), together with equation~(30) in its proof, shows that
the pole at $u=0$ is simple and its leading coefficient is a nonzero
multiple of $\delta_0$. Here the Laurent coordinate used in that paper
is $1-q_F^{-u}$, which has derivative $\log q_F\ne0$ at $u=0$.
Thus
\[
 R:=\operatorname*{Res}_{s=-2}Z_X(s)=c\delta_0,\qquad c\ne0.
\]
For the chosen normalization, the lattice formula in the proof of
\cite[Proposition~7.2]{HongSun2016} gives
\[
 I(s):=\int_{M_2(\mathcal O_F)}|\det X|^s\,dX
 =\frac{(1-q_F^{-1})(1-q_F^{-2})}
        {(1-q_F^{-s-1})(1-q_F^{-s-2})},
 \qquad
 c=\frac{(1-q_F^{-1})(1-q_F^{-2})}
        {(1-q_F)\log q_F}.
\]
Because $\tau_S$ fixes the origin, $Z_Y(s)$ has the same residue $R$.

For $g=(a,t)\in\Gamma$, the map $X\mapsto aXt^{-1}$ multiplies
$|\det X|$ by $\lambda(g)$ and has absolute Jacobian $\lambda(g)^2$.
Using the pushforward convention \eqref{m1:eq:distribution-conventions}
and the fact that conjugation by $\tau_S$ inverts $\lambda$, we obtain
\[
 gZ_X(s)=\lambda(g)^{-(s+2)}Z_X(s),\qquad
 gZ_Y(s)=\lambda(g)^{s+2}Z_Y(s).
\]
Define the finite parts using the coordinate $s+2$:
\[
 Z_X(s)=\frac{R}{s+2}+P_X+O(s+2),\qquad
 Z_Y(s)=\frac{R}{s+2}+P_Y+O(s+2).
\]
Then $P_Y=\tau_S P_X$, and comparison of constant terms gives
\[
 gP_X-P_X=-\log\lambda(g)\,R,\qquad
 gP_Y-P_Y=+\log\lambda(g)\,R,
\]
which is \eqref{m1:eq:matrix-finite-parts}.
On each rank-two open axis these finite parts restrict to the nonzero
invariant measure obtained by setting $s=-2$ in the defining integral.

\emph{Finite-part comparison.} Let $D\in\mathcal D(B)^\Gamma$, and let $c_X,c_Y$ be the coefficients of its
restrictions to the two rank-two open axes.  Put \(Q=D-c_XP_X-c_YP_Y\). Then $Q$ is supported on the
union of the two rank-at-most-one axes and
\[
    gQ-Q=(c_X-c_Y)\log\lambda(g)\,R.
\]
The right-hand side is supported at the origin.  Hence the restriction of
$Q$ to either rank-one stratum is invariant, and the rank-one
calculation above forces these restrictions to vanish.  Thus $Q$ is supported
at the origin.  A distribution supported at one point on an $l$-space is a
scalar multiple of the Dirac distribution, so $gQ-Q=0$.  Choosing $g$ with
$\lambda(g)\neq1$ gives $c_X=c_Y$.  Since $\tau P_X=P_Y$ and $\delta_0$ is
fixed by transpose, $D$ is fixed by $\tau$.
\end{proof}

\subsubsection{Singular rank-one factorizations}\label{m1:subsec:rankone-factorization}

\begin{lemma}[Rank-one factorization]\label{m1:lem:rank-one-factorization}
Let $A:W\to W$ have rank one. Suppose that $R$ is an invertible
symmetric matrix satisfying $RA=A^tR$. Put
\[
 \mathcal F_A=\{(X,Y)\in \Vlegs:YX=A\},\qquad
 \tau_R(X,Y)=(Y^tR,R^{-1}X^t).
\]
Then every $\GL(V)\times F^\times$-invariant distribution on
$\mathcal F_A$, with the second factor acting as scalar matrices on $W$,
is invariant under $\tau_R$.
\end{lemma}

\begin{proof}
The relation $RA=A^tR$ makes $\tau_R$ preserve $\mathcal F_A$;
$R^t=R$ makes it involutive. Set $Y'=RY$. Then $Y'X=RA$ is symmetric
of rank one, while the group action is still
$(X,Y')\mapsto(aXt^{-1},tY'a^{-1})$ because $t$ is scalar.
The involution is now $(X,Y')\mapsto((Y')^t,X^t)$.
Choose $C\in\GL_2(F)$ such that
\[
 C^tRA C=\begin{pmatrix}\lambda&0\\0&0\end{pmatrix},\qquad\lambda\ne0.
\]
This is congruence diagonalization of a symmetric rank-one form in
characteristic different from two. The paired changes
$(X,Y')\mapsto(XC,C^tY')$ preserve both the scalar two-leg action and
the ordinary transpose exchange. Relabel these new legs $(X,Y)$.
This congruence argument applies also to nilpotent $A$.

Write
\[
    X=(x_1,x_2),
    \qquad
    Y=\begin{pmatrix}y_1\\y_2\end{pmatrix},
\]
with $x_i\in V$ and $y_i\in V^\vee$. In these paired coordinates, the
defining equation of $\mathcal F_A$ is $YX=\operatorname{diag}(\lambda,0)$,
equivalently
\[
    y_1(x_1)=\lambda,
    \qquad
    y_1(x_2)=0,
    \qquad
    y_2(x_1)=0,
    \qquad
    y_2(x_2)=0.
\]
Consider
\[
    p:\mathcal F_A\longrightarrow
    \Omega_\lambda,
    \qquad
    p(X,Y)=(x_1,y_1),
\]
where \(\Omega_\lambda= \{(v,\ell)\in V\times V^\vee:\ell(v)=\lambda\}\). The group $\GL(V)\times
F^\times$ acts transitively on $\Omega_\lambda$.
Frobenius descent at $(e_1,\lambda e_1^\vee)$ leaves the stabilizer
\[
    H_A=
    \left\{
    \left(
    \begin{pmatrix}t&0\\0&c\end{pmatrix},
    tI_W
    \right):
    t,c\in F^\times
    \right\}
    \simeq(F^\times)^2.
\]
The fiber over $(e_1,\lambda e_1^\vee)$ is
\[
    x_2=xe_2,
    \qquad
    y_2=ye_2^\vee,
    \qquad
    xy=0.
\]
Thus the descended fiber is the scalar two-axis node \(N=\{(x,y)\in F^2:xy=0\}\). The effective quotient
of $H_A$ acts by \(q\cdot(x,y)=(qx,q^{-1}y)\), \(q\in F^\times\). All groups involved are
unimodular. Ordinary transpose sends the base
$(e_1,\lambda e_1^\vee)$ to $(\lambda e_1,e_1^\vee)$. Correct it by
$(a_0,1)$ with $a_0=\operatorname{diag}(\lambda^{-1},1)$. This fixes the
base and sends the fiber coordinates $(x,y)$ to $(y,x)$.
The correcting element composed with transpose is an involution and has
trivial modular character. Prasad's Two-Axis Lemma therefore gives the required
transpose invariance on $N$, and Frobenius descent gives the assertion on
$\mathcal F_A$.
\end{proof}

\subsubsection{The zero-product projective descent}\label{m1:subsec:zero-product-descent}
\begin{lemma}\label{m1:lem:common-U11}
Let $\mathcal F_0=\{(X,Y):YX=0\}$ and
$U_{11}=\{(X,Y)\in\mathcal F_0:\rk X=\rk Y=1\}$.
The decomposition $\mathcal F_0=U_{11}\sqcup B$ has $U_{11}$ open and
$B$ the closed matrix axes. Descent along
$(X,Y)\mapsto\operatorname{im}X=\ker Y\in\mathbf P(V)$ identifies the
ordinary invariant-distribution problem on $U_{11}$ with the
$\eta$-equivariant problem for $P\times T$ on
$(W^\vee\setminus\{0\})\times(W\setminus\{0\})$, where
\begin{equation}\label{m1:eq:common-U11-action}
 a=\begin{pmatrix}\alpha&b\\0&\beta\end{pmatrix},\qquad
 (c,d)\longmapsto(\alpha ct^{-1},\beta^{-1}td),\qquad
 \eta(a,t)=|\beta/\alpha|.
\end{equation}
The corrected odd operator is
\begin{equation}\label{m1:eq:common-U11-odd}
 (c,d)\longmapsto(d^tS,S^{-1}c^t).
\end{equation}
It squares to one, preserves the character $\eta$, and has relative
modular character one on the odd component.
\end{lemma}
\begin{proof}
When both legs have rank one, $YX=0$ forces $\operatorname{im}X=\ker Y$.
The map to $\mathbf P(V)$ is equivariant, and at $Fe_1$ the fiber
is $X=e_1c$, $Y=de_2^\vee$. Matrix multiplication gives
\eqref{m1:eq:common-U11-action}; the upper unipotent entry acts trivially.
The ordinary relative modulus is $\Delta_P^{-1}=|\beta/\alpha|$.
Writing $j=\left(\begin{smallmatrix}0&1\\1&0\end{smallmatrix}\right)$
on $V$, the operator $(j,1)\tau_S$ fixes the base line and gives
\eqref{m1:eq:common-U11-odd}. It has square one and sends the diagonal
parameters $(\alpha,\beta)$ to $(\beta^{-1},\alpha^{-1})$, preserving
$\eta$. The positive relative modular character is one on this
order-two operator, so Proposition~\ref{m1:prop:odd-Frobenius} gives the required invariance.
\end{proof}

\subsubsection{Nonsplit fibers}\label{m1:subsec:two-leg-nonsplit}

Assume that $W=K$, where $K/F$ is a quadratic field extension. We identify \(T=K^\times\subset
\GL_F(K)\) through multiplication and take $S$ to be the trace-form Gram matrix above.
Write $\bar w$ for the nontrivial
Galois conjugate of $w\in K$, and let $\iota(w)=\bar w$.

Every $F$-linear endomorphism of $K$ has a unique expression
\[
    A=m_c+m_d\iota,
    \qquad c,d\in K,
\]
where $m_c$ denotes multiplication by $c$.  For $A=YX$, define
\[
    \pi_{\mathrm{ns}}:\Vlegs\longrightarrow K\oplus F,
    \qquad
    \pi_{\mathrm{ns}}(X,Y)=\bigl(c,N_{K/F}(d)\bigr).
\]
Conjugation by $t\in K^\times$ gives
\[
    t(m_c+m_d\iota)t^{-1}
    =m_c+m_{(t/\bar t)d}\iota,
\]
so $\pi_{\mathrm{ns}}$ is $\Gamma$-invariant.

With respect to the trace pairing \(\langle x,y\rangle_K=\operatorname{Tr}_{K/F}(xy)\), the adjoint of
$A$ is \(A^\ast=m_c+m_{\bar d}\iota\). Hence $\pi_{\mathrm{ns}}$ is fixed by $\tau_S$. Values outside the
norm image have empty fiber. If $N(d)\neq0$, Hilbert 90
shows that the set of $d'$ with $N(d')=N(d)$ is one rational $T$-orbit, with
stabilizer $F^\times$.  If $N(d)=0$, then $d=0$, so the base fiber is a single
point with stabilizer $K^\times$.  These stabilizers are unimodular.

Fix $A=m_c+m_d\iota$ in one such base orbit and put \(T_A=Z_T(A)\). The operators $A$ and $A^\ast$ lie in
the same rational $T$-orbit.  If
$d\neq0$, choose $t_A\in K^\times$ such that \(\frac{t_A}{\bar t_A}=\frac{d}{\bar d}\); if $d=0$, take
$t_A=1$.  Then \(t_AA^\ast t_A^{-1}=A\). Thus the descended transpose \(\tau_A=t_A\tau\) preserves the
factorization fiber \(\mathcal F_A=\{(X,Y):YX=A\}\). Since $t_A$ is self-adjoint, $\tau_A^2=1$. More
explicitly, with \(R=S t_A^{-1}\) one has $R^t=R$, $RA=A^tR$, and
$\tau_A(X,Y)=(Y^tR,R^{-1}X^t)$. These identities follow respectively
from $St_A=t_A^tS$ and $t_AA^*t_A^{-1}=A$.
Different choices of $t_A$ differ by an element of $T_A$ and therefore induce
the same operator on $(\GL(V)\times T_A)$-invariant distributions.

Frobenius descent along the rational $T$-orbit of $A$ reduces the
fiber of $\pi_{\mathrm{ns}}$ to $\mathcal F_A$, with acting group
$\GL(V)\times T_A$ and descended transpose $\tau_A$.

If $\operatorname{rank}A=2$, then $X$ and $Y$ are isomorphisms.  Fix an
isomorphism $X_0:W\xrightarrow{\sim}V$.  Every point of $\mathcal F_A$ is
uniquely of the form
\[
    (X,Y)=\bigl(gX_0,AX_0^{-1}g^{-1}\bigr),
    \qquad g\in\GL(V).
\]
Thus $\mathcal F_A$ is one $\GL(V)$-orbit.  Its stabilizer in
$\GL(V)\times T_A$ is a diagonal copy of $T_A$, hence unimodular.  The operator $\tau_A$ preserves this rational orbit and squares to one.
Lemma~\ref{m1:lem:orbital-measures} gives its invariant-distribution statement.

Suppose next that $\operatorname{rank}A=1$.  Then $d\neq0$, and
$T_A=F^\times I_W$. With $R=S t_A^{-1}$ as above, Lemma~\ref{m1:lem:rank-one-factorization} applies and gives
transpose invariance on the rank-one factorization fiber.

It remains to treat \(\mathcal F_0=\{(X,Y)\in \Vlegs:YX=0\}\). The possible rank pairs are
\((2,0),\allowbreak(1,0),\allowbreak(0,0),\allowbreak(0,1),\allowbreak(0,2),\allowbreak(1,1)\). Let
\(U_{11}=\{(X,Y)\in\mathcal F_0: \operatorname{rank}X=\operatorname{rank}Y=1\}\). Then \(\mathcal
F_0=U_{11}\sqcup B\), with $U_{11}$ open and $B$ closed.  Lemma~\ref{m1:lem:matrix-boundary} treats
$B$.

Lemma~\ref{m1:lem:common-U11} gives the required descent. In the
trace-form model, $c(w)=\operatorname{Tr}_{K/F}(c_0w)$ for a unique
$c_0\in K^\times$.
The descended action and odd operator become
\[
 (c_0,d)\longmapsto(\alpha c_0t^{-1},\beta^{-1}td),\qquad
 (c_0,d)\longmapsto(d,c_0),
\]
with ordinary character $\eta=|\beta/\alpha|$ and odd modulus one.

Consider
\[
    \pi_{11,\mathrm{ns}}(c_0,d)=[c_0d]
    \in K^\times/F^\times.
\]
This map is invariant under the descended action.  Every fiber is one
rational orbit: if $[c_0d]=[c_0'd']$, then
$c_0'd'=r c_0d$ for some $r\in F^\times$, and the descended action above provides
$t\in K^\times$ and $\alpha,\beta\in F^\times$ carrying $(c_0,d)$ to
$(c_0',d')$.  At a representative $(1,z)$ the stabilizer is
\[
    H_z=
    \left\{
    \left(
        \begin{pmatrix}s&b\\0&s\end{pmatrix},s
    \right):
    s\in F^\times,\ b\in F
    \right\}.
\]
It is unimodular, while $\Delta_{P\times T}|_{H_z}=1$ and
$\eta|_{H_z}=1$, where $\eta=|\beta/\alpha|$. Thus the second Frobenius
compatibility character is trivial. The swap of
Lemma~\ref{m1:lem:common-U11} preserves $\eta$ and has square one.
It preserves $[c_0d]$. Lemma~\ref{m1:lem:orbital-measures} therefore fixes the positive
equivariant orbital measure.  Bernstein localization on $\pi_{11,\mathrm{ns}}$ therefore gives
\[
    \mathcal D(U_{11})^\Gamma=\mathcal D(U_{11})^{\Gamma,\tau}.
\]
An anti-invariant $\Gamma$-invariant distribution on $\mathcal F_0$ restricts to
zero on $U_{11}$ and is therefore supported on $B$, where it vanishes by
Lemma~\ref{m1:lem:matrix-boundary}.  Hence
\begin{equation}\label{m1:eq:nonsplit-zero-fiber}
    \mathcal D(\mathcal F_0)^\Gamma
    =
    \mathcal D(\mathcal F_0)^{\Gamma,\tau}.
\end{equation}

A direct two-dimensional calculation gives
\[
    \det{_F}(m_c+m_d\iota)=N_{K/F}(c)-N_{K/F}(d).
\]
Consequently,
\[
\begin{array}{c|c}
N(c)\neq N(d) & \operatorname{rank}A=2,\\
N(c)=N(d),\ (c,d)\neq(0,0) & \operatorname{rank}A=1,\\
c=d=0 & A=0.
\end{array}
\]
The three cases above show that every fiber of
$\pi_{\mathrm{ns}}$ has zero $(\widetilde \Gamma,\varepsilon)$-equivariant distribution
space.  Bernstein localization (Proposition \ref{m1:prp:localization}) gives
\[
    \mathcal D(\Vlegs)^{\widetilde \Gamma,\varepsilon}=0.
\]
This proves Theorem~\ref{m1:thm:two-leg} for the nonsplit torus.

\subsubsection{Split nonzero fibers}\label{m1:subsec:two-leg-split-nonzero}

In the split case, we take $S=I_2$ and
\(T=\{\operatorname{diag}(t_1,t_2):t_i\in F^\times\}\). We write
\[
    X=(x_1,x_2),
    \qquad
    Y=\begin{pmatrix}y_1\\y_2\end{pmatrix},
\]
where $x_i\in V$ and $y_i\in V^\vee$.  Put
\[
    u=y_1(x_1),
    \qquad
    v=y_2(x_2),
    \qquad
    \alpha=y_1(x_2),
    \qquad
    \beta=y_2(x_1),
    \qquad
    r=\alpha\beta.
\]
Then
\[
    YX=
    \begin{pmatrix}
        u&\alpha\\
        \beta&v
    \end{pmatrix},
\]
and the map
\[
    \pi_{\mathrm{sp}}:\Vlegs\longrightarrow F^3,
    \qquad
    \pi_{\mathrm{sp}}(X,Y)=(u,v,r),
\]
is $\Gamma$-invariant and fixed by transpose.

We first treat the fibers \(\pi_{\mathrm{sp}}^{-1}(u,v,r)\), \((u,v,r)\neq(0,0,0)\). Since $\det(YX)=uv-r$, the fibers divide as follows:
\[
\begin{array}{c|c}
r\neq0,\ uv-r\neq0&\text{invertible product, nonzero hyperbola}\\
r\neq0,\ uv-r=0&\text{rank-one product}\\
r=0,\ uv\neq0&\text{invertible product, complete node}\\
r=0,\ \text{exactly one of }u,v\neq0&\text{two complete nodes}\\
(u,v,r)=(0,0,0)&\text{root-component intersection and gluing}
\end{array}
\]
The last row will be treated separately.

\emph{The fiber $r\neq0$, $uv-r\neq0$.}
Then $A=YX$ is invertible, hence both $X$ and $Y$ are invertible.  For fixed
$A$, every factorization is \((X,Y)=(g,Ag^{-1})\), \(g\in\GL(V)\), so the factorization
fiber is one $\GL(V)$-orbit.  The effective split torus
acts on $(\alpha,\beta)$ by \((\alpha,\beta)\longmapsto (q\alpha,q^{-1}\beta)\), where
\(q=t_1/t_2\). Since $\alpha\beta=r\neq0$, the set of such pairs is one rational torus
orbit.  Therefore the $\pi_{\mathrm{sp}}$-fiber is a single
homogeneous $\Gamma$-orbit.  Its stabilizer is unimodular, and transpose preserves
the orbit.  The positive invariant orbital distribution is therefore fixed
by transpose, and
\[
    \mathcal D\bigl(\pi_{\mathrm{sp}}^{-1}(u,v,r)\bigr)^{\widetilde \Gamma,\varepsilon}=0.
\]

\emph{The fiber $r\neq0$, $uv-r=0$.}
Now $A=YX$ has rank one and $\alpha,\beta\neq0$.  For fixed $(u,v,r)$ the
matrices
\[
    A=\begin{pmatrix}u&\alpha\\\beta&v\end{pmatrix},
    \qquad \alpha\beta=r,
\]
form one rational $T$-orbit.  Fix one such $A_0$.  Since both off-diagonal
entries of $A_0$ are nonzero, \(T_{A_0}=F^\times I_2\). Frobenius descent along the $T$-orbit of $A_0$
reduces the fiber to
$\mathcal F_{A_0}$. If $A_0=\begin{pmatrix}u&\alpha\\\beta&v\end{pmatrix}$,
take $t_A=\operatorname{diag}(\alpha/\beta,1)$. Then
$t_AA_0^tt_A^{-1}=A_0$, and $R=t_A^{-1}$ is symmetric with
$RA_0=A_0^tR$. The descended involution is exactly $\tau_R$.
Lemma~\ref{m1:lem:rank-one-factorization} applies and gives
\[
    \mathcal D\bigl(\pi_{\mathrm{sp}}^{-1}(u,v,r)\bigr)^{\widetilde \Gamma,\varepsilon}=0.
\]

\emph{The fiber $r=0$, $uv\neq0$.}
Again $A=YX$ is invertible. In the coordinates
$(X,Y)=(g,Ag^{-1})$, the action is
$(a,t)(g,A)=(agt^{-1},tAt^{-1})$. Thus every
$\GL(V)$-invariant distribution has the form $dg\otimes D_N$,
where $dg$ is Haar measure on $\GL(V)$ and $D_N$ is a distribution
on the node $N=\{(\alpha,\beta):\alpha\beta=0\}$.
Right invariance of $dg$ identifies $T$-invariance with invariance
of $D_N$ under $(\alpha,\beta)\mapsto(q\alpha,q^{-1}\beta)$,
where $q=t_1/t_2$.
Transpose is $(g,A)\mapsto(g^{-t}A^t,A^t)$. For each fixed $A$,
inversion, transpose, and right translation preserve $dg$, so
integration in $g$ gives
$\tau_*(dg\otimes D_N)=dg\otimes(A\mapsto A^t)_*D_N$.
The latter operator exchanges $\alpha$ and $\beta$ and fixes $u,v$.
Prasad's Two-Axis Lemma therefore gives
\[
    \mathcal D\bigl(\pi_{\mathrm{sp}}^{-1}(u,v,0)\bigr)^{\widetilde \Gamma,\varepsilon}=0,
    \qquad uv\neq0.
\]

\emph{The fibers with exactly one nonzero diagonal entry.}
By symmetry assume $u\neq0$ and $v=0$.  Consider the equivariant map
\[
    p:\pi_{\mathrm{sp}}^{-1}(u,0,0)
    \longrightarrow
    \Omega_u,
    \qquad
    p(X,Y)=(x_1,y_1),
\]
where \(\Omega_u=\{(x,\ell)\in V\times V^\vee:\ell(x)=u\}\). The group $\Gamma$ acts transitively on
$\Omega_u$. We apply Frobenius descent at
$(e_1,ue_1^\vee)$. The stabilizer is
\[
    S=
    \left\{
    \left(
        \begin{pmatrix}t_1&0\\0&c\end{pmatrix},
        \operatorname{diag}(t_1,t_2)
    \right):
    t_1,t_2,c\in F^\times
    \right\}
    \simeq(F^\times)^3.
\]
Write \(x_2=xe_1+ze_2\) and \(y_2=ye_1^\vee+we_2^\vee\). The equations $v=0$ and $r=0$ become
exactly \(xy=0\) and \(zw=0\). Thus the descended fiber is \(N=N_1\times N_2\), where
\(N_1=\{(x,y):xy=0\}\) and \(N_2=\{(z,w):zw=0\}\). The two opposite-weight characters are
\(\varepsilon_1=\frac{t_1}{t_2}\) and \(\varepsilon_2=\frac{c}{t_2}\), and
\((\varepsilon_1,\varepsilon_2):S\longrightarrow(F^\times)^2\) is jointly surjective: take $t_2=1$, $t_1=q_1$, $c=q_2$ for a
prescribed pair $(q_1,q_2)$. Lemma~\ref{m1:lem:passive-products} on the
product of the two nodes gives invariance under each swap.
To track the involution, correct transpose at the chosen base by
$(\operatorname{diag}(u^{-1},1),I)$. On the fiber it acts by
\[
 (x,y,z,w)\longmapsto(u^{-1}y,ux,w,z).
\]
This is the product of the two swaps followed by the first effective
torus factor $u^{-1}$, with the second factor equal to one; both factors exist by the joint surjectivity of
$(\varepsilon_1,\varepsilon_2)$ established above. It therefore fixes every
$S$-invariant distribution. Hence
\[
    \mathcal D\bigl(\pi_{\mathrm{sp}}^{-1}(u,0,0)\bigr)^{\widetilde \Gamma,\varepsilon}=0.
\]
The case $u=0$, $v\neq0$ is obtained by interchanging the two eigenlines.

The four cases exhaust every nonzero value of $\pi_{\mathrm{sp}}$.  Therefore
\begin{equation}\label{m1:eq:split-nonzero-fibers}
    \mathcal D\bigl(\pi_{\mathrm{sp}}^{-1}(u,v,r)\bigr)^{\widetilde \Gamma,\varepsilon}=0
    \qquad
    \text{for }(u,v,r)\neq(0,0,0).
\end{equation}

\subsubsection{The split zero-product intersection}\label{m1:subsec:split-zero-intersection}

It remains to treat \(Z=\pi_{\mathrm{sp}}^{-1}(0,0,0)\). Equivalently,
\[
    Z=
    \{y_1(x_1)=y_2(x_2)=0,
      \ y_1(x_2)y_2(x_1)=0\}.
\]
Put \(Z_+=\{\beta=0\}\) and \(Z_-=\{\alpha=0\}\). Then
\[
    Z=Z_+\cup Z_-,
    \qquad
    F_0:=Z_+\cap Z_-=\{YX=0\},
\]
and transpose exchanges $Z_+$ and $Z_-$.

\begin{lemma}\label{m1:lem:split-F0}
Every $\Gamma$-invariant distribution on $F_0$ is fixed by transpose:
\[
    \mathcal D(F_0)^\Gamma=\mathcal D(F_0)^{\Gamma,\tau}.
\]
\end{lemma}

\begin{proof}
The rank pairs in $F_0$ are
\((2,0),\allowbreak(1,0),\allowbreak(0,0),\allowbreak(0,1),\allowbreak(0,2),\allowbreak(1,1)\). Let
\(U_{11}=\{(X,Y)\in F_0: \operatorname{rank}X=\operatorname{rank}Y=1\}\). Then \(F_0=U_{11}\sqcup B\),
with $U_{11}$ open and $B$ closed.  Lemma~\ref{m1:lem:matrix-boundary} treats
$B$.

By Lemma~\ref{m1:lem:common-U11}, the descended fiber has coordinates
$c=(c_1,c_2)\neq0$ and $d=(d_1,d_2)^t\neq0$, action
\eqref{m1:eq:common-U11-action}, character $\eta=|\beta/\alpha|$, and
swap $(c,d)\mapsto(d^t,c^t)$. Its odd relative modulus is one.

On the locus where \((c_1d_1,c_2d_2)\neq(0,0)\), consider
\[
    \pi_{11,\mathrm{sp}}(c,d)
    =[c_1d_1:c_2d_2]\in\mathbf P^1(F).
\]
This is invariant under the descended action and fixed by transpose.
On a fiber with both projective coordinates nonzero, if two points have
the same projective value, then for some
$\kappa\in F^\times$ one has \(c_i'd_i'=\kappa c_id_i\) for \(i=1,2\). Choose $\alpha/\beta=\kappa$ and
then choose $t_i$ so that
$c_i'=\alpha c_i/t_i$; the equations for $d_i'$ follow automatically.
Thus every fiber with both projective coordinates nonzero is a single orbit.
Its stabilizer has $\alpha=\beta$ and $t_i=\alpha$ (with $b$ arbitrary),
so it is unimodular, and both $\Delta_{P\times T}$ and $\eta$ restrict
trivially. The positive equivariant orbital distribution is fixed by the
descended involution: its action is a positive scalar of square one.

For the special fiber $[1:0]$, normalize $c_1=d_1=1$.  The remaining
variables satisfy \(c_2d_2=0\), and the residual torus acts by the standard opposite-weight action
\((c_2,d_2)\longmapsto(qc_2,q^{-1}d_2)\). The residual stabilizer has $\alpha=\beta=t_1$, with $t_2$
and $b$
free. It is unimodular, $\Delta_{P\times T}$ restricts trivially, and the
descended character is trivial. The involution fixes $c_1=d_1=1$ and
swaps $(c_2,d_2)$. Thus Prasad's
Two-Axis Lemma gives transpose invariance on this fiber.  The fiber
over $[0:1]$ is identical after interchanging the two eigenlines.

The locus omitted from $\pi_{11,\mathrm{sp}}$ consists of the two
disjoint-support orbits represented by \((c,d)=(e_1^\vee,e_2)\) and \((c,d)=(e_2^\vee,e_1)\).
For the first orbit the stabilizer is
\[
    H_z=
    \left\{
    \left(
        \begin{pmatrix}\alpha&b\\0&\beta\end{pmatrix},
        \operatorname{diag}(\alpha,\beta)
    \right)
    \right\}.
\]
The second Frobenius relative modular character is trivial, while
$\eta|_{H_z}=|\beta/\alpha|$ is nontrivial. Indeed $H_z\simeq P$ and
$\Delta_{H_z}=\Delta_{P\times T}|_{H_z}$, so there is no additional modular
factor. Hence this orbit carries no
required equivariant distribution; the transposed orbit is identical.
First apply Bernstein localization on the open domain of
$\pi_{11,\mathrm{sp}}$, using its generic and special fibers.
An odd $\eta$-equivariant distribution then has support in the omitted
closed locus. Its two disjoint-support orbits are open and closed in
that locus and both have zero required distribution space by the
character calculation above. Frobenius descent therefore gives
\[
    \mathcal D(U_{11})^\Gamma=\mathcal D(U_{11})^{\Gamma,\tau}.
\]

Now let $D\in\mathcal D(F_0)^\Gamma$ satisfy $\tau D=-D$.  Its restriction to
$U_{11}$ vanishes, so $\operatorname{supp}D\subset B$.  By
Lemma~\ref{m1:lem:matrix-boundary}, there is no $\Gamma$-invariant anti-invariant
distribution supported on $B$.  Thus $D=0$.
\end{proof}

The following consequence of the same stratification will be used in the
gluing argument.

\begin{corollary}[Boundary anomaly on $F_0$]\label{m1:cor:boundary-anomaly}
Let $Q$ be a $\tau$-anti-invariant distribution supported on $F_0$.  Suppose
that \(\operatorname{supp}(gQ-Q)\subset\{0\}\) for \(g\in \Gamma\). Then there exists $C\in\mathbf C$ such that
\begin{equation}\label{m1:eq:boundary-anomaly}
    gQ-Q=C\log\lambda(g)\,\delta_0,
\end{equation}
where, for $t=\operatorname{diag}(t_1,t_2)$, \(\lambda(a,t) =\left|\frac{\det a}{t_1t_2}\right|\). \end{corollary}

\begin{proof}
On $U_{11}$ the defect vanishes, so $Q|_{U_{11}}$ is invariant and
anti-invariant. The proof of Lemma~\ref{m1:lem:split-F0} makes this restriction
zero. Thus $Q$ is supported on $B$. Its restrictions to the two rank-two
open axes are multiples of their invariant measures. With the measures
normalized by $P_Y=\tau P_X$, anti-invariance makes their coefficients
$a,-a$. After subtracting $a(P_X-P_Y)$, the remainder is supported
on the rank-at-most-one axes and still has defect supported at zero. Hence its
restrictions to the rank-one strata are invariant, and the first part of
Lemma~\ref{m1:lem:matrix-boundary} makes them zero. The remainder is supported
at zero and is therefore a multiple of $\delta_0$, with zero defect.
Consequently
\[
 gQ-Q=-2ac\log\lambda(g)\,\delta_0,
\]
as required.
\end{proof}

\subsubsection{Root-component gluing}\label{m1:subsec:root-gluing}

Let \(D\in\mathcal D(Z)^\Gamma\) with \(\tau D=-D\). We prove that $D=0$.  Put \(U_+=Z_+\setminus
F_0\) and \(U_-=Z_-\setminus F_0\). Transpose exchanges $U_+$ and $U_-$.

\emph{Rank-mismatch orbits in $U_+$.}
On $U_+$ one has \(\beta=y_2(x_1)=0\) and \(\alpha=y_1(x_2)\neq0\). Consider
\[
    p:U_+\longrightarrow
    \Omega,
    \qquad
    p(X,Y)=(x_2,y_1),
\]
where \(\Omega=\{(x,\ell)\in V\times V^\vee:\ell(x)\neq0\}\). The group $\Gamma$ acts transitively on
$\Omega$. We apply Frobenius descent at
$(e_1,e_1^\vee)$ and write \(x_1=ae_1+de_2\) and \(y_2=be_1^\vee+ee_2^\vee\). The equations
defining $Z$ give \(a=b=0\) and \(de=0\). Hence the Frobenius fiber is \(N=\{(d,e)\in
F^2:de=0\}\). Its stabilizer is
\[
    H=
    \left\{
    \left(\operatorname{diag}(s,c),sI_2\right):
    s,c\in F^\times
    \right\},
\]
and the effective quotient acts by \(q\cdot(d,e)=(qd,q^{-1}e)\), where \(q=c/s\). The Frobenius
relative modular character is trivial.  The three strata of
$N$ correspond to
\[
\begin{array}{c|c}
(d,e) & (\operatorname{rank}X,\operatorname{rank}Y)\\ \hline
 d\neq0,\ e=0 & (2,1),\\
 d=0,\ e\neq0 & (1,2),\\
 d=e=0 & (1,1).
\end{array}
\]
The origin represents the rank-$(1,1)$ orbit in $U_+$.

Let $c_{21}^+$ and $c_{12}^+$ be the coefficients of the invariant orbital
distributions on the two punctured branches, measured relative to any
fixed positive invariant measures on the global orbits. Frobenius descent
with positive Haar measures sends positive measures to positive measures.
Thus Corollary~\ref{m1:cor:two-axis-coefficients} gives a constant $r_1>0$ with
\begin{equation}\label{m1:eq:Uplus-branch}
    c_{21}^+=r_1c_{12}^+.
\end{equation}

\emph{Comparison of the two root components.}
Consider the $\Gamma$-invariant open subset \(U_X=\{(X,Y)\in Z:\det X\neq0\}\). The projection
$U_X\to\GL_2(F)$, $(X,Y)\mapsto X$, has homogeneous base.
Frobenius descent allows us to fix $X=I_2$.  The stabilizer is the diagonal
copy of $T$, and the fiber is
\[
    N_X=
    \left\{
    \begin{pmatrix}0&\alpha\\\beta&0\end{pmatrix}:
    \alpha\beta=0
    \right\}.
\]
The effective torus acts by \((\alpha,\beta)\longmapsto(q\alpha,q^{-1}\beta)\). The two punctured
branches correspond to the rank-$(2,1)$ orbits in $U_+$
and $U_-$, respectively. The same positivity argument applied to this
chart gives a constant $r_2>0$ with
\begin{equation}\label{m1:eq:Xchart-branch}
 c_{21}^+=r_2c_{21}^-.
\end{equation}
Transpose sends the rank-$(2,1)$ orbit in $U_-$ to the
rank-$(1,2)$ orbit in $U_+$. Its pushforward preserves positivity, so anti-invariance gives
some $r_3>0$ with
\begin{equation}\label{m1:eq:transpose-branch}
 c_{21}^-=-r_3c_{12}^+.
\end{equation}
Combining \eqref{m1:eq:Uplus-branch}, \eqref{m1:eq:Xchart-branch}, and
\eqref{m1:eq:transpose-branch} gives $(r_1+r_2r_3)c_{12}^+=0$. Since all three constants are positive, this
implies, also using transpose for the fourth orbit,
\[
    c_{21}^+=c_{12}^+=c_{21}^-=c_{12}^-=0.
\]
Thus all rank-mismatch contributions on $Z\setminus F_0$ vanish.

\emph{The remaining rank-one coefficient.}
Define the closed linear root components \(S_+=\{x_1=0,\ y_2=0\}\) and \(S_-=\{x_2=0,\
y_1=0\}\). They are exchanged by transpose and satisfy \(S_+\cap S_-=\{0\}\). The open set $S_+\cap U_+$
is precisely the rank-$(1,1)$ orbit in $U_+$,
and $S_-\cap U_-$ is its transpose.

On $S_+$ write a point as $(u,\varphi)\in V\times V^\vee$ and put
\[
    q_+(u,\varphi)=\varphi(u),
    \qquad
    \rho(a,t)=\left|\frac{t_1}{t_2}\right|,
    \quad
    t=\operatorname{diag}(t_1,t_2).
\]
For $\Re(s)\gg0$, define
\[
    Z_+(s)(f)=
    \int_{S_+}f(u,\varphi)|\varphi(u)|^s\,du\,d\varphi,
\]
and define \(Z_-(s)=\tau Z_+(s)\). We view these distributions on $Z$ through the closed embeddings
$S_\pm\hookrightarrow Z$.

The determinant family in Lemma~\ref{m1:lem:matrix-boundary} also supplies
the analytic properties needed here. Indeed, the linear isomorphism
\begin{equation}\label{m1:eq:determinant-pairing-isomorphism}
 J(u,\varphi)=\begin{pmatrix}u_1&u_2\\-\varphi_2&\varphi_1\end{pmatrix},
 \qquad \det J(u,\varphi)=\varphi(u),
\end{equation}
preserves normalized product Haar measure. It therefore identifies the
test-function integrals for $Z_+(s)$ with determinant zeta integrals.
Consequently $Z_+(s)$ has a simple pole at $s=-2$, with residue
$R=c\delta_0$ for the same $c\ne0$, and is regular there away from the
origin. Since $Z_-(s)=\tau Z_+(s)$, its residue is also $R$.
This identifies the meromorphic families as distributions; the
transformation character for the action on $S_+$ is computed next.

The action of $g=(a,t)\in \Gamma$ on $S_+$ is
\[
    (u,\varphi)
    \longmapsto
    \bigl(t_2^{-1}au,\,t_1\varphi a^{-1}\bigr).
\]
Thus $q_+$ is multiplied by $t_1/t_2$ and the linear map has
absolute determinant $\rho(g)^2$. Under the pushforward convention
the product Haar measure therefore transforms by $\rho(g)^{-2}$.  Consequently
\[
    gZ_+(s)=\rho(g)^{-(s+2)}Z_+(s).
\]
Write
\[
    Z_+(s)=\frac{R}{s+2}+P_++O(s+2).
\]
Then
\begin{equation}\label{m1:eq:Pplus-defect}
    gP_+-P_+=-\log\rho(g)\,R.
\end{equation}
Choose $P_-=\tau P_+$.  Since conjugation by $\tau$ sends $(a,t)$ to $(a^{-t},t^{-1})$ and therefore
inverts $\rho$,
\[
    gP_--P_-=+\log\rho(g)\,R.
\]
Hence
\begin{equation}\label{m1:eq:Ppm-defect}
    g(P_+-P_-)-(P_+-P_-)
    =-2\log\rho(g)\,R.
\end{equation}

Let \(O_{11}^+=S_+\cap U_+\) and \(O_{11}^-=S_-\cap U_-\). On $O_{11}^+$ one has $q_+\neq0$, so
$Z_+(s)$ is regular there at $s=-2$.
Thus \(\nu_+:=P_+|_{O_{11}^+}\) is a nonzero invariant orbital distribution.  Put \(\nu_-:=\tau\nu_+\);
then $\nu_-=P_-|_{O_{11}^-}$.

Since all rank-mismatch contributions vanish, there exists
$a\in\mathbf C$ such that \(D|_{U_+}=a\nu_+\) and \(D|_{U_-}=-a\nu_-\). Define
\(Q=D-a(P_+-P_-)\). Then \(\operatorname{supp}Q\subset F_0\). Since $D$ is $\Gamma$-invariant,
\eqref{m1:eq:Ppm-defect} gives
\begin{equation}\label{m1:eq:Q-rho-defect}
    gQ-Q=2a\log\rho(g)\,R
    =2ac\log\rho(g)\,\delta_0.
\end{equation}
On the other hand, Corollary~\ref{m1:cor:boundary-anomaly} gives
\[
    gQ-Q=C\log\lambda(g)\,\delta_0
\]
for some $C\in\mathbf C$.  Thus
\[
    2ac\log\rho(g)=C\log\lambda(g)
    \qquad(g\in \Gamma).
\]
The characters $\rho$ and $\lambda$ are independent.  Indeed, let $\varpi$
be a uniformizer of $F$ and take
\[
    g=\bigl(I_2,\operatorname{diag}(\varpi,\varpi^{-1})\bigr).
\]
Then
\[
    \lambda(g)=1,
    \qquad
    \rho(g)=|\varpi|^2\neq1.
\]
Since $c\neq0$, it follows that $a=0$.

Therefore $D|_{Z\setminus F_0}=0$, so $\operatorname{supp}D\subset F_0$.
By Lemma~\ref{m1:lem:split-F0}, every $\Gamma$-invariant distribution supported on
$F_0$ is fixed by transpose.  Since $\tau D=-D$, we conclude that $D=0$.
Hence
\begin{equation}\label{m1:eq:split-zero-fiber}
    \mathcal D(Z)^{\widetilde \Gamma,\varepsilon}=0.
\end{equation}

By \eqref{m1:eq:split-nonzero-fibers} and \eqref{m1:eq:split-zero-fiber}, every
fiber of $\pi_{\mathrm{sp}}$ has zero $(\widetilde \Gamma,\varepsilon)$-equivariant
distribution space.  The map $\pi_{\mathrm{sp}}$ is continuous and
$\widetilde \Gamma$ preserves every fiber.  Bernstein localization
(Proposition~\ref{m1:prp:localization}) therefore gives
\[
    \mathcal D(\Vlegs)^{\widetilde \Gamma,\varepsilon}=0.
\]
Together with the nonsplit calculation, this proves
Theorem~\ref{m1:thm:two-leg}.

\subsection{The linear identity descendant}\label{m1:sec:identity}
\subsubsection{The representation and nonzero quadratic fibers}
\label{m1:subsec:identity-fibers}
Let $\Gamma=\GL_2(F)\times\GL_2(F)$ act on
$\mathcal N=M_2(F)^2\oplus\mathfrak{gl}_2(F)$ by
\begin{equation}\label{m1:eq:identity-action}
 (a,b)(X,Y,Z)=(aXb^{-1},bYa^{-1},bZb^{-1}),\qquad
 \tau(X,Y,Z)=(Y^t,X^t,Z^t).
\end{equation}
The involution normalizes $\Gamma$ by $(a,b)\mapsto(a^{-t},b^{-t})$;
write $\widetilde\Gamma=\Gamma\rtimes\langle\tau\rangle$ and use the
sign character $\varepsilon$.

\begin{theorem}[Identity Descendant]\label{m1:thm:identity-descendant}
Every $\Gamma$-invariant distribution on $\mathcal N$ is fixed by $\tau$; equivalently,
$\mathcal D(\mathcal N)^{\widetilde\Gamma,\varepsilon}=0$.
\end{theorem}
The proof first localizes the sign space to a quadratic null cone. The
same localization applies after partial Fourier transform, giving the
second support condition needed for odd-dimensional homogeneity.

Since $F$ has characteristic zero, write
$Z=zI+Z_0$, where $z=\tr Z/2$ and $Z_0\in W:=\mathfrak{sl}_2(F)$.
Thus $\mathcal N=M_2(F)^2\times F\times W$. The scalar $z$ is fixed
by both actions, whereas the legs and traceless variable retain their
actions. On $W$ use
\[
 B(U,V)=\tr(UV),\qquad q(U)=B(U,U).
\]
The symmetric trace pairing is nondegenerate: for
$U=\left(\begin{smallmatrix}a&b\\c&-a\end{smallmatrix}\right)$ it pairs
as $2aa'+bc'+cb'$. Both $B$ and $q$ are invariant under conjugation
and transpose.

\begin{proposition}\label{m1:prop:identity-nonzero-fibers}
For every $z\in F$ and $c\in F^\times$, the fiber
\[
 \mathcal N_{z,c}=M_2(F)^2\times\{z\}\times\mathcal O_c,
 \qquad\mathcal O_c=\{Z_0\in W:q(Z_0)=c\},
\]
has zero sign-equivariant distribution space.
\end{proposition}
\begin{proof}
Every $Z_0\in\mathcal O_c$ has characteristic polynomial $T^2-c/2$.
For $c\neq0$ this polynomial is separable, so $Z_0$ is regular
semisimple. Rational canonical form shows that $\mathcal O_c$ is one
rational $\GL_2(F)$-orbit. Its centralizer
$T_{Z_0}=F[Z_0]^\times$ is the split or nonsplit maximal torus according
to this polynomial. Frobenius descent to a chosen $Z_0$ leaves the
two-leg representation
\[
 \GL_2(F)\times T_{Z_0}\curvearrowright M_2(F)^2,
 \qquad(a,t)(X,Y)=(aXt^{-1},tYa^{-1}).
\]
Choose the rational symmetric form $S$ from
Section~\ref{m1:subsec:two-leg-statement}: the standard form in a split
eigenbasis or the trace pairing on $F[Z_0]$ in the nonsplit case.
Then $SZ_0=Z_0^tS$ and $St=t^tS$ for all $t\in T_{Z_0}$. The
correcting element $(1,S^{-1})$ sends $Z_0^t$ to $Z_0$, and the induced
operator on the fiber is
\[
 \tau_{Z_0}(X,Y)=(Y^tS,S^{-1}X^t)=\tau_S(X,Y).
\]
By \eqref{m1:eq:two-leg-adjoint}, it has square one and normalizes the
two-leg action by $(a,t)\mapsto(a^{-t},t^{-1})$. Both the ambient ordinary group and
$\GL_2(F)\times T_{Z_0}$ are unimodular. Their finite extensions are
unimodular, so the relative modulus is trivial on the odd component
as well. The descended character is the sign character, and
Theorem~\ref{m1:thm:two-leg} gives the assertion.
\end{proof}

\begin{corollary}\label{m1:cor:identity-support}
Every $T\in\mathcal D(\mathcal N)^{\widetilde\Gamma,\varepsilon}$
satisfies
\begin{equation}\label{m1:eq:identity-first-support}
 \supp T\subset M_2(F)^2\times F\times Z(B),
 \qquad Z(B)=\{U\in W:B(U,U)=0\}.
\end{equation}
\end{corollary}
\begin{proof}
On the invariant open set $q(Z_0)\neq0$, the continuous invariant map
$(X,Y,z,Z_0)\mapsto(z,q(Z_0))$ has target $F\times F^\times$.
Proposition~\ref{m1:prop:identity-nonzero-fibers} treats every
fiber, so Bernstein localization makes the restriction of $T$ to this
open set zero. The closed complement gives the stated support.
\end{proof}
For traceless $2\times2$ matrices, $q(U)=-2\det U$. Thus $Z(B)$ is
the nilpotent cone in $W$.

\subsubsection{Fourier equivariance and homogeneity}
\label{m1:subsec:identity-Fourier}
Let $dV$ be the self-dual Haar measure on $W$ for $B$ and the fixed
additive character $\psi$. For
$f\in C_c^\infty(M_2(F)^2\times F\times W)$, set
\begin{equation}\label{m1:eq:partial-Fourier}
 (\mathcal F_W f)(X,Y,z,U)
 =\int_W f(X,Y,z,V)\psi(B(U,V))\,dV,
 \qquad (\mathcal F_WT)(f)=T(\mathcal F_Wf).
\end{equation}

\begin{lemma}\label{m1:lem:identity-Fourier}
Partial Fourier transform preserves
$\mathcal D(\mathcal N)^{\widetilde\Gamma,\varepsilon}$. In particular
every distribution in this space also satisfies
\begin{equation}\label{m1:eq:identity-second-support}
 \supp(\mathcal F_WT)\subset M_2(F)^2\times F\times Z(B).
\end{equation}
\end{lemma}
\begin{proof}
Conjugation and transpose preserve the pairing on $W$:
\[
 \tr(gUg^{-1}gVg^{-1})=\tr(UV),\qquad
 \tr(U^tV^t)=\tr(VU)=\tr(UV).
\]
An isometry of a nondegenerate symmetric form has determinant of
absolute value one and preserves the self-dual measure. Changing
variables in \eqref{m1:eq:partial-Fourier} therefore proves that these
operators commute with $\mathcal F_W$ on functions and distributions.
The action on the remaining variables is independent of the $W$
integration. Thus $\mathcal F_W$ commutes with the extended group action.
Corollary~\ref{m1:cor:identity-support}, applied to $\mathcal F_WT$, now gives
\eqref{m1:eq:identity-second-support}.
\end{proof}
\begin{proof}[Proof of Theorem~\ref{m1:thm:identity-descendant}]
We apply the common homogeneity theorem, Theorem~\ref{m2:tools:homogeneity}.

Here $W=\mathfrak{sl}_2(F)$ with the nondegenerate trace form and
$M=M_2(F)^2\times F$, an affine analytic manifold. Equations~\eqref{m1:eq:identity-first-support} and
\eqref{m1:eq:identity-second-support} give both support hypotheses, with
the Fourier convention \eqref{m1:eq:partial-Fourier}. Since $\dim_F W=3$,
the theorem gives $T=0$.
\end{proof}

\subsection{Invariance of distributions on the linear normal spaces}\label{m1:sec:descendants}
\label{m1:subsec:dense-descendants}\label{m1:subsec:rankzero-descendants}
\label{m1:subsec:rankone-descendants}\label{m1:subsec:exhaustive-descendants}

\begin{proposition}\label{m1:prop:all-descendants}
\label{m1:thm:dense-descendants}\label{m1:thm:rank00-descendants}
\label{m1:thm:decomposable-descendant}\label{m1:thm:stable-system-descendant}
For every closed rational orbit $E\in\cX(F)$, let
$N_E=T_E\cX/T_E(LE)$, and let $\widetilde L_E$ be the extended
stabilizer with its restricted sign character. Then
\[
 \mathcal D(N_E)^{\widetilde L_E,\varepsilon_E}=0.
\]
\end{proposition}
\begin{proof}
Theorem~\ref{m1:thm:closed-orbits} and
Section~\ref{m1:subsec:normal-interface} give the nine normal
representations with their corrected odd operators. Throughout,
Lemma~\ref{m1:lem:passive-products} allows us to retain the summands
fixed by the extended stabilizer as parameters.

\emph{Dense orbits.}
For the absolutely irreducible type the odd operator is the identity.
The split toral type is an opposite-weight plane with rationally
surjective character $t_1/t_2$, so
Proposition~\ref{m1:prop:opposite-weight-plane} applies.
For the nonsplit type, use the rational root module
\eqref{m1:eq:dense-nonsplit-quotient} with action
\eqref{m1:eq:dense-nonsplit-action}. Hilbert~90 identifies the effective
action of $K^\times/F^\times$ with all of $K^1$. Thus each nonzero
norm fiber is a single $K^\times$-orbit with stabilizer $F^\times$;
both groups are unimodular. The odd map is $d\mapsto-\bar d$. Since $-1\in K^1$,
it acts on invariant distributions as conjugation does, and
Lemma~\ref{m2:elem:norm} applies to the root module.
For the scalar pair, Lemma~\ref{m1:lem:matrix-pair-transpose} applies
to the normal space $M_2(F)^2$.

\emph{Rank-zero orbits.}
For regular $D$, split or nonsplit,
Section~\ref{m1:subsec:rankzero-normal} identifies the nontrivial part
of the normal representation and its odd operator with the adjoint
Two-Leg representation; the two
centralizer directions are fixed. Theorem~\ref{m1:thm:two-leg}
therefore applies. For $D=dI$, $d\ne0$, the action
\eqref{m1:eq:rankzero-scalar} is the
Identity Descendant of Theorem~\ref{m1:thm:identity-descendant}.

\emph{Rank-one orbits.}
For a decomposable rank-one orbit, the normal-space formulas
\eqref{m1:eq:decomposable-normal}--\eqref{m1:eq:decomposable-transpose}
show that the characters $(\alpha,\beta)$ are jointly rationally surjective.
The invariant map on the root coordinates is
\begin{equation}\label{m1:eq:decomposable-fiber-map}
 \pi(x,y,u,v)=(xy,uv),\qquad
 \pi^{-1}(c_1,c_2)=\{xy=c_1\}\times\{uv=c_2\}.
\end{equation}
Each factor is a hyperbola or the two-axis node. Lemma~\ref{m1:lem:passive-products} and Bernstein
localization exclude sign-equivariant distributions on the normal space. For a stable-system rank-one orbit, the stabilizer character $\eta$
is rationally surjective
by \eqref{m1:eq:stable-stabilizer};
\eqref{m1:eq:stable-ext-table} and
Lemma~\ref{m1:lem:stable-parameter-basis} give one opposite-weight
plane with an ordinary swap and four fixed directions.
Proposition~\ref{m1:prop:opposite-weight-plane} therefore applies
throughout $b_1=1-ad+sb_0$, $\det H\ne0$, including nilpotent and
repeated-root values of $D$.

These cases exhaust the closed rational orbits. Conjugation transports
the normal representations and extended stabilizers from the chosen
representatives to every point of each rational orbit.
\end{proof}

\subsection{Proof of the linear multiplicity-one theorem}\label{m1:sec:completion}
\begin{proof}[Proof of Theorem~\ref{m1:thm:main}]
The residual variety $\cX\simeq\GL_4$ is smooth affine and the action
is algebraic. Corollary~\ref{m1:cor:transpose-stable} gives rational
transpose stability of its closed orbits, and
Proposition~\ref{m1:prop:all-descendants} excludes sign-equivariant
distributions on their normal spaces. Corollary~\ref{m1:cor:HC-invariance} therefore yields
\[
 \mathcal D(\cX)^{\widetilde L,\varepsilon}=0.
\]
Thus every residual invariant distribution is transpose-invariant.

The unimodular, involution-compatible Frobenius isomorphism
\eqref{m1:eq:residual-Frobenius} carries this conclusion back to all
$H_{\mathrm L}\times H_{\mathrm L}$-invariant distributions on
$G_{\mathrm L}=\GL_4(F)\times\GL_2(F)$. Proposition~\ref{m1:prop:GK} gives, for every
irreducible admissible smooth complex representation $\pi$,
\begin{equation}\label{m1:eq:final-product-bound}
 m_{H_{\mathrm L}}(\pi)m_{H_{\mathrm L}}(\pi^\vee)\leq1,
 \qquad m_{H_{\mathrm L}}(\pi)=\dim\Hom_{H_{\mathrm L}}(\pi,\Co).
\end{equation}

To identify the two multiplicities, let
$\iota_{G_{\mathrm L}}(g,z)=(g^{-t},z^{-t})$. The contragredient theorem for general
linear groups~\cite[Theorem~7.3]{BZ76}, applied to each factor, gives
$\pi^\vee\simeq\pi\circ\iota_{G_{\mathrm L}}$. Moreover
\[
 \iota_{G_{\mathrm L}}(\diag(a,b),b)=(\diag(a^{-t},b^{-t}),b^{-t}),
\]
so $\iota_{G_{\mathrm L}}$ preserves the embedded $H_{\mathrm L}$ and restricts to an automorphism
of it. The $H_{\mathrm L}$-invariant linear forms on $\pi\circ\iota_{G_{\mathrm L}}$ are exactly
those on $\pi$, and hence
\begin{equation}\label{m1:eq:final-multiplicity-equality}
 m_{H_{\mathrm L}}(\pi^\vee)=m_{H_{\mathrm L}}(\pi).
\end{equation}
Combining \eqref{m1:eq:final-product-bound} and
\eqref{m1:eq:final-multiplicity-equality} gives $m_{H_{\mathrm L}}(\pi)^2\leq1$, and
therefore $m_{H_{\mathrm L}}(\pi)\leq1$.
\end{proof}

\section{Unitary rational geometry and normal representations}\label{sec:unitary-geometry}
We use the determinant-open quiver of Section~\ref{sec:quiver}
to classify the closed unitary orbits after scalar extension.
We then determine which of these geometric orbits have $F$-rational
points, establish rational orbit uniqueness and transpose stability,
and compute the corresponding normal representations over $F$.
These results provide the input for the distributional arguments in
Section~\ref{sec:unitary-descendants}.

\subsection{The unitary residual action}\label{m2:sec:reduction}

\subsubsection{The group, subgroup, and coordinates}\label{m2:red:pair}

Let $J_{2m}$ be the anti-diagonal matrix with all anti-diagonal
entries equal to one. Write
\[
 G_{2m}(F)=\{g\in\GL_{2m}(E):g^\dagger J_{2m}g
       =\lambda(g)J_{2m},\ \lambda(g)\in F^\times\},
 \qquad g^\dagger=\bar g^{\,t}.
\]
These are the quasi-split unitary similitude groups, with
$G_2=\GU_{1,1}$ and $G_4=\GU_{2,2}$. Put $J=J_2$ and
$\U_2=U_{1,1}=\ker(\lambda:G_2\to\Gm)$ and
$X^*=JX^\dagger J$ for a $2$ by $2$ matrix. The defining equation
of $G_2$ is $h^*h=\lambda(h)I$. We also write
$E^-=\{z\in E:\bar z=-z\}$.
For the Lie algebras, write $\mathfrak{gu}_n=\Lie(G_n)$ and
$\mathfrak u_n=\ker(d\lam:\mathfrak{gu}_n\to F)$.

For $u=(u_{ij})\in G_2$ and $h\in G_2$, define
\[
 \rho(h,u)=\begin{pmatrix}
 u_{11}&0&0&u_{12}\\
 0&h_{11}&h_{12}&0\\
 0&h_{21}&h_{22}&0\\
 u_{21}&0&0&u_{22}
 \end{pmatrix}.
\]
In the notation of the introduction, the unitary pair is
\begin{equation}\label{m2:red:embedding}
 \begin{gathered}
 G_{\mathrm U}=G_4(F)\times G_2(F),\qquad
 H_{\mathrm{eq}}=\{(h,u)\in G_2^2:\lambda(h)=\lambda(u)\},\\
 j(h,u)=(\rho(h,u),h),\qquad H_{\mathrm U}=j(H_{\mathrm{eq}}(F)).
 \end{gathered}
\end{equation}
Order the coordinates as $(1,4,2,3)$, so the outer and middle
planes have form $J\oplus J$ and $\rho(h,u)=u\oplus h$.
This rational permutation preserves ordinary transpose. All subsequent
blocks $x=\left(\begin{smallmatrix}A&B\\C&D\end{smallmatrix}\right)$
use this order; the first plane is the outer plane.

The multiplier map on $G_2(F)$ is surjective: set
\begin{equation}\label{m2:red:amul}
 a_\mu=\diag(1,\mu),\qquad a_\mu^*a_\mu=\mu I
 \quad(\mu\in F^\times).
\end{equation}
The kernel of $H_{\mathrm{eq}}\to\Gm$ is $\U_2^2$.
Thus $H_{\mathrm{eq}}$, like $G_4$ and $G_2$, is connected reductive.
In particular $G_{\mathrm U}$ and $H_{\mathrm U}$ are unimodular,
and $G_{\mathrm U}/H_{\mathrm U}$ has the invariant measure required
in \cref{m2:tools:gk}.

\begin{lemma}\label{m2:red:transpose}
Ordinary transpose preserves $G_2$, $G_4$, and their multipliers.
The anti-involution $\tau(g,z)=(g^t,z^t)$ of $G_{\mathrm U}$ preserves $H_{\mathrm U}$.
\end{lemma}

\begin{proof}
For $J'=J$ or $J\oplus J$, the equation
$g^\dagger J'g=\lambda J'$ implies $gJ'g^\dagger=\lambda J'$.
Conjugating the latter gives $\bar gJ'g^t=\lambda J'$, the defining
equation for $g^t$. Also $\rho(h,u)^t=\rho(h^t,u^t)$, with unchanged
multipliers.
\end{proof}

\subsubsection{The homogeneous base and its fiber}

Let $L=H_{\mathrm{eq}}^2$ act by left and right translation through $j$.
Projection $p(g,z)=z$ has transitive base action
$((h,u),(h',v))z=hzh'^{-1}$, with the continuous section
$z\mapsto((z,a_{\lambda(z)}),(1,1))$ of the orbit map at $1$.
Its fiber is $G_4(F)$; identifying the stabilizer $h=h'$ with
$(u,v,h)$ gives
\begin{equation}\label{m2:red:K}
 K=\{(u,v,h)\in G_2^3:
         \lambda(u)=\lambda(v)=\lambda(h)\}.
\end{equation}
Its action on the fiber is
\begin{equation}\label{m2:red:action}
 (u,v,h)\cdot x=(u\oplus h)x(v\oplus h)^{-1},
 \quad
 (A,B,C,D)\longmapsto
 (uAv^{-1},uBh^{-1},hCv^{-1},hDh^{-1}).
\end{equation}
The multiplier $\lambda(x)$ is unchanged. The kernel of
$K\to\Gm$ is $\U_2^3$, so $K$ is also connected reductive.

\begin{proposition}\label{m2:red:frobenius}
There is a natural isomorphism, compatible with transpose,
\begin{equation}\label{m2:red:distiso}
 \Ddist(G_{\mathrm U})^{H_{\mathrm U}\times H_{\mathrm U}}
   \simeq\Ddist(G_4(F))^{K(F)}.
\end{equation}
\end{proposition}

\begin{proof}
The continuous section constructed above gives the quotient topology on $G_2(F)$.
Since $L$ and $K$ are reductive, their relative modulus is one;
\cref{m2:tools:frobenius} gives \eqref{m2:red:distiso}.
Transpose fixes the base point $1$ and acts on its stabilizer by
\begin{equation}\label{m2:red:alpha}
 \alpha_K(u,v,h)=(v^{-t},u^{-t},h^{-t}),
 \qquad ((u,v,h)\cdot x)^t=\alpha_K(u,v,h)\cdot x^t.
\end{equation}
Its pushforward preserves the positive Haar ray on the base; the
scalar is one because transpose squares to one. Thus Frobenius descent
intertwines total and fiber transpose, as in
\cref{m1:prop:odd-Frobenius}.
\end{proof}

Set
\begin{equation}\label{m2:red:extension}
 \Kt=K\rtimes_{\alpha_K}C_2,
 \qquad \sigma x=x^t,\qquad
 \eps|_K=1,\quad\eps(\sigma)=-1.
\end{equation}
We shall prove
\begin{equation}\label{m2:red:target}
 \Ddist(G_4(F))^{\Kt(F),\eps}=0.
\end{equation}
Taking the odd part of a $K(F)$-invariant distribution shows that
\eqref{m2:red:target} and \cref{m2:red:frobenius} imply the
bi-$H_{\mathrm U}$ distribution hypothesis of \cref{m2:tools:gk}.

\subsubsection{Contragredient duality}

Put $\iota(g,z)=(g^{-t},z^{-t})$, the automorphism associated with $\tau$.

\begin{proposition}\label{m2:red:duality}
For every irreducible admissible smooth complex representation
$\pi$ of $G_{\mathrm U}$, $\pi\circ\iota\simeq\pi^\vee$. Moreover
$\iota(H_{\mathrm U})=H_{\mathrm U}$, and hence
\[
 \dim\Hom_{H_{\mathrm U}}(\pi,\mathbb C)
 =\dim\Hom_{H_{\mathrm U}}(\pi^\vee,\mathbb C).
\]
\end{proposition}

\begin{proof}
Coordinate conjugation $c(v)=\bar v$ is an $F$-linear anti-unitary
involution for our rational symmetric form matrices. Roche--Vinroot's
Theorem~B and its character argument in \S7.2~\cite{RocheVinroot18}
therefore make $\theta_m(g)=\lambda(g)^{-1}\bar g$ a dualizing
involution on $G_{2m}(F)$: $\theta_m(g)$ is conjugate to $g^{-1}$.
The unitary equation gives $g^{-t}=J'\theta_m(g)J'$, where
$J'=J$ or $J\oplus J$ has multiplier one. Thus inverse transpose
is an inner modification of this dualizing involution. The conjugacy
and character argument applies componentwise to the reductive product,
giving $\pi\circ\iota\simeq\pi^\vee$.

Finally $\rho(h,u)^{-t}=\rho(h^{-t},u^{-t})$ and
$\lambda(h^{-t})=\lambda(h)^{-1}$. Hence $\iota$ preserves
$H_{\mathrm U}$, so the invariant linear forms on $\pi\circ\iota$
are exactly those on $\pi$. The contragredient isomorphism gives the
equality of multiplicities.
\end{proof}

\subsubsection{Rational linear algebra of the hyperbolic plane}
\label{m2:red:unitary-pair}

The following rational coordinates will give the transporters and
normal representations for dense orbits.
Choose $\delta\in E^\times$ with $\bar\delta=-\delta$, and set
\[
 d=\delta^2\in F^\times,\qquad P=\diag(1,\delta),\qquad
 J_0=\begin{pmatrix}0&1\\-1&0\end{pmatrix},\qquad
 R=P^{-2}=\diag(1,d^{-1}).
\]
Since $P^\dagger JP=\delta J_0$, the unitary equation becomes an
alternating-form equation. This rational factorization is distinct from
the scalar-extension comparison in \cref{m2:geo:split-bridge}.

\begin{lemma}\label{m2:simconj:factorization}
Every $h\in G_2(F)$ has the form
\begin{equation}\label{m2:simconj:factor}
 h=zPgP^{-1},\qquad z\in E^\times,\quad g\in\GL_2(F),
 \qquad \lambda(h)=N_{E/F}(z)\det g.
\end{equation}
The overlap of two such factorizations is precisely
$(z,g)\sim(z/a,ag)$ for $a\in F^\times$.
Writing $P^{-1}YP=A+\delta B$ with $A,B\in M_2(F)$ defines
an $F$-linear isomorphism
\begin{equation}\label{m2:simconj:Psi}
 \Psi:M_2(E)\longrightarrow M_2(F)^2,
 \quad Y\longmapsto(A,B),\qquad
 \Psi(hYh^{-1})=(gAg^{-1},gBg^{-1}).
\end{equation}
\end{lemma}

\begin{proof}
The identity $g^tJ_0g=(\det g)J_0$ for $g\in\GL_2(F)$ and
$P^\dagger JP=\delta J_0$ show that $PgP^{-1}$ has multiplier
$\det g$; multiplying by $z$ contributes $N(z)$.
Conversely, for $H_0=P^{-1}hP$, comparison of
\[
 \bar H_0^{\,t}J_0H_0=\lambda(h)J_0,
 \qquad H_0^tJ_0H_0=(\det H_0)J_0
\]
gives $\bar H_0=cH_0$, with $c=\lambda(h)/\det H_0$ and
$c\bar c=1$. Hilbert~90 supplies $z$ with $c=\bar z/z$;
then $g=z^{-1}H_0$ is rational. If $zg=z'g'$ with rational
$g,g'$, then $g'g^{-1}=(z/z')I$ forces $z/z'\in F^\times$,
proving the exact overlap.

The $1,\delta$ decomposition is unique, and conjugation by
$h=zPgP^{-1}$ acts on its coefficients by simultaneous $g$-conjugation.
Every $g\in\GL_2(F)$ lifts through $PgP^{-1}$, so $\Psi$
identifies the invariant distribution spaces. The scalar centers
act trivially.
\end{proof}

Since $P^t=P$, the odd coordinate formula is
\begin{align}
 P^{-1}Y^tP
 &=P^{-2}(P^{-1}YP)^tP^2
   =R(A^t+\delta B^t)R^{-1},\notag\\
 \Psi(Y^t)&=(RA^tR^{-1},RB^tR^{-1}).
 \label{m2:simconj:transpose}
\end{align}
Here $R$ belongs to the effective even group $\GL_2(F)$.
The corresponding adjoint identities are as follows.

\begin{lemma}\label{m2:rat:real-pair}
Write $P^{-1}DP=A+\delta B_0$, with
$A,B_0\in\Mat_2(F)$.  Then
\begin{align}
 P^{-1}D^*P&=\operatorname{adj}(A)
                      -\delta\operatorname{adj}(B_0),
               \label{m2:rat:real-adjoint}\\
 P^{-1}D^tP&=R(A^t+\delta B_0^t)R^{-1}.
               \label{m2:rat:real-transpose}
\end{align}
Under this coordinate change
$E\langle D,D^*\rangle=E\langle A,B_0\rangle$.
Every $g\in\GL_2(F)$ gives $PgP^{-1}\in G_2(F)$
with multiplier $\det g$.
\end{lemma}

\begin{proof}
The transported adjoint is $Z\mapsto J_0^{-1}\bar Z^tJ_0$.
The identity
$J_0^{-1}X^tJ_0=\operatorname{adj}(X)=\tr(X)I-X$, linear for
two-by-two matrices, gives \eqref{m2:rat:real-adjoint};
\eqref{m2:simconj:transpose} gives \eqref{m2:rat:real-transpose}.
Adding and subtracting $A+\delta B_0$ and its adjoint, and removing
the scalar trace terms, recovers $A,B_0$ over $E$. This proves the
algebra equality. The multiplier assertion is
\cref{m2:simconj:factorization} with $z=1$.
\end{proof}

\subsection{Geometric orbits and rational existence}\label{m2:geo:section}

We transfer the geometric classification from Section~\ref{sec:quiver},
then determine which parameters admit points in the prescribed Hermitian spaces.

\subsubsection{The block action and the unitary equations}

We use the block coordinates of \cref{m2:red:pair}, with
$\mathcal J=J\oplus J$ and $\mu=\lambda(x)$.
The residual action is \eqref{m2:red:action}; it preserves $\mu$
because its two translation factors have the same multiplier.

\begin{proposition}\label{m2:geo:block-equations}
An invertible block matrix $x$ belongs to $G_4(F)$ with multiplier
$\mu\in F^\times$ if and only if
\begin{align}
 A^\dagger JA+C^\dagger JC&=\mu J,&
 A^\dagger JB+C^\dagger JD&=0,\label{m2:geo:column-equations}\\
 B^\dagger JA+D^\dagger JC&=0,&
 B^\dagger JB+D^\dagger JD&=\mu J.\notag
\end{align}
For such a matrix the equivalent row equations are
\begin{align}
 AJA^\dagger+BJB^\dagger&=\mu J,&
 AJC^\dagger+BJD^\dagger&=0,\label{m2:geo:row-equations}\\
 CJA^\dagger+DJB^\dagger&=0,&
 CJC^\dagger+DJD^\dagger&=\mu J.\notag
\end{align}
At every $F$-point one has $\rk_E B=\rk_E C$.
\end{proposition}

\begin{proof}
The column equations expand $x^\dagger\mathcal Jx=\mu\mathcal J$;
inverting gives $x\mathcal Jx^\dagger=\mu\mathcal J$, hence the row
equations. For the orthogonal block planes $U,W$, restriction of $x$
identifies $\ker B$ with $xW\cap W$ and $\ker C$ with $xU\cap U$.
Since $xW=(xU)^\perp$, one has
$(xW\cap W)^\perp=xU+U$. Comparing dimensions in the four-space
gives $\dim_E(xW\cap W)=\dim_E(xU\cap U)$, proving equal ranks.
\end{proof}

\subsubsection{Passage to the split quiver action}

Fix $\mu\in F^\times$ and write $X_\mu=\lambda^{-1}(\mu)$.
Choose a splitting field $k$ of $E/F$, for example $k=E$ or
$k=\barF$, and order the two factors of
$E\otimes_F k\simeq k\times k$. Conjugation interchanges these
factors; subscripts $1,2$ will denote their matrix components.

\begin{proposition}[Comparison with the split quiver]\label{m2:geo:split-bridge}
There are algebraic isomorphisms over $k$
\[
 (G_4)_k\simeq\GL_4\times\Gm,\qquad
 K_k\simeq\GL_2^3\times\Gm.
\]
Projection to the first matrix component identifies $(X_\mu)_k$
with the determinant-open variety
\[
 \mathscr X=\left\{(A,B,C,D):
 f(A,B,C,D)=\det\begin{pmatrix}A&B\\C&D\end{pmatrix}\ne0\right\}.
\]
The action factors through the three-vertex action
\begin{equation}\label{m2:geo:split-action}
 (p,q,r)\cdot(A,B,C,D)
 =(pAr^{-1},pBq^{-1},qCr^{-1},qDq^{-1}),
 \qquad(p,q,r)=(u_1,h_1,v_1).
\end{equation}
The multiplier torus acts trivially and factors into every geometric
stabilizer. The three-vertex action has its diagonal scalar kernel,
so its effective quotient and determinant character
$\chi(p,q,r)=\det p/\det r$ agree with the linear model. Transpose becomes
\begin{equation}\label{m2:geo:split-transpose}
 (A,B,C,D)\longmapsto(A^t,C^t,B^t,D^t),
\end{equation}
with group automorphism
$(p,q,r,s)\mapsto(r^{-t},q^{-t},p^{-t},s^{-1})$.
\end{proposition}

\begin{proof}
The two-component unitary equation is
$x_2^t\mathcal Jx_1=\mu\mathcal J$. Consequently
\begin{equation}\label{m2:geo:second-component}
 x_2=\mu\mathcal Jx_1^{-t}\mathcal J.
\end{equation}
Thus arbitrary invertible $x_1$ and nonzero $\mu$ determine $x_2$
uniquely and regularly; fixing $\mu$ imposes no condition on
$\det x_1$. Similarly, $u_2=sJu_1^{-t}J$ and the analogous
formulas for $v,h$ identify $K_k$ with $\GL_2^3\times\Gm$.
Substituting their unrestricted first components gives
\eqref{m2:geo:split-action}, independently of $s$, and taking determinants
gives the stated character. The first component of $x^t$ is $x_1^t$;
\eqref{m2:red:alpha} gives the group automorphism, including inversion
of $s$. These are algebraic action and transpose identifications.
\end{proof}

The vertex correspondence is $(U,Q,P)=(\mathcal U,\mathcal W,\mathcal P)$.
The arrows, determinant open, character, and transpose are therefore
those of \cref{m1:eq:X,m1:eq:det-character,m1:eq:residual-action}.
Over $\bar F$, apply the balanced-category and polystability results
\cref{m1:lem:theta-nonnegative,m1:lem:A0-abelian,m1:prop:polystability},
then \cref{m1:lem:stable-factors,m1:prop:rank-exclusion}.
General linear rational normal forms below are taken over $E$;
their Hermitian realization is a separate question.

\subsubsection{Closed families and transpose}

On the dense locus the general linear gauge of \eqref{m1:eq:dense-gauge} is
\begin{equation}\label{m2:geo:dense-gauge}
 \begin{pmatrix}X&I\\I&Y\end{pmatrix},\qquad
 X=B^{-1}AC^{-1},\quad Y=D,
 \qquad(p,q,r)=(B^{-1},I,C).
\end{equation}
The residual action is simultaneous conjugation, with domain
$\det\left(\begin{smallmatrix}X&I\\I&Y\end{smallmatrix}\right)\ne0$.
These coordinates are retained for the Hermitian existence problem.

By \cref{m1:prop:dense-closed}, closedness in the entire determinant
open is geometric semisimplicity of this pair. For rank zero,
\cref{m1:prop:rankzero} gives
\begin{equation}\label{m2:geo:rankzero}
 \begin{pmatrix}I&0\\0&D\end{pmatrix},\qquad
 D\in\GL_2,\qquad\text{closed exactly when $D$ is semisimple.}
\end{equation}
The rational unitary normalization will instead be
$a_\mu\oplus D$, with $D\in G_2(F)$ and $\lambda(D)=\mu$.

\paragraph{Rank-one coordinates.}

For arrow ranks $(1,1)$, \cref{m1:lem:rank11-split,m1:prop:rank11}
give the $E$-rational splitting $\mathrm E\oplus\mathrm H$ of a
closed point, with dimensions $(1,0,1)$ and $(1,2,1)$.
The incidence lines $\im c,\ker b$ are either both nonstable under
$D$, making $\mathrm H$ stable, or both stable and distinct,
giving $\mathrm H=\mathrm R\oplus\mathrm Q$; all other cases are
nonclosed. Hermitian orthogonality will follow from
\cref{m2:rat:selfadjoint-idempotents}.
The Hermitian calculations use the following normal forms.
The decomposable representative \eqref{m1:eq:decomposable-rep}, with
$(x,d,e)=(\zeta,\nu,\omega)$, is
\begin{equation}\label{m2:geo:decomposable-form}
 A=\diag(\zeta,1),\quad B=C=E_{11},\quad
 D=\diag(\nu,\omega),\qquad\omega\ne0,\quad\zeta\nu\ne1.
\end{equation}
Its factors are the identity on the second end coordinates,
$\left(\begin{smallmatrix}\zeta&1\\1&\nu\end{smallmatrix}\right)$
on $(1,1,1)$, and the middle-line loop $\omega$.

For a stable system the cyclic basis $(c,Dc)$ gives
\begin{equation}\label{m2:geo:stable-form}
 A=\diag(a,1),\quad
 B=\begin{pmatrix}b_0&b_1\\0&0\end{pmatrix},\quad C=E_{11},\quad
 D=\begin{pmatrix}0&-d_0\\1&s_0\end{pmatrix},\quad
 b_1=1-ad_0+s_0b_0.
\end{equation}
As in \eqref{m1:eq:b1-normalized}, scaling $\mathcal P_1$ normalizes
$ad_0-s_0b_0+b_1$ to one and scales $a,b_0,b_1$ together.
The target Hermitian form below incorporates this change of basis. Stability is
\begin{equation}\label{m2:geo:stable-pairing}
 H_0=\begin{pmatrix}b_0&b_1\\b_1&s_0b_1-d_0b_0\end{pmatrix}
     =\begin{pmatrix}b\\bD\end{pmatrix},\qquad\det H_0\ne0,
 \quad b=(b_0,b_1).
\end{equation}
These are \eqref{m1:eq:stable-rep-general}--\eqref{m1:eq:H} with
$(d,s)=(d_0,s_0)$, including nilpotent and repeated-root loops;
$d_0$ is independent of $d=\delta^2$ above.

\begin{theorem}[Seven coarse geometric types]\label{m2:geo:classification}
On every fixed similitude fiber the closed geometric $K$-orbits
are precisely the rank-zero semisimple family
\eqref{m2:geo:rankzero}, the dense semisimple-pair family
\eqref{m2:geo:dense-gauge}, and the two rank-one families
\eqref{m2:geo:decomposable-form} and \eqref{m2:geo:stable-form},
subject to the determinant and stability conditions in those formulas.
There are seven coarse geometric types: two rank-zero types
(scalar and distinct-eigenvalue semisimple), three dense types
(absolutely irreducible, two distinct characters, and two equal
characters), and two rank-one types. Each carries continuous
parameters. The nine rational rows arise only after distinguishing
the admissible $F$-forms.
\end{theorem}

\begin{proof}
Over $\bar F$, the comparison in \cref{m2:geo:split-bridge}
identifies the closed orbits with those of the split quiver.
The classification in \cref{m1:thm:closed-orbits} therefore yields
exactly two rank-zero, three dense, and two rank-one types.
\end{proof}

\paragraph{Geometric transpose.}

\begin{lemma}[Symmetric intertwiner]\label{m2:geo:pair-transpose}
Let $(X,Y)$ be a geometrically semisimple pair of two-by-two
matrices over $k_0=F$ or $E$. There is a symmetric invertible
$H\in\GL_2(k_0)$ such that $XH=HX^t$ and $YH=HY^t$.
\end{lemma}

\begin{proof}
Choose $s\in k_0$ with $X+sI$ invertible, then choose $t$ so that
$\det\left(\begin{smallmatrix}X+sI&I\\I&Y+tI\end{smallmatrix}\right)\ne0$:
as a polynomial in $t$ it has degree two and leading coefficient
$\det(X+sI)\ne0$. Scalar translations preserve the generated algebra
and semisimplicity. Thus \cref{m1:lem:dense-transpose} applies to this
closed dense point. Subtracting the scalar terms from its intertwiner
equations and taking $H=H_{\mathrm{lin}}^{-1}$ gives the required
symmetric invertible intertwiner in our convention.
\end{proof}

\begin{proposition}[Geometric transpose stability]
\label{m2:geo:transpose}
Every individual closed geometric orbit in
\cref{m2:geo:classification} is preserved by ordinary transpose.
A split-quiver point over $E$ with closed geometric orbit has
an $E$-linear general linear transporter from its transpose to itself.
\end{proposition}

\begin{proof}
The row formulas of \cref{m1:cor:transpose-stable}, transported
through \cref{m2:geo:split-bridge}, apply over $E$ to the general
linear transporter and over $\bar F$ to geometric orbit stability. The dense
transporter is $(H,H,H)$; in rank zero apply the symmetric
intertwiner to $(D,I)$ and take $(I,H,I)$. The decomposable
representative is symmetric, while the stable system has
\begin{equation}\label{m2:geo:stable-transporter}
 (I,H_0^{-1},I)\cdot(A^t,C^t,B^t,D^t)=(A,B,C,D),
\end{equation}
in the $(p,q,r)$ order, with $H_0$ as in
\eqref{m2:geo:stable-pairing}.
The general linear normalizations transport this construction to
an $E$-linear transporter at every point. Rational unitary transpose
stability, including common-multiplier compatibility, is proved in
\cref{m2:rat:transpose}.
\end{proof}

\subsubsection{Hermitian planes and rank-zero existence}

Write
\[
 Q_F=F^\times/\N_{E/F}(E^\times).
\]
This group has order two.  Every element has square one, because
$s^2=\N_{E/F}(s)$ for $s\in F^\times$.  For a nondegenerate
two-dimensional Hermitian matrix $H$ define its relative
discriminant, relative to the hyperbolic matrix $J$, by
\begin{equation}\label{m2:rat:discriminant}
 \operatorname{disc}_J(H)=[-\det H]\in Q_F.
\end{equation}
The minus sign divides by $\det J=-1$.

\begin{lemma}\label{m2:rat:hermitian-planes}
A nondegenerate Hermitian plane over $E/F$ is hyperbolic if and
only if $-\det H\in\N_{E/F}(E^\times)$.  Two Hermitian planes
over the local field $F$ are isometric if and only if their
determinants have the same class in $Q_F$.  Scaling a plane by
an element of $F^\times$ does not change this class.
\end{lemma}

\begin{proof}
Hermitian orthogonalization diagonalizes the form as
$\diag(a,b)$, with $a,b\in F^\times$.  It has a nonzero
isotropic vector exactly when $-b/a$ is a norm.  This is
equivalent to $-ab$ being a norm, since $a^2$ is a norm.
An isotropic nondegenerate plane is hyperbolic: choose an
isotropic vector $e$, choose $f$ with $\langle e,f\rangle=1$,
and replace $f$ by
$f-\frac12\langle f,f\rangle e$ to obtain an isotropic
partner. This proves the first assertion over any quadratic
extension of characteristic different from two.

For the second assertion, put $d=-b/a$ and form the quaternion
algebra $\mathscr B=E\oplus Ej$, with $j^2=d$ and
$jz=\bar z j$ for $z\in E$. Its reduced norm satisfies
\[
 a\operatorname{Nrd}(x+yj)
 =a\N_{E/F}(x)+b\N_{E/F}(y).
\]
Over a non-Archimedean local field the reduced norm
$\operatorname{Nrd}:\mathscr B^\times\to F^\times$ is
surjective~\cite[Lemma~13.4.9]{Voight21}. Thus the plane represents $1$.
Orthogonalizing at such a vector writes it as $\diag(1,c)$,
where $[c]=[\det H]$ in $Q_F$. Two such forms are isometric
when their second coefficients differ by a norm, and the
converse follows from the change-of-basis formula for
determinants. This proves the plane case of
\cite[Theorem~3.1]{Jacobowitz62}. Finally
$\det(sH)=s^2\det H$ and $s^2=\N_{E/F}(s)$.
\end{proof}

We express the hyperbolicity conditions using \eqref{m2:rat:discriminant}.

Fix $\mu\in F^\times$ and put $a_\mu=\diag(1,\mu)$.
If $B=C=0$, the unitary equations say $A,D\in G_2(F)$ with
multiplier $\mu$.  The element $(a_\mu A^{-1},1,1)\in K(F)$
normalizes $A$ to $a_\mu$.  Hence all rank-zero representatives
are
\begin{equation}\label{m2:rat:rankzero-representative}
 x_{\mu,D}=a_\mu\oplus D,
 \qquad D\in G_2(F),\quad\lam(D)=\mu.
\end{equation}
Two such representatives lie in the same $K(F)$-orbit exactly
when their $D$-blocks are conjugate in $G_2(F)$: one
implication follows from the last block of
\eqref{m2:red:action}. Conversely, if $hDh^{-1}=D'$, take
$v=a_{\lambda(h)}$ and $u=a_\mu va_\mu^{-1}$. Both have
multiplier $\lambda(h)$, and $(u,v,h)$ realizes the required
conjugacy in $K(F)$. The stabilizer embedding is
\begin{equation}\label{m2:rat:rankzero-stabilizer}
 K_{x_{\mu,D}}
 =\{(a_\mu v a_\mu^{-1},v,h):
       hD=Dh,\ \lam(v)=\lam(h)\}
 \simeq G_2\times_{\lam}Z_{G_2}(D).
\end{equation}
Closedness is equivalent to semisimplicity of $D$ by
\cref{m2:geo:classification}.  If $D=zI$, admissibility is
$\N_{E/F}(z)=\mu$ and the stabilizer is
$G_2\times_\lam G_2$.  If $D$ is regular semisimple,
$M=E[D]$ is quadratic etale over $E$.  Since
$D^*=\mu D^{-1}$, it is stable under $*$, and
\[
 Z_{G_2}(D)=T_{M,*}.
\]
The three algebra forms in \cref{m2:rat:torus-cohomology}
distinguish the split and the two anisotropic possibilities
for its quotient by the center.

\subsubsection{The dense Hermitian domain}

For a matrix $D\in\Mat_2(E)$ and a fixed $\mu\in F^\times$
define
\begin{align}
 M_D&=\mu I-D^*D,& Q_D&=JM_D,\label{m2:rat:dense-Gram}\\
 X_D&=-M_D^{-1}D^*,&
 R_D&=J+X_D^\dagger Q_DX_D,\notag
\end{align}
whenever $M_D$ is invertible.  The matrices $Q_D$ and $R_D$
are Hermitian.  Let $\Omega_\mu$ consist of those $D$ for
which $M_D$ and $R_D$ are invertible and
\begin{equation}\label{m2:rat:dense-domain}
 [-\det Q_D]=[-\det R_D]=1\quad\hbox{in }Q_F.
\end{equation}
These conditions say that $Q_D$ is isometric to $J$
and $R_D$ is isometric to $\mu J$; the latter is again
hyperbolic.  Thus \eqref{m2:rat:dense-domain} includes all
Hermitian-form obstructions to the choices below.

\begin{proposition}\label{m2:rat:dense-param}
For $D\in\Omega_\mu$, choose invertible $B,C$ with
\begin{equation}\label{m2:rat:dense-factorizations}
 B^\dagger JB=Q_D,\qquad C^\dagger R_DC=\mu J.
\end{equation}
Then
\begin{equation}\label{m2:rat:dense-representative}
 g_{\mu,D;B,C}=
 \begin{pmatrix}BX_DC&B\\C&D\end{pmatrix}
\end{equation}
belongs to $G_4(F)$ with multiplier $\mu$ and has rank
$(2,2)$.  Every rank-$(2,2)$ point is obtained in this way.
Its orbit depends precisely on the $G_2(F)$-conjugacy class
of $D$ in $\Omega_\mu$.  Its stabilizer is
\begin{equation}\label{m2:rat:dense-stabilizer}
 K_g=\{(BhB^{-1},C^{-1}hC,h):h\in Z_{G_2}(D)\}
 \simeq Z_{G_2}(D).
\end{equation}
\end{proposition}

\begin{proof}
The definitions give
\[
 Q_DX_D=-JD^*,\qquad X_D^\dagger Q_D=-JD.
\]
Hence the off-diagonal equation is
\[
 A^\dagger JB+C^\dagger JD
 =C^\dagger(X_D^\dagger Q_D+JD)=0.
\]
The last diagonal equation is
$B^\dagger JB+D^\dagger JD=\mu J$ by the definition of
$Q_D$, and the first is
\[
 A^\dagger JA+C^\dagger JC
 =C^\dagger(X_D^\dagger Q_DX_D+J)C=\mu J.
\]
These verify all equations of \cref{m2:geo:block-equations}.
Since the right-hand side is nondegenerate, the resulting
matrix is invertible.

Conversely, let the point have invertible $B,C$ and set
$X=B^{-1}AC^{-1}$.  The last diagonal equation gives
$B^\dagger JB=Q_D$, so $M_D$ is invertible.  The
off-diagonal equation gives $X^\dagger Q_D=-JD$, whence
$X=X_D$.  The first diagonal equation then gives the second
factorization in \eqref{m2:rat:dense-factorizations}.  Thus
$D\in\Omega_\mu$ and the point has the claimed form.

For the same $D$, let $(B',C')$ be another pair of choices.
Then $u=B'B^{-1}$ and $v=C'^{-1}C$ are unitary: the first
identity follows from their common Gram matrix $Q_D$, and
the second from their common $R_D$-Gram matrix $\mu J$.
The element $(u,v,1)$ carries the original representative
to the new one.  Under an arbitrary residual element, $D$
changes by $h$-conjugation, and the construction is equivariant.
This proves the orbit assertion.

The block equations for a stabilizing triple force
$u=BhB^{-1}$, $v=C^{-1}hC$, and $hD=Dh$.
Conversely, if $h\in G_2$ commutes with $D$, taking adjoints
shows that it commutes with $D^*$ as well.  It therefore
commutes with $X_D$ and preserves the two Gram forms up to
its multiplier.  For example,
$h^\dagger Q_Dh=\lam(h)Q_D$.  This shows that
$BhB^{-1}$ belongs to $G_2$ with multiplier $\lam(h)$.
Similarly $h^\dagger R_Dh=\lam(h)R_D$, and the second
factorization proves the same assertion for $C^{-1}hC$.
The remaining $A$-equation follows from $hX_D=X_Dh$.
\end{proof}

The general linear gauge of \eqref{m2:geo:dense-gauge} at this
point is $(X_D,D)$.  Besides $X_D\in E[D,D^*]$, one has
\begin{equation}\label{m2:rat:dense-algebra}
 D^*=-\mu X_D(I-DX_D)^{-1},
 \qquad E\langle X_D,D\rangle=E\langle D,D^*\rangle.
\end{equation}
For the first identity use
$X_D=-D^*(\mu I-DD^*)^{-1}$ and
$I-DX_D=\mu(\mu I-DD^*)^{-1}$.
All inverses in this calculation exist because
$\mu I-D^*D$ is invertible.  Thus
\cref{m1:prop:dense-closed} says that the dense orbit is closed
exactly when the natural module for $E\langle D,D^*\rangle$
is semisimple.  Characteristic zero ensures that
semisimplicity is preserved by algebraic extension.

The resulting stabilizers are as follows.  If the algebra is
$\Mat_2(E)$, its commutant is $E$, so the stabilizer is
$\Res_{E/F}\Gm$.  If $D$ is scalar, the stabilizer is
$G_2$.  Otherwise the semisimple commutative algebra is
$M=E[D]$, a quadratic etale $E$-algebra stable under $*$,
and the stabilizer is $T_{M,*}$. The algebra forms in
\cref{m2:rat:torus-cohomology} give its split or anisotropic quotient
by the center.

\subsubsection{Rank-one orthogonality and existence}

The two rank-one normal forms of \cref{m2:geo:classification} are
general linear normal forms over $E$.  We next describe
exactly which Hermitian structures on them give rational
points in the prescribed hyperbolic spaces.  The description
also proves the stabilizer assertions uniformly in their
parameters.

For an $E$-quiver representation with block isomorphism $g$,
let $\mathscr A_g$ be its endomorphism algebra.  In the
$(u,v,h)$ order it consists of triples satisfying
\begin{equation}\label{m2:rat:endo-equations}
 uA=Av,\qquad uB=Bh,\qquad hC=Cv,\qquad hD=Dh.
\end{equation}
Invertibility is not imposed on these endomorphisms.
For a decomposable closed object, the stable-factor
decomposition is
\[
 \mathrm E\oplus\mathrm R\oplus\mathrm Q,
 \qquad
 (1,0,1),\ (1,1,1),\ (0,1,0),
\]
and $\mathscr A_g=E^3$.  For a stable system it is
$\mathrm E\oplus\mathrm H$, of dimensions
$(1,0,1),(1,2,1)$, and $\mathscr A_g=E^2$.
These scalar endomorphism algebras and their vertex embeddings
are the stable-factor data of \cref{m1:prop:rank11,m1:eq:decomposable-stabilizer,m1:eq:stable-stabilizer}, applied over $E$.
These are the factor and endomorphism structures over $E$. The next
step incorporates the Hermitian adjoint and orthogonality to determine
the unitary stabilizers over $F$.

\begin{lemma}\label{m2:rat:selfadjoint-idempotents}
Suppose such a rank-one object is equipped with nondegenerate
Hermitian forms at its vertices and its block map is a
similitude of multiplier $\mu$.  Its endomorphism algebra
is stable under the componentwise Hermitian adjoint, and
every primitive stable-factor idempotent is selfadjoint.
Consequently the stable-factor decompositions at each
vertex are orthogonal decompositions into nondegenerate
spaces.
\end{lemma}

\begin{proof}
Write an endomorphism as $Lg=gR$, with
$L=\diag(u,h)$ and $R=\diag(v,h)$.  The adjoints use the
appropriate target and source vertex forms.  Since
$g^*=\mu g^{-1}$, taking adjoints and multiplying by $g$
gives
\[
 R^*g^{-1}=g^{-1}L^*,\qquad L^*g=gR^*.
\]
Thus $(u^*,v^*,h^*)$ is again an endomorphism.  The induced
semilinear involution of $E^3$ or $E^2$ permutes primitive
idempotents.  Adjoint preserves the rank at every vertex,
whereas the dimension vectors of these idempotents are
distinct.  It can therefore permute none of them.

For a selfadjoint idempotent $e$, its image is orthogonal
to its kernel. Their direct sum is the entire
nondegenerate vertex space, so both restrictions are
nondegenerate. Every primitive idempotent has this property.
\end{proof}

\begin{proposition}\label{m2:rat:rankone-stabilizers}
Put
\[
 T_m=\{(z_0,\ldots,z_{m-1})\in(E^\times)^m:
          \N(z_0)=\cdots=\N(z_{m-1})\}.
\]
The rational stabilizer of every decomposable rank-one
point is $T_3$, and that of every stable-system point is
$T_2$.  In the stable-factor coordinates their vertex
representations are
\begin{equation}\label{m2:rat:rankone-embeddings}
 \begin{array}{c|ccc}
 &u\text{ on }\mathcal P&v\text{ on }\mathcal U
                       &h\text{ on }\mathcal W\\ \hline
 T_3:(r,a_1,b_2)&\diag(r,a_1)&\diag(r,a_1)&\diag(r,b_2)\\
 T_2:(r,a_1)&\diag(r,a_1)&\diag(r,a_1)&rI_2.
 \end{array}
\end{equation}
The multiplier is the common norm.  All these stabilizers
are connected tori, of dimensions four and three respectively.
\end{proposition}

\begin{proof}
The automorphisms are exactly the invertible scalar
endomorphisms of the stable factors.  By
\cref{m2:rat:selfadjoint-idempotents} each factor is
nondegenerate at every vertex where it occurs.  Scalar
multiplication on it has multiplier the norm of that
scalar.  The requirement that all three vertex
similitudes have a common multiplier therefore says
exactly that all the scalar norms agree, with the
embeddings in \eqref{m2:rat:rankone-embeddings}.
Over $\barF$ the equations are equality of the products
of two independent coordinates for each scalar.  Choosing
their common product and one coordinate in each pair
identifies $T_m$ with a split torus of dimension $m+1$
after base change.  In particular there are no additional
finite components.
\end{proof}

\paragraph{Hermitian parameter domains.}

We describe the rank-one parameter domains by equations and
norm conditions.  Let $g_0$ be one of the
general linear normal forms of
\eqref{m2:geo:decomposable-form} or \eqref{m2:geo:stable-form},
with entries in $E$.  Choose $\mu\in F^\times$ and
Hermitian matrices $\mathsf H_U,\mathsf H_P,\mathsf H_W$ satisfying
\begin{equation}\label{m2:rat:rankone-Hermitian}
 g_0^\dagger\diag(\mathsf H_P,\mathsf H_W)g_0
       =\mu\diag(\mathsf H_U,\mathsf H_W),
 \qquad
 \det \mathsf H_i\ne0,\quad[-\det \mathsf H_i]=1\text{ in }Q_F.
\end{equation}
The same form $\mathsf H_W$ occurs on the source and target
middle vertex.  The final three conditions specify that
all vertex planes are the prescribed hyperbolic planes.

For the decomposable domain $\mathscr R_{\mathrm{dec}}$,
take $g_0$ as in \eqref{m2:geo:decomposable-form} and write
\begin{equation}\label{m2:rat:decomposable-forms}
 \mathsf H_U=\diag(u_r,u_a),\qquad
 \mathsf H_P=\diag(p_r,p_a),\qquad
 \mathsf H_W=\diag(w_r,w_b),
\end{equation}
where all six entries belong to $F^\times$.  Orthogonality
in \cref{m2:rat:selfadjoint-idempotents} shows that these
diagonal forms lose no solutions.  The block matrix
equation in \eqref{m2:rat:rankone-Hermitian} is exactly
\begin{align}
 p_a&=\mu u_a,& \N(\omega)&=\mu,
       \label{m2:rat:decomposable-domain}\\
 p_r\N(\zeta)+w_r&=\mu u_r,&
 p_r\bar\zeta+w_r\nu&=0,\notag\\
 p_r+w_r\N(\nu)&=\mu w_r.&&\notag
\end{align}
To verify this, restrict to the $(1,1,1)$ stable factor.
Its matrix $\left(\begin{smallmatrix}\zeta&1\\1&\nu\end{smallmatrix}\right)$
is a similitude from $\diag(u_r,w_r)$ to
$\diag(p_r,w_r)$, giving the last three equations and
their conjugates.  The other two factors give the first
two equations.  The hyperbolicity requirements are
\begin{equation}\label{m2:rat:decomposable-norms}
 [-u_ru_a]=[-p_rp_a]=[-w_rw_b]=1\quad\text{in }Q_F.
\end{equation}
Together with $\omega\ne0$ and $\zeta\nu\ne1$,
\eqref{m2:rat:decomposable-domain} and
\eqref{m2:rat:decomposable-norms} define the whole domain.
In particular a decomposable rational point can occur
only when $\mu$ is a norm.  This restriction follows
from the scalar factor on the same middle line in both
source and target.

For the stable-system domain $\mathscr R_{\mathrm{st}}$,
take \eqref{m2:geo:stable-form}, require $\det H_0\ne0$,
and write
\begin{equation}\label{m2:rat:stable-forms}
 \mathsf H_U=\diag(u_r,u_a),\qquad \mathsf H_P=\diag(p_r,p_a),
 \qquad
 \mathsf H_W=\begin{pmatrix}w_{11}&w_{12}\\
                    \bar w_{12}&w_{22}\end{pmatrix}.
\end{equation}
Here $u_r,u_a,p_r,p_a\in F^\times$,
$w_{11},w_{22}\in F$, and $w_{12}\in E$.  The middle
factor is the stable two-dimensional space.  With
$b=(b_0,b_1)$, equation \eqref{m2:rat:rankone-Hermitian}
becomes
\begin{align}
 p_a&=\mu u_a,\label{m2:rat:stable-domain}\\
 p_r\N(a)+w_{11}&=\mu u_r,\notag\\
 \bar a\,p_rb+e_1^\dagger \mathsf H_WD&=0,\notag\\
 p_rb^\dagger b+D^\dagger \mathsf H_WD&=\mu \mathsf H_W.\notag
\end{align}
The remaining off-diagonal block is the adjoint of the
third equation.  These equations, the nonzero determinant
and norm conditions of \eqref{m2:rat:rankone-Hermitian},
and the stability inequality $\det H_0\ne0$ are the
complete restrictions.  The determinant normalization
of the $(1,2,1)$ block also gives the useful consequence
$p_r=\mu^3u_r$.

\begin{proposition}\label{m2:rat:rankone-param}
For either of these domains choose $E$-linear isometries
$b_U,b_P,b_W$ such that $\mathsf H_i=b_i^\dagger Jb_i$ and put
\begin{equation}\label{m2:rat:rankone-transport-forms}
 g=(b_P\oplus b_W)g_0(b_U\oplus b_W)^{-1}.
\end{equation}
This gives a rank-one closed point of $G_4(F)$ with
multiplier $\mu$.  Conversely, every rank-one closed
$F$-point is obtained from one of these domains.
\end{proposition}

\begin{proof}
The isometries exist by the hyperbolicity conditions.
Substitution of $\mathsf H_i=b_i^\dagger Jb_i$ into
\eqref{m2:rat:rankone-Hermitian} proves
$g^\dagger\mathcal Jg=\mu\mathcal J$ for
\eqref{m2:rat:rankone-transport-forms}.  Changes of basis
preserve the arrow ranks and the geometric closedness
of the normal form.

Conversely, regard a rational unitary point as a quiver
representation over $E$.  The classification and the
linear retraction in \cref{m1:lem:rank11-split} give
an $E$-change of basis to the appropriate general linear
normal form. The original vertex forms $J$, transported along
these changes of basis, remain hyperbolic and satisfy
\eqref{m2:rat:rankone-Hermitian}.
\Cref{m2:rat:selfadjoint-idempotents} makes their stable-factor
decompositions orthogonal, yielding
\eqref{m2:rat:decomposable-forms} or \eqref{m2:rat:stable-forms}.
Expanding the block equation gives exactly the listed
domain equations.  In the stable case, the determinant-one normalization at
$\mathcal P_1$ is a change of basis accounted for by the transported
$\mathsf H_P$.
\end{proof}

These conditions describe the existence domains;
Section~\ref{m2:rat:section} establishes uniqueness of rational orbits
within each admitted geometric orbit.

\subsection{Rational orbit uniqueness and transpose}\label{m2:rat:section}

The existence domains in Section~\ref{m2:geo:section} supply the rational base points. We now compute the stabilizer-cohomology kernels at every admitted point.

\begin{theorem}[Rational closed-orbit uniqueness]\label{m2:rat:classification}
For every closed geometric $K$-orbit $\mathscr O$ containing an
$F$-point,
\begin{equation}\label{m2:rat:unique-orbit}
 K(F)\backslash\mathscr O(F)=\{*\}.
\end{equation}
The analytically closed rational orbits form nine families:
three rank-zero families (scalar, regular split-adjoint, and regular
anisotropic-adjoint), four dense families (absolutely irreducible,
scalar, split-adjoint toral, and anisotropic-adjoint toral), and the
decomposable and stable-system rank-one families. Their stabilizers
are given by \eqref{m2:rat:rankzero-stabilizer},
\eqref{m2:rat:dense-stabilizer}, and
\cref{m2:rat:rankone-stabilizers}; the full list is recorded with the
normal representations in \cref{m2:nor:classification}.
The parameter domains are
the rank-zero condition \eqref{m2:rat:rankzero-representative} with
$D$ semisimple, the dense domain \eqref{m2:rat:dense-domain} with
$E\langle D,D^*\rangle$ acting semisimply, and the complete
decomposable and stable-system domains of
\cref{m2:rat:rankone-param}. All their determinant, multiplier, and
hyperbolicity conditions are retained.
\end{theorem}

\subsubsection{Orbit cocycles and ambient cohomology}

\begin{proposition}[The pointed-set orbit correspondence]
\label{m2:rat:orbit-correspondence}
Fix a closed geometric orbit $\mathscr O$ containing
$x\in G_4(F)$, and let $K_x$ be its algebraic stabilizer over $F$.
Then there is a natural pointed-set bijection
\begin{equation}\label{m2:rat:orbit-kernel}
 K(F)\backslash\mathscr O(F)
 \simeq\ker\bigl[H^1(F,K_x)\longrightarrow H^1(F,K)\bigr],
\end{equation}
with the orbit of $x$ corresponding to the neutral class.
\end{proposition}

\begin{proof}
Put $\Gamma=\Gal(\barF/F)$. For $y=kx\in\mathscr O(F)$,
the elements $c_\gamma=k^{-1}\gamma(k)$ lie in $K_x(\barF)$
and satisfy $c_{\gamma\eta}=c_\gamma\gamma(c_\eta)$.
Their continuous cocycle class has neutral image in $H^1(F,K)$.
Changing the lift to $kt$, with $t\in K_x(\barF)$, replaces it by
$t^{-1}c_\gamma\gamma(t)$; left multiplication by $K(F)$ leaves
it unchanged. This defines the map on rational orbits.

Conversely, a stabilizer cocycle with neutral ambient image is
$k^{-1}\gamma(k)$, and $\gamma(kx)=kc_\gamma x=kx$ gives a
rational point. If lifts $k,l$ give cohomologous stabilizer cocycles,
choose $t\in K_x(\barF)$ with
$l^{-1}\gamma(l)=(kt)^{-1}\gamma(kt)$. Then $l(kt)^{-1}\in K(F)$
sends $kx$ to $lx$, proving injectivity. Finally $k=1$ corresponds
to the neutral class.
\end{proof}

The domains of Section~\ref{m2:geo:section} decide nonemptiness.
Every kernel calculation below starts with an admitted rational point.

We use nonabelian Galois cohomology with its distinguished neutral
class, as in \cite{SerreGC}.  Descent of the defining Hermitian
space identifies $H^1(F,U_{1,1})$ with the isometry classes of
nondegenerate Hermitian planes: after forgetting the form, a
cocycle is trivial by matrix Hilbert~90, and its descended form
is well defined up to an $E$-linear isometry. Conversely every
such form becomes isometric to $J$ over $\bar F$ and hence
defines a unitary cocycle. By \cref{m2:rat:hermitian-planes},
relative discriminant therefore gives
\[
 H^1(F,U_{1,1})\simeq Q_F.
\]
The multiplier sequence
\begin{equation}\label{m2:rat:ambient-sequence}
 1\longrightarrow U_{1,1}^3\longrightarrow K
 \xrightarrow{\lam}\Gm\longrightarrow1
\end{equation}
has the rational section
$s\mapsto(a_s,a_s,a_s)$, $a_s=\diag(1,s)$.
To see the resulting quotient directly, the multiplier of a
$K$-cocycle is a $\Gm$-cocycle. Hilbert~90 and the section
give a cohomologous cocycle in $U_{1,1}^3$. If two such cocycles
are cohomologous by $b\in K(\bar F)$, then $\lambda(b)$ is
Galois-fixed. Writing $b$ as a unitary triple times the section
of $s=\lambda(b)\in F^\times$ shows that their classes in
$H^1(F,U_{1,1}^3)$ differ precisely by the $F^\times$-action
induced by this section. Thus $H^1(F,K)$ is the quotient of
$Q_F^3$ by that action. On a Hermitian form it is scalar
similarity, up to change of basis.  Its determinant contribution in
each two-dimensional factor is $s^2$, hence a norm.  The action
on all three discriminants is therefore trivial, and
\begin{equation}\label{m2:rat:H1K}
 H^1(F,K)\simeq Q_F^3.
\end{equation}
The same argument applied to
\[
 1\longrightarrow U_{1,1}\longrightarrow G_2
 \xrightarrow{\lambda}\Gm\longrightarrow1,
 \qquad s\longmapsto a_s,
\]
gives $H^1(F,G_2)\simeq Q_F$.

\subsubsection{Regular involutive tori}

The same tori occur in the regular rank-zero and dense toral
families.  Let $M$ be a quadratic etale $E$-algebra with a
semilinear involution $*$ extending conjugation on $E$, and put
\[
 L=M^{*=1},\qquad
 T_{M,*}=\{r\in M^\times:r^*r\in F^\times\}.
\]
Then $L$ is a quadratic etale $F$-algebra and $M=L\otimes_FE$.
The multiplier of $r$ is $r^*r$.  These definitions refer to
algebraic groups over $F$, as well as to their rational points.

There are three possible involutive-algebra forms:
\begin{enumerate}
 \item $M=E\times E$ and
 $(r_1,r_2)^*=(\bar r_2,\bar r_1)$; then $L\simeq E$;
 \item $M=E\times E$ and
 $(r_1,r_2)^*=(\bar r_1,\bar r_2)$; then $L=F\times F$;
 \item $M$ is a field; then $L/F$ is quadratic and distinct
 from $E$, and $M=EL$ is biquadratic over $F$.
\end{enumerate}
To see the last assertion, the fixed field of $*$ has degree
two over $F$, is not $E$, and together with $E$ generates $M$.
The split cases follow by asking whether $*$ preserves or
interchanges the primitive idempotents.

\begin{proposition}\label{m2:rat:torus-cohomology}
There is a natural identification
\begin{equation}\label{m2:rat:H1T}
 H^1(F,T_{M,*})
 \simeq L^\times/\bigl(F^\times\N_{M/L}(M^\times)\bigr).
\end{equation}
Suppose $T_{M,*}$ is embedded in $G_2$ through a rank-one
Hermitian $M$-module whose underlying $E/F$-Hermitian plane
is hyperbolic. The map to $H^1(F,G_2)$ forgets the rank-one
$M$-module structure and retains the underlying Hermitian plane
up to scalar similarity. By \cref{m2:rat:hermitian-planes}, scalar
similarity does not change its isometry class in this dimension.
The induced map is the relative determinant map
\begin{equation}\label{m2:rat:Delta}
 \Delta:[c]\longmapsto[\N_{L/F}(c)]\in Q_F.
\end{equation}
For each of the three involutive-algebra forms, $\Delta$ is
injective.
\end{proposition}

\begin{proof}
The norm morphism gives an exact sequence of tori
\[
 1\longrightarrow T_{M,*}
 \longrightarrow\Res_{M/F}\Gm
 \xrightarrow{\ r\mapsto r^*r\bmod F^\times\ }
 (\Res_{L/F}\Gm)/\Gm\longrightarrow1.
\]
Exactness follows after base change from the corresponding
sequence of products of multiplicative groups. Hilbert~90
annihilates the first cohomology of the middle torus.
Moreover
\[
 ((\Res_{L/F}\Gm)/\Gm)(F)=L^\times/F^\times,
\]
because $H^1(F,\Gm)=1$.  The connecting map consequently gives
\eqref{m2:rat:H1T}.

To calculate the map induced by the two-dimensional
representation, identify its rank-one $M$-module with $M$
and write the given Hermitian form as
\[
 \langle x,y\rangle_{c_0}
   =\Tr_{M/E}(c_0xy^*),\qquad c_0\in L^\times.
\]
Indeed $\ell(z)=\langle z,1\rangle$ is $E$-linear, and
the adjoint identity gives
$\langle x,y\rangle=\ell(xy^*)$. Nondegeneracy of the
trace pairing writes $\ell(z)=\Tr_{M/E}(c_0z)$ uniquely.
Hermitian symmetry then gives $c_0^*=c_0$, and
nondegeneracy of the form gives $c_0\in L^\times$.
A cohomology class represented by
$c$ changes this coefficient to $c_0c$, modulo a common
$F$-scalar and an $M/L$ norm.  The determinant ratio is
$\N_{M/E}(c)=\N_{L/F}(c)$.  This proves
\eqref{m2:rat:Delta}; in particular the map is independent of
the initial coefficient $c_0$.

We evaluate all three involutive-algebra cases explicitly.
For the componentwise involution on $M=E\times E$, one has
$L=F\times F$ and
\[
 H^1(F,T_{M,*})=Q_F^2/\diag(Q_F),\qquad
 \Delta([c_1,c_2])=[c_1c_2].
\]
Since $Q_F$ has exponent two, the product map has exactly the
diagonal kernel and induces an isomorphism to $Q_F$.

If $L/F$ is a field distinct from $E$, write
$E=F(\sqrt a)$.  Compatibility of quadratic norm characters
with the field norm, or equivalently the projection formula
for Hilbert symbols, is
\begin{equation}\label{m2:rat:norm-projection}
 (c,a)_L=(\N_{L/F}(c),a)_F.
\end{equation}
Here $(c,a)_L=1$ means exactly that $c$ is a norm from
$L(\sqrt a)=M$.  Formula \eqref{m2:rat:norm-projection} follows
from the norm functoriality of local reciprocity
\cite[Proposition~4.3.11]{KedlayaCFT}: under
$\Gal(EL/L)\simeq\Gal(E/F)$, the restriction to $E$
of $\phi_F(\N_{L/F}c)$ equals the restriction to
$EL$ of $\phi_L(c)$.  Evaluation at the nontrivial
quadratic character yields \eqref{m2:rat:norm-projection}.
Here $\phi_F$ and $\phi_L$ denote the respective
local reciprocity maps.
In particular it gives
\[
 c\in\N_{M/L}(M^\times)
 \quad\Longleftrightarrow\quad
 \N_{L/F}(c)\in\N_{E/F}(E^\times).
\]
For $c\in F^\times$, its norm from $L$ is $c^2$, so every
$F$-scalar already lies in $\N_{M/L}(M^\times)$.  Thus the
kernel of \eqref{m2:rat:Delta} is trivial.  The map is also
onto $Q_F$: the norm subgroups for the two distinct quadratic
extensions $E/F$ and $L/F$ are distinct subgroups of index
two by local reciprocity
\cite[Theorem~4.1.2]{KedlayaCFT}, so
$\N_{L/F}(L^\times)$ is not contained in
$\N_{E/F}(E^\times)$.

Finally, for the exchanged involution on $M=E\times E$,
one has $L\simeq E$. The algebra $M/L$ is split, its norm
is surjective, and the group in \eqref{m2:rat:H1T} is trivial.
This completes all cases.
\end{proof}

\paragraph{The effective root torus.}

For the distributional arguments we also need the quotient of
$T_{M,*}$ by its scalar subgroup.

It is useful to distinguish $L$ from the quadratic algebra that
splits the effective root torus.  Let $\iota$ be the nontrivial
$E$-automorphism of $M$ and set
$K_T=M^{\iota *=1}$.

\begin{lemma}\label{m2:rat:adjoint-torus}
There is an exact sequence of algebraic $F$-tori
\begin{equation}\label{m2:rat:root-sequence}
 1\longrightarrow\Res_{E/F}\Gm
 \longrightarrow T_{M,*}
 \xrightarrow{\ r\mapsto r/\iota(r)\ }K_T^1
 \longrightarrow1,
\end{equation}
where $K_T^1$ is the kernel of the norm from $K_T$ to $F$.
The map on $F$-points is surjective.  The quotient is split in
the first of the three cases above and anisotropic in the other
two.  In the field case $K_T$ is the third quadratic subfield
   of $M$.
\end{lemma}

\begin{proof}
For $t=r/\iota(r)$ and $r^*r=\lambda\in F^\times$,
\[
 \iota(t)=t^{-1},\qquad
 t^*=\frac{\lambda/r}{\lambda/\iota(r)}=t^{-1}.
\]
Thus $\iota *(t)=t$ and $t$ has norm one in $K_T$.
The kernel consists of the elements fixed by $\iota$, namely
$E^\times$.  Over $\barF$ this is the surjective quotient
of a three-dimensional split torus by its two-dimensional
central torus; equivalently, the ratios of the two primitive
coordinates give every norm-one pair.  This proves algebraic
exactness.  Hilbert~90 gives
$H^1(F,\Res_{E/F}\Gm)=1$, so exactness on rational points
gives surjectivity there as well.  In the exchanged-involution
case $\iota *$ is componentwise conjugation, whose fixed
algebra is $F\times F$.  In the componentwise-involution
case it has fixed algebra isomorphic to $E$.  In the field
case its fixed field is the third quadratic subfield.
\end{proof}

\subsubsection{The splitting kernels at admitted points}

We now fix points in the Hermitian domains just constructed.
The following maps measure splitting inside their geometric orbits.

\begin{proposition}[Rank-zero splitting]\label{m2:rat:rankzero-kernel}
Let $x=x_{\mu,D}$ be an admitted rank-zero $F$-point with closed
geometric orbit. The kernel in \eqref{m2:rat:orbit-kernel} is trivial.
For regular $D$ it is identified with $\ker\Delta$ for
$T_{E[D],*}\subset G_2$; for scalar $D$ it is the scalar
conjugacy problem.
\end{proposition}

\begin{proof}
The rational normalization in \eqref{m2:rat:rankzero-representative}
identifies rational orbits with $G_2(F)$-conjugacy classes of $D$.
Over $\barF$ the same normalization identifies the geometric
problem with conjugacy of its first loop component, with fixed
$\mu$. Thus two rational points in this geometric orbit reduce
to geometrically conjugate $D$-blocks. For regular $D$, apply
the cocycle construction of \cref{m2:rat:orbit-correspondence} to
the conjugation action of $G_2$; its stabilizer is $T_{E[D],*}$.
The splitting set is therefore
\begin{equation}\label{m2:rat:rankzero-regular-kernel}
 \ker[H^1(F,T_{E[D],*})\longrightarrow H^1(F,G_2)]
 =\ker\Delta=\{1\}
\end{equation}
by all three cases of \cref{m2:rat:torus-cohomology}.
Conversely, suppose
$hDh^{-1}=D'$ with $h\in G_2(F)$. Choose
$v=a_{\lambda(h)}$ and $u=a_\mu v a_\mu^{-1}$. Then
$\lambda(u)=\lambda(v)=\lambda(h)$ and
$(u,v,h)$ carries $x_{\mu,D}$ to $x_{\mu,D'}$.
This proves injectivity for the fiber-product stabilizer.

If $D=zI$, its conjugacy orbit is a point, so the same normalization
already proves uniqueness. Equivalently, the multiplier sequence
for $G_2\times_\lambda G_2$ has kernel $U_{1,1}^2$ and the
section $(a_s,a_s)$. The scalar-similarity argument proving
\eqref{m2:rat:H1K} gives $H^1(F,G_2\times_\lambda G_2)=Q_F^2$.
The three vertex forms of \eqref{m2:rat:rankzero-stabilizer} then give
\begin{equation}\label{m2:rat:rankzero-scalar-H1}
 (q_1,q_2)\longmapsto(q_1,q_1,q_2),
\end{equation}
an injective map. Conjugation of the first representation by
$a_\mu$ changes only the initial form, not its relative determinant.
\end{proof}

\begin{proposition}[Dense splitting]\label{m2:rat:dense-kernel}
Let $g=g_{\mu,D;B,C}$ be an admitted dense $F$-point with
closed geometric orbit. Its rational splitting kernel is trivial
in each of the absolutely irreducible, scalar, split-adjoint
toral, and anisotropic-adjoint toral cases.
\end{proposition}

\begin{proof}
The stabilizer embedding \eqref{m2:rat:dense-stabilizer}
allows us to treat the cases separately. For the absolutely
irreducible case its stabilizer is
$\Res_{E/F}\Gm$, whose $H^1$ is trivial by Hilbert~90.
If $D$ is scalar, the three two-dimensional representations are
conjugate copies of the defining $G_2$-representation. Pulling
back their forms gives scalar multiples of the hyperbolic form;
twisting changes the three relative determinant classes by the
same $q$. Thus
\begin{equation}\label{m2:rat:dense-scalar-H1}
 H^1(F,G_2)=Q_F\longrightarrow Q_F^3,
 \qquad q\longmapsto(q,q,q)
\end{equation}
has trivial kernel.

For regular toral $D$, the three representations are rank-one
$M$-modules, $M=E[D]$, with the same involution $*$: the
commuting Gram operators in \eqref{m2:rat:dense-Gram} do not change
the adjoint on $M$. Their trace-form coefficients may differ.
Twisting by $c\in L^\times$, however, changes each determinant
by the same ratio $\N_{L/F}(c)$, independently of that initial
coefficient, as proved in \cref{m2:rat:torus-cohomology}. Hence
\begin{equation}\label{m2:rat:dense-toral-H1}
 H^1(F,T_{M,*})\longrightarrow Q_F^3,
 \qquad[c]\longmapsto(\Delta[c],\Delta[c],\Delta[c]).
\end{equation}
For the exchanged-product involution the source is trivial.
For the componentwise-product involution and the field case,
the injectivity of $\Delta$ gives a trivial kernel. This covers
both anisotropic-adjoint forms as well as the split-adjoint form.
\end{proof}

\paragraph{Rank-one kernels.}

The factor tori of \cref{m2:rat:rankone-stabilizers} fit into the
algebraic norm-ratio sequence
\begin{equation}\label{m2:rat:Tm-sequence}
 1\longrightarrow T_m\longrightarrow(\Res_{E/F}\Gm)^m
 \xrightarrow{\psi}\Gm^{m-1}\longrightarrow1,
 \quad\psi(z)=\left(\frac{\N(z_1)}{\N(z_0)},\ldots,
                    \frac{\N(z_{m-1})}{\N(z_0)}\right).
\end{equation}
After splitting, independent primitive coordinates realize
every norm ratio, proving algebraic exactness. On $F$-points
the image is $\N(E^\times)^{m-1}$: each ratio is a norm, and
each tuple of norms is obtained with $z_0=1$. Hilbert~90
annihilates $H^1$ of the middle torus, so its rational cokernel is
\begin{equation}\label{m2:rat:H1Tm}
 H^1(F,T_m)=Q_F^{m-1}\simeq Q_F^m/\diag(Q_F).
\end{equation}
The last quotient describes changing the Hermitian coefficient
of each stable factor modulo norms and common scalar similarity.

The common determinant rule is as follows. A factor with
coefficient class $[c_i]$ and dimension $d_{vi}$ at vertex $v$
contributes $[c_i]^{d_{vi}}$ to that vertex's determinant ratio.
Orthogonality from \cref{m2:rat:selfadjoint-idempotents} makes the
whole ratio the product $\prod_i[c_i]^{d_{vi}}$.
Since every vertex has dimension two, simultaneous scaling of
all coefficients multiplies this product by a square and thus
does not change it in $Q_F$. The rule therefore descends to
the quotient in \eqref{m2:rat:H1Tm}.

\begin{proposition}[The two rank-one splitting maps]
\label{m2:rat:rankone-kernel}
At every admitted rank-one $F$-point with closed geometric
orbit, the stabilizer map to $H^1(F,K)=Q_F^3$ is
\begin{align}
 H^1(F,T_3)&\longrightarrow Q_F^3,&
 (\alpha,\beta)&\longmapsto(\alpha,\alpha,\beta),
       \label{m2:rat:T3-H1}\\
 H^1(F,T_2)&\longrightarrow Q_F^3,&
 \alpha&\longmapsto(\alpha,\alpha,1).
       \label{m2:rat:T2-H1}
\end{align}
Both maps have trivial kernel.
\end{proposition}

\begin{proof}
For $T_3$ use coefficient classes $(c_r,c_a,c_b)$ on
$(\mathrm R,\mathrm E,\mathrm Q)$ and set
$\alpha=[c_rc_a]$, $\beta=[c_rc_b]$. These are coordinates
on $Q_F^3/\diag(Q_F)$: because $Q_F$ has exponent two,
the product classes equal the ratio classes, and simultaneous
rescaling can make $[c_r]=1$. The vertex factor contents in
the $(u,v,h)$ order are $(r,a),(r,a),(r,b)$ by
\eqref{m2:rat:rankone-embeddings}. The determinant rule gives
\eqref{m2:rat:T3-H1}. A neutral image forces
$\alpha=\beta=1$, which is the neutral quotient class.

For $T_2$ the factors are $(\mathrm H,\mathrm E)$ with
classes $(c_r,c_a)$, and $\alpha=[c_rc_a]$ is the single
quotient coordinate. The contents are $(r,a),(r,a),(r,r)$.
The middle space belongs twice to the \emph{same} stable
factor, so its determinant ratio is $[c_r^2]=1$. This gives
\eqref{m2:rat:T2-H1}. A neutral image forces $\alpha=1$.
The argument uses only the factor dimensions and their
selfadjoint idempotents, so includes nilpotent and repeated-root
stable loops and every parameter in both Hermitian domains.
\end{proof}

\subsubsection{Closedness, uniqueness, and transpose}

The group $K$ is connected reductive: its kernel in
\eqref{m2:rat:ambient-sequence} is a product of connected
unitary groups, and the quotient is a torus.  Its
unipotent radical maps trivially to that torus and would
be a connected normal unipotent subgroup of the kernel,
so it is trivial.  The variety $G_4$ is affine.
The closed-orbit comparison theorem for reductive
actions over local fields
\cite[Theorem~2.3.8]{AGDescent} consequently identifies
analytically closed $K(F)$-orbits with geometrically
closed orbits at their rational points.

\begin{proof}[Proof of \cref{m2:rat:classification}]
Fix $x\in\mathscr O(F)$. The geometric classification
\cref{m2:geo:classification} exhausts the possibilities. The
rank-zero normalization, \cref{m2:rat:dense-param}, and
\cref{m2:rat:rankone-param} recover every admitted rational point;
their converses establish precisely the existence domains in
the statement. We apply \cref{m2:rat:orbit-correspondence} only
after this rational base point has been fixed. The stabilizers
are \eqref{m2:rat:rankzero-stabilizer},
\eqref{m2:rat:dense-stabilizer}, and
\eqref{m2:rat:rankone-embeddings}.
For the rank-zero, dense, and rank-one families, respectively,
\cref{m2:rat:rankzero-kernel,m2:rat:dense-kernel,m2:rat:rankone-kernel}
prove that the stabilizer kernel is trivial throughout these
domains. Thus \eqref{m2:rat:orbit-kernel} gives
\eqref{m2:rat:unique-orbit} in every row. The reductivity and
analytic closedness comparison just established gives the classification
of analytically closed $K(F)$-orbits.
\end{proof}

\begin{theorem}\label{m2:rat:transpose}
Every closed rational $K(F)$-orbit on $G_4(F)$ is
preserved by ordinary transpose.  For each of its
points there is $k_x\in K(F)$ with $k_x\cdot x^t=x$.
Its extended orbit under $\Kt(F)$ is the same
single $K(F)$-orbit.  If $\sigma_x=k_x\sigma$, then
\begin{equation}\label{m2:rat:extended-stabilizer}
 \Kt_x(F)=K_x(F)\sqcup K_x(F)\sigma_x,
 \qquad \sigma_x^2=k_x\alpha_K(k_x)\in K_x(F).
\end{equation}
\end{theorem}

\begin{proof}
\Cref{m2:geo:transpose} places $x$ and $x^t$ in the
same closed geometric orbit.  Both are rational,
and \eqref{m2:rat:unique-orbit} shows that this
geometric orbit has only one rational orbit.
Both also have the same multiplier by \cref{m2:red:transpose}.
This proves rational transpose stability. The explicit transporters in
Appendix~\ref{app:transporters} also satisfy the common-multiplier
condition and provide the choices used in the normal-space calculations.

The extended orbit is the union of $K(F)x$ and
$K(F)x^t$, so it is the same orbit.  Every odd
element fixing $x$ differs from $\sigma_x$ by
an element of $K_x(F)$, proving the coset
description.  Since $\sigma_x$ fixes $x$, its
square does as well, and its square belongs to the
even stabilizer, as expressed in \eqref{m2:rat:extended-stabilizer}.
Changing $k_x$ changes $\sigma_x$ by an even
stabilizer element.  This observation will allow
the normal calculations to use convenient
transporters without changing the invariant
distribution problem.
\end{proof}

\subsection{The normal representations}\label{m2:nor:section}

All normal spaces are over $F$. For a closed rational representative $x$,
write $S_x=K_x$ and choose $k_x=(u_x,v_x,h_x)\in K(F)$ with
$k_x\cdot x^t=x$ as in \cref{m2:rat:transpose}; put $\sigma_x=k_x\sigma$.
The odd action on $N_x$ is induced by this corrected transpose
$\sigma_x$, rather than by ordinary transpose alone. Another choice
of $k_x$ changes the action by an even stabilizer element.
We continue to write $\mathcal J=J\oplus J$.

\subsubsection{The tangent, quiver Ext, and corrected transpose}

For $\mu=\lam(x)$, left trivialize tangent vectors as $x\xi$.
Write $N_x$ for the normal space to $Kx$ in $G_4$, and $N_x^0$
for the reduced normal space in the fixed-multiplier fiber $X_\mu=\lam^{-1}(\mu)$.

\begin{proposition}[Normal spaces and scalar extension]\label{m2:nor:full}
For every rational closed-orbit representative $x$, the orbit differential is
\begin{equation}\label{m2:nor:beta}
 \beta_x(U,V,H)=\Ad(x^{-1})(U\oplus H)-(V\oplus H),
 \qquad d\lam(U)=d\lam(V)=d\lam(H).
\end{equation}
It has image in $\mathfrak u_4$, and
\begin{equation}\label{m2:nor:quotients}
 N_x^0=\frac{T_xX_\mu}{T_x(Kx)}
       =\frac{\mathfrak u_4}{\im\beta_x},\qquad
 N_x=\frac{\mathfrak{gu}_4}{\im\beta_x}
       =N_x^0\oplus F[I_4/2].
\end{equation}
Tangent spaces, $\Lie(K)$, and this differential commute with every
field extension $k/F$. If $k$ splits $E/F$, the comparison in
\cref{m2:geo:split-bridge} gives an equivariant identification
\begin{equation}\label{m2:nor:base-change}
 N_x^0\otimes_F k\simeq N^{\mathrm{lin}}_{x_1},
\end{equation}
where the right side is the normal space in the linear model's
determinant-open quiver variety. For $(u,v,h)\in S_x(F)$ and for
$\sigma_x$, respectively, the actions on \eqref{m2:nor:quotients} are induced by
\begin{align}
 \xi&\longmapsto\Ad(v\oplus h)\xi,\label{m2:nor:even}\\
 \Theta_x(\xi)&=
 \Ad\bigl((v_x\oplus h_x)(x^t)^{-1}\bigr)(\xi^t).
 \label{m2:nor:odd}
\end{align}
Both fix $I_4/2$. The odd square is the action of
$k_x\alpha_K(k_x)\in S_x(F)$; changing $k_x$ multiplies the
odd operator by an even stabilizer action.
\end{proposition}

\begin{proof}
Over $F[t]/(t^2)$ the equation for the Lie algebra is
$\xi^*+\xi=d\lam(\xi)I$, and
$\lam(x(1+t\xi))=\mu(1+t\,d\lam(\xi))$.
Thus $T_xX_\mu=x\mathfrak u_4$. Differentiating the orbit
map gives $(U\oplus H)x-x(V\oplus H)$, proving
\eqref{m2:nor:beta}. Its two terms have the same multiplier differential.
Its kernel is the differentiated stabilizer equation, hence $\Lie(S_x)$.
Algebraic groups and their stabilizers are smooth in characteristic zero;
the orbit map is separable, so its differential image is $T_x(Kx)$.

The tangent and Lie spaces are kernels of these $F$-linear matrix
equations. Tensoring gives their equations over $E\otimes_F k$, hence
\[
 T_xX_\mu\otimes_F k\simeq T_{x_k}(X_\mu)_k,\qquad
 \Lie(K)\otimes_F k\simeq\Lie(K_k).
\]
The differential tensors to $\beta_{x_k}$; flatness preserves its
kernel, image, and cokernel. Stabilizers commute with extension by
their equalizer equations.

When $k$ splits $E/F$, send $\xi$ to the additive tangent $x_1\xi_1$.
The equation $\xi^*+\xi=0$ determines its second component uniquely
from the unrestricted first component. Under this isomorphism
\eqref{m2:nor:beta} becomes
\begin{equation}\label{m2:nor:split-orbit-map}
 d\rho(L,Q_0,V)=
 (LA-AV,\ LB-BQ_0,\ Q_0C-CV,\ Q_0D-DQ_0),
\end{equation}
where $(L,Q_0,V)=(U_1,H_1,V_1)$ are unrestricted. The extra multiplier
Lie direction acts trivially. This is the quiver differential of
\cref{m1:lem:hom-ext,m1:eq:orbit-differential}. Fixing $\mu$ imposes
no condition on $\det x_1$, and the determinant open has the full
arrow-space tangent. Its cokernel proves \eqref{m2:nor:base-change}
equivariantly.

For an even element use $(u\oplus h)x=x(v\oplus h)$ to obtain
\eqref{m2:nor:even}. Transpose takes $x(1+t\xi)$ to
$(1+t\xi^t)x^t$. The transporter identity
$(u_x\oplus h_x)x^t=x(v_x\oplus h_x)$ then gives
\[
 (u_x\oplus h_x)(1+t\xi^t)x^t(v_x\oplus h_x)^{-1}
 =x(1+t\Theta_x(\xi)).
\]
This proves \eqref{m2:nor:odd}; the factor $(x^t)^{-1}$ converts
back to the left trivialization. It also verifies the split odd
compatibility explicitly: writing $y=x_1\xi_1$, and $a_1,b_1$
for the first components of $u_x\oplus h_x,v_x\oplus h_x$, one has
\begin{equation}\label{m2:nor:split-odd}
 x_1(\Theta_x\xi)_1=a_1y^tb_1^{-1}.
\end{equation}
Thus the split operator is the corrected point-fixing additive
transpose of the linear model and preserves the orbit tangent.
The semidirect product gives its square. If $k'_x$ is another
transporter, $s=k'_xk_x^{-1}\in S_x(F)$ gives
$\Theta'_x=\rho(s)\Theta_x$. After moving to a geometric representative,
the same comparison holds with its canonical linear transporter;
the extra multiplier torus acts trivially.

Finally $d\lam(I_4/2)=1$, while both \eqref{m2:nor:even} and
\eqref{m2:nor:odd} fix this central symmetric matrix. It splits the
multiplier differential equivariantly, proving the normal-space
splitting.
\end{proof}

The Hom--Ext complex of \cref{m1:lem:hom-ext}, applied to the
relation-free quiver with its loop, combines with
\eqref{m2:nor:base-change} to give
\begin{equation}\label{m2:nor:hom-ext}
 N_x^0\otimes_F k\simeq
 \operatorname{Ext}^1_{\mathscr Q}(V_{x_1},V_{x_1}).
\end{equation}
After splitting, \eqref{m1:eq:multiplicity-spaces} gives the
automorphism action on a decomposition into pairwise nonisomorphic
absolutely stable factors: each cross-Ext has target scalar divided
by source scalar. All arrows and the loop may deform.
As a dimension check,
\begin{equation}\label{m2:nor:dimension-check}
 \dim_F N_x^0=16-(13-\dim S_x)=3+\dim S_x.
\end{equation}

\subsubsection{Normal representations for dense orbits}

\begin{lemma}\label{m2:rat:dense-smooth}
On a fixed multiplier fiber, the map from the dense locus to
the $D$-coordinate is a smooth submersion onto its admissible
domain.  Its fibers are orbits of the subgroup of $K$ with
$h=1$.  Consequently its differential identifies the normal-space
quotient to a $K$-orbit with
\[
 \Mat_2(E)/\{[Z,D]:Z\in\Lie(G_2)\}.
\]
\end{lemma}

\begin{proof}
The fiber assertion follows from \cref{m2:rat:dense-param}. The
admissible domain is open on $F$-points: nondegeneracy is open and
determinant norm classes are locally constant on $F^\times$.
We verify the submersion and its tangent kernel explicitly.
For any nondegenerate Hermitian form $H$, the Gram derivative at $b$
on the increment $bZ$ is
\[
 Z^\dagger(b^\dagger Hb)+(b^\dagger Hb)Z.
\]
For a prescribed Hermitian increment $Y$, the choice
$Z=\tfrac12(b^\dagger Hb)^{-1}Y$ gives exactly $Y$.
Prescribe any $\dot D\in M_2(E)$ and differentiate
\eqref{m2:rat:dense-Gram}. The first equation in
\eqref{m2:rat:dense-factorizations} is solved by
$\dot B=\tfrac12 BQ_D^{-1}\dot Q_D$.
For the second, the right side of the $C$ Gram derivative is
$-C^\dagger\dot R_DC$, a Hermitian matrix; explicitly one may take
\[
 \dot C=\tfrac12 C(\mu J)^{-1}(-C^\dagger\dot R_DC).
\]
Define $\dot A$ by differentiating $A=BX_DC$. This solves all
linearized block equations, and the implicit function theorem over
$F$ gives local smooth Gram choices as $D$ varies.

If $\dot D=0$, then $\dot Q_D=\dot X_D=\dot R_D=0$.
The first linearized Gram equation is equivalent to
$\dot B=UB$ with a unique $U\in\mathfrak u_2$; the second is
equivalent to $\dot C=-CV$ with a unique $V\in\mathfrak u_2$.
Thus the tangent kernel is
\begin{equation}\label{m2:rat:dense-tangent-kernel}
 (\dot A,\dot B,\dot C,\dot D)=(UA-AV,UB,-CV,0),
 \qquad U,V\in\mathfrak u_2,
\end{equation}
the tangent of the $(u,v,1)$ gauge orbit. After quotienting by it,
the residual differential is $\dot D=[H,D]$ with
$H\in\mathfrak{gu}_2$. Every such $H$ is allowed in $\Lie(K)$,
by choosing $U,V$ with its multiplier differential. This proves
the asserted normal quotient.
\end{proof}

\begin{lemma}\label{m2:nor:dense-quotient}
On the entire admitted dense domain of \cref{m2:rat:dense-param},
\begin{equation}\label{m2:nor:dense-normal}
 N_x^0\simeq M_2(E)/[\mathfrak{gu}_2,D].
\end{equation}
The even action is $[Y]\mapsto[hYh^{-1}]$ for
$h\in Z_{G_2}(D)$, and the corrected odd map is
\begin{equation}\label{m2:nor:dense-odd}
 [Y]\longmapsto[h_xY^th_x^{-1}].
\end{equation}
\end{lemma}
\begin{proof}
\Cref{m2:rat:dense-smooth} proves the unitary Gram submersion and
identifies its tangent kernel with the $(u,v,1)$ gauge tangent.
The remaining differential is $[H,D]$, giving the quotient.
The embedding \eqref{m2:rat:dense-stabilizer} gives the even action.
Transpose and its transporter send the lower-right additive block
to $h_xD^th_x^{-1}$, whose derivative is \eqref{m2:nor:dense-odd}.
The norm-class conditions define open subsets, so the tangent
calculation applies throughout the admitted domain.
\end{proof}

With $\delta,P,R$ and $\Psi$ as in \cref{m2:red:unitary-pair}, we write
\begin{equation}\label{m2:nor:pair-coordinate}
 P^{-1}DP=A_0+\delta B_0,\qquad A_0,B_0\in M_2(F).
\end{equation}
The factorization $h=zPgP^{-1}$ gives the equivalence
\[
 hD=Dh\quad\Longleftrightarrow\quad
 gA_0=A_0g\ \hbox{and}\ gB_0=B_0g.
\]
Every such $g\in\GL_2(F)$ lifts, so these matrices give the
effective rational stabilizer action.
The Lie equation in the $P$ basis is
$\bar Z^tJ_0+J_0Z=cJ_0$. Writing $Z=X+\delta Y$ with real matrices
gives $\tr X=c$ and $Y=(\tr Y)I/2$; hence
$P^{-1}\mathfrak{gu}_2P=\mathfrak{gl}_2(F)+EI$.
Its central part has zero commutator. Therefore $\Psi$ induces
\begin{equation}\label{m2:nor:pair-normal}
 N_x^0\simeq
 \frac{M_2(F)^2}
 {\{([Z,A_0],[Z,B_0]):Z\in\mathfrak{gl}_2(F)\}}.
\end{equation}
Using $\mathfrak{sl}_2(F)$ gives the same image. The algebra identities
\cref{m2:rat:real-pair} and \eqref{m2:rat:dense-algebra} show that closedness
of $x$ makes this real pair geometrically semisimple.

Choose the symmetric intertwiner $H$ in
\cref{m2:rat:dense-transporter}, so $A_0H=HA_0^t$, $B_0H=HB_0^t$,
and $h_x=P(HR^{-1})P^{-1}$. The uncorrected transpose in
\eqref{m2:simconj:transpose} contributes $R$; conjugation by $h_x$
cancels it. Thus the corrected pair operator is
\begin{equation}\label{m2:nor:pair-odd}
 (Y_1,Y_2)\longmapsto(HY_1^tH^{-1},HY_2^tH^{-1}).
\end{equation}
Other rational transporters add the stabilizer action just identified.

Choose scalar translations $(A_0+sI,B_0+tI)$ in the linear
determinant open as in \cref{m2:geo:pair-transpose}. They preserve the
generated algebra, centralizer, commutator image, and intertwiner
equations. Their differential is the identity and commutes with
\eqref{m2:nor:pair-odd}. Thus \cref{m1:eq:dense-normal} identifies the
quotient, even action, and corrected odd operator with the normal
representation for a dense orbit in the linear model, using $H=H_0^{-1}$ in \cref{m1:lem:dense-transpose}.

\paragraph{Absolute and scalar pairs}

For an absolutely irreducible pair, \cref{m1:prop:abs-irred-normal}
identifies the quotient by its five invariant differentials and
proves pointwise odd fixation. Hence
\begin{equation}\label{m2:nor:dense-simple}
 S_x=\Res_{E/F}\Gm,\qquad N_x^0=F^5_{\triv},\qquad\Theta_x=1.
\end{equation}

Equivalently, use the differentials of
\begin{equation}\label{m2:nor:five-invariants}
 \tr A_0,\quad\tr B_0,\quad\tr(A_0^2),\quad
 \tr(A_0B_0),\quad\tr(B_0^2).
\end{equation}
The identity $\tr(C^2)=(\tr C)^2-2\det C$ gives an invertible change
from the trace/determinant differentials. For a scalar pair,
\cref{m1:subsec:dense-normals} gives $M_2(F)^2$ with simultaneous
conjugation and componentwise transpose, hence $N_x^0=M_2(E)$ with
$G_2$-conjugation. For $D=d_cI$, \cref{m2:rat:dense-transporter} permits
$h_x=1$, $u_x=BC^{-t}$, $v_x=C^{-1}B^t$, all of multiplier one.
Then \eqref{m2:nor:dense-odd} is $Y\mapsto Y^t$; its factor $R$ in
pair coordinates is an even conjugation.

\paragraph{Toral pairs and the root field}

For a toral pair let $L_0=F[A_0,B_0]$, a quadratic etale algebra.
Its centralizer is $L_0^\times$. In the transported algebra
$M=E\otimes_F L_0$, the $E$-linear involution $\iota$ restricts
to the nontrivial automorphism of $L_0$. The adjoint formula in
\cref{m2:rat:real-pair} restricts to
$* =\overline{\phantom z}_E\otimes\iota$:
on a quadratic matrix algebra, adjugation is its nontrivial
automorphism. Consequently
\begin{equation}\label{m2:nor:root-field}
 K_T=M^{\iota *}=L_0,
\end{equation}
as embedded $F$-algebras. This identifies the real-pair root field
with \cref{m2:rat:adjoint-torus}.

If $L_0=F\times F$, the normal-space and transpose formulas
for the split case in \cref{m1:eq:dense-split-normal} give
\begin{equation}\label{m2:nor:dense-split}
 N_x^0=F^4_{\triv}\oplus F_\alpha\oplus F_{\alpha^{-1}},
 \qquad\Theta_x(z,u,v)=(z,v,u).
\end{equation}
The four self-Ext directions are odd-fixed by that result. In a
real simultaneous eigenbasis write $g=\diag(t_1,t_2)$ in the
factorization of $h$. The character is
$\alpha(h)=t_1/t_2$; it is independent of the common-scalar overlap.
The lift of $\diag(q,1)$ realizes every $q\in F^\times$.

For a field $L_0$, use the quotient and action in
\cref{m1:eq:dense-nonsplit-quotient,m1:eq:dense-nonsplit-action}:
if $A_0=m_a$, $B_0=m_b$, put $r_a=\bar a-a$, $r_b=\bar b-b$.
At least one is nonzero, and the root quotient is
\begin{equation}\label{m2:nor:anis-root-coordinate}
 L_0^2/L_0(r_a,r_b)\xrightarrow{\ \ell\ }L_0,
 \qquad\ell(y_1,y_2)=r_by_1-r_ay_2.
\end{equation}
The cited calculation gives multiplication $w\mapsto(t/\bar t)w$
and corrected transpose $w\mapsto-\bar w$, fixing the other four
coordinates. Bars here mean the $L_0/F$ involution. For nonzero
trace-zero $\omega\in L_0$, the coordinate $w=\omega w'$ commutes
with root multiplication and satisfies
$-\overline{\omega w'}=\omega\bar w'$. Hence
\begin{equation}\label{m2:nor:dense-anis}
 N_x^0=F^4_{\triv}\oplus K_T,\qquad
 \Theta_x(z,w')=(z,\bar w').
\end{equation}
Every $t\in L_0^\times$ lifts through $PgP^{-1}$, and Hilbert~90
makes $t/\bar t$ range over $K_T^1$, as in
\cref{m2:rat:adjoint-torus}. This includes $K_T\simeq E$ and the
biquadratic case, where $K_T$ is the third quadratic field.

\subsubsection{Rank-one normal representations and rational descent}

We use \eqref{m2:nor:hom-ext} over a splitting field $k$. For the
decomposable factors $(\mathsf E,\mathsf R,\mathsf Q)$, of dimensions
$(1,0,1),(1,1,1),(0,1,0)$, denote their scalars by $(a,r,b)$.
The vertex actions are $p=v=\diag(r,a)$, $q=\diag(r,b)$.
The linear Ext table \eqref{m1:eq:decomposable-ext-table} gives
\begin{equation}\label{m2:nor:dec-geometric}
 (N_x^0)_k=k^3_{\triv}\oplus k_\alpha\oplus k_{\alpha^{-1}}
                         \oplus k_\beta\oplus k_{\beta^{-1}},
 \qquad \alpha=r/b,\quad\beta=a/b.
\end{equation}
The positive lines are $\Ext^1(\mathsf Q,\mathsf R)$ and
$\Ext^1(\mathsf Q,\mathsf E)$; the trivial part consists of two
$\mathsf R$ self-Ext classes and one $\mathsf Q$ self-Ext class.
By \cref{m1:eq:decomposable-fixed-coordinates,m1:eq:decomposable-transpose},
the canonical corrected operator $\tau_0$ fixes these three classes,
swaps each opposite pair, inverts the effective torus, and squares
to one. This holds on the entire domain $\omega\ne0$, $\zeta\nu\ne1$,
including admitted $\zeta=0$ or $\nu=\omega$.

For the stable system $\mathsf E\oplus\mathsf H$ with scalars
$(a_s,r)$, the vertex actions are $p=v=\diag(r,a_s)$, $q=rI$.
The table \eqref{m1:eq:stable-ext-table} gives
\begin{equation}\label{m2:nor:stable-ext-table}
 (N_x^0)_k=k^4_{\triv}\oplus k_\eta\oplus k_{\eta^{-1}},
 \qquad\eta=a_s/r.
\end{equation}
Here the $\eta$ line is $\Ext^1(\mathsf H,\mathsf E)$ and the
trivial part is $\Ext^1(\mathsf H,\mathsf H)$. The canonical
transporter $(I,H_0^{-1},I)$ in vertex order $(p,q,v)$ induces
\begin{equation}\label{m2:nor:stable-odd-geometric}
 (\dot A,\dot B,\dot C,\dot D)\longmapsto
 (\dot A^t,\dot C^tH_0,H_0^{-1}\dot B^t,H_0^{-1}\dot D^tH_0).
\end{equation}
The parameter derivatives in $a,b_0,d_0,s_0$ form a basis of the
self-Ext part and are fixed modulo orbit tangents by
\cref{m1:lem:stable-parameter-basis,m1:eq:stable-parameter-derivatives,m1:eq:stable-ambient-transpose}:
differentiating the parameter-dependent transporter contributes an
orbit tangent. The same operator squares to one, inverts $\eta$,
and exchanges its root lines for every $\det H_0\ne0$, including
nilpotent and repeated-root loops. Compatibility with the rational normal representation
follows from \eqref{m2:nor:split-orbit-map}--\eqref{m2:nor:split-odd}.

\begin{lemma}[Descent of a multiplicity-one opposite pair]
\label{m2:nor:root-pair-descent}
Let $T$ be an $F$-torus split over a quadratic field $L/F$. Suppose an
$F$-representation $W$ has
$W_L=L_\chi\oplus L_{\chi^{-1}}$, two distinct one-dimensional
weight spaces of multiplicity one, exchanged by Galois. Let an
$F$-linear odd operator exchange these lines and suppose its square
is the identity on this effective module. Then $W$ is the underlying
$F$-plane of $L$, with its corresponding norm-one multiplication
character. The odd map is $z\mapsto c\bar z$, with $c\bar c=1$,
and an $L$-scalar coordinate change makes it ordinary conjugation.
\end{lemma}
\begin{proof}
Choose $e_\chi$ and its Galois conjugate $\bar e_\chi$ as the two
bases. The rational vectors are
$ze_\chi+\bar z\bar e_\chi$. Galois sends $\chi$ to $\chi^{-1}$,
so $\chi(s)\overline{\chi(s)}=1$ for $s\in T(F)$ and the action
in this coordinate is multiplication by $\chi(s)$. Exchanging the
lines and commuting with descent makes the odd map $c\bar z$.
Its assumed square one gives $c\bar c=1$ before Hilbert~90 is
used. Write $c=b/\bar b$ and $z=bw$; the map becomes $w\mapsto\bar w$.
This change commutes with the multiplication action.
\end{proof}
The algebraic stabilizers, with embeddings
\eqref{m2:rat:rankone-embeddings}, are
\begin{equation}\label{m2:nor:rational-rankone-tori}
 \begin{split}
 T_3&=\{(r,a,b)\in(\Res_{E/F}\Gm)^3:N(r)=N(a)=N(b)\},\\
 T_2&=\{(r,a)\in(\Res_{E/F}\Gm)^2:N(r)=N(a)\}.
 \end{split}
\end{equation}
They split over $E$. Define the norm-one characters
\begin{equation}\label{m2:nor:rational-characters}
 \chi_\alpha(r,a,b)=r/b,\qquad
 \chi_\beta(r,a,b)=a/b,\qquad\eta(r,a)=a/r.
\end{equation}
For such a character, $E_\chi$ denotes the underlying $F$-plane
of $E$ with multiplication by $\chi$.

\begin{proposition}\label{m2:nor:rankone-rational}
Throughout the admitted rational domains,
\begin{equation}\label{m2:nor:rankone-forms}
 \begin{array}{c|c|c}
 \text{family}&N_x^0&\Theta_x\\\hline
 \text{decomposable}&F^3_{\triv}\oplus E_{\chi_\alpha}
                 \oplus E_{\chi_\beta}&(z,u,v)\mapsto(z,\bar u,\bar v)\\
 \text{stable system}&F^4_{\triv}\oplus E_\eta
                 &(z,u)\mapsto(z,\bar u)
 \end{array}
\end{equation}
in suitable $F$-linear coordinates. The effective maps
$T_3(F)\to(E^1)^2$ and $T_2(F)\to E^1$ given by
\eqref{m2:nor:rational-characters} are surjective.
\end{proposition}
\begin{proof}
The common-norm equations give
$\overline{r/b}=b/r$, $\overline{a/b}=b/a$, $\overline{a/r}=r/a$.
Thus Galois exchanges each opposite pair of the characters in
\eqref{m2:nor:dec-geometric} and \eqref{m2:nor:stable-ext-table}.
The decomposable pairs are distinct and cannot be mixed by descent
or the odd map.

It remains to verify the odd hypotheses over $F$. By
\eqref{m2:nor:split-odd}, the rational transporter gives, after passage
to the geometric representative, $\rho(s)\tau_0$ for a stabilizer
element $s$. The preceding $\tau_0$ squares to one and inverts the
effective torus. Hence on the entire effective normal module
\[
 (\rho(s)\tau_0)^2=\rho(s)\rho(s^{-1})=1.
\]
The common scalar and the additional split multiplier act trivially.
Thus the rational operator is involutive and exchanges each pair.
With $L=E$, \cref{m2:nor:root-pair-descent} identifies each root plane
with $E_\chi$ and shows that its odd coefficient has norm one;
Hilbert~90 normalizes it to conjugation. The two decomposable
coordinate changes can be made independently.

The three or four even-trivial coordinates are odd-fixed by the
linear operator calculations above; stabilizer corrections act
trivially there, so this descends over $F$. Finally
$(r,a,b)=(q_1,q_2,1)$ realizes any $(q_1,q_2)\in(E^1)^2$ jointly,
and $(r,a)=(1,q)$ realizes any $q\in E^1$, all with common norm one.
\end{proof}

\subsubsection{The three rank-zero families}

We take the representative $x=a_\mu\oplus D$ from
\cref{m2:rat:classification}, where $a_\mu=\diag(1,\mu)$,
$\lam(D)=\mu$, and $D$ is semisimple.  Its stabilizer is
\[
 S_x=\{(a_\mu v a_\mu^{-1},v,h):
       (v,h)\in G_2\times_\lam Z_{G_2}(D)\}.
\]
For $Q\in M_2(E)$ put $Q^*=JQ^\dagger J$.  Every vector in
$\mathfrak u_4$ has the unique form
\begin{equation}\label{m2:nor:unitary-blocks}
 \xi=\begin{pmatrix}P_0&Q\\-Q^*&S_0\end{pmatrix},
 \qquad P_0,S_0\in\mathfrak u_2,
 \quad Q\in M_2(E).
\end{equation}

\begin{proposition}\label{m2:nor:zero-normal}
For the rank-zero representative above,
\[
 N_x^0\simeq M_2(E)\oplus\mathfrak c_D,
 \qquad
 \mathfrak c_D=
 \begin{cases}
 \mathfrak u_2,&D\text{ scalar},\\
 \mathfrak u_2\cap E[D],&D\text{ regular semisimple}.
 \end{cases}
\]
With $S_x$ identified as above, the stabilizer acts by
\begin{equation}\label{m2:nor:zero-action}
 (v,h)\cdot(Q,S_0)=(vQh^{-1},hS_0h^{-1}).
\end{equation}
Put
\[
 c_U=v_xa_\mu^{-t},\qquad c_W=h_xD^{-t}.
\]
Then
\begin{equation}\label{m2:nor:zero-odd}
 \Theta_x(Q,S_0)=
 \bigl(-c_UJ\bar QJc_W^{-1},\ c_WS_0^tc_W^{-1}\bigr),
 \qquad
 \lam(c_U)=\lam(c_W)=\frac{\lam(k_x)}\mu,
 \qquad c_WD^tc_W^{-1}=D.
\end{equation}
For regular $D$ the two-dimensional space $\mathfrak c_D$ is
pointwise fixed by both actions.  For scalar $D=d_cI$, the odd
operator may be taken, up to an even action, to be
\[
 (Q,S_0)\longmapsto(-J\bar QJ,S_0^t).
\]
\end{proposition}

\begin{proof}
At $a_\mu\oplus D$, the image in \eqref{m2:nor:beta} is block
diagonal, with blocks
\[
 a_\mu^{-1}Ua_\mu-V,
 \qquad D^{-1}HD-H.
\]
The first block fills $\mathfrak u_2$, by taking $V=H=0$ and
$U=a_\mu P_0a_\mu^{-1}$.  Conversely, after subtracting a common
central infinitesimal scalar from $U,V,H$, we may suppose that all
three have multiplier differential zero; that subtraction does
not change either of the two blocks just computed.  Consequently the image is exactly
\[
 \mathfrak u_2\oplus
 (\Ad(D^{-1})-1)\mathfrak u_2
\]
in the two diagonal blocks, and it has no off-diagonal part.

Semisimplicity of $D$ gives, after a splitting extension and then
descent to $F$,
\[
 \mathfrak u_2=
 \ker(\Ad(D)-1)\oplus\im(\Ad(D)-1).
\]
The kernel is $\mathfrak u_2$ for scalar $D$, and
$\mathfrak u_2\cap E[D]$ for regular $D$, since the matrix
centralizer is $E[D]$. In the latter case it is the Lie algebra of
the two-dimensional unitary centralizer. This proves the quotient;
\eqref{m2:nor:even} gives its even action.

Formula \eqref{m2:nor:odd} now has block-diagonal conjugating matrix
$c_U\oplus c_W$.  The top-right block of $\xi^t$ is
$-(Q^*)^t=-J\bar QJ$, giving \eqref{m2:nor:zero-odd}.
The multiplier equality follows from
$\lam(v_x)=\lam(h_x)$ and
$\lam(a_\mu)=\lam(D)=\mu$.  The lower-right transporter
equation is $h_xD^th_x^{-1}=D$; multiplying $h_x$ by $D^{-t}$
does not change this conjugation, so
$c_WD^tc_W^{-1}=D$ as claimed.  If $S_0\in E[D]$, write it
as $f(D)$ with coefficients in $E$.  Transpose leaves those scalar
coefficients unchanged, and therefore
$c_WS_0^tc_W^{-1}=f(D)=S_0$.  The centralizer also fixes this
space by conjugation.

If $D=d_cI$, the representative is symmetric and we can take
$k_x=1$.  Thus $c_U=a_\mu^{-1}$ and $c_W=d_c^{-1}I$.
The resulting leg map is $(d_ca_\mu^{-1})(-J\bar QJ)$.
Since $N_{E/F}(d_c)=\mu$, the matrix $d_ca_\mu^{-1}$ has
multiplier one.  Its left action, with right factor $1$, is an even
stabilizer action.  Removing it leaves the stated odd operator and
does not alter $S_0^t$.
\end{proof}

The scalar row has dimensions $(\dim S_x,\dim N_x^0)=(9,12)$,
and each regular row has $(7,10)$. In the regular rows the two
central directions are fixed, while the adjoint-torus type affects
the leg representation with the conditions in \eqref{m2:nor:zero-odd}.

\subsubsection{Classification of the rational normal representations}

\begin{theorem}\label{m2:nor:classification}
The nine rational closed-orbit families of \cref{m2:rat:classification}
have the algebraic $F$-stabilizers and reduced normal spaces
in the table below, with their indicated parameter-domain references.
Here $T=T_{M,*}$ is the torus of \cref{m2:rat:torus-cohomology}; its
split or anisotropic adjoint type and $K_T$ are those of
\cref{m2:rat:adjoint-torus}. The stabilizer embeddings in the table are \eqref{m2:rat:rankzero-stabilizer}, \eqref{m2:rat:dense-stabilizer},
and \eqref{m2:rat:rankone-embeddings}. Put
\[
 \kappa(Q)=-J\bar QJ,\qquad
 \kappa_x(Q)=c_U\kappa(Q)c_W^{-1},\qquad
 c_U=v_xa_\mu^{-t},\quad c_W=h_xD^{-t}.
\]
In both regular rank-zero rows the transporter conditions are
\begin{equation}\label{m2:nor:table-transport-conditions}
 \lam(c_U)=\lam(c_W)=\lam(k_x)/\mu,\qquad
 c_WD^tc_W^{-1}=D.
\end{equation}
Every $F^j$-summand in the table is fixed by both actions. In the action column,
only the other coordinates are written; the first line is even
and the second is odd. The last column gives
$(\dim S_x,\dim_F N_x^0)$.

{\small
\begin{longtable}{@{}L{.20\textwidth}L{.31\textwidth}L{.32\textwidth}c@{}}
\toprule
Family and domain & $S_x$; reduced normal $N_x^0$
 & Even action; corrected odd action & Dims.\\
\midrule
\endfirsthead
\toprule
Family and domain & $S_x$; reduced normal $N_x^0$
 & Even action; corrected odd action & Dims.\\
\midrule
\endhead
\bottomrule
\endfoot
Zero scalar,
$D=d_cI$, $N(d_c)=\mu$;
\eqref{m2:rat:rankzero-representative}
 & $G_2\times_\lam G_2$;
   $M_2(E)\oplus\mathfrak u_2$
 & $(Q,Z)\mapsto(vQh^{-1},hZh^{-1})$;
   $(Q,Z)\mapsto(a_\mu^{-1}\kappa(Q)d_c,Z^t)$
 & $9,12$\\[5pt]
Zero regular split-adjoint;
\eqref{m2:rat:rankzero-representative}
 & $G_2\times_\lam T$;
   $M_2(E)\oplus F^2_{\triv}$
 & $Q\mapsto vQh^{-1}$;
   $Q\mapsto\kappa_x(Q)$
 & $7,10$\\[5pt]
Zero regular anisotropic-adjoint;
\eqref{m2:rat:rankzero-representative}
 & $G_2\times_\lam T$;
   $M_2(E)\oplus F^2_{\triv}$
 & $Q\mapsto vQh^{-1}$;
   $Q\mapsto\kappa_x(Q)$
 & $7,10$\\[5pt]
Dense absolute;
\eqref{m2:rat:dense-domain}
 & $\Res_{E/F}\Gm$;
   $F^5_{\triv}$
 & identity; identity
 & $2,5$\\[5pt]
Dense scalar;
\eqref{m2:rat:dense-domain}
 & $G_2$;
   $M_2(E)$
 & $Y\mapsto hYh^{-1}$;
   $Y\mapsto Y^t$ ($h_x=1$)
 & $5,8$\\[5pt]
Dense split-adjoint toral;
\eqref{m2:rat:dense-domain}
 & $T$;
   $F^4_{\triv}\oplus F_\alpha\oplus F_{\alpha^{-1}}$
 & $(u,v)\mapsto(\alpha(h)u,\alpha(h)^{-1}v)$;
   $(u,v)\mapsto(v,u)$
 & $3,6$\\[5pt]
Dense anisotropic-adjoint toral;
\eqref{m2:rat:dense-domain}
 & $T$;
   $F^4_{\triv}\oplus K_T$
 & $w\mapsto\chi_T(h)w$;
   $w\mapsto\bar w$ over $K_T/F$
 & $3,6$\\[5pt]
One decomposable;
\eqref{m2:rat:decomposable-domain},
\eqref{m2:rat:decomposable-norms}
 & $T_3$;
   $F^3_{\triv}\oplus E_{\chi_\alpha}\oplus E_{\chi_\beta}$
 & $(u,v)\mapsto(\chi_\alpha u,\chi_\beta v)$;
   $(u,v)\mapsto(\bar u,\bar v)$
 & $4,7$\\[5pt]
One stable system;
\eqref{m2:rat:stable-domain},
\cref{m2:rat:rankone-param}
 & $T_2$;
   $F^4_{\triv}\oplus E_\eta$
 & $u\mapsto\eta u$;
   $u\mapsto\bar u$
 & $3,6$
\end{longtable}}

Here $\alpha(T(F))=F^\times$, and
$\chi_T:T(F)\twoheadrightarrow K_T^1$, $h\mapsto h/\iota(h)$,
is the root ratio in \cref{m2:rat:adjoint-torus}. The rank-one maps
\eqref{m2:nor:rational-characters} are surjective, jointly for two
roots. Both allowed anisotropic involutive-algebra forms are included.
The scalar zero-row operator uses $k_x=1$; removing the even action
$(d_ca_\mu^{-1},1)$ gives $(\kappa(Q),Z^t)$ and retains the entire
moving $\mathfrak u_2$. Other rows use the coordinates and rescalings
proved above.

For every row the normal space and extended stabilizer are
\begin{equation}\label{m2:nor:full-classification}
 N_x=N_x^0\oplus F[I_4/2],\qquad
 \widetilde S_x(F)=S_x(F)\sqcup S_x(F)k_x\sigma.
\end{equation}
The multiplier line is pointwise fixed. The sign $\eps$ is $1$
on $S_x(F)$ and $-1$ on its odd coset. A new transporter adds an
even stabilizer action; a new point of the same rational orbit
transports the entire representation by its rational group element.
The absence of sign-equivariant distributions on $N_x^0$
therefore implies the same on $N_x$.
\end{theorem}
\begin{proof}
\Cref{m2:nor:zero-normal} gives the three unitary rank-zero rows with
their transporter conditions. \Cref{m2:nor:dense-quotient} and
\eqref{m2:nor:dense-simple}--\eqref{m2:nor:dense-anis} give the four dense
rows through the rational matrix-pair comparison. \Cref{m2:nor:rankone-rational}
gives the remaining two. Their stabilizers and admitted domains are
the unchanged ones in \cref{m2:rat:classification}, which proves
exhaustiveness. All nine dimensions agree with
\eqref{m2:nor:dimension-check}. Finally \cref{m2:nor:full} proves the
normal-space splitting and transporter comparison. Projection to its
fixed multiplier line has complete fibers isomorphic to $N_x^0$;
\cref{m2:tools:localization} proves the last assertion.
\end{proof}

\section{Unitary descendants and multiplicity one}\label{sec:unitary-descendants}
We prove that invariant distributions on the rational normal spaces
classified in Section~\ref{sec:unitary-geometry} are fixed by the
corresponding corrected involutions. For regular rank-zero orbits,
we either reduce to the linear two-leg theorem or use Hermitian
Gram fibers. The scalar rank-zero case is treated by partial Fourier
transform and homogeneity. Together with generalized descent,
these arguments prove Theorem~\ref{thm:main-unitary}.

\subsection{Elementary descendants and simultaneous conjugation}
\label{m2:sec:elementary}

We first prove the elementary vanishing statements for the
representations in \cref{m2:nor:classification}. An odd
operator is always the corrected operator of the relevant
stabilizer. Multiplication of it by an even group element does not
change its action on even-invariant distributions.

We apply \cref{m2:elem:split,m2:elem:norm,m2:elem:products}
through the rational images of these stabilizers.

\subsubsection{Rational root maps on normal spaces}

For a regular toral dense orbit, let $T$ be its rational
stabilizer and $T_{\mathrm{ad}}$ its effective one-dimensional
torus. The root action computed in \cref{m2:nor:classification}
fits into
\begin{equation}\label{m2:elem:rootsequence}
 1\longrightarrow\Res_{E/F}\Gm\longrightarrow T
 \xrightarrow{\chi}T_{\mathrm{ad}}\longrightarrow1.
\end{equation}
By \cref{m2:rat:adjoint-torus}, $\chi$ is the eigenvalue ratio,
with scalar kernel, and is surjective on rational points;
$T_{\mathrm{ad}}$ is $\Gm$ in the split case and $K_T^1$ otherwise.
For the rank-one families, \cref{m2:nor:rankone-rational}
proves that the characters \eqref{m2:nor:rational-characters}
give surjections $T_3(F)\to(E^1)^2$ and $T_2(F)\to E^1$.
Thus the effective actions are the opposite-weight and
norm-one actions of \cref{m2:elem:split,m2:elem:norm,m2:elem:products}.

\begin{proposition}\label{m2:elem:vanishing}
The spaces of sign-equivariant distributions on the reduced normal
spaces vanish for the dense absolutely
irreducible, dense split toral, dense anisotropic toral,
rank-one decomposable, and rank-one stable-system families.
\end{proposition}

\begin{proof}
The calculations in \cref{m2:nor:classification} identify the
following extended representations on the reduced normal spaces.
In the absolutely irreducible case, the normal space is $F^5$ and
the odd operator acts trivially, so sign equivariance forces zero.
The split toral case is $F^4\oplus F^2$ with the last two
coordinates of opposite weights and swapped by the odd map.
The anisotropic toral case is $F^4\oplus K_T$, with
$K_T^1$ acting on the last coordinate and conjugation odd.
The rank-one cases are $F^3\oplus E^2$ and $F^4\oplus E$,
with the coordinate norm actions and conjugations just described.
The indicated copies of $F$ are fixed by the entire extended
stabilizer. The fixed coordinates are passive by \cref{m2:tools:passive}.
The rational root surjections then permit the application of
\cref{m2:elem:split,m2:elem:norm,m2:elem:products} to the remaining
coordinates, proving all five assertions.
\end{proof}

\subsubsection{Simultaneous-conjugation vanishing}\label{m2:sec:simconj}

We use the rational factorization, the coordinates $\Psi$, and the
matrices $J_0,R$ of \cref{m2:red:unitary-pair}, together with
\eqref{m2:simconj:transpose}.

\begin{theorem}\label{m2:simconj:vanishing}
Every $G_2(F)$-conjugation-invariant distribution on $M_2(E)$
is invariant under ordinary transpose.
\end{theorem}

\begin{proof}
Let $T'=\Psi_*T$. By \cref{m2:simconj:factorization}, it is invariant
under simultaneous $\GL_2(F)$-conjugation on $M_2(F)^2$.
The common matrix-pair result \cref{m1:lem:matrix-pair-transpose}
therefore fixes $T'$ under componentwise transpose. The odd
coordinate formula \eqref{m2:simconj:transpose} adds conjugation by
$R\in\GL_2(F)$, which also fixes $T'$. Pulling back through $\Psi$
proves the assertion.
\end{proof}

\subsection{Regular unitary two-leg representations}\label{m2:leg:section}\label{m2:sec:twoleg}

Let $D\in G_2(F)$ be regular semisimple, and put
\[
 T=Z_{G_2}(D),\qquad S=G_2\times_{\lambda}T,
 \qquad \mathcal L=\Mat_2(E).
\]
The even action is
\begin{equation}\label{m2:leg:action}
 (v,h)\cdot Q=vQh^{-1}.
\end{equation}
Throughout the regular-descendant argument the odd operator is specified by matrices
$c_U,c_W\in G_2(F)$ satisfying
\begin{equation}\label{m2:leg:transport-condition}
 \lambda(c_U)=\lambda(c_W),\qquad c_WD^tc_W^{-1}=D,
\end{equation}
and is
\begin{equation}\label{m2:leg:odd}
 \Theta(Q)=c_U\kappa(Q)c_W^{-1},
 \qquad \kappa(Q)=-J\bar QJ.
\end{equation}
These are precisely the conditions furnished by the regular rank-zero
normal representation in \cref{m2:nor:zero-normal,m2:nor:zero-odd}.
We first check the extended action, then prove its distribution theorem.

\begin{lemma}\label{m2:leg:extended-action}
The operator $\Theta$ normalizes the action of $S(F)$, and $\Theta^2$
belongs to that even action. In particular it makes sense to require a
distribution to be $S(F)$-invariant and $\Theta$-anti-invariant.
Moreover,
\begin{equation}\label{m2:leg:gram-odd}
 \Theta(Q)^*\Theta(Q)=c_W(Q^*Q)^tc_W^{-1}.
\end{equation}
\end{lemma}

\begin{proof}
For $g\in G_2(F)$, the unitary equation gives
$J\bar gJ=\lambda(g)g^{-t}$. Hence on an equal-multiplier pair,
\[
 \kappa(vQh^{-1})=v^{-t}\kappa(Q)h^t.
\]
Every $h\in T(F)$ belongs to $E[D]$, because the centralizer of a
regular two-by-two matrix is $E[D]$.
Condition \eqref{m2:leg:transport-condition} consequently gives
\begin{equation}\label{m2:leg:transpose-centralizer}
 c_Wh^tc_W^{-1}=h\quad(h\in T(F)),\qquad c_WT^tc_W^{-1}=T.
\end{equation}
Thus $\Theta$ conjugates $(v,h)$ to
$(c_Uv^{-t}c_U^{-1},c_Wh^{-t}c_W^{-1})\in S(F)$.
Their multipliers are both $\lambda(v)^{-1}$.
Since $\kappa^2=1$, a second application yields
\begin{equation}\label{m2:leg:odd-square}
 \Theta^2(Q)=(c_Uc_U^{-t})Q(c_Wc_W^{-t})^{-1}.
\end{equation}
Transposing the second condition in \eqref{m2:leg:transport-condition}
shows that $c_Wc_W^{-t}$ commutes with $D$. Both factors in
\eqref{m2:leg:odd-square} have multiplier one, proving the assertion.

For the Gram formula, $\kappa(Q)^*=-Q^t$, so
\[
 \kappa(Q)^*\kappa(Q)=Q^tJ\bar QJ=(Q^*Q)^t.
\]
Using $c_U^*c_U=\lambda(c_U)I$ and
$(c_W^{-1})^*=\lambda(c_W)^{-1}c_W$ now proves
\eqref{m2:leg:gram-odd}. The equality of multipliers is essential in
this cancellation.
\end{proof}

We write $\Ddist(\mathcal L)^{S,\Theta=-1}$ for the distributions
that are $S(F)$-invariant and $\Theta$-anti-invariant.

\begin{theorem}\label{m2:leg:regular}
For every regular semisimple $D\in G_2(F)$ and every pair
$c_U,c_W$ satisfying \eqref{m2:leg:transport-condition},
\[
 \Ddist(\Mat_2(E))^{S,\Theta=-1}=0.
\]
Equivalently, every distribution invariant under \eqref{m2:leg:action}
is fixed by \eqref{m2:leg:odd}.
\end{theorem}

The proof divides according to whether the adjoint torus $T/Z(G_2)$
is split or anisotropic.

\subsubsection{Rational factorization and the effective action}

As in \cref{m2:red:unitary-pair}, we have $\bar\delta=-\delta$,
$d=\delta^2$, $P=\diag(1,\delta)$ and $R=P^{-2}$.
By \cref{m2:simconj:factorization}, write
\[
 D=z_DPg_DP^{-1},\qquad z_D\in E^\times,
 \quad g_D\in\GL_2(F).
\]
The matrix $g_D$ is regular semisimple. Set
\[
 T_0(F)=Z_{\GL_2(F)}(g_D)=F[g_D]^\times.
\]
If $h=zPgP^{-1}$, then $hD=Dh$ is equivalent to $gg_D=g_Dg$.
It follows directly that
\begin{equation}\label{m2:leg:rational-torus}
 T(F)=E^\times PT_0(F)P^{-1}.
\end{equation}
The common scalar overlap is exactly $F^\times$:
$E^\times I\cap P\GL_2(F)P^{-1}=F^\times I$.
In particular the quadratic \emph{étale} $F$-algebra $F[g_D]$
determines the split or anisotropic type of the adjoint torus.

The kernel of \eqref{m2:leg:action} is the diagonal scalar group
\begin{equation}\label{m2:leg:kernel}
 \{(zI,zI):z\in E^\times\}.
\end{equation}
Indeed an element in the kernel has $v=h$ by evaluation at $Q=I$;
that matrix then commutes with every $Q\in\Mat_2(E)$ and is scalar.
Conversely, every pair $(zI,zI)$ in \eqref{m2:leg:kernel} has equal multipliers and
acts trivially.

Write $v=zPgP^{-1}$ and $h=wPtP^{-1}$, with
$g\in\GL_2(F)$ and $t\in T_0(F)$. Put $r=z/w$.
Since $\lambda(PgP^{-1})=\det g$, the multiplier condition becomes
\begin{equation}\label{m2:leg:effective-triples}
 N_{E/F}(r)\det g=\det t.
\end{equation}
Changing the two scalar lifts produces exactly
\[
 (r,g,t)\sim((b/a)r,ag,bt),\qquad a,b\in F^\times.
\]
If $P^{-1}QP=A+\delta B$ and $r=a_0+\delta b_0$, the action is
\begin{equation}\label{m2:leg:mixing}
 (A,B)\longmapsto
 \bigl(g(a_0A+db_0B)t^{-1},\ g(b_0A+a_0B)t^{-1}\bigr).
\end{equation}
This follows by multiplying $r g(A+\delta B)t^{-1}$.
The two rational matrix components are therefore mixed by a general
even element. We use a different linear map in the split case.

\subsubsection{Reduction for split adjoint tori}

Suppose $F[g_D]\simeq F\times F$. A change of basis by
$b=PfP^{-1}\in G_2(F)$ diagonalizes $D$. Such a change is compatible
with all the specified data: replace $D,Q$ by $bDb^{-1},bQb^{-1}$,
replace the even factors by their conjugates, and replace
\[
 c_U,c_W\quad\hbox{by}\quad bc_Ub^t,bc_Wb^t.
\]
Substitution in \eqref{m2:leg:odd} gives exactly the new odd operator.
Both new multipliers are multiplied by $\lambda(b)^2$, and the
transpose-transporter condition is preserved. We may thus assume
that $D$ and $T$ are diagonal in the hyperbolic basis with Gram matrix
$J$. An element $h=\diag(h_1,h_2)$ of $T(F)$ satisfies
$\bar h_1h_2\in F^\times$.

\begin{lemma}\label{m2:leg:split-effective}
The map
\begin{equation}\label{m2:leg:split-group-map}
 (g,s)\longmapsto
 \left(PgP^{-1},\ h_{g,s}\right),\qquad
 h_{g,s}=\diag\left(s,\frac{\det g}{\bar s}\right),
\end{equation}
identifies the effective action of $S(F)$ with that of
$(\GL_2(F)\times E^\times)/\diag F^\times$.
\end{lemma}

\begin{proof}
Both matrices in \eqref{m2:leg:split-group-map} have multiplier $\det g$,
and the formula is multiplicative. For surjectivity on effective
actions, factor an arbitrary $v=zPgP^{-1}$ and divide the pair
$(v,h)$ by $(zI,zI)$. The remaining right matrix is diagonal, of
multiplier $\det g$, so its second diagonal entry is
$\det g/\bar s$ if its first entry is $s$.
Finally the pair acts trivially precisely when
$PgP^{-1}=h_{g,s}$ is scalar. This means $g=aI$, $s=a$,
$a\in F^\times$. This proves the asserted kernel as well.
\end{proof}

Regard $E$ as an $F$-vector space with ordered basis $(1,\delta)$, and
let $m_s$ be multiplication by $s$ on this space. Define
\[
 \mathcal T(x)(u)=
 \begin{pmatrix}\Tr_{E/F}(x_1u)\\\Tr_{E/F}(x_2u)\end{pmatrix},
 \qquad x=(x_1,x_2)^t\in E^2.
\]
The trace pairing is nondegenerate. More explicitly, if
$x_i=x_{i0}+\delta x_{i1}$, then
\begin{equation}\label{m2:leg:trace-matrix}
 \mathcal T(x)=2
 \begin{pmatrix}x_{10}&d x_{11}\\x_{20}&d x_{21}\end{pmatrix}.
\end{equation}
For a two-by-two matrix over $F$ write
$Y^\#=\tr(Y)I-Y$; this is the adjugate, and is an $F$-linear
involution. Write $P^{-1}QP=(u,v)$ by columns and define
\begin{equation}\label{m2:leg:linear-map}
 \mathcal I(Q)=(X,Y),\qquad X=\mathcal T(u),
 \quad Y^\#=\mathcal T(v).
\end{equation}
It is an $F$-linear homeomorphism
\[
 \Mat_2(E)\longrightarrow
 \Hom_F(E,F^2)\oplus\Hom_F(F^2,E).
\]

\begin{lemma}\label{m2:leg:split-intertwining}
Under \eqref{m2:leg:split-group-map} and \eqref{m2:leg:linear-map}, the
even action becomes
\begin{equation}\label{m2:leg:model-action}
 (g,s)(X,Y)=(gXm_s^{-1},m_sYg^{-1}).
\end{equation}
Let $S_E=\diag(2,2d)$, the Gram matrix of
$\Tr_{E/F}(uu')$, and put
\begin{equation}\label{m2:leg:model-odd}
 \tau_E(X,Y)=(Y^tS_E,S_E^{-1}X^t).
\end{equation}
Then the odd intertwining formula is
\begin{equation}\label{m2:leg:odd-intertwining}
 \mathcal I\kappa\mathcal I^{-1}
   =(R,m_{2\delta})\cdot\tau_E.
\end{equation}
Here the pair on the right is an element of the even group in
\eqref{m2:leg:model-action}.
\end{lemma}

\begin{proof}
The two $E$-columns transform by
\[
 u\longmapsto gu s^{-1},\qquad
 v\longmapsto(\det g)^{-1}gv\bar s.
\]
We have $\mathcal T(xs)=\mathcal T(x)m_s$. Also adjugation reverses
products and
\[
 (g^{-1})^\#=(\det g)^{-1}g,
 \qquad m_s^\#=N_{E/F}(s)m_s^{-1}=m_{\bar s}.
\]
Consequently
\[
 (m_sYg^{-1})^\#=(\det g)^{-1}gY^\#m_{\bar s},
\]
which proves \eqref{m2:leg:model-action} for both legs.

For the odd formula set
$L=\left(\begin{smallmatrix}0&-\delta\\\delta^{-1}&0\end{smallmatrix}\right)$.
Direct substitution of $Q=P(u,v)P^{-1}$ gives
\[
 P^{-1}\kappa(Q)P=L\overline{(u,v)}L.
\]
Its columns are
\[
 \begin{pmatrix}-\bar v_2\\d^{-1}\bar v_1\end{pmatrix},
 \qquad
 \begin{pmatrix}d\bar u_2\\-\bar u_1\end{pmatrix}.
\]
Using \eqref{m2:leg:trace-matrix} and
$Y^\#=J_0Y^tJ_0^{-1}$ yields
\begin{equation}\label{m2:leg:odd-coordinates}
 \mathcal I\kappa\mathcal I^{-1}(X,Y)
       =(RY^tJ,JX^tR^{-1}).
\end{equation}
Since
\[
 m_{2\delta}=\begin{pmatrix}0&2d\\2&0\end{pmatrix}=JS_E,
 \qquad S_Em_{2\delta}^{-1}=J,
\]
applying $(R,m_{2\delta})$ to \eqref{m2:leg:model-odd} gives exactly
\eqref{m2:leg:odd-coordinates}.
\end{proof}

We apply the linear two-leg theorem \cref{m1:thm:two-leg}, using
its trace-form adjoint \eqref{m1:eq:two-leg-adjoint}.

In our application
$V=F^2$, $W=E$, $T'=E^\times$ acts by multiplication, and
$S_W=S_E$. The trace identity
$\Tr((su)u')=\Tr(u(su'))$ gives
$S_Em_s=m_s^tS_E$. Thus this is the nonsplit case of that theorem,
even though the unitary adjoint torus considered here is split.
The scalar subgroup $(aI,aI)$ acts trivially; passage to its
effective quotient changes neither the distribution space nor the
sign character. The identification is on $F$-points and on the entire
linear representation.

For singular legs, the nonsplit proof is given in \cref{m1:subsec:two-leg-nonsplit}, using
\cref{m1:lem:matrix-boundary,m1:lem:rank-one-factorization,m1:lem:common-U11}. It localizes by the complete composition invariants
of $YX$, then treats ranks two, one, and zero of that composition.
For rank one, its factorization lemma includes singular individual
legs. For zero composition, it first treats the open locus
$\rk X=\rk Y=1$, and only then the closed union of matrix axes.
On a rank-one matrix axis, descent over $\mathbf P(V)$ leaves the
upper triangular group
\[
 P_V=\left\{\begin{pmatrix}\alpha&b\\0&\beta\end{pmatrix}\right\}
\]
with required Frobenius character $|\beta/\alpha|$; the element
$\diag(1,\beta)$ acts trivially on the resulting geometric fiber.
This is the modular obstruction on that stratum. On the two rank-two
axes the finite-part distributions have opposite logarithmic defects
and the same nonzero residue at the origin. 

We can finish the split case of \cref{m2:leg:regular}. For diagonal
regular $D$, \eqref{m2:leg:transport-condition} says that $c_W$ commutes
with $D$, so $c_W\in T(F)$ and $(c_U,c_W)\in S(F)$.
Thus $\Theta$ differs from $\kappa$ by an even action. By
\cref{m2:leg:split-intertwining}, $\kappa$ in turn differs from
$\tau_E$ by the explicit even pair $(R,m_{2\delta})$.
Thus \cref{m1:thm:two-leg} makes every even-invariant distribution
$\Theta$-invariant.

\subsubsection{Anisotropic tori and the self-adjoint Gram space}

Assume now that $T/Z(G_2)$ is anisotropic. Let $W=E^2$ be the right
Hermitian space, let $M=E[D]\subset\End_E(W)$, and denote the
Hermitian adjoint on $M$ by $*$. Since
$D^*=\lambda(D)D^{-1}$, this algebra is stable under $*$.
Put $L=M^{*=1}$. There are two forms to consider:
\begin{align}
 &M=E\times E,\qquad
   (r_1,r_2)^*=(\bar r_1,\bar r_2);\label{m2:leg:an-product}\\
 &M=EL\text{ a biquadratic field},\qquad
   \langle *\rangle=\Gal(M/L).\label{m2:leg:an-field}
\end{align}
These are the anisotropic cases of \cref{m2:rat:torus-cohomology}.
Let $\iota$ denote the nontrivial $E$-automorphism of $M$.
If $K_T$ is the splitting field of the adjoint torus, then
\begin{equation}\label{m2:leg:splitting-field}
 K_T=E\quad\hbox{in \eqref{m2:leg:an-product}},\qquad
 K_T=M^{*\iota}\quad\hbox{in \eqref{m2:leg:an-field}}.
\end{equation}
These identifications are those of \cref{m2:rat:adjoint-torus}.
In both cases
\[
 T(F)=\{r\in M^\times:r^*r\in F^\times\}.
\]

Since $D$ is cyclic, a cyclic vector identifies the rank-one
$M$-module $W$ with $M$. The trace-form calculation in
\cref{m2:rat:torus-cohomology} writes its Hermitian form as
\begin{equation}\label{m2:leg:trace-hermitian}
 \langle u,v\rangle=\Tr_{M/E}(c uv^*),
 \qquad c\in(M^{*=1})^\times.
\end{equation}
These coordinates carry the given hyperbolic form on $W$.

Put
\[
 \mathcal H=\{A\in\End_E(W):A^*=A\},\qquad
 \Gamma(Q)=Q^*Q.
\]
The action of $S$ on Gram matrices is
\begin{equation}\label{m2:leg:gram-even}
 \Gamma(vQh^{-1})=h\Gamma(Q)h^{-1}.
\end{equation}
For example this follows from
$(h^{-1})^*=\lambda(h)^{-1}h$ and the equality
$\lambda(v)=\lambda(h)$. We compute the rational $T$-orbits of
$\mathcal H$ before analyzing the Gram factorization fibers.

\begin{lemma}\label{m2:leg:root-image}
There is a $T$-stable decomposition over $F$
\[
 \mathcal H=C\oplus W_{\mathrm{rt}},\qquad
 \dim_F C=2,\qquad W_{\mathrm{rt}}\simeq K_T,
\]
where $C$ is pointwise fixed and $K_T^1$ acts on the
root plane by multiplication. The rational root maps have kernel
$E^\times$ and are surjective.
\end{lemma}

\begin{proof}
In \eqref{m2:leg:an-product}, write $c=(c_1,c_2)$ with $c_i\in F^\times$.
Relative to the two idempotents, the self-adjoint matrices are
\begin{equation}\label{m2:leg:product-selfadjoint}
 A=\begin{pmatrix}a&b\\(c_1/c_2)\bar b&e\end{pmatrix},
 \qquad a,e\in F,\quad b\in E.
\end{equation}
The fixed coordinates are $a$ and $e$; the root coordinate is $b$.
Conjugation by $(r_1,r_2)\in T$ multiplies $b$ by $r_1/r_2$.
This is the root ratio in \cref{m2:rat:adjoint-torus}.

In \eqref{m2:leg:an-field}, $c\in L^\times$ and every $E$-endomorphism
of $M$ has a unique expression $m_x+m_y\iota$, $x,y\in M$.
The two summands are linearly independent: otherwise $\iota$ would
be multiplication by an element of $M$, contrary to its action on $1$
and on an element outside $E$. Dimensions then give the assertion.
Taking the adjoint in \eqref{m2:leg:trace-hermitian} shows
\[
 (m_y\iota)^*=m_{c^{-1}\iota(c)\iota(y^*)}\iota.
\]
Consequently the self-adjoint matrices are precisely
\begin{equation}\label{m2:leg:field-selfadjoint}
 A=m_x+m_{c^{-1}z}\iota,\qquad x\in L,\quad z\in K_T.
\end{equation}
Here $y=c^{-1}z$ satisfies the adjoint identity above precisely
when $z=\iota(z^*)$. Multiplication operators $m_x$ give $C$, and
the $z$-plane is the root plane.

Conjugation by $r\in T$ multiplies $z$ by $r/\iota(r)$.
This is again the root ratio of \cref{m2:rat:adjoint-torus},
which gives its rational surjectivity and scalar kernel in both cases.
\end{proof}

Choose any $F$-basis of $C$ and write a self-adjoint matrix as
$(a_1,a_2,w)\in F^2\oplus K_T$. Define
\begin{equation}\label{m2:leg:an-invariants}
 \pi(Q)=\bigl(a_1(\Gamma(Q)),a_2(\Gamma(Q)),
                     N_{K_T/F}(w(\Gamma(Q)))\bigr)\in F^3.
\end{equation}

\begin{lemma}\label{m2:leg:gram-orbits}
Every nonempty fiber of $(a_1,a_2,N(w))$ on $\mathcal H$ is one
rational $T(F)$-orbit. Its points all have the same rank.
The map \eqref{m2:leg:an-invariants} is invariant under both $S$ and
$\Theta$.
\end{lemma}

\begin{proof}
A nonzero norm fiber in the field $K_T$ is one $K_T^1$-orbit:
if $N(w')=N(w)\ne0$, then $w'/w\in K_T^1$.
The zero norm fiber consists only of zero. Rational surjectivity in
\cref{m2:leg:root-image} therefore gives the orbit assertion.
In product coordinates \eqref{m2:leg:product-selfadjoint},
\[
 \det_E A=ae-(c_1/c_2)N_{E/F}(b).
\]
In field coordinates \eqref{m2:leg:field-selfadjoint},
\begin{equation}\label{m2:leg:field-determinant}
 \det_E A=N_{L/F}(x)-N_{L/F}(c)^{-1}N_{K_T/F}(z).
\end{equation}
For the latter identity, after extending scalars to split $M/E$,
$m_x+m_y\iota$ has matrix
$\left(\begin{smallmatrix}x&y\\\iota(y)&\iota(x)\end{smallmatrix}\right)$;
its determinant is $x\iota(x)-y\iota(y)$, giving the formula over $F$.
These formulas determine the determinant on every invariant fiber.
They also distinguish the zero matrix: in the product case all its
coordinates vanish; in the field case $z=0,x\ne0$ gives an invertible
multiplication operator. Rank is consequently constant.

Even invariance follows from \eqref{m2:leg:gram-even}.
For odd invariance put $\vartheta(A)=c_WA^tc_W^{-1}$.
Transpose preserves $\mathcal H$: if $A^*=A$, then
$A^t=J\bar AJ$ and $(A^t)^*=J\bar AJ=A^t$.
Conjugation by a unitary similitude also preserves self-adjointness,
so $\vartheta$ preserves $\mathcal H$.
It fixes $C$ pointwise: every element of $C$ is in $E[D]$, and
$c_Wf(D)^tc_W^{-1}=f(D)$.
Furthermore \eqref{m2:leg:transpose-centralizer} gives
\[
 \vartheta(hAh^{-1})=h^{-1}\vartheta(A)h.
\]
Thus on the nontrivial root summand it intertwines multiplication
by $q$ with multiplication by $q^{-1}$, $q\in K_T^1$.
It follows that it has the form $w\mapsto\eta\bar w$:
composing with field conjugation gives an $F$-linear map commuting
with $K_T^1$, hence with the algebra $K_T$ which that infinite group
spans. Such a map is multiplication by an element $\eta\in K_T^\times$.
Since $\vartheta$ preserves determinants and fixes the $C$-coordinates,
the determinant formulas force $N(\eta)=1$.
It therefore fixes $N(w)$, proving odd invariance of \eqref{m2:leg:an-invariants}.
\end{proof}

\subsubsection{Nonzero Gram fibers}

For a represented $A\in\mathcal H$ let
\[
 \mathcal F_A=\{Q:Q^*Q=A\},\qquad T_A=Z_T(A),
 \qquad S_A=G_2\times_\lambda T_A.
\]
The base $T$-orbit of $A$ is homogeneous for $S$, with stabilizer $S_A$,
and its Gram fiber is $\mathcal F_A$.
By \cref{m2:leg:root-image}, $T_A$ is either $T$ or the scalar torus
$E^\times$; in particular $S_A$ is reductive and unimodular.
Thus Frobenius descent over this base has relative modulus one.

\begin{lemma}\label{m2:leg:gram-factorization}
If $A$ has rank two, $U_{1,1}(F)$ acts simply transitively on
$\mathcal F_A$, and the $S$-stabilizer of $Q\in\mathcal F_A$ is
\begin{equation}\label{m2:leg:rank-two-stabilizer}
 S_Q=\{(QhQ^{-1},h):h\in T_A\}.
\end{equation}
If $A$ has rank one and is nonzero, $\mathcal F_A$ is one rational
$U_{1,1}(F)$-orbit. Its $S$-stabilizer is the torus
\begin{equation}\label{m2:leg:rank-one-stabilizer}
 S_Q\simeq\{(h,s)\in T_A\times E^\times:
                                      N_{E/F}(s)=\lambda(h)\}.
\end{equation}
In both cases $\Delta_S|_{S_Q}/\Delta_{S_Q}=1$.
\end{lemma}

\begin{proof}
For rank two, $Q$ is invertible. If $Q_1^*Q_1=Q_2^*Q_2$, then
$Q_2Q_1^{-1}$ has adjoint times itself equal to $I$. Its action is
the unique unitary element taking $Q_1$ to $Q_2$.
The stabilizer equation $vQ=Qh$ forces $v=QhQ^{-1}$, and the
similitude condition for $v$ is equivalent to $hA=Ah$.
This proves \eqref{m2:leg:rank-two-stabilizer}.

For rank-one nonzero $A$, $Q$ has rank one: it cannot be invertible,
and it cannot be zero. Its image cannot be an isotropic line,
since that would give $Q^*Q=0$. Thus its image is a nondegenerate
Hermitian line. Since $\ker Q\subset\ker A$ and both are lines,
$\ker Q=\ker A$. For two factorizations of $A$ the rule
$Q_1w\mapsto Q_2w$ is therefore well-defined. It is an isometry,
because
\[
 \langle Q_iw,Q_iw'\rangle=\langle w,Aw'\rangle.
\]
To construct the unitary extension, choose $w_0$ whose image under either factorization
is nonzero and put $e_i=Q_iw_0$. Their common norm
$a=\langle w_0,Aw_0\rangle$ is nonzero. Choose a nonzero
$f_i\in(Ee_i)^\perp$ and write $b_i=\langle f_i,f_i\rangle$.
The complementary lines are nondegenerate, so $b_i\ne0$.
Let $B_i$ be the $E$-matrix with columns $e_i,f_i$ in the original
hyperbolic basis. Then
\[
 B_i^\dagger JB_i=\diag(a,b_i),\qquad
 ab_i=-N_{E/F}(\det B_i).
\]
Thus $b_1/b_2=N_{E/F}(\det B_1/\det B_2)$.
The map sending $e_1$ to $e_2$ and $f_1$ to
$(\det B_1/\det B_2)f_2$ is an isometry of both orthogonal lines.
It is therefore the required rational unitary extension on the
whole plane. This proves transitivity and supplies explicitly the
two-dimensional Hermitian Witt extension used here.

If $h$ commutes with $A$, it preserves the Hermitian form
$\langle w,Aw'\rangle$ with multiplier $\lambda(h)$.
On its nondegenerate one-dimensional quotient it therefore has an
eigenvalue of norm $\lambda(h)$. The equation $vQ=Qh$ uniquely
determines $v$ on $\im Q$ with this eigenvalue. On its orthogonal
complement $v$ is an arbitrary scalar $s$ of norm $\lambda(h)$.
This is \eqref{m2:leg:rank-one-stabilizer}; conversely each such choice
is a stabilizing similitude. The stabilizers in \eqref{m2:leg:rank-two-stabilizer} and
\eqref{m2:leg:rank-one-stabilizer} are tori, so they and $S$ are unimodular.
\end{proof}

Consider a nonempty fiber of \eqref{m2:leg:an-invariants} whose Gram
orbit is nonzero. It is one $S(F)$-orbit. Indeed, for two points first
choose $h\in T(F)$ conjugating their Gram matrices. Choose
$v_0\in G_2(F)$ of multiplier $\lambda(h)$, for example
$\diag(1,\lambda(h))$. Then $v_0Q_1h^{-1}$ and $Q_2$ have equal
Gram matrices, and \cref{m2:leg:gram-factorization} connects them by
a unitary left factor.
The preceding lemma supplies its unimodular stabilizer, so its
invariant distribution space is generated by a positive orbital
measure. The complete fiber is $\Theta$-stable by
\cref{m2:leg:gram-orbits}. The pushforward of this positive measure is
a positive scalar multiple of itself, and \eqref{m2:leg:odd-square}
forces that scalar to square to one. It is therefore one.
Consequently every nonzero complete Gram-invariant fiber has zero
sign-distribution space. 

\subsubsection{The zero Gram fiber and its nonreductive stabilizers}

It remains to prove vanishing on the entire closed fiber
\[
 \mathcal N=\Gamma^{-1}(0),\qquad
 \mathcal N^\times=\mathcal N\setminus\{0\}.
\]
For $Q\in\mathcal N$, its image is totally isotropic. The target is
a hyperbolic Hermitian plane, so every nonzero such $Q$ has rank one
and can be written
\begin{equation}\label{m2:leg:null-factorization}
 Q=e\ell,\qquad e\ne0\text{ isotropic},\quad
 \ell\in W^\vee\setminus\{0\}.
\end{equation}
Its kernel is an $E$-line. Two such maps lie in the same $S$-orbit
exactly when their kernel lines lie in the same $T$-orbit.
Necessity follows from \eqref{m2:leg:action}. For sufficiency, use $h$
to align the kernel lines. The remaining maps have proportional
nonzero functionals, and their isotropic image vectors can be
matched by a similitude of any prescribed multiplier.
Explicitly, in hyperbolic bases beginning with the respective
image vectors, an element with first eigenvalue $a\in E^\times$
and prescribed multiplier $\lambda$ has second eigenvalue
$\lambda/\bar a$. Taking $\lambda=\lambda(h)$ gives the required
left factor, with no norm restriction on $a$.

Under the identification $W\simeq M$ of
\eqref{m2:leg:trace-hermitian}, the open locus of projective lines
with an invertible representative is
$M^\times/E^\times$. Its $T$-orbit space is
\begin{equation}\label{m2:leg:null-quotient}
 \Xi=M^\times/(E^\times T(F))=M^\times/T(F),
\end{equation}
since $E$-scalars already lie in $T$.
For a field $M$, every nonzero vector is invertible, so this is the
whole projective-line quotient. For the product algebra the missing
lines are its two idempotent eigenlines, which we treat below.

The group $T(F)$ is closed in $M^\times$, and its projective action
factors through the compact group $T(F)/E^\times\simeq K_T^1$.
The quotient of the compact projective $E$-line by this compact action
is Hausdorff and totally disconnected. For the latter assertion,
given an orbit and an invariant open neighborhood of it, compactness
of the orbit gives a clopen set $C$ containing the entire orbit and
contained in that neighborhood. Compactness and continuity show
that $C$ has open stabilizer in the compact acting group, so its
translates form a finite family. The intersection of this family
still contains the entire orbit and is invariant and clopen.
The images of these neighborhoods provide
a clopen basis in the quotient. The image of the invariant
open invertible locus is open, and identifies with \eqref{m2:leg:null-quotient}.
In particular $\Xi$ is a separated locally compact totally
disconnected space. The kernel-line map on the rank-one locus is
continuous. Its composition with \eqref{m2:leg:null-quotient} is
$S$-invariant, with each fiber a whole rational $S$-orbit, by the
preceding orbit calculation.

\begin{lemma}\label{m2:leg:null-stabilizer}
Let $Q=e\ell\in\mathcal N^\times$, let $T_{\ker\ell}$ be the
stabilizer of its kernel line, and define $a(h)$ by
$\ell h=a(h)\ell$. In a hyperbolic basis beginning with $e$ its
stabilizer consists of
\begin{equation}\label{m2:leg:null-matrices}
 h\in T_{\ker\ell},\qquad
 v=\begin{pmatrix}
 a(h)&a(h)u\\0&\lambda(h)/\overline{a(h)}
 \end{pmatrix},\quad u\in E^-.
\end{equation}
It has the form $E^-\rtimes T_{\ker\ell}$, and
\begin{equation}\label{m2:leg:null-modulus}
 \Delta_{S_Q}(h)=
       \left|N_{E/F}(a(h))/\lambda(h)\right|_F.
\end{equation}
For every kernel line in the present anisotropic cases this modulus
equals one.
\end{lemma}

\begin{proof}
The stabilizer equation is $v e\ell=e\ell h$, so it requires $h$
to fix the kernel and $ve=a(h)e$. The unitary-similitude equation
with multiplier $\lambda(h)$ then forces exactly
\eqref{m2:leg:null-matrices}; its upper-right condition is
$u+\bar u=0$. Conjugation by the diagonal part of $v$ acts on
the $F$-line $E^-$ by
\[
 u\longmapsto\frac{N(a(h))}{\lambda(h)}u.
\]
The torus has modulus one, and this determinant on the unipotent
radical gives \eqref{m2:leg:null-modulus}, with the convention of
\cref{m2:tools:frobenius}.

For an invertible line $Ew$, the relation $hw=aw$ forces $h=a\in E^\times$,
so $T_{Ew}=E^\times$. In this case $N(a(h))=\lambda(h)$.
For either exceptional product eigenline, its stabilizer is all $T$;
$a(h)$ is one of the two components of $h$, and
\eqref{m2:leg:an-product} gives $N(h_1)=N(h_2)=\lambda(h)$.
Thus the quotient is one in every case under consideration.
\end{proof}

We next compute the odd action on the rational orbits, using
\eqref{m2:leg:transport-condition}.

\begin{lemma}\label{m2:leg:null-odd}
The operator $\Theta$ preserves each $S(F)$-orbit in $\mathcal N^\times$.
\end{lemma}

\begin{proof}
On kernel lines the induced anti-linear map is
\[
 j(w)=c_WJ\bar w.
\]
Indeed $\ker\kappa(Q)=J\overline{\ker Q}$ and the right factor
$c_W^{-1}$ in \eqref{m2:leg:odd} then multiplies this kernel by $c_W$.
Put $\mu=\lambda(D)$. The unitary equation and the transporter give
\[
 jDj^{-1}=c_WJ\bar DJc_W^{-1}
          =\mu c_WD^{-t}c_W^{-1}
          =\mu D^{-1}=D^*.
\]
Since $j$ conjugates $E$-scalars by the nontrivial field automorphism,
it follows that $j(rw)=r^*j(w)$ for every $r\in M$.
Writing $W$ as $M$ and putting $a=j(1)$ therefore gives
\begin{equation}\label{m2:leg:null-antilinear}
 j(w)=aw^*,\qquad a\in M^\times.
\end{equation}
Invertibility of $a$ follows from invertibility of $j$.

The map $J\bar{\phantom w}$ is anti-unitary for the original form
with matrix $J$. Thus $j$ is an anti-similitude of multiplier
$\lambda(c_W)$: $\langle j(u),j(v)\rangle
=\lambda(c_W)\overline{\langle u,v\rangle}$.
In the trace coordinates \eqref{m2:leg:trace-hermitian}, this says
\[
 \Tr_{M/E}(c a^*a u^*v)
       =\lambda(c_W)\Tr_{M/E}(c u^*v)\qquad(u,v\in M).
\]
Nondegeneracy of the trace pairing and $c\in M^\times$ imply
$a^*a=\lambda(c_W)\in F^\times$, so $a\in T(F)$.
For invertible $w$,
\[
 (w^*/w)^*(w^*/w)=1,
\]
and hence $w^*/w\in T(F)$. Formula
\eqref{m2:leg:null-antilinear} therefore fixes the class of $w$ in
$\Xi$. In the componentwise product case $*$ fixes each idempotent;
multiplication by $a$ also preserves each eigenline. Those two
exceptional projective lines are consequently fixed individually.
The kernel-orbit classification proves the lemma.
\end{proof}

Every orbit just considered has the unimodular stabilizer
of \cref{m2:leg:null-stabilizer}. By Frobenius descent its invariant
distribution space is one-dimensional, generated by a positive
orbital measure. The odd operator preserves its orbit by
\cref{m2:leg:null-odd}; its positive scaling constant squares to one
by \eqref{m2:leg:odd-square}. Thus that measure is fixed.

By \cref{m2:tools:localization}, vanishing on the open
invertible-kernel locus in $\mathcal N^\times$ reduces to the complete
fibers of the map into $\Xi$. Every such fiber is an orbit already
shown to have zero sign space, so vanishing holds on the whole open locus. In the product case its
complement in $\mathcal N^\times$ is the union of the two closed
exceptional orbits. Each is odd-stable and has the same vanishing;
\cref{m2:tools:support} completes the passage to
$\mathcal N^\times$. There is no exceptional locus in the field case.
It remains to pass from $\mathcal N^\times$ to $\mathcal N$.
Any remaining sign distribution is supported at zero. Such a
distribution on this non-Archimedean space is a scalar multiple of
$\delta_0$, which is fixed by $\Theta$. Therefore
\begin{equation}\label{m2:leg:null-vanishing}
 \Ddist(\mathcal N)^{S,\Theta=-1}=0.
\end{equation}

\subsubsection{Conclusion of the regular theorem}

For the anisotropic case, every complete fiber of
\eqref{m2:leg:an-invariants} is empty, a single orbit with nonzero Gram
matrix, or the zero Gram fiber. The first case contributes nothing,
the second was treated by \cref{m2:leg:gram-factorization} and the
positive-measure argument, and the third is
\eqref{m2:leg:null-vanishing}. The map is continuous and invariant
under the whole extended action by \cref{m2:leg:gram-orbits}.
Bernstein localization therefore proves \cref{m2:leg:regular} in both
anisotropic forms. The split case was proved using
\cref{m2:leg:split-intertwining,m1:thm:two-leg}. This completes the theorem.

The normal representation for a regular rank-zero orbit also
contains the two-dimensional summand
$\mathfrak u_2\cap E[D]$. It is even-fixed, and
$c_Wf(D)^tc_W^{-1}=f(D)$ makes it odd-fixed as well.
Consequently \cref{m2:tools:passive} applies to this summand and to the
multiplier line in $N_x$.

\subsection{The unitary identity descendant}\label{m2:id:section}\label{m2:sec:identity}

The normal representation for a scalar rank-zero orbit has
an additional moving Lie-algebra summand. Let
\[
 S=G_2\times_\lambda G_2,\qquad
 \mathcal L_{\mathrm{id}}=\Mat_2(E)\oplus\mathfrak u_2,
\]
with action
\begin{equation}\label{m2:id:action}
 (v,h)(Q,Z)=(vQh^{-1},hZh^{-1}).
\end{equation}
For $\mu\in F^\times$, put $A=a_\mu=\diag(1,\mu)$,
and let $D=\zeta I$ with $N_{E/F}(\zeta)=\mu$.
At the scalar representative $A\oplus D$, which is transpose-fixed,
\cref{m2:nor:zero-odd} gives the odd map
\begin{equation}\label{m2:id:odd}
 \Theta(Q,Z)=\bigl(A^{-1}\kappa(Q)\zeta,Z^t\bigr)
             =\bigl(-A^{-1}J\bar QJ\zeta,Z^t\bigr).
\end{equation}
Its left and right factors on $Q$ are $c_U=A^{-1}$ and
$c_W=\zeta^{-1}I$, both of multiplier $\mu^{-1}$.
The map normalizes the even action and is involutive: the two
factors are symmetric, so \eqref{m2:leg:odd-square} gives the identity
on $Q$, while transpose squares to the identity on $Z$.

\begin{theorem}\label{m2:id:vanishing}
For \eqref{m2:id:action} and \eqref{m2:id:odd},
\[
 \Ddist(\Mat_2(E)\oplus\mathfrak u_2)^{S,\Theta=-1}=0.
\]
\end{theorem}

The proof first localizes the nonzero values of a quadratic form on
the moving Lie-algebra factor. Those fibers are the regular two-leg
problems of \cref{m2:leg:regular}. A partial Fourier transform then
places both a distribution and its transform on the same
three-dimensional quadratic null cone. The odd-dimensional
homogeneity theorem finishes the argument.

\subsubsection{The passive central line and the quadratic space}

The equation $Z^*=-Z$ gives
\[
 \mathfrak u_2=
 \left\{\begin{pmatrix}a&b\\c&-\bar a\end{pmatrix}:
          a\in E,\ b,c\in E^-\right\}.
\]
Writing $a=(a-\bar a)/2+(a+\bar a)/2$ yields the direct sum
\begin{equation}\label{m2:id:lie-splitting}
 \mathfrak u_2=E^-I\oplus\mathfrak{su}_{1,1},\qquad
 \mathfrak{su}_{1,1}=
 \left\{\begin{pmatrix}x&b\\c&-x\end{pmatrix}:
          x\in F,\ b,c\in E^-\right\}.
\end{equation}
The central line is fixed pointwise by conjugation and by ordinary
transpose. Projection to it is therefore invariant under the entire
extended action. By \cref{m2:tools:passive}, it suffices to prove the
theorem on $\Mat_2(E)\oplus W$, where
$W=\mathfrak{su}_{1,1}$.

The rational linear map
\begin{equation}\label{m2:id:sl2-map}
 \Phi:W\longrightarrow\mathfrak{sl}_2(F),\qquad
                  \Phi(Z)=P^{-1}ZP
\end{equation}
is a bijection. Indeed, writing $b=\delta y$, $c=\delta z$ gives
\[
 Z=\begin{pmatrix}x&\delta y\\\delta z&-x\end{pmatrix},
 \qquad \Phi(Z)=\begin{pmatrix}x&dy\\z&-x\end{pmatrix}.
\]
All three variables range freely over $F$ and $d\ne0$.
By \cref{m2:simconj:factorization}, every $h\in G_2(F)$ has the form
$h=sPgP^{-1}$, and
\[
 \Phi(hZh^{-1})=g\Phi(Z)g^{-1}.
\]
Conversely every $g\in\GL_2(F)$ occurs through $PgP^{-1}$.
Thus rational conjugation on $W$ is the rational $\GL_2$ adjoint
action, with its scalar kernel. Ordinary transpose
is transported to
\begin{equation}\label{m2:id:transported-transpose}
 X\longmapsto RX^tR^{-1}.
\end{equation}

On $W$ define
\begin{equation}\label{m2:id:quadratic-form}
 B_u(Z,Z')=\tr(ZZ'),\qquad q_u(Z)=\tfrac12\tr(Z^2).
\end{equation}
In the coordinates of \eqref{m2:id:sl2-map},
\[
 q_u(Z)=x^2+dyz=-\det\Phi(Z),
 \qquad
 B_u(Z,Z')=2xx'+d(yz'+zy').
\]
The determinant of the Gram matrix of $B_u$ in coordinates
$(x,y,z)$ is $-2d^2$, which is nonzero in $F$. Thus $W$ is a
nondegenerate quadratic space of $F$-dimension three.
Conjugation preserves $B_u$ by the trace identity, and transpose
does also, since $\tr(Z^tZ'^t)=\tr(Z'Z)=\tr(ZZ')$.
The first $G_2$ factor acts trivially on $W$. It follows that the whole
extended action on $\Mat_2(E)\oplus W$ preserves $q_u$ in the second
variable.

\subsubsection{Nonzero quadratic levels and regular fibers}

For $c\in F^\times$ put
\[
 \mathcal O_c=\{Z\in W:q_u(Z)=c\},\qquad
 X_c=\begin{pmatrix}0&c\\1&0\end{pmatrix},\qquad
 Z_c=PX_cP^{-1}.
\]

\begin{lemma}\label{m2:id:nonzero-orbit}
The level $\mathcal O_c$ is one rational $G_2(F)$-orbit. Its
centralizer $T_c=Z_{G_2}(Z_c)$ is a regular torus, split modulo the
center if $c$ is a square and anisotropic otherwise. Its stabilizer
for the action of $S$ is $S_c=G_2\times_\lambda T_c$.
\end{lemma}

\begin{proof}
Under \eqref{m2:id:sl2-map}, every matrix on this level has
characteristic polynomial $t^2-c$. Since $c\ne0$ and $F$ has
characteristic zero, this polynomial is separable. A trace-zero
matrix with this polynomial is not scalar, so its minimal polynomial
also has degree two. A cyclic vector gives a rational basis in which
the matrix is $X_c$. The rational adjoint action identified above
therefore gives a single $G_2(F)$-orbit.
Its centralizer is
\[
 T_c(F)=E^\times P F[X_c]^\times P^{-1}.
\]
The algebra $F[X_c]\simeq F[t]/(t^2-c)$ is split exactly when $c$ is
a square. Finally $S$ acts on the base through its right factor $h$,
so fixing $Z_c$ means $h\in T_c$, with $v$ subject only to the
common-multiplier condition. This is exactly $S_c$.
\end{proof}

The projection
$\Mat_2(E)\times\mathcal O_c\longrightarrow\mathcal O_c$
is equivariant, with fiber $\Mat_2(E)$ at $Z_c$ and stabilizer
$S_c$. Both $S$ and $S_c$ are reductive and unimodular, so Frobenius
descent carries the even invariant distributions to ordinary
$S_c$-invariant distributions on this fiber, with no modular twist.
To use it for the odd character we now compute the transporter and
the induced odd operator explicitly.

From \eqref{m2:id:transported-transpose},
\[
 RX_c^tR^{-1}=\begin{pmatrix}0&d\\c/d&0\end{pmatrix}.
\]
The rational diagonal matrix $g_c=\diag(c,d)$ conjugates this
matrix to $X_c$. Therefore
\begin{equation}\label{m2:id:level-transporter}
 h_c=Pg_cP^{-1}\in G_2(F),\qquad
 h_cZ_c^th_c^{-1}=Z_c,\qquad \lambda(h_c)=cd.
\end{equation}
Choose $v_c\in G_2(F)$ with $\lambda(v_c)=cd$; one explicit choice
is $v_c=h_c$. The pair $(v_c,h_c)\in S(F)$ corrects the base
action of $\Theta$ at $Z_c$. On the $Q$-fiber its odd action is
\begin{equation}\label{m2:id:fiber-odd}
 Q\longmapsto v_cA^{-1}\kappa(Q)\zeta h_c^{-1}
                 =c_U'\kappa(Q)(c_W')^{-1},
 \qquad c_U'=v_cA^{-1},\quad c_W'=h_c\zeta^{-1}.
\end{equation}
These factors satisfy
\[
 \lambda(c_U')=\lambda(c_W')=cd/\mu,
 \qquad c_W'Z_c^t(c_W')^{-1}=Z_c.
\]
To apply \cref{m2:leg:regular}, observe that $Z_c^2=cI$ and $Z_c^*=-Z_c$,
so
\[
 Z_c^*Z_c=-cI.
\]
Thus $Z_c$ itself is an invertible regular semisimple element of
$G_2(F)$, of multiplier $-c$. Taking $D=Z_c$ in that theorem gives
the torus $T_c$, the even group $S_c$, and the odd map
\eqref{m2:id:fiber-odd}.

The homogeneous base has a positive $S$-invariant measure because
$S$ and $S_c$ are unimodular. Transpose preserves its orbit and
normalizes the action; its positive measure-scaling constant is
one because its square is the identity. The correcting even
element also preserves the base measure. Naturality of the
Frobenius construction with this measure therefore intertwines
the odd action on the total space with the corrected fiber action
\eqref{m2:id:fiber-odd}.
Thus \cref{m1:prop:odd-Frobenius} and \cref{m2:leg:regular} prove
\begin{equation}\label{m2:id:nonzero-vanishing}
 \Ddist(\Mat_2(E)\times\mathcal O_c)^{S,\Theta=-1}=0
                         \qquad(c\ne0).
\end{equation}
The sign is trivial on the even fiber stabilizer and is minus one
on its corrected odd coset.

\subsubsection{Support and partial Fourier transform}

Let
\[
 T\in\Ddist(\Mat_2(E)\oplus W)^{S,\Theta=-1}.
\]
On the invariant open set $q_u(Z)\ne0$, the continuous map
$(Q,Z)\mapsto q_u(Z)\in F^\times$ has the complete fibers
treated in \eqref{m2:id:nonzero-vanishing}. Bernstein localization
therefore makes the restriction of $T$ to this open set zero.
The supported-distribution exact sequence gives
\begin{equation}\label{m2:id:support}
 \supp(T)\subset\Mat_2(E)\times\{Z\in W:q_u(Z)=0\}.
\end{equation}

Choose a nontrivial additive character $\psi$ of $F$ and the
self-dual Haar measure for $B_u$. Let $\mathcal F_W$ be partial
Fourier transform in the $W$-variable, using the kernel
$\psi(B_u(Z,Z'))$. Every even or odd transformation on $W$ is an
isometry of $B_u$. If its matrix is $L$ and the Gram matrix of the
pairing is $H$, the identity $L^tHL=H$ implies
$|\det L|_F=1$. Changing variables in the Fourier integral
therefore introduces no absolute determinant factor and gives
\[
 \mathcal F_W L_*=L_*\mathcal F_W.
\]
The same statement holds on distributions by duality. The
$Q$-transformations in \eqref{m2:id:action} and \eqref{m2:id:odd} are
linear in $Q$ and independent of $Z$; they commute with this partial
transform. Hence
\[
 \widehat T=\mathcal F_WT
       \in\Ddist(\Mat_2(E)\oplus W)^{S,\Theta=-1}.
\]
Applying the same open localization to $\widehat T$ gives the
second support statement
\begin{equation}\label{m2:id:fourier-support}
 \supp(\widehat T)\subset
              \Mat_2(E)\times\{Z\in W:q_u(Z)=0\}.
\end{equation}
We can now invoke the non-Archimedean homogeneity theorem
\cref{m2:tools:homogeneity} (\cite[Theorem~5.1.5]{AGDescent})
with parameter manifold $\Mat_2(E)$ and nondegenerate quadratic
space $(W,B_u)$. Its quadratic null cone is
$\{Z:B_u(Z,Z)=0\}=\{Z:q_u(Z)=0\}$, since $2\ne0$ in $F$.
Its two hypotheses are exactly
\eqref{m2:id:support} and \eqref{m2:id:fourier-support}, using partial
Fourier transform in $W$. The final clause of that theorem gives
zero whenever $\dim_F W$ is odd. Here this dimension is three, so
$T=0$. Restoring the passive line in \eqref{m2:id:lie-splitting}
proves \cref{m2:id:vanishing}.

\begin{remark}\label{m2:id:nilpotent-remark}
The zero level of $q_u$ consists of the origin and the single nonzero
rational nilpotent adjoint orbit. The latter assertion follows from
the rational $\GL_2$ adjoint description in \eqref{m2:id:sl2-map}.
The support inclusions and homogeneity apply to this entire level.
\end{remark}

\subsection{Proof of the unitary multiplicity-one theorem}
\label{m2:sec:completion}

\subsubsection{Invariance under the corrected involutions}

The nine normal representations, with their rational stabilizers,
odd operators, and parameter domains, are listed in
\cref{m2:nor:classification}. Their vanishing results are:
\begin{itemize}
\item \cref{m2:id:vanishing} for the scalar rank-zero family;
\item \cref{m2:leg:regular} for both regular rank-zero families;
\item \cref{m2:simconj:vanishing} for the dense scalar family;
\item \cref{m2:elem:vanishing} for the other dense families and both
rank-one families.
\end{itemize}

\begin{proposition}\label{m2:completion:normal-vanishing}
For every closed rational $K(F)$-orbit in $G_4(F)$, with
representative $x$ and extended stabilizer $\Kt_x(F)$,
\[
 \Ddist(N_x(F))^{\Kt_x(F),\eps}=0.
\]
\end{proposition}

\begin{proof}
By \cref{m2:rat:classification,m2:rat:transpose}, the nine rational
families are exhaustive and every extended stabilizer has an odd
coset. Its normal action is the one in \cref{m2:nor:classification}.
In the two regular rank-zero rows, the transporter factors have
equal multipliers and satisfy $c_WD^tc_W^{-1}=D$. Thus
\cref{m2:leg:regular} applies to the $M_2(E)$ summand with its
corrected odd action. The remaining $F^2$ is fixed by both actions,
so \cref{m2:tools:passive} adjoins it. For the scalar rank-zero
row, \cref{m2:id:vanishing} applies to the entire reduced normal space
$M_2(E)\oplus\mathfrak u_2$. The remaining six rows are covered by
\cref{m2:elem:vanishing,m2:simconj:vanishing}.
Finally \cref{m2:nor:full} gives
$N_x=N_x^0\oplus F[I_4/2]$, with the last line fixed by the
extended stabilizer. Applying \cref{m2:tools:passive} proves the assertion.
\end{proof}

\subsubsection{Generalized descent for the residual action}

\begin{theorem}\label{m2:completion:residual}
Every $K(F)$-invariant distribution on $G_4(F)$ is fixed by
ordinary transpose. Equivalently,
\[
 \Ddist(G_4(F))^{\Kt(F),\eps}=0.
\]
\end{theorem}

\begin{proof}
We apply \cref{m1:cor:HC-invariance} to the action of $K$ on
$G_4$, taking $\tau$ to be ordinary transpose. The group $K$ is connected
reductive, and $G_4$ is smooth affine. The action
\eqref{m2:red:action}, transpose, and \eqref{m2:red:alpha} are
algebraic over $F$. By \cref{m2:rat:transpose}, every closed rational
orbit is transpose-stable, and its extended stabilizer is
$K_x(F)\sqcup K_x(F)k_x\sigma$. Its normal space and sign action
are those of \cref{m2:nor:full};
\cref{m2:completion:normal-vanishing} supplies the required vanishing.
The corollary therefore gives transpose invariance, equivalently
the required vanishing of sign-equivariant distributions.
\end{proof}

\begin{corollary}\label{m2:completion:biH}
Every bi-$H_{\mathrm U}$-invariant distribution on $G_{\mathrm U}$ is fixed by
$\tau(g,z)=(g^t,z^t)$.
\end{corollary}

\begin{proof}
The modular factor and compatibility with the invariant base measure
for the isomorphism \eqref{m2:red:distiso} were verified in
\cref{m2:red:frobenius}. Its image is fixed by fiber transpose
by \cref{m2:completion:residual}, so its inverse image is fixed by
$\tau$.
\end{proof}

\subsubsection{The product inequality and equality of multiplicities}

\begin{proof}[Proof of \cref{thm:main-unitary}]
Fix an irreducible admissible smooth complex representation
$\pi$ of $G_{\mathrm U}$, and put
$m(\pi)=\dim\Hom_{H_{\mathrm U}}(\pi,\mathbb C)$.
The groups $G_{\mathrm U}$ and $H_{\mathrm U}$ are unimodular, so $G_{\mathrm U}/H_{\mathrm U}$
has an invariant measure. The anti-involution preserves the
subgroup by \cref{m2:red:transpose}, and fixes bi-invariant
distributions by \cref{m2:completion:biH}. Thus the
non-Archimedean Gelfand--Kazhdan criterion
\cite[Lemma~4.2]{Prasad90}, in the form of \cref{m2:tools:gk},
gives
\[
 m(\pi)m(\pi^\vee)\leq1.
\]
Finally \cref{m2:red:duality}, obtained from
Roche--Vinroot's dualizing theorem and the preservation of $H_{\mathrm U}$
by inverse transpose, gives $m(\pi)=m(\pi^\vee)$. Their common
dimension is therefore at most one.
\end{proof}

\section{The linear Shalika model}\label{sec:gl-shalika}
We deduce Rallis's multiplicity-one theorem for twisted Jacquet
modules, as recorded in Prasad~\cite[Theorem~1]{PrasadShalika}, from
Theorem~\ref{thm:main-linear}. This argument also gives the theta-transfer
construction used for the unitary similitude groups below.

Write
\[
 N_{2,2}=\left\{n(X)=\begin{pmatrix}I_2&X\\0&I_2\end{pmatrix}:
                                  X\in\operatorname{M}_2(F)\right\},
 \qquad S_4=\Delta\GL_2(F)\ltimes N_{2,2},
\]
where $\Delta g=\diag(g,g)$, and put
$\psi_N(n(X))=\psi(\tr X)$. Since trace is invariant under
conjugation, $\sigma\boxtimes\psi_N$ defines a representation
of $S_4$ for every smooth representation $\sigma$ of $\GL_2(F)$.

\begin{corollary}[Rallis; see Prasad]\label{cor:gl-shalika}
For irreducible admissible representations $\Pi$ of $\GL_4(F)$
and $\sigma$ of $\GL_2(F)$,
\[
 \dim\Hom_{S_4}(\Pi,\sigma\boxtimes\psi_N)\leq1.
\]
Equivalently, $\dim\Hom_{\GL_2(F)}(\Pi_{N_{2,2},\psi_N},\sigma)\leq1$.
\end{corollary}

We use the equal-rank type-II oscillator normalization on
$\mathcal S(\operatorname{M}_4(F))$ given by
\[
 ((g,h)\phi)(x)=|\det(g^{-1}h)|_F^2\phi(g^{-1}xh).
\]
The big theta lift $\Theta(\tau)$ is characterized by the maximal
$\tau$-isotypic quotient $\tau\boxtimes\Theta(\tau)$ of this
oscillator. We use two facts about this lift. First, with
$L=\GL_2(F)\times\GL_2(F)$ embedded block diagonally,
Gan's big-lift period identity is
\begin{equation}\label{eq:gl-big-transfer}
 \Hom_{S_4}(\Theta(\tau),\sigma\boxtimes\psi_N)
 \simeq\Hom_L(\tau^\vee,\sigma\boxtimes\mathbf1).
\end{equation}
See \cite[Theorem~3.1]{GanPeriods}; the big-lift formulation
also appears in \cite[Proposition~3.5]{YangV1}.
This is an identity for the big lift of an arbitrary irreducible
admissible $\tau$. Second, equal-rank type-II correspondence gives
\begin{equation}\label{eq:gl-theta-quotient}
 q_\Pi:\Theta(\Pi^\vee)\twoheadrightarrow\Pi
\end{equation}
for every irreducible admissible $\Pi$.
This follows from M\'{i}nguez~\cite[Theorem~1]{Minguez08};
for the oscillator normalization above it is recalled in
\cite[Theorem~1.2 and the preceding paragraph]{FangSunXue}.

\begin{proof}
Precomposition with the surjection \eqref{eq:gl-theta-quotient},
followed by \eqref{eq:gl-big-transfer} with $\tau=\Pi^\vee$, gives
\begin{equation}\label{eq:gl-precomposition}
 \Hom_{S_4}(\Pi,\sigma\boxtimes\psi_N)
 \lhook\joinrel\longrightarrow
 \Hom_L(\Pi,\sigma\boxtimes\mathbf1).
\end{equation}
Here admissibility identifies $(\Pi^\vee)^\vee$ with $\Pi$.
To place the coefficient on the second block, let
$w=\left(\begin{smallmatrix}0&I_2\\I_2&0\end{smallmatrix}\right)$.
The block swap and evaluation pairing give an injection
\[
 \begin{aligned}
 \Hom_L(\Pi,\sigma\boxtimes\mathbf1)
 &\lhook\joinrel\longrightarrow
 \Hom_{H_{\mathrm L}}(\Pi\boxtimes\sigma^\vee,\mathbb C),\\
 A&\longmapsto
 \bigl[v\otimes\eta\longmapsto
          \langle A(\Pi(w)v),\eta\rangle\bigr].
 \end{aligned}
\]
Indeed, $w\diag(a,b)=\diag(b,a)w$ and the evaluation pairing is
$\GL_2(F)$-invariant. Injectivity follows because the smooth dual
of $\sigma$ separates vectors and $\Pi(w)$ is invertible.
Theorem~\ref{thm:main-linear} bounds the last Hom space by one,
since $\Pi\boxtimes\sigma^\vee$ is irreducible admissible
by~\cite[Proposition~2.16]{BZ76}.
\end{proof}

\section{Unitary similitude oscillator and Jacquet module}\label{sec:gu-transfer}
We relate the unitary spherical pair to Shalika periods with
coefficient representations through the similitude theta correspondence. For ordinary
unitary groups, the theta-period identity is given in
Gan~\cite[\S2]{GanPeriods} and Chen--Gan
\cite[\S3, Proposition~3.4 and Theorem~3.5]{ChenGan}.
Here the multiplier coordinate enters the rational orbit and the
stabilizer action. Theorem~\ref{thm:main-unitary} then bounds the
resulting Hom space.
\label{gu:sec:setup}

\subsection{The unitary similitude dual pair}
We use the notation of the introduction and regard vectors as column
vectors, with groups acting on the left. Let $J_4$ be the antidiagonal
identity matrix and set
\[
 V=E^4,\qquad h_V(v,v')=\bar v'^{\,t}J_4v.
\]
Let $W=X\oplus Y$, with $X=Y=E^2$, and put
\begin{equation}\label{gu:eq:W-form}
 h_W((x,y),(x',y'))=x^t\bar y'-y^t\bar x',
 \qquad J_W=\begin{pmatrix}0&I_2\\-I_2&0\end{pmatrix}.
\end{equation}
Thus
\[
 G_\theta=\GU(V)=\{g:g^\dagger J_4g=\lambda_{G_\theta}(g)J_4\},\qquad
 G'_\theta=\GU(W)=\{h:h^\dagger J_Wh=\lambda_{G'_\theta}(h)J_W\}.
\]
Write $G_\theta^1=U(V)$ and $(G'_\theta)^1=U(W)$ for their kernels. When the group is clear, we abbreviate its multiplier to $\lambda$. Multiplication of $h_W$ by a nonzero trace-zero element of $E$ identifies $G'_\theta$ with the split Hermitian similitude group $\GU_{2,2}(F)$, preserving its multiplier. The skew-Hermitian realization is adapted to the dual pair.

Fix a hyperbolic Hermitian form on $X$,
\[
 B(x,x')=\bar x'^{\,t}\mathsf Bx,\qquad \mathsf B^\dagger=\mathsf B,
\]
and put $K=\GU(B)$. One may take $\mathsf B=J_2$. Keeping $\mathsf B$ in the formulas allows arbitrary hyperbolic bases.

The conjugate space and its form are
\[
 X^c=X\otimes_{E,\bar{\ }}E,
 \qquad x\otimes z=\bar z x\otimes1,
 \qquad B^c(x\otimes1,x'\otimes1)=B(x',x).
\]
For $a\in\GL(X)$ define $\iota(a)(x\otimes1)=ax\otimes1$. In conjugate bases, $\iota(a)$ has matrix $\bar a$, and $B^c$ has matrix $\overline{\mathsf B}$. Consequently
\begin{equation}\label{gu:eq:iota-det}
 \iota:K\xrightarrow{\sim}\GU(B^c),\qquad
 \lambda(\iota(a))=\lambda_B(a),\qquad
 \det_{X^c}\iota(a)=\overline{\det_X a}.
\end{equation}
These are the conjugate-space conventions of~\cite[Sections~2.2--2.4]{ChenGan}.

\subsection{The Shalika character and subgroup}
The Siegel Levi and unipotent radical in $G'_\theta$ are
\begin{equation}\label{gu:eq:siegel-coordinates}
 m(a,t)=\begin{pmatrix}a&0\\0&t\bar a^{-t}\end{pmatrix},
 \qquad
 n(u)=\begin{pmatrix}I_2&u\\0&I_2\end{pmatrix},
 \quad u^\dagger=u.
\end{equation}
Here $a\in\GL_2(E)$, $t\in F^\times$, and $\lambda(m(a,t))=t$. Let $P^{\GU}$ be the Siegel parabolic with this Levi and unipotent radical $N(X)$. The form $B$ corresponds to the opposite-unipotent map $b_B:X\to Y$ with matrix
\[
 b_B=-\mathsf B.
\]
Indeed, $h_W(b_Bx,x')=-\bar x'^{\,t}b_Bx=B(x,x')$.

We use the additive character $\psi$ fixed in the introduction and
write $\tr_E$ for matrix trace over $E$. Our normalization is
\begin{equation}\label{gu:eq:psiB}
 \psi_B(n(u))
 =\psi\left(-\tfrac12\tr_E(u\mathsf B)\right)
 =\psi\left(-\tfrac14\Tr_{E/F}\tr_E(u\mathsf B)\right).
\end{equation}
Here $\tr_E(u\mathsf B)\in F$, because both matrices are Hermitian.

Write $\operatorname{Herm}_2(E)$ for the space of Hermitian $2$ by $2$ matrices. The trace pairing $(u,D)\mapsto\tr_E(uD)$ on this space
is nondegenerate: diagonal test matrices recover the diagonal entries,
and the off-diagonal entries are detected by the field-trace pairing.
Since
\begin{equation}\label{gu:jac:levi-conjugation}
 m(a,t)n(u)m(a,t)^{-1}=n(t^{-1}aua^\dagger),
\end{equation}
we have
\[
 \psi_B(m(a,t)n(u)m(a,t)^{-1})
 =\psi\left(-\tfrac12\tr_E(u\,t^{-1}a^\dagger\mathsf Ba)\right).
\]
Nondegeneracy and the nontriviality of $\psi$ show that the character
stabilizer is given by
\begin{equation}\label{gu:jac:character-stabilizer}
 a^\dagger\mathsf Ba=t\mathsf B.
\end{equation}
Thus
\begin{equation}\label{gu:eq:shalika-embedding}
 m_B(a)=m(a,\lambda_B(a)),\qquad
 S_B^{\GU}=m_B(K)N(X)\simeq K\ltimes N(X).
\end{equation}
On $\operatorname{Herm}_2(E)$ the map $u\mapsto aua^\dagger$ has absolute
determinant $|\det a|_E^2$. For diagonal $a$, this follows from its two
norm-scaled diagonal coordinates and its off-diagonal $E$-coordinate;
elementary unipotents have determinant one. Consequently
\begin{equation}\label{gu:jac:siegel-modulus}
 \delta_{P^{\GU}}(m(a,t))=|\det a|_E^2|t|_F^{-4}.
\end{equation}
Taking determinants in \eqref{gu:jac:character-stabilizer} gives
\begin{equation}\label{gu:jac:plane-determinant}
 N_{E/F}(\det a)=t^2,\qquad |\det a|_E=|t|_F^2.
\end{equation}
Hence $\delta_{P^{\GU}}|_K=1$, the group $K\ltimes N(X)$ is unimodular,
and normalized and unnormalized Jacquet actions agree on $K$.

For a smooth representation $\sigma$ of $K$, we use $\sigma\boxtimes\psi_B$ with the embedding $m_B$. Since
\[
 m_B(zI_2)=m(zI_2,N_{E/F}z)=zI_4,
\]
its central character on $E^\times$ is $\omega_\sigma$. In particular,
\begin{equation}\label{gu:eq:central-necessary}
 \Hom_{S_B^{\GU}}(\Pi,\sigma\boxtimes\psi_B)\ne0
 \quad\Longrightarrow\quad\omega_\Pi=\omega_\sigma.
\end{equation}

\subsection{The Schr\"odinger model}
On the underlying $F$-space of $V\otimes_EW$ put
\begin{equation}\label{gu:eq:tensor-symplectic}
 [v\otimes w,v'\otimes w']
 =-\tfrac12\Tr_{E/F}\bigl(h_V(v,v')h_W(w,w')\bigr).
\end{equation}
This is a nondegenerate alternating form, with polarization $(V\otimes X)\oplus(V\otimes Y)$. The second Lagrangian is identified with
\begin{equation}\label{gu:eq:schrodinger-Z}
 Z=Y\otimes_EV\simeq\Hom_E(X^c,V),
 \qquad
 (y\otimes v)(x\otimes1)=\overline{h_W(x,y)}\,v.
\end{equation}
It has dimension $16$ over $F$. Define the form-valued moment map by
\begin{equation}\label{gu:eq:moment-map}
 Q(T)=T^*h_V,\qquad G_T=T^\dagger J_4T.
\end{equation}
Thus $G_T$ is the matrix of the Hermitian form $Q(T)$ on $X^c$. The corresponding form matrix on $X$ is $\overline{G_T}$. The notation $T^*h_V$ means pullback of a form. If one instead forms an adjoint using $B^c$, the equation $Q(T)=cB^c$ becomes $T^*T=cI_{X^c}$.

For a locally compact totally disconnected space $Z'$, let $\mathcal S(Z')=C_c^\infty(Z',\mathbb C)$. The ordinary oscillator representation acts on $\mathcal S(Z)$. We use Heisenberg multiplication
\[
 (z,t)(z',t')=(z+z',t+t'+\tfrac12[z,z'])
\]
and Schr\"odinger operators
\begin{equation}\label{gu:eq:heisenberg-action}
 [\rho_\psi(x+y,t)f](T)
 =\psi\bigl(t-[x,T]-\tfrac12[x,y]\bigr)f(T+y)
\end{equation}
for $x\in V\otimes X$ and $y,T\in Z$. Weil operators satisfy left covariance
\[
 \omega(s)\rho_\psi(z,t)\omega(s)^{-1}=\rho_\psi(sz,t).
\]
The choices in~\eqref{gu:eq:tensor-symplectic} and~\eqref{gu:eq:heisenberg-action} give~\eqref{gu:eq:psiB}. Initially both ordinary splitting characters are one, which is permissible because both dimensions are even. Splitting changes will be introduced after the normalized calculation.

\subsection{Induction, quotients, and the unitary multiplicity-one theorem}
For a closed subgroup $D\subset L$ and a smooth representation $\rho$ of $D$, we realize the unnormalized compact induction $\cInd_D^L\rho$ as right-smooth functions satisfying
\[
 f(d\ell)=\rho(d)f(\ell),\qquad
 [\pi(\ell_0)f](\ell)=f(\ell\ell_0),
\]
with support compact modulo $D$. Ordinary smooth induction $\Ind_D^L\rho$ uses the same covariance without a compact-support condition. For a smooth representation $M$ the twisted Jacquet module is the unnormalized quotient
\[
 M_{N,\psi_B}
 =M/\langle nv-\psi_B(n)v:n\in N,\ v\in M\rangle.
\]

For $V=V_1\perp V_2$ with hyperbolic planes $V_1,V_2$, we retain
$H_{\mathrm U}$ for the equal-multiplier subgroup
\eqref{eq:main-unitary-subgroup}, to which
Theorem~\ref{thm:main-unitary} applies.

\subsection{The matching subgroup and its normalized extension}
\label{gu:jac:oscillator-section}
The direct tensor action has multiplier $\lambda(g)\lambda(h)$. We use the
equal-multiplier convention
\[
 R=\{(g,h)\in G_\theta\times G'_\theta:\lambda(g)=\lambda(h)\},\qquad
 \alpha_W(h)=\lambda(h)^{-1}h,
\]
and the symplectic action
\begin{equation}\label{gu:jac:tensor-action}
 (g,h)\longmapsto g\otimes\alpha_W(h).
\end{equation}
The involution $\alpha_W$ identifies this with the inverse-multiplier
convention and is the identity on $(G'_\theta)^1$. The tensor kernel is
\begin{equation}\label{gu:jac:tensor-kernel}
 K_E=\{(zI_V,\bar zI_W):z\in E^\times\},
\end{equation}
since $\alpha_W(\bar zI_W)=z^{-1}I_W$.
We use the similitude oscillator formalism of~\cite[\S\S2--4]{ZhangSimilitude}
and~\cite[\S\S3, 6.1]{BGS}. The following formulas fix its normalization.

For the ordinary Weil representation $\omega$ with splitting characters
$(1,1)$, the formulas needed below are
\begin{align}
 [\omega(g,1)f](T)&=f(g^{-1}T), &&g\in U(V),
 \label{gu:jac:ordinary-g}\\
 [\omega(1,m(a,1))f](T)&=|\det a|_E^2f(T\bar a),
 &&a\in\operatorname{GL}_2(E),\label{gu:jac:ordinary-levi}\\
 [\omega(1,n(u))f](T)&=
 \psi\left(-\frac12\operatorname{tr}_E
                     (u\overline{G_T})\right)f(T),
 &&u^\dagger=u.\label{gu:jac:ordinary-n}
\end{align}
These are the ordinary unitary formulas of \cite[\S3.1]{ChenGan} in our matrix
and trace conventions.

Their phase, including its sign and conjugation, follows from
\eqref{gu:eq:heisenberg-action}; the calculation and the Fourier operators
are recorded in Appendix~\ref{gu:app:normalization}.

For $s\in F^\times$, put
\begin{equation}\label{gu:jac:sections}
 r_s=\operatorname{diag}(s,s,1,1)\in G_\theta,
 \qquad d_s=\operatorname{diag}(I_2,sI_2)\in G'_\theta,
 \qquad e_s=(r_s,d_s).
\end{equation}
Both $r_s$ and $d_s$ have multiplier $s$, and these are multiplicative
sections.  If $(g,h)\in R$ has common multiplier $s$, its unique
decomposition is
\begin{equation}\label{gu:jac:semidirect-decomposition}
 (g,h)=(gr_s^{-1},hd_s^{-1})e_s,
 \qquad (gr_s^{-1},hd_s^{-1})\in U(V)\times U(W).
\end{equation}
Define
\begin{equation}\label{gu:jac:dilation}
 [D_sf](T)=|s|_F^{-4}f(r_s^{-1}T).
\end{equation}
The map $T\mapsto r_sT$ has absolute $F$-determinant $|s|_F^8$, since
$|\det_Er_s|_E=|s|_F^4$ and $T$ has two columns.  Thus the coefficient
in \eqref{gu:jac:dilation} is its inverse positive square root.  In
particular, $D_sD_t=D_{st}$, and $D_s$ is unitary for the Schr\"odinger
$L^2$-measure.

\begin{proposition}\label{gu:jac:extension}
The formula
\begin{equation}\label{gu:jac:extension-definition}
 \omega_{R,0}(g,h)
 =\omega(gr_s^{-1},hd_s^{-1})D_s,
 \qquad s=\lambda(g)=\lambda(h),
\end{equation}
defines a smooth oscillator representation of $R$ implementing
\eqref{gu:jac:tensor-action} and extending $\omega$.  Its action on the
matching Siegel Levi is
\begin{equation}\label{gu:jac:normalized-levi}
 [\omega_{R,0}(g,m(a,s))f](T)
 =|s|_F^{-4}|\det a|_E^2f(g^{-1}T\iota(a)),
 \qquad s=\lambda(g).
\end{equation}
The subgroup $K_E$ acts trivially.
\end{proposition}

\begin{proof}
The semidirect decomposition reduces the representation law to
\begin{equation}\label{gu:jac:covariance}
 D_s\omega(g_0,h_0)D_s^{-1}
 =\omega(r_sg_0r_s^{-1},d_sh_0d_s^{-1}).
\end{equation}
The generator calculation is given in Appendix~\ref{gu:app:normalization}.
Substituting $hd_s^{-1}=m(a,1)$ in \eqref{gu:jac:extension-definition}
gives \eqref{gu:jac:normalized-levi}. Its coefficient is the positive square
root of the absolute determinant $|s|_F^{-8}|\det a|_E^4$ of
$T\mapsto g^{-1}T\bar a$. For $(zI_V,\bar zI_W)\in K_E$, the argument
is $T$ and the scalar is
$|N_{E/F}z|_F^{-4}|\det(\bar zI_2)|_E^2=1$.
\end{proof}

\subsection{The product group and its multiplier coordinate}
The quotient map and a multiplicative section are
\[
 q(g,h)=\lambda(h)/\lambda(g),\qquad q_c=(1,d_c).
\]
The decomposition $G_\theta\times G'_\theta=R\{q_c:c\in F^\times\}$ identifies
$R\backslash(G_\theta\times G'_\theta)$ homeomorphically with $F^\times$.
Both groups are unimodular: $R$ is reductive, with kernel
$G_\theta^1\times (G'_\theta)^1$ under its common multiplier. Set
\begin{equation}\label{gu:jac:full-product}
 \Omega_0^{\GU}=\cInd_R^{G_\theta\times G'_\theta}\omega_{R,0}.
\end{equation}
\begin{proposition}\label{gu:jac:full-model}
Evaluation along $q_c$ gives a realization
\[
 \Omega^{\GU}_0=\mathcal S(F^\times\times Z),
 \qquad \phi(c,T)=F(q_c)(T).
\]
For $(g,h)\in G_\theta\times G'_\theta$, with $s=\lambda(g)$ and $t=\lambda(h)$, its action is
\begin{equation}\label{gu:jac:induced-action}
 [\Omega^{\GU}_0(g,h)\phi](c,\cdot)
 =\omega_{R,0}\bigl(g,d_chd_{ct/s}^{-1}\bigr)
                   \phi(ct/s,\cdot).
\end{equation}
In particular,
\begin{align}
 [\Omega^{\GU}_0(1,n(u))\phi](c,T)
 &=\psi\left(-\frac1{2c}\operatorname{tr}_E
                            (u\overline{G_T})\right)\phi(c,T),
 \label{gu:jac:induced-n}\\
 [\Omega^{\GU}_0(g,m(a,t))\phi](c,T)
 &=|s|_F^{-4}|\det a|_E^2
                         \phi(ct/s,g^{-1}T\bar a).
 \label{gu:jac:induced-levi}
\end{align}
\end{proposition}

\begin{proof}
Evaluation and the inverse $F(rq_c)=\omega_{R,0}(r)\phi(c,\cdot)$
identify the spaces: compactly supported section values have a finite
compact open constancy partition, and smoothness supplies a common
small right stabilizer. These maps also preserve the Schwartz
inductive-limit topologies. The identity
\[
 q_c(g,h)=(g,d_chd_{ct/s}^{-1})q_{ct/s}
\]
gives \eqref{gu:jac:induced-action}. Now
$d_cn(u)d_c^{-1}=n(c^{-1}u)$ and
$d_cm(a,t)d_{ct/s}^{-1}=m(a,s)$ give the two asserted formulas.
\end{proof}
The central subgroup $K_E$ acts trivially on $\Omega_0^{\GU}$ as well.

\subsection{The scaled Gram orbit}\label{gu:jac:jacquet-section}
Fix an isometric embedding $j:(X^c,B^c)\hookrightarrow V$, and put
$L=G_\theta\times K$. Define
\begin{equation}\label{gu:jac:scaled-locus}
 \widetilde{\mathcal O}_B
 =\{(c,T)\in F^\times\times Z:Q(T)=cB^c\}.
\end{equation}
The following is the unitary similitude analogue of
\cite[Proposition~3.4]{ChenGan}.

\begin{proposition}[Similitude Jacquet module]\label{gu:jac:jacquet}
Restriction, followed by multiplication by $|c|_F^4$, gives isomorphisms
\begin{equation}\label{gu:jac:J}
 (\Omega_0^{\GU})_{N(X),\psi_B}
 \simeq\mathcal S(\widetilde{\mathcal O}_B)
 \simeq\cInd_{H_j}^{G_\theta\times K}\mathbf1,
\end{equation}
where the middle term has geometric pullback action for
\begin{equation}\label{gu:jac:orbit-action}
 (g,a)\cdot(c,T)
 =\left(c\frac{\lambda(g)}{\lambda_B(a)},\ gT\iota(a)^{-1}\right),
\end{equation}
and
\begin{equation}\label{gu:jac:stabilizer}
 H_j=\{(g,a):\lambda(g)=\lambda_B(a),\ gj=j\iota(a)\}.
\end{equation}
These are homeomorphisms for the Schwartz inductive-limit topologies.
In these coordinates, the composite is
\begin{equation}\label{gu:jac:intertwiner}
 \mathcal I_j[\phi](g,a)
 =\left|\frac{\lambda_B(a)}{\lambda(g)}\right|_F^4
 \phi\left(\frac{\lambda_B(a)}{\lambda(g)},g^{-1}j\iota(a)\right).
\end{equation}
\end{proposition}

\begin{proof}
\emph{Twisted support.}
By \eqref{gu:jac:induced-n} and the nondegenerate trace pairing from
the conventions above, the unipotent character equals $\psi_B$ exactly
when $c^{-1}\overline{G_T}=\mathsf B$, or $Q(T)=cB^c$.
The locus is closed, and every $T$ on it is injective with nondegenerate
image. Lemma~\ref{gu:jac:restriction} gives the surjectivity of
restriction, identifies its kernel with the algebraic Jacquet relation
space, and establishes the quotient topology.

For $s=\lambda(g)$ and $t=\lambda_B(a)$, equations
\eqref{gu:jac:induced-levi} and \eqref{gu:jac:plane-determinant} give the
action before changing realization:
\begin{equation}\label{gu:jac:quotient-action}
 [\Omega_0^{\GU}(g,a)\phi](c,T)
 =|t/s|_F^4\phi(ct/s,g^{-1}T\iota(a)).
\end{equation}
The Gram matrix at the new argument is
$s^{-1}a^tG_T\bar a=(ct/s)\overline{\mathsf B}$, by
\eqref{gu:jac:character-stabilizer}. Set
\begin{equation}\label{gu:jac:gauge}
 \Phi(c,T)=|c|_F^4\phi(c,T).
\end{equation}
Since $|c|_F^4|t/s|_F^4=|ct/s|_F^4$, this gives the pullback action $\pi_0$
for \eqref{gu:jac:orbit-action}; multiplication preserves support and the
Schwartz topology.

\emph{Rational orbit.}
For $(c,T)$ in the locus, $r_c^{-1}T$ is an isometric embedding.
Witt extension~\cite[Chapter~7, \S9]{Scharlau85} gives $u\in U(V)$ with $uj=r_c^{-1}T$; hence
$(r_cu,1)\cdot(1,j)=(c,T)$.
The orbit map is a submersion and induces a homeomorphism
$L/H_j\simeq\widetilde{\mathcal O}_B$, by
Lemma~\ref{gu:jac:orbit-topology} in the appendix.

\emph{Stabilizer and block transport.}
At $o=(1,j)$, equation \eqref{gu:jac:orbit-action} gives
\eqref{gu:jac:stabilizer}. Put
\begin{equation}\label{gu:jac:planes}
 V_2=j(X^c),\qquad V_1=V_2^\perp,\qquad
 t_j(a)=j\iota(a)j^{-1}.
\end{equation}
Both planes are hyperbolic, by Witt extension and orthogonal
complementation in the split four-space. The map
$t_j:K\xrightarrow{\sim}\GU(V_2)$ preserves multipliers.
A similitude preserving $V_2$ also preserves $V_1$, and therefore
\begin{equation}\label{gu:jac:block-stabilizer}
 H_j=\{(u\oplus h,t_j^{-1}(h)):
       \lambda(u)=\lambda(h),\ u\in\GU(V_1),\ h\in\GU(V_2)\}.
\end{equation}
Transport by $1\times t_j$ identifies it with the subgroup
\eqref{eq:main-unitary-subgroup}.

\emph{Fiber character and modulus.}
For a stabilizer element of common multiplier $s$, the matching Levi
operator at $c=1$ has fixed argument and scalar
\begin{equation}\label{gu:jac:trivial-fiber}
 |s|_F^{-4}|\det a|_E^2=1.
\end{equation}
The splitting characters and the multiplier normalization are those of
Proposition~\ref{gu:jac:extension}, so the fiber character is trivial.
The equal-multiplier product $H_j$ is reductive, with kernel
$U(V_1)\times U(V_2)$ under its multiplier. Thus $H_j$ and $L$ are
unimodular and
\begin{equation}\label{gu:jac:relative-modulus}
 \Delta_L|_{H_j}\Delta_{H_j}^{-1}=1.
\end{equation}
The elements
$(zI_V,\bar zI_X)$ of $H_j$ transport to the common diagonal
$E^\times$-center of the block subgroup.

\emph{Compact induction.}
Compose restriction and \eqref{gu:jac:gauge} with
$[\mathcal A\Phi](k)=\Phi(k^{-1}o)$. The quotient homeomorphism
identifies compact support with compact support modulo $H_j$, and gives
the inverse. For $h\in H_j$ and $k_0\in L$,
\[
 \mathcal A\Phi(hk)=\mathcal A\Phi(k),\qquad
 \mathcal A(\pi_0(k_0)\Phi)(k)=\mathcal A\Phi(kk_0).
\]
The composite is \eqref{gu:jac:intertwiner}.
\end{proof}

Changing $j$ to $vj$, with $v\in U(V)$, conjugates $H_j$ by $(v,1)$
and $t_j$ by $v$, giving the natural isomorphism of the induced models.
On restriction to $U(V)\times U(W)$, the coordinate $c$ remains;
its $c=1$ fiber recovers the ordinary Weil-model calculation of
\cite[Proposition~3.4]{ChenGan}.

\subsection{Splitting characters and the Hom identity}\label{gu:sec:hom}
For $\tau\in\Irr(G_\theta)$ define the big theta lift by
\begin{equation}\label{gu:eq:big-theta-zero}
 \Theta_0^{\GU}(\tau)=(\Omega_0^{\GU}\otimes\tau^\vee)_{G_\theta}.
\end{equation}
Lemma~\ref{gu:app:comparison} identifies this coinvariant construction with the matching-subgroup construction, with the forward
integral \eqref{gu:eq:forward-comparison}.

We use the family of splitting data $D=(\chi,1,1)$, where $\chi$ is a
unitary character of $E^\times/F^\times$. Let
\[
 \mathfrak i:E^\times/F^\times\xrightarrow{\sim}E^1,
 \quad z\mapsto z/\bar z,\qquad \mu_\eta=\eta\circ\mathfrak i^{-1}.
\]
For either four-dimensional similitude group set
\begin{equation}\label{gu:eq:ambient-character}
 d_4(g)=\frac{\det_Eg}{\lambda(g)^2}\in E^1,\qquad
 \mathcal A_\eta(g)=\mu_\eta(d_4(g)).
\end{equation}
The determinant relation $N_{E/F}(\det_Eg)=\lambda(g)^4$ makes this a
character on either similitude group. Its isometry restriction gives the
splitting change in~\cite[equations~(3.1)--(3.3)]{ChenGan}: on $U(V)$ it is
$\mu_\eta\circ\det_E$, and on $m(a,1)$ it is $\eta(\det_Ea)$.
It is trivial on the upper unipotent radical and on the skew-Hermitian
partial Weyl generators, whose determinants are one. Put
\begin{equation}\label{gu:eq:chi-oscillator}
 \omega_{R,\chi}=\omega_{R,0}(1\boxtimes\mathcal A_\chi)|_R,
 \qquad \Omega_\chi^{\GU}=\cInd_R^{G_\theta\times G'_\theta}\omega_{R,\chi},
 \qquad \Theta_{\chi,1,1}^{\GU}(\tau)
       =(\Omega_\chi^{\GU}\otimes\tau^\vee)_{G_\theta}.
\end{equation}
Here $\chi=1$ recovers the subscript $0$ notation. Multiplying inducing
functions by the ambient character gives
\begin{equation}\label{gu:eq:chi-twist}
 \Omega_\chi^{\GU}\simeq
 \Omega_0^{\GU}\otimes(1\boxtimes\mathcal A_\chi),\qquad
 \Theta_{\chi,1,1}^{\GU}(\tau)
 \simeq\mathcal A_\chi\otimes\Theta_0^{\GU}(\tau).
\end{equation}
The general data $(\chi_V,\chi_W,\kappa)$ are treated in
Appendix~\ref{gu:app:splitting}.

\begin{theorem}[Shalika period identity for similitude theta lifts]\label{gu:thm:hom}
For irreducible admissible $\tau$ of $G_\theta$ and $\sigma$ of $K$, and any
unitary character $\chi$ of $E^\times/F^\times$, there is an isomorphism
\begin{equation}\label{gu:eq:hom-coverage}
 \begin{aligned}
 &\Hom_{S_B^{\GU}}(\Theta_{\chi,1,1}^{\GU}(\tau),
                         \sigma\boxtimes\psi_B)\\
 &\qquad\simeq
 \Hom_{H_{\mathrm U}}\left(\tau^\vee\boxtimes
  \bigl((\sigma^j)^\vee\otimes(\chi^{-1}\circ\det_E)\bigr),
  \mathbf1\right),\qquad \sigma^j=\sigma\circ t_j^{-1}.
 \end{aligned}
\end{equation}
Its scalar normalization is fixed by a nonzero invariant measure on
$H_j\backslash(G_\theta\times K)$. In particular, the Hom space has dimension
at most one.
\end{theorem}
\begin{proof}
First take $\chi=1$. Tensor evaluation (Lemma~\ref{gu:app:tensor}),
Proposition~\ref{gu:jac:jacquet}, and closed-subgroup duality
(Lemma~\ref{gu:app:duality}) give the canonical chain
\begin{equation}\label{gu:eq:hom-zero}
 \begin{aligned}
 \Hom_{S_B^{\GU}}(\Theta_0^{\GU}(\tau),\sigma\boxtimes\psi_B)
 &\simeq\Hom_{G_\theta\times K}((\Omega_0^{\GU})_{N,\psi_B},
                                      \tau\boxtimes\sigma)\\
 &\simeq\Hom_{G_\theta\times K}(\cInd_{H_j}^{G_\theta\times K}\mathbf1,
                                      \tau\boxtimes\sigma)\\
 &\simeq\Hom_{H_j}(\tau^\vee\boxtimes\sigma^\vee,\mathbf1).
 \end{aligned}
\end{equation}
In the first identification, tensor evaluation is followed by factoring
through the $(N,\psi_B)$-quotient.
The last step applies because $H_j$ and $G_\theta\times K$ are unimodular and
$\tau\boxtimes\sigma$ is admissible by~\cite[Proposition~2.16]{BZ76}.
Transporting the coefficient by $t_j$ gives the corresponding
spherical-pair Hom space.

For general $\chi$, on the Shalika subgroup
\[
 d_4(m_B(a))=\mathfrak i(\det_Ea),\qquad
 \mathcal A_\chi(m_B(a))=\chi(\det_Ea),\qquad
 \mathcal A_\chi|_{N(X)}=1.
\]
After absorbing this target character through \eqref{gu:eq:chi-twist},
the normalized identity gives the following three changes of coefficient:
\begin{equation}\label{gu:eq:coefficient-three-steps}
 \begin{aligned}
 \sigma_0&=\sigma\otimes(\chi^{-1}\circ\det_X),\\
 \sigma_0^j&=\sigma^j\otimes(\chi\circ\det_{V_2}),\\
 (\sigma_0^j)^\vee&=(\sigma^j)^\vee
                        \otimes(\chi^{-1}\circ\det_{V_2}).
 \end{aligned}
\end{equation}
Indeed, $t_j$ conjugates determinants by \eqref{gu:eq:iota-det}, and
$\chi(\bar z)=\chi(z)^{-1}$. This proves \eqref{gu:eq:hom-coverage}.
The representation on its right is irreducible admissible on the
specified split product, so Theorem~\ref{thm:main-unitary} gives the bound.
\end{proof}

\section{Variable splitting and nonvanishing of theta lifts}\label{sec:gu-coverage}
\label{gu:cov:section}
To pass from the Shalika period identity for big theta lifts to
arbitrary irreducible admissible representations of $\GU(W)$,
we prove that each such representation is a quotient of a suitable
big theta lift after varying the splitting character. The proof
combines analytic equal-rank theta dichotomy, highly ramified
stability, and a comparison of ordinary and similitude theta
correspondence.

\subsection{Analytic equal-rank theta and stability}\label{ana:section}

\subsubsection{Forms, characters, and normalization}

We use the extension $E/F$ and its conjugation $x\mapsto\bar x$
fixed in the introduction. The associated quadratic character of
$F^\times$ is denoted by $\omega_{E/F}$. We use the fixed nontrivial
additive character $\psi$ of $F$ and choose $\delta\in E^\times$ with
$\bar\delta=-\delta$.

For $n=2,4$, let $W_n$ be a hyperbolic skew-Hermitian space of dimension
$n$ over $E$. Let $V_n^+$ be a hyperbolic Hermitian space of dimension
$n$, and let $V_n^-$ be the other Hermitian space of that dimension.
For either type of form set
\[
 \operatorname{disc}(A)=(-1)^{n(n-1)/2}\det(A),
\]
where the determinant is taken modulo $N_{E/F}(E^\times)$. Thus
$\operatorname{disc}(V)\in F^\times/N_{E/F}(E^\times)$ in the
Hermitian case and
$\operatorname{disc}(W)\in\delta^n F^\times/N_{E/F}(E^\times)$
in the skew-Hermitian case; see \cite[\S2.1]{GanIchino}. Define
\[
 \epsilon(V)=\omega_{E/F}(\operatorname{disc}(V)),\qquad
 \epsilon(W)=\omega_{E/F}(\delta^{-n}\operatorname{disc}(W)).
\]
Thus $\epsilon(W_n)=\epsilon(V_n^+)=1$ and
$\epsilon(V_n^-)=-1$. We also put $W_0=V_0^+=0$ and interpret their
unitary groups as the trivial group. There is no space $V_0^-$.

The tensor form is fixed in \eqref{gu:eq:tensor-symplectic}.
Gan--Ichino use the full-trace tensor form
\cite[\S4.1]{GanIchino}. The change of central
coordinate $(z,t)\mapsto(z,-t/2)$ identifies their convention with ours.
Accordingly, all their analytic factors in this section are evaluated
using the fixed character
\begin{equation}\label{ana:psi}
 \psi_0(t)=\psi(-t/2).
\end{equation}

Let $\chi:E^\times\to\C^\times$ be a smooth unitary character satisfying
$\chi|_{F^\times}=1$. The splitting pair $(\chi,1)$ is permissible
for all the even-dimensional spaces just specified. For the full-trace
form with $(\chi,1,\psi_0)$, equivalently our form with
$(\chi,1,\psi)$, write $\omega_{V,W}$ for the ordinary Weil
representation of $\U(V)\times\U(W)$, and put
\begin{equation}\label{ana:big-lift}
 T_{V,W}(\pi)
    =(\omega_{V,W}\otimes\pi^\vee)_{\U(W)}.
\end{equation}
We write $T^\chi_{V,W}$ when the dependence on $\chi$ is needed.
The subscript in \eqref{ana:big-lift} denotes smooth coinvariants.

\subsubsection{Conservation and the analytic sign}

\begin{theorem}[Equal-rank dichotomy and finite length]
\label{ana:dichotomy}
Let $n=2$ or $4$, and let $\pi\in\Irr(\U(W_n))$.
Exactly one of
\[
 T_{V_n^+,W_n}(\pi),\qquad T_{V_n^-,W_n}(\pi)
\]
is nonzero. Both big theta lifts have finite length as smooth
representations of their respective target groups.
\end{theorem}

\begin{proof}
Fix the splitting character $\chi$ in the two compatible even
Hermitian Witt towers. Let $m^+(\pi)$ and $m^-(\pi)$ denote the
dimensions of first occurrence. Sun--Zhu's conservation relation gives
\[
 m^+(\pi)+m^-(\pi)=2n+2;
\]
see \cite[Theorem~1.10 and Proposition~1.6]{SunZhu}.
Both indices are even. If both were at most $n$, their sum would be
at most $2n$. If both exceeded $n$, their sum would be at least
$2n+4$. Hence precisely one is at most $n$. Kudla's persistence
principle, recalled in \cite[\S1.5]{SunZhu}, now proves the assertion
about occurrence. 
Finite length here is the finiteness theorem for the ordinary
\emph{big} theta lift due to Kudla, as recalled in
\cite[Introduction]{GanTakeda} and \cite[\S4.2]{GanIchino}.
The finiteness theorem holds in all residue characteristics.
Consequently a nonzero big theta lift has an irreducible quotient. \end{proof}

Write
$\gamma_{\mathrm{LR}}(s,\pi,\eta,\psi_0)$ for the analytic doubling
factor, for a smooth character $\eta$ of $E^\times$, in the normalization
of \cite[\S10.1 and Appendix~A.5]{GanIchinoFormal}. It is defined by the
normalized local functional equation. 
\begin{proposition}[Tempered analytic criterion]\label{ana:criterion}
Let $n=2$ or $4$ and let $\pi$ be an irreducible tempered
representation of $\U(W_n)$. The analytic factor
$\gamma_{\mathrm{LR}}(s,\pi,\chi^{-1},\psi_0)$ is finite and nonzero
at $s=\tfrac12$. For $V=V_n^+$ or $V_n^-$,
\begin{equation}\label{ana:sign}
 T_{V,W_n}(\pi)\ne0
 \quad\Longleftrightarrow\quad
 \gamma_{\mathrm{LR}}\left(\frac12,\pi,\chi^{-1},\psi_0\right)
       =\omega_\pi(-1)\epsilon(V).
\end{equation}
\end{proposition}

\begin{proof}
The proof uses the analytic argument of \cite[Appendix~A.5]{GanIchinoFormal}.
Let $\mathbf W=W_n\oplus W_n^{-}$, where $W_n^{-}$ has the negative
of the form on $W_n$, and let $I(0,\chi)$ be the normalized Siegel
principal series of $\U(\mathbf W)$ at the center.
The two Siegel--Weil submodules satisfy
\begin{equation}\label{ana:siegel-weil-sum}
 I(0,\chi)=R(V_n^+)\oplus R(V_n^-).
\end{equation}
This is \cite[Proposition~4.1(ii), p.~961]{HKS96}, transported to
the conventions of \cite{GanIchinoFormal}.

Let $\beta_W$ be the nonsingular auxiliary map used to normalize the
doubling intertwiner. By \cite[Appendix~A.4, Lemma~A.3]{GanIchinoFormal},
$M^{\mathrm{LR}}_{\psi_0}(0,\chi)$ acts on $R(V)$ by
\begin{equation}\label{ana:intertwiner-eigenvalue}
 \lambda_V
  =\chi(\det\beta_W)\chi(\delta)^n
        \epsilon(V)\epsilon(W_n)^n.
\end{equation}
The change from the Harris--Kudla--Sweet normalization is computed
in that appendix; the original scalar calculation is
\cite[Proposition~6.8, p.~974]{HKS96}.

For a matrix coefficient $f$ of $\pi^\vee$, the first functional
equation in \cite[Appendix~A.5]{GanIchinoFormal} is
\begin{align}
 &Z\bigl(0,M^{\mathrm{LR}}_{\psi_0}(0,\chi)\mathcal F,f\bigr)
 \notag\\
 &\quad=
   \epsilon(W_n)^{n+1}\omega_\pi(-1)\chi(\det\beta_W)
   \gamma_{\mathrm{LR}}\left(\frac12,\pi,\chi^c,\psi_0\right)
   Z(0,\mathcal F,f),
 \label{ana:functional-equation}
\end{align}
where $\chi^c(x)=\chi(\bar x)$.
For tempered $\pi$, the integrals at the center converge absolutely;
the relevant convergence statement is
\cite[Lemma~9.5(ii)]{GanIchinoFormal}, used expressly in Appendix~A.5.
The resulting family of central zeta functionals is nonzero. Indeed,
one can take a section supported on a sufficiently small compact
neighborhood in the open doubling orbit, and a matrix coefficient
nonzero at the identity.

Finiteness of the central gamma value follows
directly from the meromorphic functional equation. Extend such a
section to a flat holomorphic section near the center. Its zeta
integral is nonzero there, and the normalized intertwiner is regular
at the center by \cite[Proposition~8.1 and Appendix~A]{GanIchinoFormal}.
The functional equation therefore makes the gamma factor holomorphic
at $\tfrac12$. Its value cannot be zero: the operator
$M^{\mathrm{LR}}_{\psi_0}(0,\chi)$ is invertible, as is also evident
from \eqref{ana:siegel-weil-sum} and its two nonzero eigenvalues
\eqref{ana:intertwiner-eigenvalue}. A zero gamma value in
\eqref{ana:functional-equation} would consequently make every
central zeta integral vanish.

By \eqref{ana:siegel-weil-sum}, the zeta functional is nonzero on
at least one $R(V)$. Its restriction to that summand yields theta
occurrence by the doubling see-saw
\cite[Proposition~3.1, p.~960]{HKS96}. The existence hypothesis
on an irreducible theta quotient in that proposition is supplied by
the big-lift finite length in Theorem~\ref{ana:dichotomy}.
Thus this $V$ is the unique nonzero tower in that theorem.
On this summand, \eqref{ana:functional-equation} forces equality
between \eqref{ana:intertwiner-eigenvalue} and the scalar on the
right side of the functional equation. The two eigenvalues have
opposite signs, so this equality also singles out that summand.

Finally, $\chi^c=\chi^{-1}$ because $x\bar x\in F^\times$,
and $\chi(\delta)^n=1$ because $\delta^2\in F^\times$ and $n$ is
even. We have $\epsilon(W_n)=1$ and $\omega_\pi(-1)^2=1$.
Canceling the common factor $\chi(\det\beta_W)$ gives
\eqref{ana:sign}.
\end{proof}

\subsubsection{Highly ramified comparison}

Let $\mathcal O_E$ be the ring of integers of $E$ and $\mathfrak p_E$
its maximal ideal. Put $U_E^0=\mathcal O_E^\times$ and
$U_E^A=1+\mathfrak p_E^A$ for $A\ge1$. For a smooth character
$\eta$ of $E^\times$, write
\[
 a(\eta)=\min\{A\ge0:\eta|_{U_E^A}=1\}.
\]
In particular, $a(\eta^{-1})=a(\eta)$.

\begin{proposition}[Analytic stability with a fixed central character]
\label{ana:stability}
Let $n=2$ or $4$, and let $\pi_1,\pi_2\in\Irr(\U(W_n))$
have the same central character on $E^1$.
For the fixed additive character $\psi_0$, there is an integer
$N=N(\pi_1,\pi_2,\psi_0)$ such that
\begin{equation}\label{ana:stable-gamma}
 \gamma_{\mathrm{LR}}(s,\pi_1,\eta,\psi_0)
   =\gamma_{\mathrm{LR}}(s,\pi_2,\eta,\psi_0)
\end{equation}
as meromorphic functions of $s$, for every smooth character
$\eta$ of $E^\times$ satisfying $a(\eta)>N$.
\end{proposition}

\begin{proof}
This is the unitary case of \cite[Proposition~5.3]{FinisLapid}.
That proposition concerns the analytic Piatetski-Shapiro--Rallis
factors explicated by Lapid--Rallis, and requires equality of central
characters in the unitary case. The group, additive character, and
normalization are fixed in the comparison. The unitary normalization
correction $\epsilon(W_n)^{n+1}$ used in \cite[\S10.1]{GanIchinoFormal}
equals one for the present source; the other common normalizing
factors cancel~\cite[Remarks~5.3--5.5]{Kakuhama19}. Thus the comparison is precisely
\eqref{ana:stable-gamma} in the notation above.
The proposition applies to all sufficiently ramified characters of
$E^\times$, so in particular it applies after imposing
$\eta|_{F^\times}=1$.
\end{proof}

\begin{lemma}[Unbounded conductors]\label{ana:characters}
The smooth unitary characters of $E^\times/F^\times$ have
unbounded conductors.
\end{lemma}

\begin{proof}
The group $Q=E^\times/F^\times$ is compact and totally
disconnected. Indeed, valuations are taken modulo
$e(E/F)\mathbf Z$, and the remaining quotient of unit groups is
compact. Hilbert's Theorem~90 identifies $Q$ with the infinite
group $E^1$ through $x\mapsto x/\bar x$.
Its finite-quotient characters therefore separate points and
form an infinite set. They are unitary, since they have finite
image. For each fixed $A$, a character of conductor at most $A$
factors through the finite group
\[
 E^\times/(F^\times U_E^A).
\]
There are only finitely many such characters. Their conductors
cannot be bounded.
\end{proof}

\subsection{Inducing blocks and the principal-series base case}
\label{blk:section}

\subsubsection{Normalized Jacquet modules}

We retain the splitting characters $(\chi,1)$ and the additive character
fixed above.  Write $G_n=\U(W_n)$ and $H_m=\U(V_m)$, where $W_n$ is
hyperbolic skew-Hermitian and $V_m$ is Hermitian.  In this section all
dimensions are even.  We abbreviate the Weil representation to
$\omega_{V_m,W_n}$ and use $T_{V_m,W_n}(\pi)$ for its ordinary big theta
lift.  All groups denote their groups of $F$-points.  We use normalized
parabolic induction and normalized Jacquet modules:
\[
 i_P^G(\rho)=\operatorname{Ind}_P^G
       (\delta_P^{1/2}\rho),
 \qquad
 r_P(A)=\delta_P^{-1/2}\otimes A_{N_P}.
\]
Thus normalized Frobenius reciprocity uses $r_P$, whereas an induced
representation in the first argument of $\Hom$ uses the opposite
Jacquet module:
\begin{align}
 \Hom_G(A,i_P^G B)&\simeq\Hom_M(r_P A,B),
       \label{blk:first-adjunction}\\
 \Hom_G(i_P^G A,B)&\simeq\Hom_M(A,r_{\overline P}B).
       \label{blk:second-adjunction}
\end{align}
Here $A$ may be smooth and $B$ admissible in
\eqref{blk:second-adjunction}.  For these hypotheses the latter identity
follows from smooth duality, ordinary Frobenius reciprocity, and
Casselman's Jacquet pairing:
\[
\begin{aligned}
 \Hom_G(i_P A,B)
 &\simeq\Hom_G(B^\vee,i_P(A^\vee))\\
 &\simeq\Hom_M(r_P(B^\vee),A^\vee)\\
 &\simeq\Hom_M\bigl(A,(r_P(B^\vee))^\vee\bigr)
  \simeq\Hom_M(A,r_{\overline P}B).
\end{aligned}
\]
The duality for induction applies to smooth $A$ because $P\backslash G$
is compact, and the admissible modules in this calculation are reflexive
\cite[Proposition~3.1.2, Theorems~3.2.4 and~3.3.1,
Corollary~4.2.5]{Casselman}.
For representations of a product group, these identities apply with
admissibility for the product representation.

Let $P_k$ be the parabolic of $G_n$ stabilizing a $k$-dimensional
isotropic subspace $X_k$, with Levi
\[
 L_k=\GL(X_k)\times G_{n-2k}.
\]
If $Y_a\subset V_m$ is isotropic of dimension $a$, let $Q'_a$ be its
stabilizer in $H_m$.  Write
$X_k=X_{k-a}\oplus X'_a$, and denote by $Q_{k-a,a}$ the parabolic of
$\GL(X_k)$ stabilizing $X_{k-a}$.  A subscript of dimension zero denotes
the trivial group or zero space, as appropriate.

\begin{lemma}[The equal-rank boundary exponents]
\label{blk:kudla}
For $g\in\GL_k(E)$, $c\in\GL_a(E)$, and
$(d,b)\in\GL_{k-a}(E)\times\GL_a(E)$, the relevant modulus characters
are
\begin{align}
 \delta_{P_k}(g)&=\abs{\det g}_E^{n-k},
       \label{blk:source-modulus}\\
 \delta_{Q'_a}(c)&=\abs{\det c}_E^{m-a},
       \label{blk:target-modulus}\\
 \delta_{Q_{k-a,a}}(d,b)
   &=\abs{\det d}_E^a\abs{\det b}_E^{-(k-a)}.
       \label{blk:gl-modulus}
\end{align}
The normalized Jacquet module $r_{P_k}(\omega_{V_m,W_n})$ has a
decreasing filtration
\[
 R^0\supset R^1\supset\cdots\supset R^k\supset R^{k+1}=0,
 \qquad J^a=R^a/R^{a+1}.
\]
If $V_m$ has no isotropic $a$-space, then $J^a=0$.  Otherwise $J^a$ is
normalized induction from
\[
 (Q_{k-a,a}\times G_{n-2k})\times Q'_a
   \ \subset\ L_k\times H_m
\]
of a representation with underlying space
\[
 \cS\bigl(\operatorname{Isom}_{E,c}(X'_a,Y_a)\bigr)
       \otimes\omega_{V_{m-2a},W_{n-2k}}.
\]
Here $\operatorname{Isom}_{E,c}$ denotes conjugate-linear $E$-isomorphisms;
thus $xb$ means composition with $b$ on the source.
On this coefficient space the unmatched factor $\GL(X_{k-a})$ acts
by
\begin{equation}
 \chi(\det d)\abs{\det d}_E^{(m-n+k-a)/2}.
 \label{blk:boundary-character}
\end{equation}
The matched source and target factors act on the first tensor factor by
\[
 ((b,c)f)(x)=\chi(\det b)f(c^{-1}xb),
\]
and the two remaining unitary groups act on the smaller Weil
representation.  In particular, for $m=n$ every boundary term $a<k$
has strictly positive unmatched real exponent $(k-a)/2$.
\end{lemma}

\begin{proof}
The off-diagonal part of the unipotent radical of $P_k$ contributes
$\abs{\det g}_E^{n-2k}$ to its modulus.  The Hermitian or skew-Hermitian
$k$-by-$k$ part contributes $\abs{\det g}_E^k$, giving
\eqref{blk:source-modulus}.  The same calculation gives
\eqref{blk:target-modulus}; the action $x\mapsto dxb^{-1}$ on
$\operatorname{Hom}_E(X'_a,X_{k-a})$ gives
\eqref{blk:gl-modulus}.  These are determinant identities over $F$ and
apply in residue characteristic two as well.

At rank zero, the mixed
Schr\"odinger model has scalar factor
$\chi(\det g)\abs{\det g}_E^{m/2}$ before taking normalized Jacquet
modules.  Dividing by $\delta_{P_k}^{1/2}$ gives
\[
 \chi(\det g)\abs{\det g}_E^{m/2-(n-k)/2}
   =\chi(\det g)\abs{\det g}_E^{(m-n+k)/2}.
\]
The source-space filtration is \cite[Lemma~3.1]{GanTakeda}; compare
\cite[Appendix~C.1, Lemma~C.2]{GanIchinoFormal} with the two unitary spaces
interchanged.  The parameter $\epsilon_0$ there is zero for a unitary
pair.  Normalized induction from these parabolics then gives
\eqref{blk:boundary-character}.  Setting $m=n$ proves the last assertion.
\end{proof}

Only the following boundary terms occur in our argument:
\[
\begin{array}{c|c|c|c}
 n&k&a&\text{unmatched real exponent}\\ \hline
 2&1&0&\tfrac12\\
 4&1&0&\tfrac12\\
 4&2&0&1\\
 4&2&1&\tfrac12
\end{array}
\]
We will repeatedly use a simple consequence of the decreasing
filtration.  If $\Hom(J^a,B)=0$ for all $a<k$, then restriction is
injective:
\begin{equation}
 \Hom(R^0,B)\longrightarrow\Hom(J^k,B).
 \label{blk:restriction-injective}
\end{equation}
Indeed, a map vanishing on $J^k=R^k$ factors through $R^0/R^k$, and its
successive restrictions vanish along the finite filtration of that
quotient.  In particular, if $J^k=0$, the first Hom space is zero.

\subsubsection{The zero-exponent principal series}

\begin{proposition}
\label{blk:principal-split}
\label{blk:zero}
Let $n=2\ell\in\{2,4\}$, and let $\alpha_1,\ldots,\alpha_\ell$ be
smooth unitary characters of $E^\times$.  Every irreducible constituent
$\pi$ of
\[
 I=i_B^{G_n}(\alpha_1\boxtimes\cdots\boxtimes\alpha_\ell)
\]
satisfies
\[
 T_{V_n^-,W_n}(\pi)=0,
 \qquad T_{V_n^+,W_n}(\pi)\ne0
\]
for every permissible unitary $\chi$.
\end{proposition}

\begin{proof}
The representation $I$ is admissible, of finite length, and unitary
\cite[Proposition~2.3 and Theorem~2.8]{BZ77}
\cite[Proposition~3.1.4]{Casselman}.  It is therefore semisimple.
Here is the elementary argument relevant to this use.  If $A\subset I$
is invariant and $v\in I^K$ for a compact open subgroup $K$, decompose
$v=a+b$ orthogonally in the finite-dimensional space $I^K$, with
$a\in A^K$ and $b\perp A^K$.  Averaging over $K$ shows that
\[
 \langle b,u\rangle=\langle b,e_Ku\rangle=0
       \qquad(u\in A).
\]
Thus $I=A\oplus A^\perp$.  Finite length now gives a finite direct sum
of irreducibles; see also
\cite[Propositions~2.1.5 and~2.1.14]{Casselman}.  In particular,
$\pi\hookrightarrow I$, even when $I$ is reducible.

Group the characters in the source Siegel Levi.  Normalized
transitivity gives
\begin{equation}
 \rho=i_{B_\ell}^{\GL_\ell(E)}
       (\alpha_1\boxtimes\cdots\boxtimes\alpha_\ell),
 \qquad I=i_{P_\ell}^{G_n}(\rho).
 \label{blk:grouped-principal}
\end{equation}
The center of $\GL_\ell(E)$ acts unitarily on $\rho$.  For $\ell=2$
the geometric lemma gives
\begin{equation}
 r_{\overline B_2}(\rho)^{\mathrm{ss}}
  =(\alpha_1\boxtimes\alpha_2)
       \oplus(\alpha_2\boxtimes\alpha_1).
 \label{blk:opposite-principal}
\end{equation}
This statement about composition factors applies when the characters
coincide as well \cite[Lemma~2.12]{BZ77}.

Suppose $T_{V_n^-,W_n}(\pi)\ne0$.  By
Theorem~\ref{ana:dichotomy} it has finite length, so it has an
irreducible quotient $\sigma$. Composing the corresponding nonzero
map $\omega_{V_n^-,W_n}\to\pi\boxtimes\sigma$ with the embedding
of $\pi$ in \eqref{blk:grouped-principal} gives, by
\eqref{blk:first-adjunction},
\begin{equation}
 \Hom_{\GL_\ell(E)\times H_n}
   \bigl(r_{P_\ell}\omega_{V_n^-,W_n},
               \rho\boxtimes\sigma\bigr)\ne0.
 \label{blk:principal-hom}
\end{equation}

If $\ell=1$, the rank-zero source character is
$\chi\abs{\cdot}_E^{1/2}$, whereas $\rho=\alpha_1$ is unitary.
Consequently the Hom space from $J^0$ is zero.

If $\ell=2$, the source $\GL_2$ acts on $J^0$ by
$\chi(\det)\abs{\det}_E$.  On a scalar $zI_2$ its absolute value is
$\abs{z}_E^2$, whereas the center acts unitarily on $\rho$.
Thus the rank-zero Hom space again vanishes.  For rank one apply
\eqref{blk:second-adjunction} to the source $B_2$ and target $Q'_1$
occurring in Lemma~\ref{blk:kudla}.  The unmatched source $E^\times$
acts in the coefficient by $\chi\abs{\cdot}_E^{1/2}$, while the first
characters of \eqref{blk:opposite-principal} are unitary.  This excludes
the rank-one term.  Possible extensions between the composition factors
in \eqref{blk:opposite-principal} do not affect the argument: a character
with positive real exponent cannot occur in a finite-length module all
of whose characters have real exponent zero.

By \eqref{blk:restriction-injective}, the nonzero Hom space in
\eqref{blk:principal-hom} must be nonzero on $J^\ell$.  But $V_n^-$
has Witt index $\ell-1$, so $J^\ell=0$.  This contradiction proves
vanishing in the nonsplit tower.  Equal-rank dichotomy
(Theorem~\ref{ana:dichotomy}) then gives nonvanishing in the split tower.
\end{proof}

\subsubsection{Positive inducing exponents}

\begin{lemma}[Removal of inducing blocks in dimensions two and four]
\label{blk:remove}
Let $n\in\{2,4\}$.  Suppose that an irreducible representation $\pi$ of
$G_n$ embeds in
\begin{equation}
 i_P^{G_n}\bigl(
   \tau_1\abs{\det}_E^{-s_1}\boxtimes\cdots\boxtimes
   \tau_r\abs{\det}_E^{-s_r}\boxtimes\pi_0\bigr),
 \label{blk:anti-standard}
\end{equation}
where $s_1\ge\cdots\ge s_r>0$, each $\tau_i$ is a unitary discrete
series representation of $\GL_{k_i}(E)$, and $\pi_0$ is irreducible
tempered on $G_{n_0}$, with $n_0=n-2\sum_i k_i$.  If
$T_{V_n,W_n}(\pi)\ne0$, then
\[
 V_n\simeq\mathbb H^{\sum_i k_i}\perp V_{n_0},
 \qquad T_{V_{n_0},W_{n_0}}(\pi_0)\ne0.
\]
Here $\mathbb H$ is a Hermitian hyperbolic plane, and
$T_{0,0}(\mathbf1)=\C$.
\end{lemma}

\begin{proof}
If $r=0$, the assertion is immediate.  Put $k=\sum_i k_i>0$ and group
the general linear blocks into a representation $\rho$ of $\GL_k(E)$.
Normalized transitivity turns \eqref{blk:anti-standard} into an
embedding
\[
 \pi\hookrightarrow i_{P_k}^{G_n}(\rho\boxtimes\pi_0).
\]
Choose an irreducible quotient $\sigma$ of the nonzero big lift.
Composing as in the proof of Proposition~\ref{blk:principal-split}
gives a nonzero map
\begin{equation}
 r_{P_k}\omega_{V_n,W_n}\longrightarrow
         \rho\boxtimes\pi_0\boxtimes\sigma.
 \label{blk:removal-hom}
\end{equation}
We first show that its restriction to $J^k$ is nonzero.

If $k=1$, then $\rho=\alpha\abs{\cdot}_E^{-s}$ with $s>0$ and
$\alpha$ unitary.  The only boundary term has source character
$\chi\abs{\cdot}_E^{1/2}$.  Its Hom space to the target of
\eqref{blk:removal-hom} is zero because the real exponents differ.

If $k=2$, necessarily $n=4$ and $n_0=0$.  There are two cases.  First,
suppose
\[
 \rho=\tau\abs{\det}_E^{-s},\qquad s>0,
\]
with $\tau$ unitary discrete series on $\GL_2(E)$.  The center acts on
$\rho$ with absolute value $\abs{z}_E^{-2s}$ at $zI_2$, whereas it acts
on $J^0$ with absolute value $\abs{z}_E^2$.  This excludes rank zero.
For rank one, second adjointness uses $r_{\overline B_2}(\tau)$.
The real exponent on its first $E^\times$ factor is negative whenever
this Jacquet module is nonzero, by Casselman's square-integrability
criterion; this is precisely the opposite-parabolic use of the
criterion in \cite[Appendix~C.1, equations~(C.1)--(C.2)]{GanIchinoFormal}.
After twisting by $\abs{\det}_E^{-s}$ the exponent remains negative,
so it cannot equal the boundary exponent $1/2$.  For example, a
unitary Steinberg twist has
\[
 r_{\overline B_2}\bigl(\operatorname{St}_2\otimes(\alpha\circ\det)\bigr)
 =\alpha\abs{\cdot}_E^{-1/2}
       \boxtimes\alpha\abs{\cdot}_E^{1/2},
\]
which gives the unmatched exponent $-s-1/2$.

Second, suppose
\[
 \rho=i_{B_2}^{\GL_2(E)}
       \bigl(\alpha_1\abs{\cdot}_E^{-s_1}
             \boxtimes\alpha_2\abs{\cdot}_E^{-s_2}\bigr).
\]
The center of $\rho$ has absolute value
$\abs{z}_E^{-(s_1+s_2)}$, excluding $J^0$.  Its opposite Jacquet
composition factors are the two ordered inducing characters and their
permutation.  Thus the unmatched first exponent is $-s_1$ or $-s_2$,
both negative, excluding $J^1$.  This argument uses the entire induced
representation $\rho$ and remains valid when $\rho$ is reducible.

These cases exhaust $k\le2$.  Equation~\eqref{blk:restriction-injective}
now shows that \eqref{blk:removal-hom} is nonzero on $J^k$.  Consequently
$V_n$ contains an isotropic $k$-space and has a decomposition
\[
 V_n\simeq\mathbb H^k\perp V_{n_0}.
\]
If $n_0=0$, this proves the assertion, including the stated nonvanishing
for the zero-dimensional pair.

It remains to consider $n_0=2$.  Then $n=4$ and $k=1$.  The coefficient
of the bottom term is
\[
 \mathcal A^1=
 \cS\bigl(\operatorname{Isom}_{E,c}(X_1,Y_1)\bigr)
       \otimes\omega_{V_2,W_2},
\]
with the actions of Lemma~\ref{blk:kudla}, and $J^1$ is its normalized
induction in the target from $Q'_1$.  Applying
\eqref{blk:second-adjunction} in that target gives a nonzero map
\begin{equation}
 \mathcal A^1\longrightarrow
   (\alpha\abs{\cdot}_E^{-s})\boxtimes\pi_0
                   \boxtimes r_{\overline Q'_1}(\sigma).
 \label{blk:bottom-map}
\end{equation}
This is a map of representations of
\[
 E^\times\times G_2\times E^\times\times\U(V_2).
\]
The Jacquet module in \eqref{blk:bottom-map} is admissible of finite
length.  Choose its composition series, with factors
$\nu_j\boxtimes\sigma_{0,j}$.  Tensoring this series with the fixed
irreducible $(\alpha\abs{\cdot}_E^{-s})\boxtimes\pi_0$ gives a
composition series of the target of \eqref{blk:bottom-map}; the factors
are irreducible by the tensor-product theorem for representations of a
product group \cite[\S1.6]{BZ77}.

Let $C$ be the nonzero image of \eqref{blk:bottom-map}, and choose a
simple quotient $C\twoheadrightarrow S$.  By Jordan--H\"older,
\[
 S\simeq
 (\alpha\abs{\cdot}_E^{-s})\boxtimes\pi_0
                      \boxtimes\nu\boxtimes\sigma_0
\]
for some constituent $\nu\boxtimes\sigma_0$ of the chosen composition series.  Composing
$\mathcal A^1\twoheadrightarrow C\twoheadrightarrow S$ gives a
nonzero map with this irreducible target.  Choose
$f\in\cS(\operatorname{Isom}_{E,c}(X_1,Y_1))$ on which the resulting map is
nonzero.  Fixing $f$ and identifying the two one-dimensional character
factors with $\C$ yields a nonzero
$G_2\times\U(V_2)$-equivariant map
\[
 \omega_{V_2,W_2}\longrightarrow\pi_0\boxtimes\sigma_0.
\]
Therefore $T_{V_2,W_2}(\pi_0)\ne0$, as required.
\end{proof}

\subsubsection{An irreducible principal series with prescribed center}

\begin{lemma}
\label{blk:prescribed-center}
\label{blk:principal}
For every smooth unitary character $\zeta$ of $E^1$ and every
$n\in\{2,4\}$, there exists an irreducible unitary tempered principal
series $\pi^*$ of $G_n$ with central character $\zeta$.  For every
permissible unitary $\chi$,
\[
 T_{V_n^+,W_n}(\pi^*)\ne0.
\]
\end{lemma}

\begin{proof}
Since $E^1\subset\mathcal O_E^\times$, choose $A\ge1$ such that
$\zeta$ is trivial on $E^1\cap(1+\mathfrak p_E^A)$.  The induced
character of the image of $E^1$ in the finite abelian group
$\mathcal O_E^\times/(1+\mathfrak p_E^A)$ extends to that entire group.
Indeed, characters of subgroups of finite abelian groups extend by
successively choosing roots in $\C^\times$ on cyclic generators.
The extension is unitary.  Giving a uniformizer a unitary value then
defines a smooth unitary character $\alpha$ of $E^\times$ extending
$\zeta$.

Put $\ell=n/2$ and initially choose
$(\alpha_1,\ldots,\alpha_\ell)=(\alpha,1,\ldots,1)$. We may twist
the $\alpha_i$ independently by unramified unitary characters $\nu_i$.
These twists are trivial on $E^1$ and preserve
\[
 (\alpha_1\nu_1\cdots\alpha_\ell\nu_\ell)|_{E^1}=\zeta.
\]
They can be chosen so that the resulting character $\beta$ of the
Borel torus is Weyl-regular.  To see this explicitly, write
$\nu_i(\varpi_E)=z_i\in S^1$.  The relative Weyl group acts by signed
permutations, where a sign replaces a character by its inverse
conjugate.  For any nonidentity Weyl element, fixing $\beta$ either is
impossible on $\mathcal O_E^\times$, or imposes a nontrivial equation
\[
 z_i=c z_j,\qquad z_i=c z_j^{-1},\qquad\text{or}\qquad z_i^2=c
\]
with $c\in S^1$ fixed.  Each stabilizer locus is a proper closed subset
of $(S^1)^\ell$ with empty interior.  Their finite union does not cover
the torus.  Choose the $\nu_i$ in its complement.

Let $T$ be the Levi torus of $B$, and let $I=i_B^{G_n}(\beta)$.
The geometric lemma gives one composition
factor $w\beta$ in $r_BI$ for each Weyl element $w$.  Regularity implies
that $\beta$ occurs exactly once.  Hence
\[
 \dim\operatorname{End}_{G_n}(I)
   =\dim\Hom_T(r_BI,\beta)\le1.
\]
The identity endomorphism gives equality.  Since $I$ is unitary,
admissible, and of finite length, it is semisimple by the argument in
Proposition~\ref{blk:principal-split}; the equality forces $I$ to be
irreducible.  This is also
\cite[Theorem~6.6.1]{Casselman}.  Normalized induction preserves
temperedness \cite[Lemma~III.2.3]{Waldspurger03}, so $I$ is tempered.

A scalar $z\in E^1$ corresponds to $(z,\ldots,z)$ in the Borel torus,
and $\delta_B(z)=1$.  Its action on $I$ is consequently
\[
 \prod_{i=1}^{\ell}\alpha_i(z)\nu_i(z)=\zeta(z).
\]
Thus $\pi^*=I$ has the prescribed central character.
Proposition~\ref{blk:principal-split} gives the asserted split
nonvanishing for every permissible $\chi$.
\end{proof}

\subsection{Simultaneous ordinary nonvanishing}
\label{cov:ordinary-section}

Konno's Langlands classification
\cite[Theorem~3.5 and Corollary~3.2]{Konno03}, applied to $\pi^\vee$
and followed by smooth duality, embeds $\pi$ into a normalized
anti-standard induction from an irreducible tempered Levi representation
with strictly negative real exponents on its general linear factors.
For the hyperbolic $\U(W_n)$, $n=2,4$, such a Levi is
\[
 \GL_{k_1}(E)\times\cdots\times\GL_{k_r}(E)
                  \times\U(W_{n_0}),\qquad
 n_0=n-2\sum_i k_i\in\{0,2,4\}.
\]
If a general linear tempered factor is not discrete series, refine it
using the direct-summand statement
\cite[Proposition 2.2]{Konno03}, which is
\cite[Proposition III.4.1]{Waldspurger03}. Exactness and transitivity of
normalized induction preserve the embedding. In these dimensions the
only possible refinement replaces a tempered $\GL_2$ factor by an
induction of two unitary $\GL_1$ characters, both carrying the same
strictly negative twist. Thus we have the embedding required by
Lemma~\ref{blk:remove}, with unitary discrete-series blocks and a
tempered $\pi_0$. The zero-dimensional core is the trivial representation of the trivial
group.

\begin{theorem}[Every sufficiently ramified permissible character works]
\label{cov:ordinary}
Fix the additive character and the hyperbolic source and target spaces.
Let $\pi_1,\ldots,\pi_r$ be a finite family of irreducible admissible
smooth representations of $\U(W_4)$. There is an integer $N$ such that
for every smooth character $\chi:E^\times/F^\times\to\C^\times$ with
$a(\chi)>N$,
\[
 T^{\chi}_{V_4^+,W_4}(\pi_i)\ne0
             \qquad(1\le i\le r).
\]
The assertion also holds for a finite family drawn from both
$\U(W_2)$ and $\U(W_4)$, with each target of the corresponding dimension.
Such characters exist with arbitrarily large conductor.
\end{theorem}

\begin{proof}
Choose a tempered core $\pi_{i,0}$ and an anti-standard embedding as above
for each $\pi_i$. Each positive-dimensional core has a smooth unitary
central character on the compact group $E^1$. By
Lemma~\ref{blk:principal}, choose an irreducible tempered unitary principal
series $\pi_{i,0}^*$ of the \emph{same} unitary group with
\[
 \omega_{\pi_{i,0}^*}=\omega_{\pi_{i,0}}
                  \quad\hbox{on all of }E^1,
\]
whose equal-rank lift is in the hyperbolic tower for every permissible
$\chi$. All these comparators are chosen before imposing a conductor
bound on $\chi$.

For each of these finitely many pairs, the stability theorem
\cite[Proposition~5.3]{FinisLapid} provides a conductor bound with
the same additive character $\psi_0$ and the normalization fixed
in the analytic section. Let $N$ be the maximum of these bounds (and take $N=0$ if there are no positive-dimensional
cores). For every permissible $\chi$ with $a(\chi)>N$, inversion preserves
conductor, so stability applies to $\eta=\chi^{-1}$. The tempered analytic
criterion then gives
\begin{align*}
 \gamma_{\mathrm{LR}}\!\left(\tfrac12,\pi_{i,0},\chi^{-1},
                         \psi_0\right)
 &=\gamma_{\mathrm{LR}}\!\left(\tfrac12,\pi_{i,0}^*,\chi^{-1},
                         \psi_0\right)\\
 &=\omega_{\pi_{i,0}^*}(-1)
  =\omega_{\pi_{i,0}}(-1).
\end{align*}
These central values are finite and nonzero. The same criterion places
each positive-dimensional core in its hyperbolic equal-rank tower.
For a zero-dimensional core the only target is the zero space, and its
lift is nonzero.

For the fixed character $\chi$ under consideration, equal-rank dichotomy
(Theorem~\ref{ana:dichotomy}) provides a nonzero lift of each original
$\pi_i$ in one of the two towers. Lemma~\ref{blk:remove} forces
the target of that lift to be
\[
 \mathbb H^{\sum_j k_{ij}}\perp V_{n_{i,0}},
 \qquad T^{\chi}_{V_{n_{i,0}},W_{n_{i,0}}}(\pi_{i,0})\ne0.
\]
Uniqueness of the core's nonzero equal-rank tower makes $V_{n_{i,0}}$
hyperbolic. Adjoining the hyperbolic summand $\mathbb H^{\sum_j k_{ij}}$ therefore gives
the prescribed $V_4^+$, or $V_2^+$ when the original source has dimension two. In
dimension zero, the target is the sum of the hyperbolic planes
in the decomposition above.
Unbounded conductors of characters of $E^\times/F^\times$ give the last
assertion.
\end{proof}

\subsection{Restriction, reverse theta, and the quotient argument}
\label{fun:gu-passage}

Let $G=G_\theta$, $H=G'_\theta$, $G_1=G_\theta^1$, and
$H_1=(G'_\theta)^1$. The matching group $R$ is as in
Section~\ref{gu:jac:oscillator-section}. We fix permissible splitting data
$D=(\chi_V,\chi_W,\kappa)$ as in Appendix~\ref{gu:app:splitting},
and write $\varpi_D=\omega_{R,D}$ on $\cS=\mathcal S(Z)$. Define
\[
 \Omega_D=\cInd_R^{G\times H}\varpi_D,\qquad
 \Theta_D^{\GU}(\tau)=(\Omega_D\otimes\tau^\vee)_G.
\]
For $D=(\chi,1,1)$, the restriction to $G_1\times H_1$ is the ordinary
Weil representation used above. Ordinary big theta lifts have finite
length by Theorem~\ref{ana:dichotomy}; their irreducible subquotients
are admissible by \cite[Theorem~2.5]{BZ77}. Exactness of compact-open
invariants makes the big lifts admissible as well.

\subsubsection{Restriction and smooth duality}

\begin{lemma}\label{fun:restriction}
For \(\Pi\in\Irr(H)\), there are pairwise inequivalent
irreducible admissible representations \(\pi_i\) of \(H_1\) such that
\[
 \Pi|_{H_1}\simeq\bigoplus_{i=1}^r\pi_i,\qquad
 \Pi^\vee|_{H_1}\simeq\bigoplus_{i=1}^r\pi_i^\vee,\qquad r\in\{1,2\}.
\]
Their isomorphism classes form one \(H\)-conjugacy orbit.
\end{lemma}

\begin{proof}
The scalar center is \(Z=E^\times\). The section
\(d_s=\operatorname{diag}(I_2,sI_2)\) has multiplier \(s\), so
\[
 H_0=ZH_1=\lambda^{-1}(N_{E/F}E^\times),\qquad H/H_0\simeq C_2.
\]
The subgroup \(H_0\) is normal and open. The center acts on \(\Pi\)
by a smooth character \(\omega_\Pi\).

Finite-index Clifford theory makes \(\Pi|_{H_0}\) semisimple of
finite length. \(\Pi\) is generated by one vector
over \(H\), hence is finitely generated over \(H_0\), and has an
irreducible \(H_0\)-quotient \(\rho\). Frobenius reciprocity for the
open subgroup \(H_0\) embeds \(\Pi\) in
\(\operatorname{Ind}_{H_0}^H\rho\). Restricted to \(H_0\), this
induction is the direct sum of the two conjugates of \(\rho\).

Because \(Z\) acts by scalars, the \(H_1\)-submodules and
\(H_0\)-submodules of \(\Pi\) are identical. Irreducibility makes
the restriction constituents one conjugacy orbit. Their common
multiplicity is the dimension of an irreducible projective
representation of an inertia quotient in \(C_2\). Such a
representation is one-dimensional: a matrix representing a
generator has square \(cI\), and rescaling by a square root of
\(c\ne0\) reduces to an ordinary cyclic representation.

The restriction is admissible. If \(J\subset H_1\) is compact open,
choose a compact open \(Z_0\subset Z\) in \(\ker\omega_\Pi\).
The central-product map \(Z\times H_1\to H_0\) is open, and
\(Z_0J\) is compact open in \(H\). Thus
\(\Pi^J=\Pi^{Z_0J}\) is finite-dimensional.

If a functional on \(\Pi\) is \(H_1\)-smooth, a subgroup \(J\) as
above fixes it, and so does \(Z_0J\). It is therefore \(H\)-smooth.
The two smooth subspaces of the algebraic dual coincide, proving
the assertion about contragredients. These are the elementary statements of~\cite[Lemmas~5.2--5.3, \S9]{BGS}.
\end{proof}

\subsubsection{The reverse module}

Define
\[
 M_D(\Pi)=(\cS\otimes\Pi^\vee)_{H_1}.
\]
The \(R\)-action on the tensor product descends across the normal
subgroup \(1\times H_1\), so \(M_D(\Pi)\) carries an action of
\(R/(1\times H_1)\simeq G\). Namely,
\begin{equation}\label{fun:matching-action}
 g[v\otimes\ell]=[\varpi_D(g,h)v\otimes\Pi^\vee(h)\ell],
 \qquad \lambda(h)=\lambda(g).
\end{equation}
The formula is independent of the matching \(h\). The action is
smooth, since the quotient projection has the continuous section
\(g\mapsto(g,d_{\lambda(g)})\).

Lemma~\ref{gu:app:comparison} gives the $G$-equivariant isomorphism
\[
 (\Omega_D\otimes\Pi^\vee)_H\simeq M_D(\Pi),
\]
with quotient map \eqref{gu:asm:quotient-integral}.

\subsubsection{The fixed-data equivalence}

\begin{proposition}\label{fun:fixed-data}\label{gu:asm:assembly}
For \(\Pi\in\Irr(H)\) and fixed \(D\), the following
are equivalent:
\begin{enumerate}
\item Some irreducible admissible \(\tau\) of \(G\) has a quotient
\(\Theta_D^{\GU}(\tau)\twoheadrightarrow\Pi\).
\item \(M_D(\Pi)\ne0\).
\item Some constituent of \(\Pi|_{H_1}\) has nonzero ordinary big
theta lift to \(G_1\) for \(D\).
\end{enumerate}
The reverse module is of finite length and admissible. Any
irreducible quotient of a nonzero \(M_D(\Pi)\) gives the
surjection in (1).
\end{proposition}

\begin{proof}
Lemma~\ref{fun:restriction} and \eqref{fun:matching-action} give
\[
 M_D(\Pi)|_{G_1}
 \simeq\bigoplus_i(\cS\otimes\pi_i^\vee)_{H_1}.
\]
The summands are ordinary big theta lifts. Thus (2) and (3) are
equivalent. Their finite length gives finite length of \(M_D(\Pi)\),
because every strict chain of \(G\)-submodules is a strict chain
of \(G_1\)-submodules.

The reverse module is also admissible. For compact open \(K\subset G\),
\(K\cap G_1\) is compact open in \(G_1\), and
\[
 M_D(\Pi)^K\subset
 \bigl(M_D(\Pi)|_{G_1}\bigr)^{K\cap G_1}.
\]
The right side is finite-dimensional. Compact averaging makes
\(K\)-invariants exact, so every quotient of \(M_D(\Pi)\) is
admissible.

Suppose \(M_D(\Pi)\ne0\), and choose an irreducible quotient
\(p:M_D(\Pi)\twoheadrightarrow\tau\). Lemmas~\ref{fun:tensor}
and \ref{fun:reverse-comparison} give
\[
\begin{aligned}
 \Hom_G(M_D(\Pi),\tau)
 &\simeq\Hom_{G\times H}
       (\Omega_D,\tau\boxtimes\Pi)\\
 &\simeq\Hom_H
       (\Theta_D^{\GU}(\tau),\Pi).
\end{aligned}
\]
The two tensor-evaluation steps use admissibility of \(\Pi\)
and \(\tau\), respectively. The map corresponding to \(p\) is
nonzero and hence surjective. Conversely, a nonzero map in (1)
gives a nonzero element of the left Hom space, proving (2).
\end{proof}

For these groups, some constituent in (3) can be replaced by
every constituent. They form one \(H\)-orbit, and each \(h\in H\)
has a matching \(g\in G\). The operator \(\varpi_D(g,h)\) transports
the associated ordinary theta modules up to conjugation of the
\(G_1\)-action. Their nonvanishing is therefore identical. In
particular, the existence part needs a splitting choice for only
one restriction constituent.

The scalar and multiplier actions for general $D$ are recorded in
Appendix~\ref{gu:app:splitting}, together with the possible multiplier
self-twists.

\subsection{Coverage and Shalika uniqueness}
\label{cov:gu-section}

\begin{theorem}[Finite-family variable-splitting coverage]
\label{cov:gu}\label{gu:asm:coverage}
The spaces $V_4^+$ and $W_4$, the additive character, and the
normalized matching oscillator are fixed throughout. For every finite family
\[
 \Pi_1,\ldots,\Pi_r\in\Irr(\GU(W_4))
\]
there is an integer $N$ such that for every smooth character
$\chi:E^\times/F^\times\to\C^\times$ with $a(\chi)>N$, there are
irreducible admissible smooth representations
$\tau_j\in\Irr(\GU(V_4^+))$ and surjections
\[
 \Theta^{\GU}_{\chi,1,1}(\tau_j)\twoheadrightarrow\Pi_j
                    \qquad(1\le j\le r).
\]
All these characters $\chi$ are unitary. They have unbounded conductor.
In particular, varying only $\chi_V=\chi$ suffices: $\chi_W=1$ and
$\kappa=1$ may remain fixed. 
\end{theorem}

\begin{proof}
Restriction of each $\Pi_j$ to $\U(W_4)$ is a sum of one or two
inequivalent irreducible admissible representations, by the restriction
lemma in the preceding section. Applying Theorem~\ref{cov:ordinary}
to all these constituents yields a common bound $N$ and simultaneous
ordinary nonvanishing for every $\chi$ above it. For each such $\chi$ and each $j$, the fixed-data equivalence
of Proposition~\ref{fun:fixed-data} supplies the required quotient.
Compactness of $E^\times/F^\times$ makes every smooth complex character
of this quotient unitary. 
\end{proof}

\begin{proof}[Proof of Theorem~\ref{thm:main-gu}]
Choose $q:\Theta^{\GU}_{\chi,1,1}(\tau)\twoheadrightarrow\Pi$ by
Theorem~\ref{cov:gu}. Precomposition by a surjective map is injective,
so
\[
 \Hom_{S_B^{\GU}}(\Pi,\sigma\boxtimes\psi_B)
 \lhook\joinrel\longrightarrow
 \Hom_{S_B^{\GU}}(\Theta^{\GU}_{\chi,1,1}(\tau),
                                  \sigma\boxtimes\psi_B).
\]
By Theorem~\ref{gu:thm:hom}, the space on the right is
\[
 \Hom_{H_{\mathrm U}}\!\left(
   \tau^\vee\boxtimes
   \bigl((\sigma^j)^\vee\otimes(\chi^{-1}\circ\det_E)\bigr),
   \mathbf1\right).
\]
Both $\tau^\vee$ and the indicated twist of $(\sigma^j)^\vee$ are irreducible admissible: smooth
contragredients preserve irreducibility in the admissible category,
and character twists preserve it. Theorem~\ref{thm:main-unitary} bounds this Hom space by one.
\end{proof}

\Needspace{8\baselineskip}
\appendix
\section{Functional analysis, tensor evaluation, and induction duality}\label{app:functional}
\label{gu:app:functional}
\subsection{Tensor evaluation and closed-subgroup duality}
\begin{lemma}[Tensor evaluation]\label{gu:app:tensor}\label{fun:tensor}
Let $A$ and $B$ be locally compact totally disconnected groups,
$\mathscr E$ a smooth $A\times B$-module, $\pi$ an admissible smooth
$A$-module, and $V$ a smooth $B$-module. Contraction gives a natural
isomorphism
\begin{equation}\label{gu:eq:tensor-lemma}
 \Hom_B((\mathscr E\otimes\pi^\vee)_A,V)
 \simeq\Hom_{A\times B}(\mathscr E,\pi\boxtimes V).
\end{equation}
All tensor products here are algebraic.
\end{lemma}
\begin{proof}
Let $b:\mathscr E\times\pi^\vee\to V$ represent a map on the left.
For $e\in\mathscr E$, choose a compact open subgroup $J\subset A$
fixing $e$. Diagonal invariance makes $b(e,-)$ factor through
$(\pi^\vee)_J=(\pi^\vee)^J$, where the equality is given by compact
averaging. Its pairing with $\pi^J$ is perfect: any functional on
$\pi^J$ extends to a smooth functional on $\pi$ by composition with
the averaging projection. Admissibility makes both spaces finite
dimensional. For dual bases $e_i\in\pi^J$ and
$e_i^\vee\in(\pi^\vee)^J$, the corresponding tensor is
\[
 F(e)=\sum_i e_i\otimes b(e,e_i^\vee)\in\pi\otimes V.
\]
It is characterized by $(\ell\otimes1)F(e)=b(e,\ell)$ for all
$\ell\in\pi^\vee$. This proves independence of $J$ and the bases,
linearity, and equivariance. In particular,
$b(ae,\ell)=b(e,\pi^\vee(a^{-1})\ell)$ gives the $\pi(a)$
action. Contraction is the inverse. 
\end{proof}

For irreducible \(\pi\), this identifies the coinvariant construction
with the usual big theta lift. To see the maximal-isotypic
statement explicitly, put \(T=(\mathscr E\otimes\pi^\vee)_A\).
The identity of \(T\) gives a map \(q:\mathscr E\to\pi\otimes T\).
The target is semisimple and isotypic, so its submodules have the
form \(\pi\otimes T_0\). If this is the image of \(q\), all
contractions \([e\otimes\ell]\) belong to \(T_0\). Those classes
span \(T\), proving surjectivity. Lemma~\ref{gu:app:tensor} supplies
the universal property.

\begin{lemma}[Closed-subgroup duality]\label{gu:app:duality}
Let $D$ be a closed unimodular subgroup of a unimodular locally compact
totally disconnected group $L$, and let $V$ be an admissible smooth
$L$-module. A nonzero invariant measure on $D\backslash L$ gives
\begin{equation}\label{gu:eq:closed-duality}
 \Hom_L(\cInd_D^L\mathbf1,V)
 \simeq\Hom_D(V^\vee,\mathbf1).
\end{equation}
\end{lemma}
\begin{proof}
Integration against the quotient measure identifies
\begin{equation}\label{gu:eq:compact-dual}
 (\cInd_D^L\mathbf1)^\vee=\Ind_D^L\mathbf1.
\end{equation}
Indeed, for a compact open $J\subset L$, the $J$-orbits on
$D\backslash L$ are compact open. A $J$-invariant distribution is a
scalar multiple of invariant measure on each orbit. Smooth
distributions are therefore right-smooth locally constant densities,
with no compact-support condition; every such density pairs with
compactly supported functions.

Smooth reflexivity of $V$ and ordinary smooth-induction reciprocity give
\[
 \begin{aligned}
 \Hom_L(\cInd_D^L\mathbf1,V)
 &\simeq\Hom_L(V^\vee,\Ind_D^L\mathbf1)\\
 &\simeq\Hom_D(V^\vee,\mathbf1).
 \end{aligned}
\]
For inverse transposition, evaluating on a compactly induced vector
fixed by a compact open subgroup gives a smooth functional on
$V^\vee$, hence a vector of $V$. The last map is evaluation at $1$; its inverse sends a
$D$-invariant functional $\ell$ to
$v^\vee\mapsto[k\mapsto\ell(kv^\vee)]$.
\end{proof}

\subsection{Fiber relations and integration}
\begin{lemma}\label{gu:app:fibers}
Let $N$ act smoothly on a vector space $E_0$ by $T$, and let
$\alpha_c$ be automorphisms of $N$, for $c\in F^\times$, such that
$c\mapsto T_{\alpha_c(n)}v$ is locally constant for fixed $n,v$.
Give $C_c^\infty(F^\times,E_0)$ the action
$[nF](c)=T_{\alpha_c(n)}F(c)$. Then pointwise projection induces
\[
 C_c^\infty(F^\times,E_0)_N
 \simeq C_c^\infty(F^\times,(E_0)_N).
\]
If the remaining multiplier action is pure translation
$F(c)\mapsto F(ct)$ or $F(c/t)$, its coinvariant quotient is
$(E_0)_N$, with quotient map integration against $d^\times c$.
\end{lemma}
\begin{proof}
Surjectivity follows by lifting the finitely many values of a locally
constant compactly supported function. For the kernel, partition the
support into finitely many compact open pieces and express each value
as a finite sum $(T_w-1)v$. At a point $c_0$ put
$n=\alpha_{c_0}^{-1}(w)$. Local constancy gives
$T_{\alpha_c(n)}v=T_wv$ on a neighborhood of $c_0$. Multiplication by
its characteristic function realizes this as a global $n$-relation.
A finite cover and disjoint refinement give the whole kernel; the
reverse inclusion is pointwise.

For translation, write $F=\sum_i1_{c_iJ}m_i$ for a common compact open
subgroup $J\subset F^\times$. Its integral is
$\operatorname{vol}(J)\sum_i m_i$. If this is zero, subtracting one
fixed coset expresses $F$ as a sum of translation differences.
Surjectivity follows from
$m\mapsto\operatorname{vol}(J)^{-1}1_Jm$. Both orientations have the
same relation span. All integrals are finite vector-valued sums.
\end{proof}

\subsection{Matching and product coinvariants}
\begin{lemma}\label{gu:app:comparison}\label{fun:reverse-comparison}\label{fun:forward-comparison}
Let $\varpi$ be a smooth representation of the matching group $R$ on
$\mathscr S$, and let $\Omega=\cInd_R^{G_\theta\times G'_\theta}\varpi$, realized
by section values $f(c)=f(q_c)$. For irreducible admissible $\tau$ of
$G_\theta$ and $\Pi$ of $G'_\theta$, put
\[
 V_\tau=(\mathscr S\otimes\tau^\vee)_{G_\theta^1},\qquad
 M_\Pi=(\mathscr S\otimes\Pi^\vee)_{(G'_\theta)^1}.
\]
Their matching actions are
\begin{align}
 h[v\otimes\ell_\tau]
  &=[\varpi(r_t,h)v\otimes\tau^\vee(r_t)\ell_\tau],
     &&t=\lambda(h),\label{gu:eq:forward-matching-action}\\
 g[v\otimes\ell]
  &=[\varpi(g,d_s)v\otimes\Pi^\vee(d_s)\ell],
     &&s=\lambda(g).\label{gu:eq:reverse-matching-action}
\end{align}
Integration gives $G'_\theta$- and $G_\theta$-equivariant isomorphisms:
\begin{align}
 (\Omega\otimes\tau^\vee)_{G_\theta}&\simeq V_\tau,
 &f\otimes\ell_\tau&\longmapsto
 \int[\varpi(e_c)^{-1}f(c)\otimes\tau^\vee(r_c^{-1})\ell_\tau]\,d^\times c,
 \label{gu:eq:forward-comparison}\\
 (\Omega\otimes\Pi^\vee)_{G'_\theta}&\simeq M_\Pi,
 &f\otimes\ell&\longmapsto
 \int[f(c)\otimes\Pi^\vee(d_c)\ell]\,d^\times c.
 \label{gu:asm:quotient-integral}
\end{align}
For $\varpi=\omega_{R,0}$, the first integrand uses $D_c^{-1}$.
\end{lemma}
\begin{proof}
The matching actions are independent of the chosen matching elements:
two choices differ by the group whose coinvariants are being taken.
Normality preserves the relations, and the multiplicative sections give
the group law and smoothness.

For the forward map, $\gamma(g)=\varpi(g,d_{\lambda(g)})$ is a
representation of $G_\theta$, and its induced action is
$f(c)\mapsto\gamma(g)f(c/\lambda(g))$. Lemma~\ref{gu:app:fibers}, initially
with constant automorphism, gives $G_\theta^1$-coinvariants
$C_c^\infty(F^\times,V_\tau)$. The residual action has fiber operator
$\varpi(e_s)\otimes\tau^\vee(r_s)$. Applying its inverse at $c$
turns this into pure translation. The same lemma proves that the
forward integral has exactly the $G_\theta$-relation space as kernel and is
onto. For $G'_\theta$-equivariance, substitute $q=ct$ and use
\[
 \varpi(e_{q/t})^{-1}\varpi(1,d_{q/t}hd_q^{-1})
 =\varpi(r_t,h)\varpi(e_q)^{-1},\qquad
 r_{q/t}^{-1}=r_t r_q^{-1}.
\]
This transports the action to \eqref{gu:eq:forward-matching-action}.

For the reverse map, identify the tensor space with
$C_c^\infty(F^\times,E_0)$, where $E_0=\mathscr S\otimes\Pi^\vee$,
and change fibers by
\begin{equation}\label{gu:asm:fiber-change}
 (\mathcal BF)(c)=(1\otimes\Pi^\vee(d_c))F(c).
\end{equation}
This and its inverse preserve local constancy and compact support.
For $u\in (G'_\theta)^1$ the new fiber action is
$T_{d_cud_c^{-1}}$, where $T_u=\varpi(1,u)\otimes\Pi^\vee(u)$.
Lemma~\ref{gu:app:fibers} applies with $\alpha_c=\operatorname{Ad}(d_c)$.
After taking these relations, $d_t$ acts by $F(c)\mapsto F(ct)$.
The lemma therefore proves that the reverse integral is onto with
the $G'_\theta$-relation space as kernel.

The equivariance directions can also be read directly. Put
$u(c,h)=d_chd_{ct}^{-1}\in (G'_\theta)^1$. Diagonal $h$-translation changes the
reverse integrand to
\[
 [\varpi(1,u(c,h))f(ct)\otimes\Pi^\vee(d_ch)\ell]
 =[f(ct)\otimes\Pi^\vee(d_{ct})\ell],
\]
using the coinvariant relation. Multiplicative Haar invariance proves
$G'_\theta$-invariance. For $g\in G_\theta$, the identity
$q_c(g,1)=(g,d_s)q_{c/s}$ and substitution $c=ys$ give
\[
 \int[\varpi(g,d_s)f(y)\otimes
             \Pi^\vee(d_s)\Pi^\vee(d_y)\ell]\,d^\times y,
\]
which is the action in \eqref{gu:eq:reverse-matching-action}.
Unimodularity of $R$ and $G_\theta\times G'_\theta$ makes the quotient
measure multiplicative Haar measure.
\end{proof}

\section{Normalization and orbit topology}\label{app:normalization}
\label{gu:app:normalization}
\subsection{The phase and Fourier covariance}
We prove the covariance used in Proposition~\ref{gu:jac:extension}, with
the conventions of Section~\ref{gu:sec:setup}. An upper shear
$(x,y)\mapsto(x+Ay,y)$ acts by multiplication by
$\psi(\tfrac12[AT,T])$: conjugating \eqref{gu:eq:heisenberg-action} uses
\[
 [A(T+y),T+y]-[AT,T]=2[Ay,T]+[Ay,y].
\]
For the shear $n(u)$, the tensor pairing gives
\[
 [(1\otimes u)T,T]
 =-\tfrac12\Tr_{E/F}\tr_E(u\overline{G_T}).
\]
Its half-phase is \eqref{gu:jac:ordinary-n}, since the matrix trace lies
in $F$. The conjugate of $G_T$ transports the form from $X^c$ to $X$.

For the Weyl element
\[
 w=\begin{pmatrix}0&I_2\\-I_2&0\end{pmatrix},
\]
the operator is the Fourier transform
\begin{equation}\label{gu:jac:fourier}
 [\mathcal F_\beta f](T)
 =\int_Zf(S)\psi(\beta(T,S))\,dS,
 \qquad
 \beta(T,S)=-\frac12\operatorname{Tr}_{E/F}
                         \operatorname{tr}_E(T^\dagger J_4S),
\end{equation}
where $dS$ is self-dual for $\beta$ and $\psi$.  The partial Weyl elements
act by the corresponding partial Fourier transforms.  The Weil-index
coefficients in these formulas are $1$: the quadratic $F$-space underlying
the split Hermitian space $V$ is hyperbolic, a hyperbolic plane has Weil
index $1$, and the index is multiplicative under orthogonal direct sum.
The same assertion holds after multiplying the quadratic form or the
additive character by any nonzero scalar.

To prove \eqref{gu:jac:covariance}, the $U(V)$ formula is immediate;
left multiplication by $r_s$ commutes with the Levi's right action.
For the upper radical use
\[
 G_{r_s^{-1}T}=s^{-1}G_T,\qquad
 d_sn(u)d_s^{-1}=n(s^{-1}u).
\]
It remains to verify covariance at the Fourier operator.  The bilinear
form in \eqref{gu:jac:fourier} satisfies
\[
 \beta(r_s^{-1}T,r_s^{-1}S)=s^{-1}\beta(T,S).
\]
Using $U=r_sS$ in the Fourier integral therefore gives
\begin{align*}
 [D_s\mathcal F_\beta D_s^{-1}f](T)
 &=\int_Z f(r_sS)\psi(\beta(r_s^{-1}T,S))\,dS\\
 &=|s|_F^{-8}\int_Z f(U)\psi(s^{-1}\beta(T,U))\,dU\\
 &=|s|_F^{-8}[\mathcal F_\beta f](s^{-1}T).
\end{align*}
On the group side,
\[
 d_swd_s^{-1}=m(s^{-1}I_2,1)w,
 \qquad |\det(s^{-1}I_2)|_E^2=|s|_F^{-8}.
\]
The final integral is consequently the operator on the right
side of \eqref{gu:jac:covariance}. The Siegel parabolic and $w$ generate
both Siegel radicals; in each hyperbolic
coordinate plane, a product of upper, lower, and upper one-parameter
unipotents gives the partial Weyl element.  The Siegel Bruhat
decomposition then gives all of $U(W)$.  This proves
\eqref{gu:jac:covariance} on the whole isometry subgroup.

Equations \eqref{gu:jac:semidirect-decomposition} and \eqref{gu:jac:covariance},
together with $D_sD_t=D_{st}$, prove the representation law.
The operator $D_s$ is the canonical polarized Heisenberg intertwiner for
the symplectic transformation \eqref{gu:jac:tensor-action} associated with
$e_s$.  Thus \eqref{gu:jac:extension-definition} has Heisenberg
covariance.  Smoothness follows from the operator formulas above and the
smoothness of $\omega$.

\subsection{Schwartz restriction and the orbit quotient}\label{gu:app:topology}
\begin{lemma}\label{gu:jac:restriction}
Restriction from $M=F^\times\times Z$ to $C=\widetilde{\mathcal O}_B$
is surjective, its kernel is the algebraic $(N(X),\psi_B)$-relation
space, and the induced map is a homeomorphism for the Schwartz
inductive-limit topologies. In particular, the relation space is closed.
\end{lemma}
\begin{proof}
By \eqref{gu:jac:induced-n}, the unipotent character on $C$ is
$\psi_B$, so every relation vanishes there. If $f|_C=0$, local constancy makes the support
of $f$ a compact subset of $M\setminus C$. At every point of this support
some scalar character difference
\[
 d_n(c,T)=\psi\left(-\tfrac1{2c}\tr_E(u\overline{G_T})\right)
                  -\psi_B(n(u))
\]
is nonzero. A finite disjoint compact open refinement of these
neighborhoods writes $f=\sum_i f_i$, with $d_{n_i}$ nonzero on the
support of $f_i$. Extending $f_i/d_{n_i}$ by zero gives a Schwartz
function $g_i$ and
$f_i=(\Omega_0^{\GU}(n_i)-\psi_B(n_i))g_i$.
Conversely, a finite constancy partition of a Schwartz function on $C$
extends to compact open neighborhoods in $M$, chosen to avoid the other
values and the zero part of $C$. Disjoint refinement proves surjectivity.

The spaces under consideration are second-countable locally compact
totally disconnected spaces.  Their Schwartz spaces are inductive limits
of finite-dimensional spaces obtained by bounding the support and fixing
a finite compact open constancy partition.  Every finite-dimensional
subspace is contained in one such stage.  They therefore have the finest
locally convex topology on their underlying countable-dimensional vector
spaces.  Restriction is continuous, its kernel is closed, and the
algebraic quotient just constructed is a homeomorphism for these
topologies.  Equivalently, the extension argument supplies a continuous
lift on each finite-dimensional stage.
\end{proof}

\begin{lemma}\label{gu:jac:orbit-topology}
The transitive orbit map for \eqref{gu:jac:orbit-action} is a submersion.
Its induced bijection $L/H_j\to\widetilde{\mathcal O}_B$ is a
homeomorphism.
\end{lemma}
\begin{proof}
We verify the quotient topology. At a point of the scaled Gram locus,
$G_T=c\overline{\mathsf B}$ is invertible because $c\ne0$ and $B$ is
nondegenerate. The differential of the Gram map is
\[
 dG_T(S)=S^\dagger J_4T+T^\dagger J_4S.
\]
For any Hermitian matrix $D$, taking
$S=\tfrac12TG_T^{-1}D$ gives $dG_T(S)=D$.  Thus the equations defining
$\widetilde{\mathcal O}_B$ make it a smooth $F$-analytic manifold of
dimension $1+16-4=13$.

The orbit map from $G_\theta$ is a submersion.  It suffices to check this at
$(1,j)$.  A tangent vector $(\dot c,S)$ there satisfies
$dG_j(S)=\dot c\overline{\mathsf B}$. We use the decomposition
$V=j(X^c)\oplus j(X^c)^\perp$ and write $S=(b,C)$, identifying the first
summand by $j$.  With $*$ denoting adjoint for the forms on these two
nondegenerate summands, the tangent equation is $b^*+b=\dot c I$.
The block endomorphism
\[
 A=\begin{pmatrix}b&-C^*\\C&\dot c I/2\end{pmatrix}
\]
satisfies $A^*+A=\dot c I$ and $Aj=S$.  It is a tangent vector of $G_\theta$
with image $(\dot c,S)$.  Translation gives surjectivity of the
differential at every point.  The nonarchimedean inverse-function
theorem now gives local continuous analytic sections of the orbit map.
In particular, the orbit map is open.  Its induced bijection from the
quotient is therefore a homeomorphism.
\end{proof}

\section{General splitting data}\label{app:splitting}
\label{gu:app:splitting}
The preceding arguments use $D=(\chi,1,1)$. Here we give the twisting
formulas and the corresponding theta quotient statement for general
fixed data.

\subsection{Extensions and character twists}
Any oscillator extension implementing the same Heisenberg
automorphisms and agreeing with $\omega$ on $G_\theta^1\times (G'_\theta)^1$ differs
from $\omega_{R,0}$ by scalars, by the Stone--von Neumann theorem.
The representation law makes these scalars a character of
$R/(G_\theta^1\times (G'_\theta)^1)\simeq F^\times$. Thus
\begin{equation}\label{gu:jac:extension-twist}
 \omega_{R,\kappa}=\omega_{R,0}(\kappa\circ\lambda_R).
\end{equation}
Here $\lambda_R$ is the common multiplier. The extensions trivial
on $K_E$ are those with $\kappa|_{N_{E/F}E^\times}=1$. The matching Levi
normalization \eqref{gu:jac:normalized-levi} selects $\kappa=1$.

Let $D=(\chi_V,\chi_W,\kappa)$, with $\chi_V,\chi_W$ unitary
characters trivial on $F^\times$ and $\kappa$ smooth on $F^\times$.
Using \eqref{gu:eq:ambient-character}, put
\begin{equation}\label{gu:eq:A-C}
 A(g)=\mathcal A_{\chi_W}(g)\kappa(\lambda_{G_\theta}(g)),\qquad
 C(h)=\mathcal A_{\chi_V}(h),\qquad
 \omega_{R,D}=\omega_{R,0}(A\boxtimes C)|_R.
\end{equation}
Its ordinary restriction has splitting characters $(\chi_V,\chi_W)$.
Define $\Omega_D^{\GU}=\cInd_R^{G_\theta\times G'_\theta}\omega_{R,D}$ and
$\Theta_D^{\GU}(\tau)=(\Omega_D^{\GU}\otimes\tau^\vee)_{G_\theta}$.
Multiplication of inducing functions by the ambient character preserves
support and gives
\begin{equation}\label{gu:eq:ambient-twist}
 \Omega_D^{\GU}\simeq\Omega_0^{\GU}\otimes(A\boxtimes C),\qquad
 \Theta_D^{\GU}(\tau)
 \simeq C\otimes\Theta_0^{\GU}(\tau\otimes A^{-1}).
\end{equation}
Write $\Omega_\kappa^{\GU}=\Omega_{(1,1,\kappa)}^{\GU}$.
For the pure multiplier case $D=(1,1,\kappa)$, specialize
\eqref{gu:eq:ambient-twist} to $A=\kappa\circ\lambda_{G_\theta}$ and $C=1$.
Writing $\lambda_{H_j}$ for the common multiplier on $H_j$,
Proposition~\ref{gu:jac:jacquet} gives
\[
 (\Omega_\kappa^{\GU})_{N,\psi_B}
 \simeq\cInd_{H_j}^{G_\theta\times K}(\kappa\circ\lambda_{H_j}),
\]
with intertwiner $\kappa(\lambda(g))\mathcal I_j$ in the same coordinates.

The scalar formula $\mathcal A_\eta(zI_4)=\eta(z)^2$ gives the
character of $K_E$:
\begin{equation}\label{gu:eq:kernel-character}
 k_D(z)=\chi_W(z)^2\chi_V(z)^{-2}\kappa(N_{E/F}z).
\end{equation}
Every nonzero quotient $\tau\boxtimes\Pi$ of $\Omega_D^{\GU}$ thus has
\begin{equation}\label{gu:eq:theta-central}
 \omega_\tau(z)\omega_\Pi(\bar z)=k_D(z).
\end{equation}
For normalized data the scalar character on the nonzero big lift
is $\omega_{\Theta_0^{\GU}(\tau)}(z)=\omega_\tau(\bar z)^{-1}$,
directly from the tensor kernel.

\begin{proposition}[Shalika period identity for general splitting data]
\label{gu:thm:hom-general}
For these data and irreducible admissible $\tau$ and $\sigma$,
\begin{equation}\label{gu:eq:hom-general}
 \begin{aligned}
 &\Hom_{S_B^{\GU}}(\Theta_D^{\GU}(\tau),\sigma\boxtimes\psi_B)\\
 &\quad\simeq\Hom_{H_{\mathrm U}}\left(
   (\tau^\vee\otimes A)\boxtimes
   ((\sigma^j)^\vee\otimes(\chi_V^{-1}\circ\det_E)),\mathbf1\right).
 \end{aligned}
\end{equation}
The dimension is at most one.
\end{proposition}
\begin{proof}
On the Shalika subgroup,
$C(m_B(a))=\chi_V(\det_Ea)$ and $C|_N=1$. We apply
\eqref{gu:eq:ambient-twist} and the normalized case of
Theorem~\ref{gu:thm:hom} to $\tau\otimes A^{-1}$ and
$\sigma\otimes(\chi_V^{-1}\circ\det_X)$. The transport and duality
in \eqref{gu:eq:coefficient-three-steps}, with $\chi_V$ in place of
$\chi$, give \eqref{gu:eq:hom-general}; Theorem~\ref{thm:main-unitary} gives
the bound.
\end{proof}
Equivalently, the period character before absorbing the twists is
$\Lambda_D(u,h)=A(u\oplus h)^{-1}\chi_V(\det_Eh)$; for a pure
multiplier twist it is $\kappa^{-1}\circ\lambda_{H_{\mathrm U}}$.
On the common scalar center the ambient representation in
\eqref{gu:eq:hom-general} has character
$\omega_\tau(z)^{-1}\omega_\sigma(\bar z)^{-1}k_D(z)$.
If $\Pi$ is a theta quotient with $\omega_\Pi=\omega_\sigma$, this
is one by \eqref{gu:eq:theta-central}.

\subsection{The reverse module for fixed general data}
For the module $M_D(\Pi)$ and matching action of
Section~\ref{fun:gu-passage}, the scalar and multiplier twists are
\begin{align}
 \omega_{M_D(\Pi)}(z)&=k_D(z)\omega_\Pi(\bar z)^{-1},\\
 M_D(\Pi\otimes(\eta\circ\lambda_{G'_\theta}))
 &\simeq M_D(\Pi)\otimes(\eta^{-1}\circ\lambda_{G_\theta}),
 \label{gu:asm:multiplier-twist}
\end{align}
where $\eta:F^\times\to\mathbb C^\times$ is smooth.
The first formula follows by taking the matching element
$h=\bar zI_W$ in \eqref{fun:matching-action}; the second follows by
twisting its second tensor factor.

Put $H=G'_\theta$, $H_1=(G'_\theta)^1$, and $H_0=E^\times H_1$.
If $\Pi\simeq\Pi\otimes(\eta\circ\lambda_H)$,
its central character forces $\eta|_{N(E^\times)}=1$. Thus $\eta$
is trivial or $\omega_{E/F}$. Write $\xi=\omega_{E/F}\circ\lambda_H$.
The index-two induction formula and Frobenius reciprocity give
\[
 r=\dim\End_{H_0}(\Pi)
   =1+\dim\Hom_H(\Pi,\Pi\otimes\xi).
\]
Here $r\in\{1,2\}$ is the restriction length in
Lemma~\ref{fun:restriction}. The nontrivial multiplier self-twist
therefore occurs precisely when $r=2$, a case included in the proof.

\begin{corollary}[Fixed general splitting data]\label{gu:app:general-assembly}
For fixed $D$ and $\Pi\in\Irr(G'_\theta)$, a constituent of
$\Pi|_{(G'_\theta)^1}$ has nonzero ordinary lift if and only if
$\Theta_D^{\GU}(\tau)\twoheadrightarrow\Pi$ for some
$\tau\in\Irr(G_\theta)$.
\end{corollary}
\begin{proof}
This is Proposition~\ref{fun:fixed-data} for the fixed data $D$.
\end{proof}

\section{Rational transpose transporters}\label{app:transporters}
We retain the notation of Section~\ref{sec:unitary-geometry} and
write $\mu=\lambda(x)$. The two rank-one families are the
decomposable and stable-system families of \cref{m2:rat:rankone-param}.

The following rational examples illustrate the two rank-one
families. We set
\[
 P_\pm=\frac12\begin{pmatrix}1&\pm1\\\pm1&1\end{pmatrix},
 \qquad
 x_{\mathrm{dec}}=\begin{pmatrix}P_-&P_+\\P_+&P_-\end{pmatrix}.
\]
The two projectors are $J$-selfadjoint and have orthogonal
images.  The matrix swaps the two $P_+$ lines and fixes
the two $P_-$ lines, so it is unitary.  It is a
decomposable point.  For $t\in F^\times$, $t^2\ne1$,
let $h_t=\diag(t,t^{-1})$ and
\begin{equation}\label{m2:rat:rankone-examples}
 x_t=x_{\mathrm{dec}}\diag(I,h_t).
\end{equation}
This is again unitary.  Its incidence lines are
$\im P_+$ and $h_t^{-1}\im P_-$, neither stable under
$D=P_-h_t$, so it is a stable system.  Both examples have multiplier one.

We use the real-pair coordinates of \cref{m2:rat:real-pair}.

\begin{lemma}\label{m2:rat:dense-transporter}
For every dense closed point
$x=\left(\begin{smallmatrix}A_x&B_x\\C_x&D\end{smallmatrix}\right)$
there is a rational transporter $k_x=(u_x,v_x,h_x)\in K(F)$
with $k_x\cdot x^t=x$.  It can be constructed by
linear equations over $F$.  If $D$ is scalar, one may
choose $h_x=1$.
\end{lemma}

\begin{proof}
In the real-pair coordinates of \cref{m2:rat:real-pair},
closedness implies that
$(A,B_0)$ is geometrically semisimple: its generated
algebra after extension to $E$, and hence to $\barF$,
is the algebra in \eqref{m2:rat:dense-algebra}.
Choose an invertible symmetric solution $H$ over $F$
of
\[
 AH=HA^t,\qquad B_0H=HB_0^t
\]
by \cref{m2:geo:pair-transpose}.  This is a linear system
on three entries; the lemma shows that an invertible
solution exists in every semisimple case.  Put
\begin{equation}\label{m2:rat:dense-h}
 h_x=P(HR^{-1})P^{-1}.
\end{equation}
By \eqref{m2:rat:real-transpose},
$h_xD^th_x^{-1}=D$, and by
\cref{m2:rat:real-pair} $h_x\in G_2(F)$.
Moreover $R^{-1}=P^2$, so $h_x=PHP$ is symmetric.
Define
\begin{equation}\label{m2:rat:dense-k}
 u_x=B_xh_xC_x^{-t},\qquad
 v_x=C_x^{-1}h_xB_x^t.
\end{equation}
The $B,C,D$ transporter equations follow immediately.
We check the similitudes as well.  Writing
$\lambda=\lam(h_x)$, the relation
$Dh_x=h_xD^t$ gives
\[
 h_x^\dagger Q_Dh_x
 =\lambda(\mu J-\bar D J D^t)
 =\lambda\bar C_x J C_x^t,
\]
where the last equality is the transpose of the last
row equation in \eqref{m2:geo:row-equations}.  Therefore
$u_x^\dagger Ju_x=\lambda J$.
Symmetry of $h_x$ gives
$v_x=C_x^{-1}h_x^tB_x^t=u_x^t$.
Ordinary transpose preserves $G_2$ and its multiplier,
so $v_x^\dagger Jv_x=\lambda J$ as well.
Thus $(u_x,v_x,h_x)\in K(F)$.
Acting by this element on $x^t$ produces a unitary
matrix with the same multiplier and the same $B,C,D$
blocks as $x$.  By the uniqueness of $A=B_xX_DC_x$
in \cref{m2:rat:dense-param}, its $A$-block also agrees.

If $D$ is scalar, it is fixed by transpose, and one
can take $h_x=1$ directly in \eqref{m2:rat:dense-k}.
Here $B_x$ and $C_x$ have the same multiplier
$\mu-\N(D)$, by the two Gram equations, so the two
remaining factors have multiplier one.  The same
uniqueness argument verifies the $A$-block.
\end{proof}

For rank zero write $P^{-1}DP=A+\delta B_0$ over $F$
as in \cref{m2:rat:real-pair}.  Since $D$ is semisimple
and $D^*=\mu D^{-1}$, the algebra
$E\langle D,D^*\rangle=E[D]$ is semisimple.
Thus the real pair $(A,B_0)$ is geometrically semisimple.
Applying \cref{m2:geo:pair-transpose} to this pair and
using \eqref{m2:rat:dense-h} gives $h_x\in G_2(F)$ with
$h_xD^th_x^{-1}=D$.  Choose any $v_x\in G_2(F)$
of multiplier $\lam(h_x)$, for example
$\diag(1,\lam(h_x))$, and set
$u_x=a_\mu v_xa_\mu^{-1}$.  This gives an explicit
transporter for \eqref{m2:rat:rankzero-representative}.

For rank one we use the scalar endomorphism algebra to construct
a transporter.

\begin{proposition}\label{m2:rat:rankone-transporter}
Let $x$ be an arbitrary rank-one closed $F$-point.
A rational transporter from $x^t$ to $x$ can be
constructed by solving an $E$-linear system followed
by quadratic local norm equations.
\end{proposition}

\begin{proof}
First solve, without imposing the unitary condition,
\begin{equation}\label{m2:rat:linear-transporter}
 u_0A^t=Av_0,\quad u_0C^t=Bh_0,\quad
 h_0B^t=Cv_0,\quad h_0D^t=Dh_0.
\end{equation}
These equations are linear in the entries of
$u_0,v_0,h_0$ over $E$.  An invertible solution exists
by \cref{m2:geo:transpose} applied to the split-quiver
representation over $E$.  After row reduction it may
be selected by avoiding the zero locus of the product
of the three determinants in the free variables.
Let $t_0=(u_0,v_0,h_0)$ be such a solution, let
$y=x^t$, and write $L_0y=xR_0$ for its block equation.

Set
\[
 a=(u_0^*u_0,v_0^*v_0,h_0^*h_0).
\]
The unitary equations for $x,y$, with the same
multiplier $\mu$, imply $L_0^*x=yR_0^*$.
It follows that
\[
 L_0^*L_0y=yR_0^*R_0.
\]
Thus $a$ is an invertible selfadjoint element of
$\mathscr A_y=E^m$, with $m=3$ or $2$.
Its decomposition along the selfadjoint primitive idempotents
of \cref{m2:rat:selfadjoint-idempotents} is
$a=\sum_i a_ie_i$, with $a_i\in F^\times$.

At every vertex the determinant of $t_0^*t_0$ is
$\N_{E/F}(\det t_0)$, hence a norm.  In the
decomposable case the first and middle vertices
therefore say that $a_ra_a$ and $a_ra_b$ are norms.
In the stable-system case the first vertex says
that $a_ra_a$ is a norm.  Since $Q_F$ has exponent
two, in either case all the $a_i$ have the same
norm class, by the incidence-determinant computation of
\cref{m2:rat:rankone-kernel} applied to $t_0^*t_0$.
Set $\lambda=a_r$ and choose
$s_i\in E^\times$ with
\[
 \N(s_i)=\lambda/a_i.
\]
These are solvable local norm equations.  For
$s=\sum_i s_ie_i\in\mathscr A_y^\times$ one has
$s^*as=\lambda I$ at every vertex.  Hence
$k_x=t_0s$ is a triple of Hermitian similitudes
with the same multiplier $\lambda$.  Since $s$
stabilizes $y$, the triple $k_x$ still transports
$y$ to $x$.  This constructs the required element
of $K(F)$.
\end{proof}

For the examples at the start of this appendix, $x_{\mathrm{dec}}$ is
symmetric and one may take $k_x=1$.  For
\eqref{m2:rat:rankone-examples}, direct block
multiplication gives the transporter
$(I,I,h_t^{-1})$.

\bibliographystyle{amsalpha}
\bibliography{references}
\end{document}